\makeatletter

\newdimen\paperwidth
\newdimen\paperheight

\def\papersize#1#2{\let\p@persize\relax\paperwidth#1\paperheight#2}

\def\Afour{\papersize{210truemm}{297truemm}}

\let\p@persize\Afour

\let\onesidestyle\@twosidefalse
\let\twosidestyle\@twosidetrue

\def\margins{\@ifnextchar[{\@margins}{\@margins[\z@]}}

\def\@margins[#1]#2#3{
  \p@persize\dimen0 #3\dimen0 .5\dimen0\normalsize%
  \oddsidemargin-1truein\advance\oddsidemargin#2%
  \evensidemargin-1truein\advance\evensidemargin#2%
  \topmargin-1truein\advance\topmargin\dimen0\headsep\dimen0\footskip\dimen0%
  \textwidth\paperwidth\advance\textwidth-#2\advance\textwidth-#2%
  \textheight\paperheight\advance\textheight-#3\advance\textheight-#3%
  \headheight\baselineskip\advance\topmargin-.5\baselineskip%
  \advance\headsep-.5\baselineskip%
  \footheight\baselineskip
  \advance\textwidth-#1\advance\oddsidemargin#1
  \if@twoside\def\@themargin%
    {\ifodd\count\z@\oddsidemargin\else\evensidemargin\fi}\fi}

\def\headlinesep#1{\advance\topmargin\headsep\advance\topmargin -#1
  \advance\topmargin.5\baselineskip\headsep #1\advance\headsep-.5\baselineskip}

\def\headline{\if@twoside\let\n@xt\h@dlin@\else\let\n@xt\h@@dlin@\fi\n@xt}
  
\def\h@dlin@#1#2{%
  \def\@oddhead{%
    {{\leftskip\z@\rightskip\z@\noindent\normalsize#1}}}
  \def\@evenhead{%
    {{\leftskip\z@\rightskip\z@\noindent\normalsize#2}}}}

\def\h@@dlin@#1{%
  \def\@oddhead{{{\leftskip\z@\rightskip\z@\noindent\normalsize#1}}}}

\def\footline{\if@twoside\let\n@xt\f@tlin@\else\let\n@xt\f@@tlin@\fi\n@xt}

\def\f@tlin@#1#2{%
  \def\@oddfoot{%
    {{\leftskip\z@\rightskip\z@\noindent\normalsize#1}}}
  \def\@evenfoot{%
    {{\leftskip\z@\rightskip\z@\noindent\normalsize#2}}}}

\def\f@@tlin@#1{%
  \def\@oddfoot{{{\leftskip\z@\rightskip\z@\noindent\normalsize#1}}}}

\def\normalpage{\global\@specialpagefalse}

\makeatother

\makeatletter

\def\ft{\@ifnextchar[{\ft@s}{\ft@}}
\def\ft@{\ft@@@s[\f@size]}
\def\ft@s[{\@ifnextchar{a}{\ft@sz[}{\ft@@s[}}
\def\ft@@s[{\@ifnextchar{s}{\ft@sz[}{\ft@@@s[}}
\def\ft@@@s[#1]{\ft@sz[at #1pt]}
\def\ft@sz[#1]#2{\font\fonttemp=#2 #1\fonttemp\ignorespaces}

\makeatother

\makeatletter

\input{epsf.sty}

\def\showfig#1#2{\epsfbox{#2}}
\def\fig@#1#2{\leavevmode{\framebox{\figstyl@\strut{ #1 }}}}

\def\figstyle#1{\def\figstyl@{#1}}

\figstyle{\shape{n}\series{m}\selectfont\normalsize}

\def\showfigurestrue{\let\fig\showfig}
\def\showfiguresfalse{\let\fig\fig@}

\makeatother

\showfiguresfalse

\def\smallcircc{\mathop{\mkern3.5mu\hbox{\raise.58ex\hbox{\ft{lcircle10}a}}}}
\def\varemptyset{{\hbox{\raise.21ex\hbox{$\not$}}\mkern.15mu\mathrm{O}\mkern.15mu}}

\let\texepsilon\epsilon  \let\epsilon\varepsilon
\let\textheta\theta      \let\theta\vartheta
          \let\phi\varphi
   \let\emptyset\varemptyset

\documentstyle[12pt]{article}

\makeatletter

\let\Larg@\Large
\let\hug@\huge

\def\usepackage#1{\input{#1.sty}}

\input{geom.sty}

\def\r@adlabel#1#2{\global\@namedef{#1@\the\@key}{#2}}

\let\Large\Larg@
\let\huge\hug@

\def\smallskip{\vskip\smallskipamount}
\def\medskip{\vskip\medskipamount}
\def\bigskip{\vskip\bigskipamount}

\def\mytrivlist{\parsep\parskip\@nmbrlistfalse
  \my@trivlist \labelwidth\z@ \leftmargin\z@
  \itemindent\z@ \def\makelabel##1{##1}}

\def\my@trivlist{\global\@newlisttrue \@outerparskip\parskip}

\def\end#1{\csname end#1\endcsname\@checkend{#1}%
  \expandafter\endgroup\if@endpe\@doendpe\fi
  \if@ignore \global\@ignorefalse \ignorespaces\fi}

\def\put{\@ifnextchar[{\@put}{\@@rput[\z@,\z@][r]}}
\def\@put[#1]{\@ifnextchar[{\@@put[#1]}{\@@@@@put[#1]}}
\def\@@put[#1][{\@ifnextchar{l}{\@@lput[#1][}{\@@@put[#1][}}
\def\@@@put[#1][{\@ifnextchar{c}{\@@cput[#1][}{\@@@@put[#1][}}
\def\@@@@put[#1][{\@ifnextchar{r}{\@@rput[#1][}{\relax}}
\def\@@@@@put[{\@ifnextchar{l}{\@@lput[\z@,\z@][}{\@@@@@@put[}}
\def\@@@@@@put[{\@ifnextchar{c}{\@@cput[\z@,\z@][}{\@@@@@@@put[}}
\def\@@@@@@@put[{\@ifnextchar{r}{\@@rput[\z@,\z@][}{\@@@@@@@@put[}}
\def\@@@@@@@@put[#1]{\@@rput[#1][r]}

\let\hm@d@\leavevmode

\long\def\@@lput[#1,#2][l]#3{\setbox0\hbox{#3}\hm@d@\raise#2\hbox to\z@{\dimen0 #1%
  \advance\dimen0-\wd0\kern\dimen0\dp0\z@\ht0\z@\wd0\z@\box0\hss}\ignorespaces}
\long\def\@@cput[#1,#2][c]#3{\setbox0\hbox{#3}\hm@d@\raise#2\hbox to\z@{\dimen0 #1%
  \advance\dimen0-.5\wd0\kern\dimen0\dp0\z@\ht0\z@\wd0\z@\box0\hss}\ignorespaces}
\long\def\@@rput[#1,#2][r]#3{\setbox0\hbox{\kern#1\raise#2\hbox{#3}}%
  \dp0\z@\ht0\z@\wd0\z@\hm@d@\box0\ignorespaces}

\def\flbox{\@ifnextchar[{\@flbox}{\@@rflbox[\z@,\z@][r]}}
\def\@flbox[#1]{\@ifnextchar[{\@@flbox[#1]}{\@@@@@flbox[#1]}}
\def\@@flbox[#1][{\@ifnextchar{l}{\@@lflbox[#1][}{\@@@flbox[#1][}}
\def\@@@flbox[#1][{\@ifnextchar{c}{\@@cflbox[#1][}{\@@@@flbox[#1][}}
\def\@@@@flbox[#1][{\@ifnextchar{r}{\@@rflbox[#1][}{\relax}}
\def\@@@@@flbox[{\@ifnextchar{l}{\@@lflbox[\z@,\z@][}{\@@@@@@flbox[}}
\def\@@@@@@flbox[{\@ifnextchar{c}{\@@cflbox[\z@,\z@][}{\@@@@@@@flbox[}}
\def\@@@@@@@flbox[{\@ifnextchar{r}{\@@rflbox[\z@,\z@][}{\@@@@@@@@flbox[}}
\def\@@@@@@@@flbox[#1]{\@@rflbox[#1][r]}
\long\def\@@lflbox[#1,#2][l]#3{\@@lput[#1,#2][l]{%
  \vtop{\leftskip\z@\parindent\z@\raggedleft\hm@d@#3}}}
\long\def\@@cflbox[#1,#2][c]#3{\@@cput[#1,#2][c]{%
  \vtop{\leftskip\z@\parindent\z@\raggedcenter\hm@d@#3}}}
\long\def\@@rflbox[#1,#2][r]#3{\@@rput[#1,#2][r]{%
  \vtop{\leftskip\z@\parindent\z@\raggedright\hm@d@#3}}}

\def\maketitle{\par
 \begingroup
 \def\thefootnote{\fnsymbol{footnote}}
 \def\@makefnmark{\hbox to 0pt{$^{\@thefnmark}$\hss}} 
 \if@twocolumn 
 \twocolumn[\@maketitle] 
 \else 
 \global\@topnum\z@ \@maketitle \fi\thispagestyle{plain}\@thanks
 \endgroup
 \setcounter{footnote}{0}
 \let\maketitle\relax
 \let\@maketitle\relax
 \gdef\@thanks{}\gdef\@author{}\gdef\@title{}\let\thanks\relax}

\def\@maketitle{ 
 \null
 \vskip 2em \begin{center}
 {\LARGE \@title \par} \vskip 1.5em {\large \lineskip .5em
\begin{tabular}[t]{c}\@author 
 \end{tabular}\par} 
 \vskip 1em {\large \@date} \end{center}
 \par
 \vskip 1.5em}
 
\def\partbeforeskip#1{\def\p@rtbeforeskip{#1}}
\def\partstyle#1{\def\p@rtstyl@{#1}}
\def\partdot#1{\def\partd@t{#1}}
\def\partafterskip#1{\def\p@rtafterskip{#1}}
\def\partintrostyle#1{\def\partintr@styl@{#1}}
\def\partintrodot#1{\def\partintr@dot{#1}}
\long\def\partintrosep#1{\long\def\partintr@sep{#1}}
\def\partnewpagetrue{\def\p@rtnewp@ge{\newpage}}
\def\partnewpagefalse{\long\def\p@rtnewp@ge{\par}}

\partbeforeskip{4ex}
\partstyle{\centering\Large\bf}
\partdot{}
\partafterskip{3ex}
\partintrostyle{\large}
\partintrodot{}
\partintrosep{\par}
\partnewpagefalse

\def\partname{Part}
\def\part{\p@rtnewp@ge\addvspace\p@rtbeforeskip\@afterindentfalse\secdef\@part\@spart}

\def\@part[#1]#2{\ifnum \c@secnumdepth >-1\relax  
        \refstepcounter{part}                     
        \def\@tempa{\addcontentsline{toc}{part}}  %
        \expandafter\@tempa\expandafter{\thepart  
          \hspace{1em}#1}\else                    
        \addcontentsline{toc}{part}{#1}\fi        
   {\p@rtstyl@                       
    \ifnum \c@secnumdepth >-1\relax        
      {\partintr@styl@\partname\ \thepart  
       \partintr@dot}\partintr@sep\nobreak 
    \fi                                    
    #2\partd@t\markboth{}{}\par}
    \nobreak                       
    \vskip\p@rtafterskip           
   \@afterheading                  
    }                              

\def\@spart#1{{\p@rtcentering\p@rtstyl@                      
    #1\partd@t\par}                 
    \nobreak                        
    \vskip\p@rtafterskip            
    \@afterheading                  
  }                                 

\newif\ifsection@ftind
\newif\ifsection@ftpar

\def\sectionbeforeskip#1{\def\s@ctbeforeskip{#1}}
\def\sectionstyle#1{\def\s@ctstyl@{#1}}
\def\sectiondot#1{\def\sectiond@t{#1}}
\def\sectionafterskip#1{\def\s@ctafterskip{#1}}
\def\sectionintrostyle#1{\def\sectionintr@styl@{#1}}
\def\sectionintro#1{\def\sectionintr@{#1}}
\def\sectionintrodot#1{\def\sectionintr@dot{#1}}
\def\sectionintrosep#1{\def\sectionintr@sep{#1}}
\def\sectionindenttrue{\def\s@ctind{\parindent}}
\def\sectionindentfalse{\def\s@ctind{\z@}}
\def\sectionafterindenttrue{\section@ftindtrue}
\def\sectionafterindentfalse{\section@ftindfalse}
\def\sectionafternewlinetrue{\section@ftpartrue}
\def\sectionafternewlinefalse{\section@ftparfalse}

\newif\ifsubsection@ftind
\newif\ifsubsection@ftpar

\def\subsectionbeforeskip#1{\def\ss@ctbeforeskip{#1}}
\def\subsectionstyle#1{\def\ss@ctstyl@{#1}}
\def\subsectiondot#1{\def\subsectiond@t{#1}}
\def\subsectionafterskip#1{\def\ss@ctafterskip{#1}}
\def\subsectionintrostyle#1{\def\subsectionintr@styl@{#1}}
\def\subsectionintro#1{\def\subsectionintr@{#1}}
\def\subsectionintrodot#1{\def\subsectionintr@dot{#1}}
\def\subsectionintrosep#1{\def\subsectionintr@sep{#1}}
\def\subsectionindenttrue{\def\ss@ctind{\parindent}}
\def\subsectionindentfalse{\def\ss@ctind{\z@}}
\def\subsectionafterindenttrue{\subsection@ftindtrue}
\def\subsectionafterindentfalse{\subsection@ftindfalse}
\def\subsectionafternewlinetrue{\subsection@ftpartrue}
\def\subsectionafternewlinefalse{\subsection@ftparfalse}

\newif\ifsubsubsection@ftind
\newif\ifsubsubsection@ftpar

\def\subsubsectionbeforeskip#1{\def\sss@ctbeforeskip{#1}}
\def\subsubsectionstyle#1{\def\sss@ctstyl@{#1}}
\def\subsubsectiondot#1{\def\subsubsectiond@t{#1}}
\def\subsubsectionafterskip#1{\def\sss@ctafterskip{#1}}
\def\subsubsectionintrostyle#1{\def\subsubsectionintr@styl@{#1}}
\def\subsubsectionintro#1{\def\subsubsectionintr@{#1}}
\def\subsubsectionintrodot#1{\def\subsubsectionintr@dot{#1}}
\def\subsubsectionintrosep#1{\def\subsubsectionintr@sep{#1}}
\def\subsubsectionindenttrue{\def\sss@ctind{\parindent}}
\def\subsubsectionindentfalse{\def\sss@ctind{\z@}}
\def\subsubsectionafterindenttrue{\subsubsection@ftindtrue}
\def\subsubsectionafterindentfalse{\subsubsection@ftindfalse}
\def\subsubsectionafternewlinetrue{\subsubsection@ftpartrue}
\def\subsubsectionafternewlinefalse{\subsubsection@ftparfalse}

\newif\ifparagraph@ftind
\newif\ifparagraph@ftpar

\def\paragraphbeforeskip#1{\def\p@rbeforeskip{#1}}
\def\paragraphstyle#1{\def\p@rstyl@{#1}}
\def\paragraphdot#1{\def\paragraphd@t{#1}}
\def\paragraphafterskip#1{\def\p@rafterskip{#1}}
\def\paragraphintrostyle#1{\def\paragraphintr@styl@{#1}}
\def\paragraphintro#1{\def\paragraphintr@{#1}}
\def\paragraphintrodot#1{\def\paragraphintr@dot{#1}}
\def\paragraphintrosep#1{\def\paragraphintr@sep{#1}}
\def\paragraphindenttrue{\def\p@rind{\parindent}}
\def\paragraphindentfalse{\def\p@rind{\z@}}
\def\paragraphafterindenttrue{\paragraph@ftindtrue}
\def\paragraphafterindentfalse{\paragraph@ftindfalse}
\def\paragraphafternewlinetrue{\paragraph@ftpartrue}
\def\paragraphafternewlinefalse{\paragraph@ftparfalse}

\newif\ifsubparagraph@ftind
\newif\ifsubparagraph@ftpar

\def\subparagraphbeforeskip#1{\def\sp@rbeforeskip{#1}}
\def\subparagraphstyle#1{\def\sp@rstyl@{#1}}
\def\subparagraphdot#1{\def\subparagraphd@t{#1}}
\def\subparagraphafterskip#1{\def\sp@rafterskip{#1}}
\def\subparagraphintrostyle#1{\def\subparagraphintr@styl@{#1}}
\def\subparagraphintro#1{\def\subparagraphintr@{#1}}
\def\subparagraphintrodot#1{\def\subparagraphintr@dot{#1}}
\def\subparagraphintrosep#1{\def\subparagraphintr@sep{#1}}
\def\subparagraphindenttrue{\def\sp@rind{\parindent}}
\def\subparagraphindentfalse{\def\sp@rind{\z@}}
\def\subparagraphafterindenttrue{\subparagraph@ftindtrue}
\def\subparagraphafterindentfalse{\subparagraph@ftindfalse}
\def\subparagraphafternewlinetrue{\subparagraph@ftpartrue}
\def\subparagraphafternewlinefalse{\subparagraph@ftparfalse}

\sectionbeforeskip{\bigskipamount}
\sectionstyle{\large\bf}
\sectiondot{}
\sectionafterskip{.5\bigskipamount}
\sectionintrostyle{}
\sectionintro{}
\sectionintrodot{.}
\sectionintrosep{1.25ex}
\sectionindentfalse
\sectionafterindenttrue
\sectionafternewlinetrue

\subsectionbeforeskip{.8\bigskipamount}
\subsectionstyle{\normalsize\bf}
\subsectiondot{}
\subsectionafterskip{.4\bigskipamount}
\subsectionintrostyle{}
\subsectionintro{}
\subsectionintrodot{.}
\subsectionintrosep{1.25ex}
\subsectionindentfalse
\subsectionafterindenttrue
\subsectionafternewlinetrue

\subsubsectionbeforeskip{.6\bigskipamount}
\subsubsectionstyle{\normalsize\bf}
\subsubsectiondot{}
\subsubsectionafterskip{.3\bigskipamount}
\subsubsectionintrostyle{}
\subsubsectionintro{}
\subsubsectionintrodot{.}
\subsubsectionintrosep{1.25ex}
\subsubsectionindentfalse
\subsubsectionafterindenttrue
\subsubsectionafternewlinetrue

\paragraphbeforeskip{.5\bigskipamount}
\paragraphstyle{\normalsize\bf}
\paragraphdot{.}
\paragraphafterskip{1.25ex}
\paragraphintrostyle{}
\paragraphintro{}
\paragraphintrodot{.}
\paragraphintrosep{1.25ex}
\paragraphindentfalse
\paragraphafterindenttrue
\paragraphafternewlinefalse

\subparagraphbeforeskip{.5\bigskipamount}
\subparagraphstyle{\normalsize\bf}
\subparagraphdot{.}
\subparagraphafterskip{1.25ex}
\subparagraphintrostyle{}
\subparagraphintro{}
\subparagraphintrodot{.}
\subparagraphintrosep{1.25ex}
\subparagraphindenttrue
\subparagraphafterindenttrue
\subparagraphafternewlinefalse

\let\@partoken\par
\long\def\@@gobble#1{}
\def\ignorepar{\@ifnextchar\@partoken{\expandafter\ignorepar\@@gobble}{\ignorespaces}}

\def\@startsection#1#2#3#4#5#6{
   \@tempskipa #4\relax
   \csname if#1@ftind\endcsname\@afterindenttrue\else\@afterindentfalse\fi
   \advance\@tempskipa by\presection
   \if@nobreak \everypar{}\else
     \addpenalty{\@secpenalty}\addvspace{\@tempskipa}%
     \allowbreak\vskip -\presection \fi \@ifstar
     {\@ssect{#1}{#2}{#3}{#4}{#5}{#6}}{\@dblarg{\@sect{#1}{#2}{#3}{#4}{#5}{#6}}}}

\def\@sect#1#2#3#4#5#6[#7]#8{\def\object@type{#1}%
   \ifnum #2>\c@secnumdepth\def\@svsec{}\def\@tempb{}%
      \else\refstepcounter{#1}\def\@svsec{{\csname #1intr@styl@\endcsname%
        {\csname #1intr@\endcsname}\csname the#1\endcsname%
        \csname #1intr@dot\endcsname\kern\csname #1intr@sep\endcsname}}%
        \edef\@tempb{\noexpand\numberline{\csname the#1\endcsname}}\fi%
   \def\@tempa{\addcontentsline{toc}{#1}}%
   \csname if#1@ftpar\endcsname%
      \begingroup #6\relax%
        \@hangfrom{\hskip #3\relax\@svsec}{\interlinepenalty \@M{#8}%
        \csname #1d@t\endcsname\par}%
      \endgroup%
      \csname #1mark\endcsname{#7}%
      \expandafter\@tempa\expandafter{\@tempb #7}%
      \ifautolabel\label*{#8}\fi%
   \else%
      \def\@svsechd{#6\hskip #3\relax%
         \@svsec{#8}\csname #1mark\endcsname{#7}%
         \expandafter\@tempa\expandafter{\@tempb #7}%
         \ifautolabel\label*{#8}\fi}\fi%
   \@xsect{#1}{#5}\ignorepar}

\def\@ssect#1#2#3#4#5#6#7{%
   \ifnum #2>\c@secnumdepth\def\@tempb{}\else \def\@tempb{\numberline{}}\fi%
     \def\@tempa{\addcontentsline{toc}{s#1}}%
     \csname if#1@ftpar\endcsname
        \begingroup #6\relax
           \@hangfrom{\hskip #3}{\interlinepenalty \@M{#7}%
           \csname #1d@t\endcsname\par}%
        \endgroup
        \csname s#1mark\endcsname{#7}%
        \ifstarredcontents\expandafter\@tempa\expandafter{\@tempb #7}\fi%
        \ifautolabel\label*{#7}\fi%
     \else%
        \def\@svsechd{#6\hskip #3\relax{#7}\csname s#1mark\endcsname{#7}%
        \ifautolabel\label*{#7}\fi}\fi
   \@xsect{#1}{#5}\ignorepar}

\def\@xsect#1#2{
   \csname if#1@ftpar\endcsname 
       \par \nobreak \vskip #2\relax \@afterheading
    \else \global\@nobreakfalse \global\@noskipsectrue
       \everypar{\if@noskipsec \global\@noskipsecfalse
                   \clubpenalty\@M \hskip -\parindent
                   \begingroup \@svsechd \endgroup \unskip
                   \hskip #2\relax  
                  \else \clubpenalty \@clubpenalty
                    \everypar{}\fi}\fi\ignorespaces}

\def\section{\@startsection{section}{1}{\s@ctind}
  {\s@ctbeforeskip}{\s@ctafterskip}{\s@ctstyl@}}
\def\subsection{\@startsection{subsection}{2}{\ss@ctind}
  {\ss@ctbeforeskip}{\ss@ctafterskip}{\ss@ctstyl@}}
\def\subsubsection{\@startsection{subsubsection}{3}{\sss@ctind}
  {\sss@ctbeforeskip}{\sss@ctafterskip}{\sss@ctstyl@}}
\def\paragraph{\@startsection{paragraph}{4}{\p@rind}
  {\p@rbeforeskip}{\p@rafterskip}{\p@rstyl@}}
\def\subparagraph{\@startsection{subparagraph}{4}{\sp@rind}
  {\sp@rbeforeskip}{\sp@rafterskip}{\sp@rstyl@}}

\def\statementabove#1{\def\th@bove{#1}}
\def\statementstyle#1{\def\thstyl@{#1}}
\def\statementbelow#1{\def\thb@low{#1}}
\def\statementindentfalse{\let\thind@nt\relax}
\def\statementindenttrue{\let\thind@nt\indent}

\def\statementintrostyle#1{\def\thintr@style{#1}}
\def\statementintrodot#1{\def\thintr@dot{#1}}
\def\statementintrosep#1{\def\thintr@sep{#1}}
\def\statementintrobrackets#1#2{\def\thintr@left{#1}\def\thintr@right{#2}}

\statementabove{\medskip}
\statementstyle{\sl}
\statementbelow{\medskip}
\statementindenttrue

\statementintrostyle{\normalshape\bf}
\statementintrodot{.}
\statementintrosep{\kern1.25ex}
\statementintrobrackets{(}{)}

\def\@thskip{\dimen0\lastskip\vskip-\dimen0%
  \th@bove\dimen1\lastskip\vskip-\dimen1%
  \ifdim\dimen0>\dimen1\else\dimen0\dimen1\fi\vskip\dimen0}

\long\def\@@newtheorem#1#2#3{%
  \newenvironment{#3}%
    {\def\object@type{#3}\par\@thskip%
     \@ifnextchar[{\@enva{#3}{\thstyl@#1{#2}}}{\@envb{#3}{\thstyl@#1{#2}}}}%
    {\end{#3@}}%
  \@ifnextchar[{\@nothm{#3}}{\@nnthm{#3}}}

\def\@nothm#1[#2]#3{%
  \@ifundefined{c@#2}{\@latexerr{No theorem environment `#2' defined}\@eha}%
  {\expandafter\@ifdefinable\csname #1@\endcsname
  {\global\@namedef{the#1}{\@nameuse{the#2}}%
   \global\@namedef{c@#1}{\@nameuse{c@#2}}
   \global\@namedef{p@#1}{\@nameuse{p@#2}}
   \global\@namedef{#1@}{\@nnnthm{#2}{#3}}%
   \global\@namedef{end#1@}{\@endtheorem}}}}

\def\@nnnthm#1#2{\refstepcounter
    {#1}\@ifnextchar[{\@ynnnthm{#1}{#2}}{\@xnnnthm{#1}{#2}}}
 
\def\@xnnnthm#1#2{\@begintheorem{#2}{\csname the#1\endcsname}\ignorespaces}
\def\@ynnnthm#1#2[#3]{\@opargbegintheorem{#2}{\csname
       the#1@\endcsname}{#3}\ignorespaces}

\def\renewtheorem{\@ifnextchar[{\@renewtheorem}{\@renewtheorem[{}{}]}}

\long\def\@renewtheorem[#1]{\@@renewtheorem#1}

\long\def\@@renewtheorem#1#2#3{%
  \expandafter\let\csname#3@\endcsname\undefined
  \renewenvironment{#3}%
    {\def\object@type{#3}\par\@thskip%
     \@ifnextchar[{\@enva{#3}{\thstyl@#1{#2}}}{\@envb{#3}{\thstyl@#1{#2}}}}%
    {\end{#3@}}%
  \@ifnextchar[{\@nothm{#3}}{\@nnthm{#3}}}

\def\@begintheorem#1#2{\@opargbegintheorem{#1}{#2}{}}

\def\@opargbegintheorem#1#2#3{%
        \def\@tempx{#1}%
        \expandafter\let\expandafter\@tempy#2
        \def\@tempz{#3}%
        \mytrivlist\item[\thind@nt\hskip\labelsep%
        {\thintr@style%
          #1\if\@tempx\@empty\else\if\@tempy\relax\else\kern1ex\fi\fi#2%
          \ifx\@tempz\@empty%
            \if\@tempx\@empty\if\@tempy\relax%
            \else\thintr@dot\thintr@sep\fi\else\thintr@dot\thintr@sep\fi%
            \else%
            \if\@tempx\@empty\if\@tempy\relax\else\kern1ex\fi\else\kern1ex\fi%
           \thintr@left{#3}\thintr@right\thintr@dot\thintr@sep\fi}%
            \hskip-\labelsep]%
        \ifautolabel\label*{#3}\fi}

\def\@endtheorem{\strut\endtrivlist\thb@low}

\def\proofabove#1{\def\pf@bove{#1}}
\def\proofstyle#1{\def\pfstyl@{#1}}
\def\proofbelow#1{\def\pfb@low{#1}}
\def\proofindentfalse{\let\pfind@nt\relax}
\def\proofindenttrue{\let\pfind@nt\indent}

\def\proofintrostyle#1{\def\pfintr@style{#1}}
\def\proofintrodot#1{\def\pfintr@dot{#1}}
\def\proofintrosep#1{\def\pfintr@sep{#1}}
\def\proofintrobrackets#1#2{\def\pfintr@left{#1}\def\pfintr@right{#2}}

\proofabove{\medskip}
\proofstyle{}
\proofbelow{\medskip}
\proofindenttrue

\proofintrostyle{\sl}
\proofintrodot{.}
\proofintrosep{\kern1.25ex}
\proofintrobrackets{of\kern1ex}{}

\def\@pfskip{\dimen0\lastskip\vskip-\dimen0%
  \pf@bove\dimen1\lastskip\vskip-\dimen1%
  \ifdim\dimen0>\dimen1\else\dimen0\dimen1\fi\vskip\dimen0}

\renewenvironment{proof}%
  {\@pfskip\mytrivlist\item[\pfind@nt]\@ifnextchar[{\pro@f}{\pro@f[\prooftag]}}
  {\ifvoid\provedbox\else\hproved\fi\endtrivlist\pfb@low}

\def\pro@f[#1]{\setbox\provedbox\hbox{\provedboxcontents{#1}}\proofintro{#1}}

\def\proofintro#1{\expandafter\def\expandafter\@tempa\expandafter{#1}%
  {\pfintr@style{Proof\ifx\@tempa\empty\else\kern1ex\pfintr@left{#1}%
  \pfintr@right\fi}\pfintr@dot\pfintr@sep}\pfstyl@\ignorespaces}

\def\provedmark#1{\def\prm@rk{#1}}
\def\provedsep#1{\def\prs@p{#1}}

\provedmark{$\square$}
\provedsep{\kern1.25ex}

\def\provedtexttrue{\def\prb@x##1{\fbox{\small##1}}}
\def\provedtextfalse{\def\prb@x##1{\prm@rk}}
\def\provedmarkrighttrue{\let\prhf@l\hfill}
\def\provedmarkrightfalse{\let\prhf@l\relax}

\provedtextfalse
\provedmarkrightfalse

\def\provedboxcontents#1{\expandafter\def\expandafter\@tempa\expandafter{#1}%
  \ifx\@tempa\empty\prm@rk\else\prb@x{#1}\fi}

\def\proved{\ifmmode\eqno{\box\provedbox}\else\hproved\fi}

\def\hproved{\unskip\nobreak\prhf@l\penalty50\prs@p\hbox{}\nobreak\prhf@l
  \box\provedbox{\finalhyphendemerits=0\par}}

\def\captionstyle#1{\def\c@ptstyl@{#1}}
\def\captionintrostyle#1{\def\c@pintr@style{#1}}
\def\captionintrodot#1{\def\c@pintr@dot{#1}}
\def\captionintrosep#1{\def\c@pintr@sep{#1}}

\captionstyle{\small\sf}
\captionintrostyle{\bf}
\captionintrodot{.}
\captionintrosep{\hskip1.25ex}

\long\def\@makecaption#1#2{%
    \vskip\captionskip
    \setbox\@tempboxa\hbox{%
      \ifproofing\@ifundefined{the@label}{}
        {\hbox to 0pt{\vbox to 0pt{\vss\hbox{\tiny\the@label}\bigskip}\hss}}\fi
      \c@ptstyl@{\c@pintr@style #1\c@pintr@dot}\ignorespaces #2}%
    \@captionwidth=\hsize \advance\@captionwidth-2\@captionmargin
    \ifdim \wd\@tempboxa >\@captionwidth {%
        \rightskip=\@captionmargin\leftskip=\@captionmargin
        \unhbox\@tempboxa\par}%
      \else
        \hbox to\hsize{\hfil\box\@tempboxa\hfil}%
    \fi}

\def\end@Float#1{%
  \expandafter\caption\expandafter[\the@title]{%
   {\c@pintr@style%
   \ifx\the@caption\empty\ifx\the@title\empty
   \else\c@pintr@sep\fi\else\c@pintr@sep\fi
    \the@title\ifx\the@caption\empty%
     \expandafter\label\expandafter*\expandafter{\the@label}%
    \else\ifx\the@title\empty%
     \expandafter\label\expandafter*\expandafter{\the@label}%
    \else\c@pintr@dot\c@pintr@sep%
     \expandafter\label\expandafter*\expandafter{\the@label}\fi\fi}%
   \ignorespaces\the@caption}%
  \end{#1}}

\@definecounter{bibenumi}

\def\thebibliography#1{%
 \section*{\refname}\vskip-\lastskip%
 \list{[\arabic{bibenumi}]}{\topsep0pt\settowidth\labelwidth{[#1]}%
 \leftmargin\labelwidth\advance\leftmargin\labelsep\usecounter{bibenumi}}%
 \def\newblock{\hskip .11em plus .33em minus .07em}%
 \sloppy\clubpenalty4000\widowpenalty4000\sfcode`\.=1000\relax}

\makeatother

\proofingfalse
\autolabelfalse

\showfiguresfalse

\newtheorem{stat}{\statname}  \unnumbered{stat}

\newtheorem{nstat}{\nstatname}[section]

\newtheorem{lemma}[nstat]{Lemma}
\newtheorem{proposition}[nstat]{Proposition}
\newtheorem{theorem}[nstat]{Theorem}
\newtheorem{corollary}[nstat]{Corollary}

\newtheorem[{\ns}{}]{remark}[nstat]{Remark}

\papersize{215.9truemm}{279.4truemm}
\margins{3.295cm}{2.12cm}

\headline{\hfill}
\footline{\small\hfill--\kern1ex\thepage\kern1ex--\hfill}

\abovedisplayskip\smallskipamount
\belowdisplayskip\smallskipamount
\abovedisplayshortskip\smallskipamount
\belowdisplayshortskip\smallskipamount

\showfigurestrue

\sectionbeforeskip{1.5\bigskipamount}
\sectionstyle{\centering\normalsize\bf}
\sectionafterskip{\bigskipamount}
\sectionafterindenttrue

\subsectionbeforeskip{\bigskipamount}
\subsectionstyle{\normalsize\bf}
\subsectionafterskip{.5\bigskipamount}
\subsectionafterindenttrue

\paragraphbeforeskip{\bigskipamount}
\paragraphstyle{\centering\normalsize\sl}
\paragraphafterskip{.5\bigskipamount}
\paragraphafternewlinetrue
\paragraphafterindenttrue
\paragraphdot{}

\statementintrostyle{\sc}
\renewtheorem[{\ns}{}]{definition}[nstat]{Definition}
\renewtheorem[{\ns}{}]{remark}[nstat]{Remark}
\newtheorem[{\ns}{}]{block}[nstat]{}

\captionstyle{\small}
\captionintrostyle{\sc}

\def\ms{\mediumseries}
\def\ns{\normalshape}

\def\({\mbox\bgroup{\ns(}\sl\aux}\gdef\nextsp{}
\def\aux#1{\def\tempx{#1}\let\next\aux%
 \if\tempx)\let\next\egroup{\/\ns)}%
 \else\if\tempx f{\nextsp f\gdef\nextsp{\kern0.2ex}}%
 \else\if\tempx i{\nextsp i\gdef\bextsp{\kern0.1ex}}%
 \else\if\tempx I{\nextsp I\gdef\nextsp{\kern0.1ex}}%
 \else\if\tempx'{\/$'$\gdef\nextsp{}}%
 \else\if\tempx-{\/\ns-\gdef\nextsp{}}%
 \else\nextsp\tempx\gdef\nextsp{}%
 \fi\fi\fi\fi\fi\fi\next}

\def\varemptyset{%
 {\text{\raise.21ex\hbox{$\not$}}\mkern.15mu\mathrm{O}\mkern.15mu}}
\def\widebar#1{%
 \text{\setbox0\hbox{$#1$}\dimen0\ht0\dimen0.25\dimen0
 \hbox to 0pt{$\kern\dimen0\overline{\kern-\dimen0\phantom{#1}}$\hss}}#1}

\def\_{{\hbox to 1.2ex{\hss\vrule width1ex height0pt depth.4pt\hss}}}

\let\emptyset\varemptyset

\let\bar\widebar
\let\hat\widehat
\let\tilde\widetilde

\let\epsilon\texepsilon
\let\theta\textheta

\newcommand{\alg}{{{\cal A}\mkern0.5mu l\mkern-1.3mu g}}
\newcommand{\cob}{{{\cal C}\mkern-0.8mu o\mkern0.5mu b}}

\newcommand{\one}{{\boldmath 1}}

\newcommand{\Int}{\mathop{\mathrm{Int}}\nolimits} 
\newcommand{\Bd}{\mathop{\mathrm{Bd}}\nolimits}

\newcommand{\Obj}{\mathop{\mathrm{Obj}}\nolimits}
\newcommand{\Mor}{\mathop{\mathrm{Mor}}\nolimits}

\newcommand{\id}{\mathop{\mathrm{id}}\nolimits}
\newcommand{\rot}{\mathop{\mathrm{rot}}\nolimits}
\newcommand{\rev}{\mathop{\mathrm{rev}}\nolimits}

\newcommand{\diam}{\mathrel{\diamond}}

\newcommand{\C}{{\cal C}}
\newcommand{\G}{{\cal G}}
\renewcommand{\H}{{\cal H}}
\newcommand{\K}{{\cal K}}

\renewcommand{\S}{{\cal S}}
\newcommand{\T}{{\cal T}}

\newcommand{\rs}{rs}

\makeatletter
\def\up{\@ifnextchar[{\@up}{{\uparrow}}}
\def\@up[#1]{{\uparrow}\text{\raise .6ex\hbox{$_#1$}}}
\def\down{\@ifnextchar[{\@down}{{\downarrow}}}
\def\@down[#1]{{\downarrow}\text{\raise .6ex\hbox{$_#1$}}}
\makeatother

\newcommand{\bs}{\text{\raise.4ex\hbox{\bf$\scriptscriptstyle\backslash$}}}

\newcommand{\mapright}[1]{\smash{\mathop{\longrightarrow}\limits^{#1}}}

\font\ftt cmtt10 at 11pt
\font\fsc cmcsc10 at 12pt
\font\fsl cmsl12 at 12pt
\font\bfs cmbxsl10 at 12pt

\def\mypagebreak{{\parfillskip0pt\relax\par}\break\noindent}

\begin{document}

\title{\large\bf
A UNIVERSAL INVARIANT OF FOUR-DIMENSIONAL 2-HANDLEBODIES AND THREE-MANIFOLDS
\label{Version 0.9(50) / \today}}
\author{\fsc I. Bobtcheva\\
\fsl Dipartimento di Scienze Matematiche\\[-3pt]
\fsl Universit\`a Politecnica delle Marche -- Italia\\
\ftt bobtchev@dipmat.univpm.it
\and
\fsc R. Piergallini\\
\fsl Dipartimento di Matematica e Informatica\\[-3pt]
\fsl Universit\`a di Camerino -- Italia\\
\ftt riccardo.piergallini@unicam.it}
\date{}

\maketitle

\begin{abstract}
\baselineskip13.5pt
\smallskip

\noindent
In \cite{BP} it is shown that up to certain set of local moves, connected  simple
coverings of $B^4$ branched over ribbon surfaces, bijectively represent connected
orientable 4-dimensional 2-handlebodies up to 2-deformations (handle slides and
creations/cancellations of handles of index $\leq 2$). We factor this bijective
correspondence through a map onto the closed morphisms in a universal braided category
$\H^r$ freely generated by a Hopf algebra object $H$. In this way we obtain a complete
algebraic description of 4-dimensional 2-handlebodies. This result is then used to
obtain an analogous description of the boundaries of such handlebodies, i.e.
3-dimensional manifolds, which resolves for closed manifolds the problem posed by
Kerler in \cite{Ke02} (cf. \cite[Problem 8-16 (1)]{O02}).



\end{abstract}

\section{Introduction}

During the last twenty years the developments in the quantum theory have build a
bridge between two distinct areas of mathematics: the topology of low dimensions (2,
3 and 4) and the theory of Hopf algebras. An important result in this direction is
the one of Shum \cite{Sh94}, which states that the category of framed tangles is
equivalent to the universal monoidal tortile category generated by a single object.
Since the category of representations of a wide class of Hopf algebras, the so
called ribbon Hopf algebras, is a tortile monoidal category (cf. \cite{RT91, L93}
etc.), the result of Shum implies that the subcategory generated by any
representation of such a Hopf algebra produces an invariant of framed tangles.
Framed tangles and links are interesting and rich object of study by themselves, but
through surgery they are also the main tool of describing 3- and 4-dimensional
manifolds. Actually, 4-dimensional 2-handlebodies modulo 2-deformations (handle
slidings and addition/deletion of canceling handles of indices $\leq 2$),
bijectively correspond to the equivalence classes of Kirby link diagrams modulo
Kirby calculus moves (cf. \cite{GS99} and Section \ref{kirby/sec} below), while
3-dimensional manifolds bijectively correspond to the equivalence classes of framed
links modulo the Fenn-Rourke move \cite{FR79}. Reshetikhin and Turaev \cite{RT91}
used this last fact to construct invariants of 3-manifolds, by showing that some
finite dimensional Hopf algebras have a finite subset of representations $S$, such
that a proper linear combination of the framed link invariants corresponding to the
elements of $S$, is also invariant under the Fenn-Rourke move. On the other hand,
Hennings \cite{He96} defined invariants of 3-manifolds starting directly with a
unimodular ribbon Hopf algebra and a trace function on it. Eventually Lyubashenko
and Kerler \cite{KL01} constructed 3-manifold invariants (2+1 topological quantum
field theory) out of a Hopf algebra in a linear abelian braided monoidal category
with certain coends, and showed that Reshetikhin-Turaev's and Hennings' invariants
are particular cases of such more general approach (cf. \cite{Ke97}).

All these results lead Kerler to define in \cite{Ke02}  a surjective functor from a
braided monoidal category $\alg$, freely generated by a Hopf algebra object, to the
category of 3-dimensional relative cobordisms $\cob^{2+1}$. Moreover he asked what
additional relations should be introduced in $\alg$, so that the functor above
defines an equivalence between the quotient of $\alg$ by these relations and
$\cob^{2+1}$ (cf. \cite[Problem 8-16 (1)]{O02}). Finding these relations would
complete the characterization of $\cob^{2+1}$ in purely algebraic terms.

\medskip

In the present work we propose two new relations such that the corresponding functor 
on the quotient induces a bijective map between the closed morphisms of the algebraic
category and the set of closed connected orientable 3-manifolds, obtaining in this way
an algebraic characterization of such manifolds.

Our approach is independent of Kerler's work. Actually, we obtain the result for
3-manifolds as a consequence of the algebraic description of 4-dimensional
2-handlebodies modulo 2-deformations. 

In particular, we construct a bijective correspondence between the set of equivalence
classes of connected 4-dimensional 2-handlebodies modulo 2-deformations and the set of
closed morphisms of a universal category $\H^r$. The objects of $\H^r$ are the elements
of the free $(\diam,\one)$-algebra on a single object $H$ and $\H^r$ is universal with
respect to the properties that:
\begin{itemize}\itemsep\smallskipamount
\item[{\sl a}\/)] \vskip-\lastskip\smallskip
$H$ is a unimodular braided Hopf algebra;
\item[{\sl b}\/)] there is a special morphism $v: H \to H$ which makes $H$ into a
ribbon algebra as defined in \ref{ribbon} (over the trivial groupoid: one object and
one relation).
\end{itemize}\vskip-\lastskip\medskip

We point out the following differences between $\H^r$ and  the algebra $\alg$ proposed
by Kerler in \cite{Ke02}.

Firstly, since we are interested in 4-dimensional 2-handlebodies, we do not require
that the Hopf copairing $\one \to H \diam H$ is non-degenerate. Such condition is
equivalent to the requirement that the algebra integral and cointegral are dual with
respect to this copairing (cf. (34) in \cite{Ke02}) and is necessary only when one
restricts to 3-manifolds (cf. \ref{defn-boundH/par}).

Secondly, we introduce the ribbon morphism as a morphism $v: H \to H$ instead of 
$v_K: \one \to H$ as it is done in \cite{Ke02}. These two ways are
equivalent, being $v$ ob\-tainable by taking the product of $v_K$ with the
identity morphism on $H$ and then composing with the multiplication. However,
viewing the ribbon morphism as a morphism from the algebra to itself allows to talk
about its propagation through the diagrams, which is a main tool in the proofs. 

Besides these cosmetic differences, the new relations with respect to the ones
given by Kerler are \(r10) and \(r11) in \ref{ribbon}. The first one describes the 
propagation of the ribbon morphism through the copairing, and the second one relates
the copairing and the braiding morphism. 

\medskip

Before stating the main results, we would like to sketch our approach to the problem.
In \cite{Ke02} Kerler constructs a functor from the algebraic category to the category
of Kirby tangles (similar map is presented also in \cite{H00}). The main difficulty in
going the other way consists in the fact that the elementary morphisms in the universal
braided Hopf algebra are local (involve only objects which are ``close'' to each
other), while handle slides in general are not. So, to be able to go naturally from
the topological category to the algebra, one would like to have an alternative
description of the topological objects in which the moves of equivalence are more
similar to the algebraic ones. For 4-dimensional 2-handlebodies such description is
given in \cite{BP} in terms of simple coverings of $B^4$ branched over ribbon
surfaces.

In particular, starting from a surgery description of a connected 4-dimensional
2-handlebody $M$ as an ordinary Kirby diagram $K$ (with a single 0-handle), in order to
obtain a description of $M$ as a simple covering of $B^4$ of degree $n \geq 4$
branched over a ribbon surface, one first stabilizes with $n - 1$ pairs of canceling
0/1-handles. This transforms $K$ into a generalized Kirby diagram $\tilde K$, which
represents the attaching maps of the 1- and 2-handles on the boundary of the $n$
0-handles. The generalized Kirby diagrams differ from the ordinary ones by the fact
that in them each side of the spanning disks of the dotted components and each
component of the framed ones, after cutting them at the intersections with those
disks, carries a label from $\{1, 2, \dots ,n\}$, indicating the 0-handle where it
lives. Examples of generalized Kirby tangles are shown in the right column of Figure
\ref{comp-funct/fig}, while the precise definition can be found in Section
\ref{kirby/sec}. By using the extra 0-handles, $\tilde K$ can be symmetrized and the
resulting diagram is interpreted as a labeled ribbon surface $F_K$ in $B^4$, that is
the branching surface of an $n$-fold covering. Moreover, if we change $M$ by
2-deformation, i.e. change $K$ by Kirby calculus moves, $F_K$ changes by a sequence of
local moves (1-isotopy and ribbon moves). The correspondence $K \mapsto F_K$ between
labeled ribbon surfaces modulo 1-isotopy and ribbon moves, and generalized Kirby
diagrams modulo Kirby calculus moves, is invertible (for the definition of the inverse
$F \mapsto K_F$, see Section 2 in \cite{BP} and Section \ref{surf-kirby/sec} below). 

\begin{Figure}[p]{comp-funct/fig}{}{}
\centerline{\fig{}{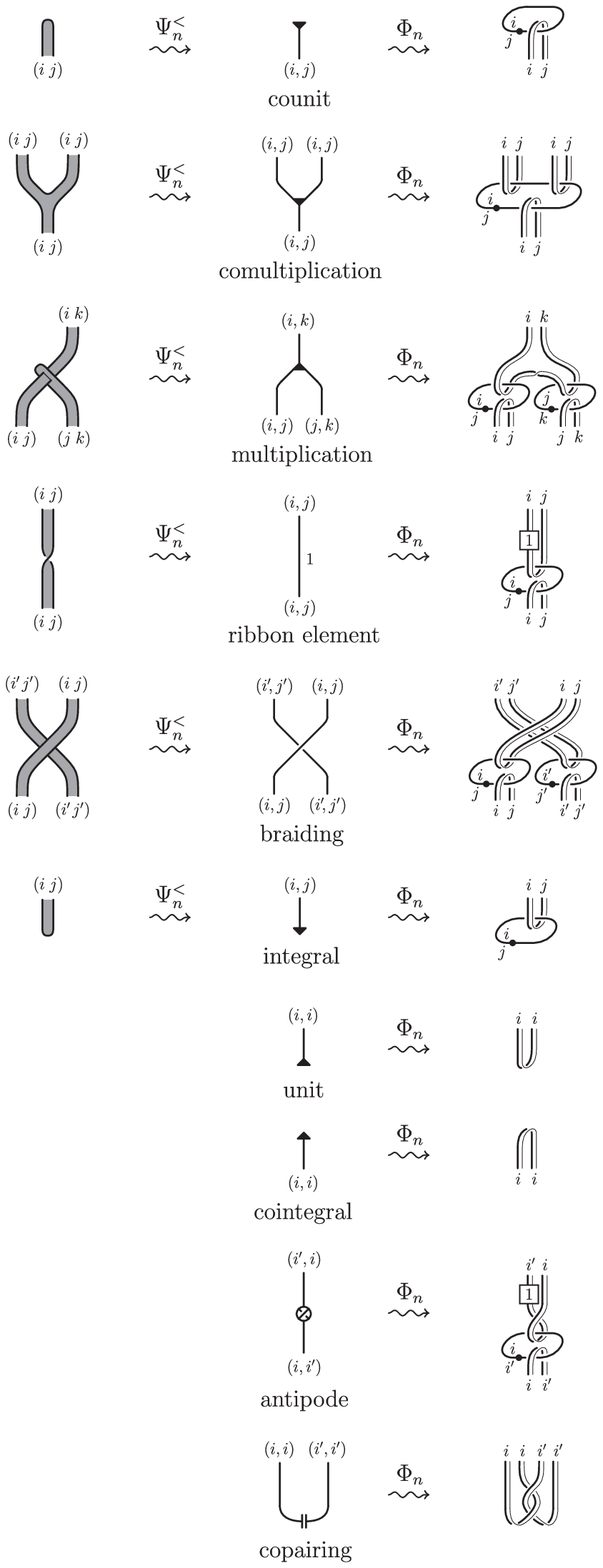}}
\end{Figure}

Therefore, in order to obtain the desired algebraic description of 4-dimensional
2-handlebodies it is enough to factor the correspondence $F \mapsto K_F$ through a map
onto the morphisms of the algebraic category and since the equivalence moves on ribbon
surfaces are local, this should be relatively easy to do. Yet, there is a price to
pay: we will need a family of ribbon Hopf algebras $\H^r_n$, which generalize $\H^r =
\H^r_1$ and describe 4-dimensional 2-handlebodies with $n$ 0-handles. Then we will
also need a ``reduction'' map, which is the algebraic analog of canceling pairs of
0/1-handles. As one may expect, the construction of this map is the main technical
difficulty to be overcome.

\medskip

We now proceed with some more details and list the main results. Given a
groupoid $\G$, we define a ribbon Hopf $\G$-algebra in a braided monoidal category
(generalizing the concept of group Hopf algebra described in \cite{V01}) and
introduce the universal category $\H^r(\G)$ freely generated by a ribbon Hopf 
$\G$-algebra. Our basic example of a braided monoidal category with a ribbon Hopf
$\G$-algebra in it, is the category of admissible generalized Kirby tangles $\K_n$
with $n$ labels (equiv. $n$ 0-handles) modulo slides of handles of index $\leq
2$ and creations/cancellations of pairs of 1/2-handles. In this case, $\G = \G_n$ is
the groupoid with objects $1,2,\dots,n$ and a single morphism, denoted by $(i,j)$
between any $i$ and $j$.
In particular, if $\H^r_n = \H^r(\G_n)$ we have the following theorem, which extends
results of Habiro \cite{H00} and Kerler \cite{Ke02}.

\begin{theorem}\label{alg-kirby/theo}
 There is a braided monoidal functor $\Phi_n: \H^r_n \to \K_n$.
\end{theorem}

The elementary morphisms in $\H^r_n$ and their images under $\Phi_n$ are presented
on the right in Figure \ref{comp-funct/fig} (one needs to add the inverses of the
braiding, the antipode and the ribbon morphism, but this is straightforward). 

\medskip

Another class of braided monoidal categories that we introduce are the categories
$\S_n$ of {\sl ribbon surface tangles} (\rs-tangles) labeled by transpositions in
the symmetric group $\Sigma_n$.

\begin{Figure}[b]{comp-funct2/fig}{}
 {($(i\;j)$ and $(k\;l)$ disjoint)}
\centerline{\fig{}{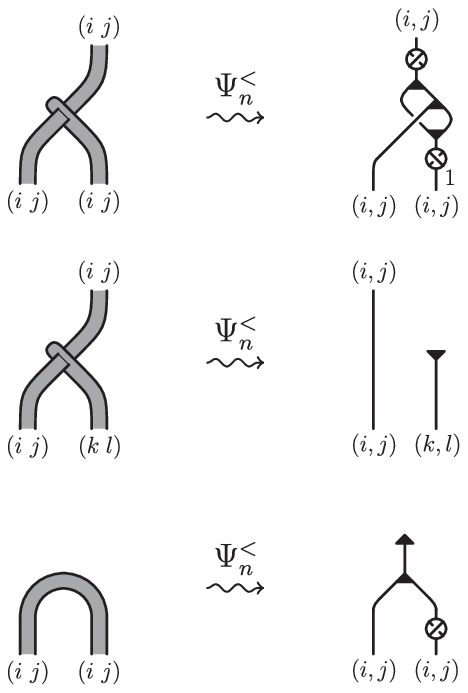}}\vskip-3pt
\end{Figure}

A smooth compact surface $F \subset B^4$ with $\Bd F \subset S^3 = R^3\cup \infty$ is
called a {\sl ribbon surface} if the Euclidean norm restricts to a Morse function on
$F$ with no local maxima in $\Int F$. Then an \rs-tangle should be thought as a slice
of a ribbon surface. We present ribbon surfaces and \rs-tangles through their
projections in $R^3$. In particular, the objects in $\S_n$ are intervals labeled by
transpositions in $\Sigma_n$, while any morphism is a composition of juxtapositions of
the elementary morphisms  presented on the left in Figures \ref{comp-funct/fig} and
\ref{comp-funct2/fig} (the transpositions $(i\;j)$ and $(k\;l)$ in Figure
\ref{comp-funct2/fig} are assumed to be disjoint, i.e. $\{i,j\}\cap \{k,l\} =
\emptyset$). The relations in the category are given by 1-isotopy moves, which are
special type of isotopy moves, and the two ribbon moves presented in Figure
\ref{ribbon-m/fig}, where $(i\;j)$ and $(k\;l)$ are disjoint. It is not known if
1-isotopy coincides with isotopy of ribbon surfaces, for some discussions of this
problem we refer the reader to \cite{BP}.

\begin{Figure}[htb]{ribbon-m/fig}{}
 {Ribbon relations in $\S_n$ ($(i\;j)$ and $(k\;l)$ disjoint).}
\centerline{\fig{}{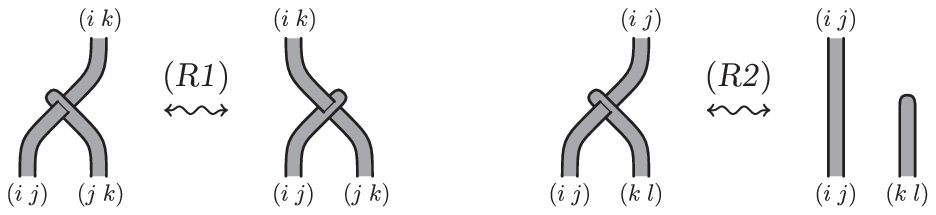}}\vskip-3pt
\end{Figure}

Observe that each morphism in the monoidal categories we are considering is a
composition of products (juxtapositions) of elementary ones and any elementary morphism
involves some set of labels, which are morphisms in a given connected groupoid $\G$
($\G=\G_n$ for $\K_n$ and $\H^r_n$, while $\G=\Sigma_n$ for $\S_n$). Then we call a
morphism {\sl complete} if the labels occurring in it together with the identities of
$\G$ generate all $\G$, i.e. any element of $\G$ which is not an identity can be
obtained as a product of those labels or their inverses.
As it will be shown later, for each one of the categories used, the notion of
completeness is well-posed, i.e it is preserved by equivalence of morphisms. In
particular, a generalized Kirby diagram is complete if and only if it describes a
connected 2-handlebody. Analogously, a labeled ribbon surface is complete if and only
if it describes a connected branched covering.

Given such a monoidal category $\C$, we denote by $\hat\C$ the set of {\sl closed}
morphisms, i.e. the ones having $\one$ (the empty object in $\H^r_n$ and $\K_n$) as
source and target, and by $\hat\C^c \subset \hat\C$ the set of closed complete
morphisms.

\medskip

For $m < n$, there is an injection $\up_m^n: \hat\K_m^c \to \hat\K_n^c$ which sends
$K \in \K_m^c$ to its disjoint union with $n - m$ dotted components of labels $(m,m+1),
(m+1,\break m+2), \dots, (n-1,n)$. We call $\up_m^n K$ the $n$-stabilization of $K$. 
In the handlebody, described by $K$, this corresponds to the creation of $n-m$ pairs 
of canceling 0/1-handles. As it is well known (cf. Section \ref{kirby/sec}), any
connected handlebody with $n$ 0-handles is equivalent (through 1- and 2-handle slides)
to the $n$-stabilization of a 2-handlebody with $m < n$ 0-handles. Therefore the
stabilization map is invertible, and its inverse $\down_m^n: \hat\K_n^c \to \hat\K_m^c$
will be called a reduction map.

The next theorem states the existence of analogous maps between the closed complete
morphisms in $\H_m^r$ and $\H_n^r$ with $m < n$. Hopefully without causing a
confusion, we use for these maps the same notation as for the corresponding maps
between generalized Kirby diagrams.

\begin{theorem}\label{reduction/theo}
 Given $m < n$, there exist a stabilization map $\up_m^n: \hat\H_m^{r,c} \to
\hat\H_n^{r,c}$ and a reduction map $\down_m^n: \hat\H_n^{r,c} \to \hat\H_m^{r,c}$
which are the inverse of each another, such that $\Phi_n \circ \up_m^n = \up_m^n
\circ \Phi_m.$
\end{theorem}

One can also define a stabilization map $\up_m^n: \hat\S_m^c \to \hat\S_n^c$ between
the corresponding sets of labeled ribbon surfaces, by sending $F\in\hat\S_m^c$ to the
disjoint union of $F$ and $n-m$ disks labeled by the permutations $(m\;m{+}1),
(m{+}1\;m{+}2), \dots, (n{-}1\;n)$, and we denote by $\up_m^n F$ the $n$-stabilization
of $F$. Proposition 4.2 in \cite{BP} states that this map is invertible for $n > m
\geq 3$ (that is  any $F' \in \hat\S_n^c$ is equivalent through 1-isotopy and ribbon
moves to the $n$-stabilization of a surface $F\in\hat\S_m^c$). The inverse is called
again a reduction map and is denoted by $\down_m^n: \hat\S_n^c \to \hat\S_m^c$.

Recall that our goal is to factor the bijective correspondence $F \mapsto K_F$ between
the set $\hat\S^c_n$ of complete labeled ribbon surfaces and the set $\hat\K^c_n$ of
complete generalized Kirby diagrams with $n \geq 4$, through a map onto the closed
complete morphisms in the universal algebraic category $\hat\H^{r,c}_n$. This is done
by defining a functor $\S_n \to \H^r_n$  and then composing it with $\Phi_n$. Observe
that the objects in $\S_n$ and $\H^r_n$ are not the same: in $\S_n$ they are intervals
labeled by simple permutations, i.e. unordered pairs of indices, while in $\H^r_n$
they are intervals labeled by morphisms in $\G_n$, i.e. ordered pairs of indices.
Then, any functor $\S_n \to \H^r_n$ requires a choice of an ordering of the indices.
However, its restriction to the closed morphisms is independent on such a choice, as
indicated in the next theorem.

\begin{theorem}\label{defnpsi/theo} 
 Let $<$ denote a strict total order on the set of objects of $\G_n$. Then there exists
a braided monoidal functor $\Psi_n^<: \S_n \to \H^r_n$. Moreover, if $<'$ is another
strict total order on the set of objects of $\G_n$, there is a natural equivalence
$\tau: \Psi_n^< \to \Psi_n^{<'}$ which is identity on the empty set. In particular,
the restriction of $\Psi_n^<$ to $\hat\S_n$ is independent on $<$ and is denoted by
$\Psi_n$.
\end{theorem}

On the left in Figures \ref{comp-funct/fig} and in Figure \ref{comp-funct2/fig} we
present the images under $\Psi_n^<$ of the elementary morphisms with some choices of
labels. We assume $i < j < k$ and $i' < j'$ in Figure \ref{comp-funct/fig}, and $i <
j$, $k < l$, $\{k,l\} \cap \{i,j\} = \emptyset$ in Figure \ref{comp-funct2/fig}. The
images under $\Psi_n^<$ of elementary morphisms with different orderings of the labels
are similar (see Section \ref{surf-alg/sec}).

The following theorem summarizes and completes the algebraic description of
4-dimensional 2-handlebodies.

\begin{theorem} \label{eq-alg-kirby/theo}
 For $n \geq 4$ there is a commutative diagram of bijective maps:
\vskip9pt\centerline{\fig{}{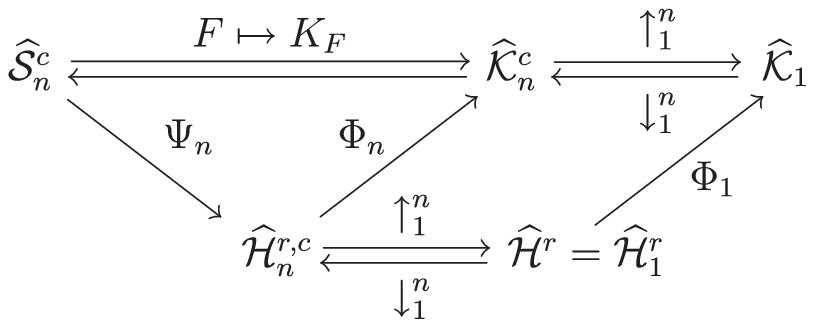}}
\end{theorem}

\vskip-3mm

Observe that, since $F \mapsto K_F$ is a bijective map, the commutativity of the
diagram implies that $\Psi_n$ is injective and $\Phi_n$ is surjective. Therefore to
prove the bijectivity of all maps in the diagram it is enough to show that $\Psi_n$ is 
also surjective when $n \geq 4$.

\medskip
 
In order to deal with 3-manifolds we basically need to quotient all the categories
involved in Theorem \ref{eq-alg-kirby/theo} by some additional relations. Indeed,
according to Kirby \cite{Ki89}, we can think of closed connected 3-manifolds as 
connected 4-dimensional 2-handlebodies modulo 1/2-handle trading and
positive/negative blowups. The quotient category of $\K_n$ modulo these extra relations
is denoted by $\partial \K_n$. Analogously, following \cite{BP}, one can introduce in
$\S_n$ the additional local relations presented in Figure \ref{defn-boundS/fig}.
Those relations preserve up to isotopy the boundary of the ribbon surface while
changing its Euler characteristic. Let $\partial \S_n$ denote the quotient of $\S_n$
modulo these additional relations. Then Theorem 2 in \cite{BP} states that $F \mapsto
K_F$ defines a bijective correspondence between $\partial\hat\S_n^c$ and
$\partial\hat\K_n^c$ for $n \geq 4$.

\begin{Figure}[htb]{defn-boundS/fig}{}
 {Additional relations in $\partial\S_n$.}
\vskip3pt\centerline{\fig{}{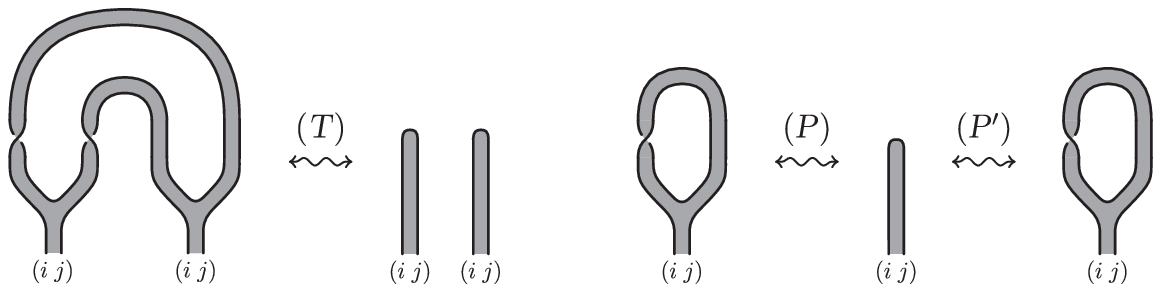}}\vskip-3pt
\end{Figure}

In the same spirit, in Paragraph \ref{defn-boundH/par} we introduce two additional
relations in $\H_n^r$: the first one states the duality of the algebra integral
andcointegral with respect to the copairing, and the second one is a normalization in
sense that a specific closed morphism is sent to the trivial morphism on $\one$. Let
$\partial\H_n^r$ denote the quotient of $\H_n$ modulo these two relations and put
$\partial\hat\H^r = \partial \hat\H_1^r$. Then the algebraic description of
3-manifolds follows from our last theorem.

\begin{theorem}\label{3-mani/theo}
 The braided monoidal functors $\Phi_n$ and $\Psi_n$ induce well defined functors
$\partial\Phi_n: \partial\H^r_n \to \partial\K_n$ and $\partial\Psi_n: \partial\S_n 
\to \partial\H_n^r$ between the quotient categories. Moreover, $\up_m^n$ and
$\down_m^n$ induce well defined bijective maps $\partial\up_m^n: \partial\hat\H_m^{r,c}
\to \partial\hat\H_n^{r,c}$ and $\partial\down_m^n: \partial\hat\H_n^{r,c} \to
\partial\hat\H_m^{r,c}$, for any $m < n$. 
 Consequently, the commutative diagram in Theorem \ref{eq-alg-kirby/theo} induces an
analogous commutative diagram of bijective maps between the corresponding quotient
sets.
\end{theorem}

\bigskip

\centerline{\bfs \hfill Comments and open questions\hfill}

\begin{itemize}\itemsep\smallskipamount
\item[{\sl a}\/)]
 The proofs of all theorems are constructive, i.e. the functors $\Phi_n, \Psi_n$ and
the map $\down_m^n$ are explicitly defined. Moreover, once we know that $\Phi_1$ is a
bijection on the set of closed morphisms, it is easy to describe its inverse, i.e. the
map which associates to each surgery description of a 4-dimensional 2-handlebody a
morphism in $\hat\H^r$, without passing through $\hat\S^c_n$ (cf. \ref{Phi-inv/par}).
\item[{\sl b}\/)] 
 The reader might wonder why we have introduced the general concept of a groupoid
ribbon Hopf algebra, even if it is being used only in the case of a specific and
very simple groupoid. The reasons are two. The first one is that working with the
general case does not make heavier the algebraic part, actually it makes it easier to
follow. The second one is that, in our believe, the group ribbon Hopf algebra (which
is another particular case of the construction) should be useful in finding an
algebraic description of other types of topological objects, for example the group
manifolds studied in \cite{V01}. 
\item[{\sl c}\/)]
 The only reason for which our result concerns closed 3-manifolds and not
cobordisms, is because it is based on the result of \cite{BP}, where the map $F
\mapsto K_F$ is defined and shown to be a bijection only for surfaces/links and not
for tangles. Nevertheless, the whole spirit of the present work (observe that the
factorization of $F \mapsto K_F$ is done through functors defined on the categories
of tangles), suggests that these functors themselves are probably equivalences of
categories, i.e. that $\partial \H^r$ is indeed the category which represents the
algebraic characterization of the 3-dimensional relative cobordisms (see \cite{Ke02}).
\item[{\sl d}\/)]
 By restricting the map $\down_1^2 \circ \Psi_2$ to double branched coverings of
$B^4$, i.e. to ribbon surfaces labeled with the single permutation $(1\;2)$, one
obtains an invariant of ribbon surfaces embedded in $R^4$ under 1-isotopy taking
values in $\H^r$.
\item[{\sl e}\/)]
 Obviously, given any particular braided (selfdual) ribbon Hopf algebra over a ring,
Theorem \ref{eq-alg-kirby/theo} (resp. Theorem \ref{3-mani/theo}) can be used to
construct particular invariants of 4-dimensional 2-handlebodies (resp. 3-manifolds).
All examples of such algebras that we know are simply braidings of (ordinary)
finite-dimensional unimodular ribbon Hopf algebras, therefore the resulting
invariants are the HKR-type invariants from \cite{BM03} (for the definition of the
braided Hopf algebra associated to an ordinary ribbon Hopf algebra, see for example
Lemma 3.7 in \cite{V01}). Nevertheless, we do not know if any finite-dimensional braided
ribbon Hopf algebra is a braiding of an ordinary one.
\end{itemize}

\section*{Contents}
\vskip-\lastskip
\begin{itemize}\itemsep0pt
\item[1.]\vskip-\lastskip Introduction 
 \hfill 1
\item[2.] The category $\K_n$ of admissible generalized Kirby tangles
 \hfill{\pageref{kirby/sec}}
\item[3.] The category $\S_n$ of labeled \rs-tangles 
 \hfill{\pageref{surfaces/sec}}
\item[4.] From labeled ribbon surfaces to Kirby diagrams 
 \hfill{\pageref{surf-kirby/sec}}
\item[5.] The universal braided Hopf algebras $\H(\G)$ and $\H^u(\G)$
 \hfill{\pageref{hopf-alg/section}}
\item[6.] The universal ribbon Hopf algebra $\H^r(\G)$ 
 \hfill{\pageref{ribbon-alg/section}}
\item[7.] From the algebra to Kirby diagrams 
 \hfill{\pageref{alg-kirby/sec}}
\item[] Proof of Theorem \ref{alg-kirby/theo} 
 \hfill{\pageref{alg-kirby/par}}
\item[8.] The reduction map
 \hfill{\pageref{reduction/sec}}
\item[] Proof of Theorem \ref{reduction/theo}
 \hfill{\pageref{reduction/par}} 
\item[9.] From surfaces to the algebra
 \hfill{\pageref{surf-alg/sec}}
\item[] Proof of Theorem \ref{defnpsi/theo}
 \hfill{\pageref{connect-surf/par}}
\item[] Proof of Theorem \ref{eq-alg-kirby/theo}
 \hfill{\pageref{4-mani/par}}
\item[] Proof of Theorem \ref{3-mani/theo}
 \hfill{\pageref{3-mani/par}}
\item[10.] Appendix: proof of Proposition \ref{surfcat/theo}
 \hfill{\pageref{appendix1/sec}}
\item[11.] Appendix: proof of some relations in $\H^r(\G)$
 \hfill{\pageref{appendix/sec}}
\item[] References
 \hfill{\pageref{references/sec}}
\end{itemize}

\vskip12pt\centerline{\fig{}{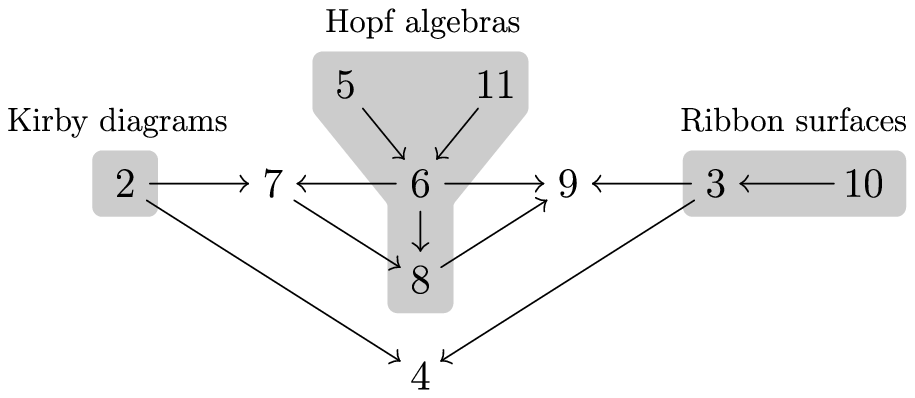}}
\centerline{\small Table of interdependence of sections}

\section{The category $\K_n$ of admissible generalized Kirby tangles%
\label{kirby/sec}}

In this section we review the concept of generalized Kirby diagram from \cite{BP}. 

\begin{block}
 A {\sl Kirby diagram} describes an orientable 4-dimensional 2-handlebody $H^0
\cup H^1_1 \cup \dots \cup H^1_p \cup H^2_1 \cup \dots \cup H^2_q$ with only one
0-handle, by encoding 1- and 2-handles in a suitable link $K \subset S^3 \cong
\Bd H^0$. Namely, $K$ has $p$ dotted components spanning disjoint flat disks which
represent the 1-handles and $q$ framed components which determine the attaching maps
of the 2-handles. We refer to \cite{Ki89} or \cite{GS99} for details and basic facts
about Kirby diagrams, limiting ourselves to recall here only the ones relevant for our
purposes.

The assumption of having only one 0-handle, is crucial in order to make a natural
convention on the framings, that allows to express them by integers fixing as zero
the homologically trivial ones.

However, we renounce this advantage on the notation for framings in favour of
more flexibility in the representation of multiple 0-handles. The reason is that an
$n$-fold covering of $B^4$ branched over a ribbon surface turns out to have a
natural handlebody structure with $n$ 0-handles.

\medskip

We call a generalized Kirby diagram a representation of an orientable
4-dimen\-sional 2-handlebody with multiple 0-handles. It is essentially defined by
identifying the boundaries of all 0-handles with $S^3$ (where the diagram takes
place), and by putting labels in the diagram in order to keep trace of the original
0-handle where each part of it is from. If there is only one 0-handle, the labels
can be omitted and we have an ordinary Kirby diagram.

\begin{Figure}[b]{diag1/fig}{}{One-handle notation.}
\vskip-6pt\centerline{\fig{}{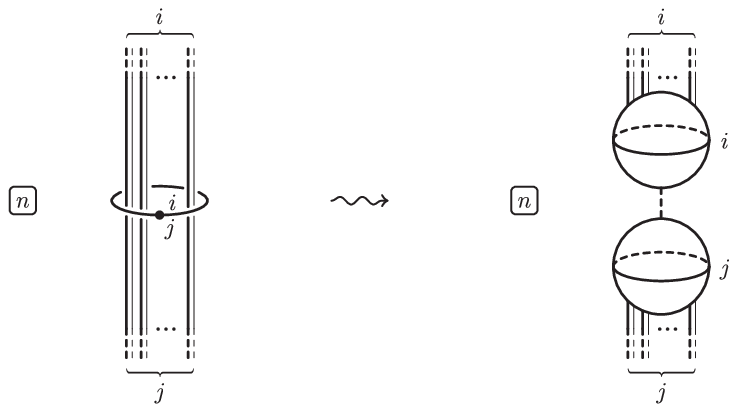}}\vskip-6pt
\end{Figure}

More precisely, a {\sl generalized Kirby diagram} representing an orientable
4-di\-mensional handleboby $H^0_1 \cup \dots \cup H^0_n \cup H^1_1 \cup \dots \cup
H^1_p \cup H^2_1 \cup \dots \cup H^2_q$ consists of the following data: a boxed
label indicating the number $n$ of 0-handles; $p$ dotted unknots spanning disjoint
flat disks, each side of which has a label from $\{1, \dots, n\}$; $q$ integer framed
disjoint knots transversal with respect to those disks, with a label from $\{1, \dots,
n\}$ for each component of the complement of the intersections with the disks. The
labeling must be admissible in the sense that all the framed arcs coming out from one
side of a disks have the same label of that side (cf. Figure \ref{diag1/fig}). This
rule makes the labeling redundant and sometimes we omit the superfluous labels.
Moreover, the framings are always drawn as parallel curves, i.e. the undotted
components represent embeddings of $S^1\times I$ in $R^3$.

To establish the relation between a generalized Kirby diagram and the handlebody it
represents, we first convert the dot notation for the 1-handles into a ball
notation, as shown in Figure \ref{diag1/fig}. Here, the two balls, together with the
relative framed arcs, are symmetric with respect to the horizontal plane containing
the disk and squeezing them vertically on the disk we get back the original diagram.
After that, we take a disjoint union of $n$ 4-balls $H^0_1 \cup \dots \cup H^0_n$,
which are going to be the 0-handles, and draw on the boundary of each $H^0_i$ the
portion of the diagram labeled with $i$. Then, we attach to $H^0_1 \cup \dots \cup
H^0_n$ a 1-handle between each pair of balls (possibly lying in different
0-handles), according to the diffeomorphism induced by the above symmetry, so that
we can join longitudinally along the handle the corresponding framed arcs. Of
course, the result turns out to be defined only up to 1-handle full twists. At this
point, we have a 1-handlebody $H^0_1 \cup \dots \cup H^0_n \cup H^1_1 \cup \dots
\cup H^1_p$ with $q$ framed loops in its boundary and we use such framed loops as
attaching instructions for the 2-handles $H^2_1, \dots, H^2_q$.

We observe that any orientable 4-dimensional 2-handlebody can be represented, up to
isotopy, by a generalized Kirby diagram. In fact, in order to reverse our
construction, we only need that the identification of the boundaries of the
0-handles with $S^3$ is injective on the attaching regions of 1- and 2-handles and
that the attaching maps of the 2-handles run longitudinally along the 1-handles.
These properties can be easily achieved by isotopy.

\begin{Figure}[b]{diag3/fig}{}
 {Labeled isotopy moves ($k \neq l$ in the last crossing change).}
\vskip-8pt\centerline{\fig{}{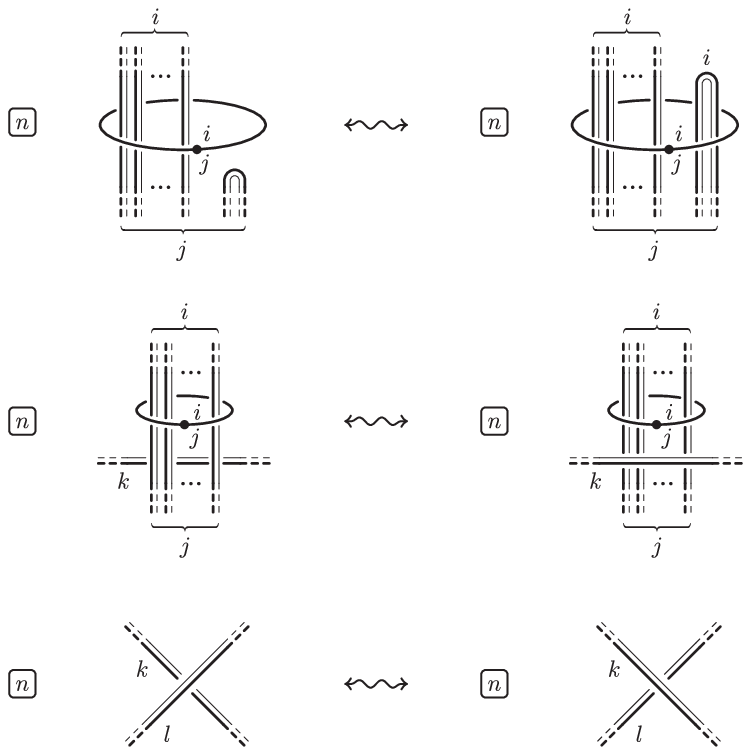}}
\end{Figure}

The above construction gives isotopic handlebody structures if and only if the
starting generalized Kirby diagrams are equivalent up to {\sl labeled isotopy},
generated by labeled diagram isotopy, preserving all the intersections between loops
and disks (as well as labels), and by the three moves described in Figure
\ref{diag3/fig}. Here, in the last move we assume $k \neq l$, so that the crossing
change at the bottom of the figure preserves the isotopy class of the framed link in
$H^0_1 \cup \dots \cup H^0_n \cup H^1_1 \cup \dots \cup H^1_p$. It is due to this
crossing change that the framing convention usually adopted for ordinary Kirby
diagrams cannot be extended to generalized Kirby diagrams.
On the contrary, the other two moves make sense whatever are $i$, $j$ and $k$. In
particular, if $i = j = k$ they reduce to the ordinary ones. Actually, this is the
only relevant case for the second move, usually referred to as ``sliding a 2-handle
over a 1-handle'', since the other cases can be obtained by crossing changes.
Moreover, even this ordinary case becomes superfluous in the context of
2-deformations, since it can be realized by addition/deletion of canceling
1/2-handles and 2-handle slides (cf. \cite{GS99}).

The following Figures \ref{diag4/fig} and \ref{diag5/fig} represent 2-deformations of
4-dimensional 2-handlebodies in terms of generalized Kirby diagrams. Namely: the moves
in Figure \ref{diag4/fig} correspond to addition/deletion of canceling 0/1-handles and
1/2-handles, where $i \leq n$ and the canceling framed component in the bottom move is
assumed to be closed; the moves in Figure \ref{diag5/fig} correspond to 1- and
2-handle slidings (in the low left corner of the figure we assume that the two
framed segments correspond to different components of the link). Except for the
addition/deletion of canceling 0/1-handles, which does not make sense for ordinary
Kirby diagrams, the rest of the moves reduce to the ordinary ones if $i = j = k$.

\begin{Figure}[htb]{diag4/fig}{}{Addition/deletion of canceling pairs of handles.}
\vskip4pt\centerline{\fig{}{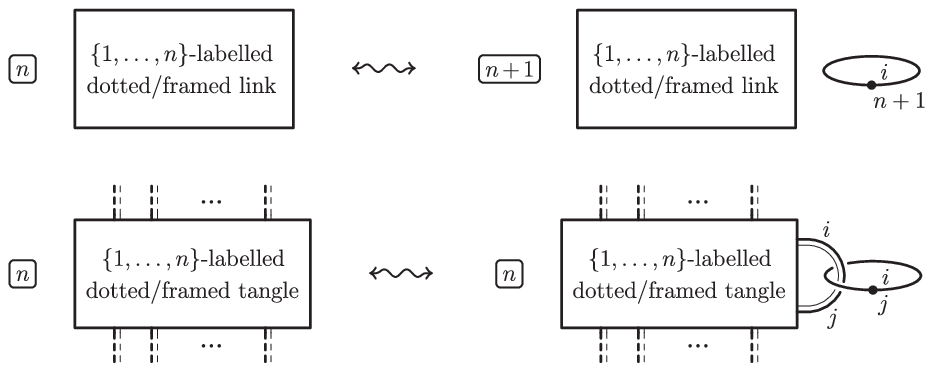}}
\end{Figure}

\begin{Figure}[htb]{diag5/fig}{}{Handle slidings.}
\vskip-9pt\centerline{\fig{}{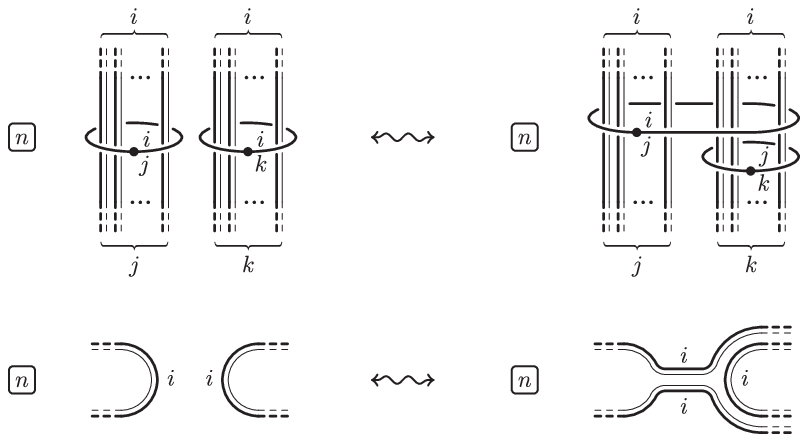}}
\end{Figure}

The 1-handle sliding is included for the sake of completeness, but it can be
generated by addition/deletion of canceling 1/2-handles and 2-handle sliding, just
like in the ordinary case (cf. \cite{GS99}).

Summing up, two generalized Kirby diagrams represent 2-equivalent 4-dimen\-sional
2-handlebodies if and only if they can be related by the first and third moves in
Figure \ref{diag3/fig} (labeled isotopy), the two moves in Figure \ref{diag4/fig}
(addition/deletion of canceling handles) and the second move in Figure
\ref{diag5/fig} (2-handle sliding).

Of these, only the first move in Figure \ref{diag3/fig} and the second ones in
Figures \ref{diag4/fig} and \ref{diag5/fig} (for $i = j = n = 1$) make sense in the
case of ordinary Kirby diagrams. Actually, such three moves suffice to realize
2-equivalence of 4-dimensional 2-handlebodies with only one 0-handle, since any
extra 0-handle occurring during a 2-deformation can be eliminated by a suitable
fusion of 0-handles.
\end{block}

\begin{block}\label{red-kirby/par}
 Given any generalized Kirby diagram $K$ representing a connected handlebody with $n$
0-handles, we can use 2-deformation moves to transform it into an ordinary one, by
reducing the number of 0-handles to 1. In fact, assuming $n > 1$,
\mypagebreak
we can eliminate the $n$-th 0-handle as follows. Since the handlebody is connected, there exists at least one dotted component such that one side of its spanning disk is labeled by $n$, and the other one is labeled by $i\neq n$. We may assume that all such disks lie in the $xy$ plane with the side labeled $n$ facing up. Moreover, by changing the crossings of framed components with different labels and by isotopy we may also assume that the diagram intersects the $xy$ plane only in these spanning disks and that all dotted and framed components above the $xy$ plane contain only the label $n$. Now we choose one of the components in the $xy$ plane, say $U$ of label $(i_0,n)$, and slide it over all dotted components in this plane changing their labels from $(i,n)$ to $(i,i_0)$
(cf. the top move in Figure \ref{diag5/fig}). Then we pull $U$ up until it becomes
disjoint from the rest of the diagram and this changes the labels of the dotted
components above the plane from $(n,n)$ to $(i_0,i_0)$. The resulting diagram is the
disjoint union of $U$ and a diagram $K^U$ which does not contain the label $n$. Finally
we cancel the $n$-th 0-handle with the one handle represented by $U$. An example of
such reduction is presented in Figure \ref{example-red/fig}, where the 1-handle $U$ is
the one in the upper right corner of the diagram. Observe that if we had chosen
another 1-handle, for example the one in the lower left corner, we would obtain the
diagram on the right in Figure \ref{example-red1/fig} which, as it is shown in the
figure (and as it should be), is equivalent to the previous one through 1-handle
slides.
\end{block}

\begin{Figure}[htb]{example-red/fig}{}{}
\centerline{\fig{}{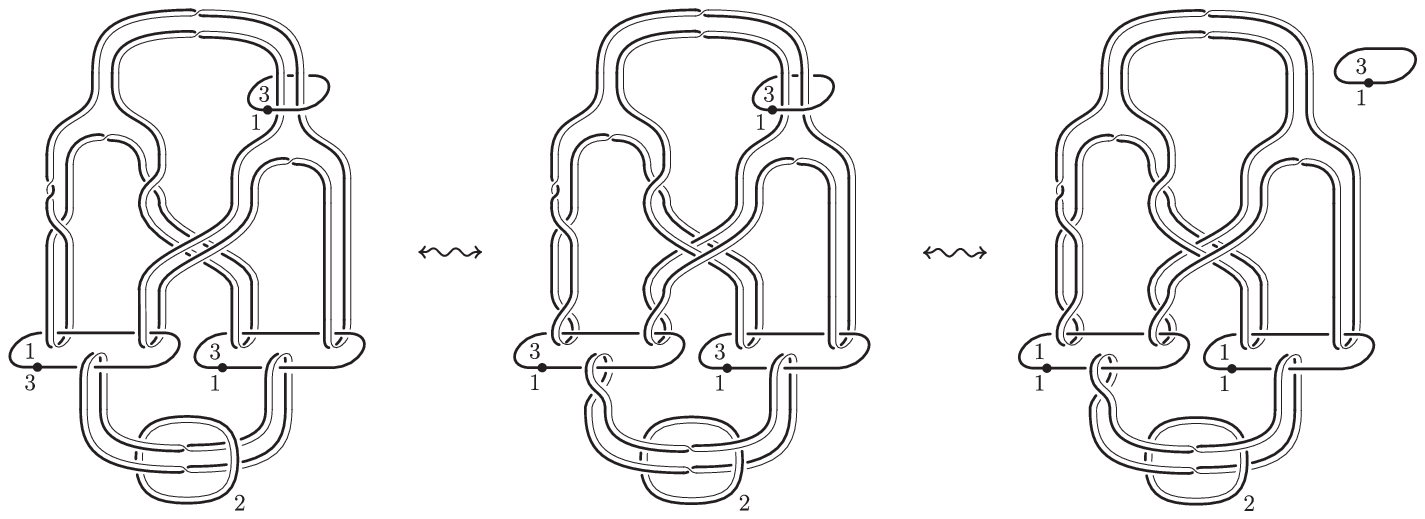}}
\end{Figure}

\begin{Figure}[htb]{example-red1/fig}{}{}
\centerline{\fig{}{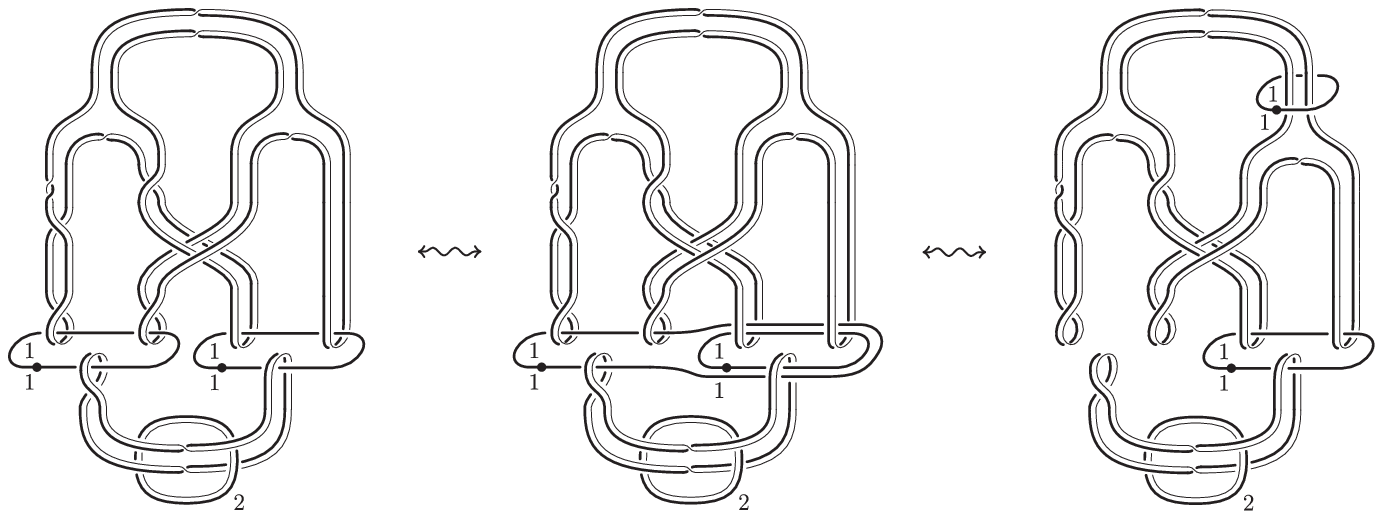}}
\end{Figure}

\begin{block}\label{monoidal/par}
 Now we introduce the notion of an admissible generalized Kirby tangle and organize
such tangles into a monoidal category. 

We recall that a {\sl braided monoidal} category $(\C,\,\diam,\, \one,\, \iota_A,\,
{}_A\iota,\, \alpha,\gamma)$ is a category $\C$ equipped with a functor ${\diam}: \C
\times \C \to \C$, an object $\one\in \C$ and natural isomorphisms:
\vskip-9pt
$$
\begin{array}{rl}
 \iota_A: &\one \diam A \to A,\\
 {}_A\iota: &A \diam \one \to A,\\
 \alpha_{A,B,C}: &(A \diam B) \diam C \to A \diam (B \diam C)
 \quad\text{(associativity)},\\
 \gamma_{A,B}: &A \diam B \to B \diam A 
 \quad\text{(commutativity or braiding)},
\end{array}
$$
\vskip-6pt\noindent
 which satisfy a well known set of axioms (see for example \cite{Sh94}). $\diam$ is
called the product of the category. In the category of generalized Kirby tangles this
product will be given by the juxtaposition of diagrams. This is possible only if the
number of the 0-handles (labels) is fixed; otherwise the equivalence move presented on
the top in Figure \ref{diag4/fig} would violate the monoidal structure. So, in what
follows we assume that the box of all diagrams contains a fixed $n$, and therefore we
omit it.

A {\sl generalized Kirby tangle} is a slice of a generalized Kirby diagram, i.e. the
intersection of a generalized Kirby diagram $K$ with $R^2 \times [0,1]$, where the
dotted components of $K$ (and their spanning disks) do not intersect $R^2 \times
\{0, 1\}$, while the framed components of $K$ intersect $R^2 \times \{0, 1\}$
transversely. Two generalized Kirby tangles are equivalent if they can be
transformed into each other through ambient isotopy of $R^2 \times [0,1]$, the first
and the third moves in Figure \ref{diag3/fig}, the lower move in Figure
\ref{diag4/fig}, and the second move in \ref{diag5/fig}, leaving $R^2 \times \{0, 1\}$
fixed. Moreover, in the move in Figure \ref{diag5/fig} we assume that the two segments
belong to different components of the tangle and that the component over which the
sliding is done is closed.
	
We define the category $\K_n$ of ($n$-labeled) {\em admissible} generalized
Kirby tangles as follows. An object of $\K_n$ is an ordered set of pairs of intervals
$\{(I_1^+,I_1^-), \ldots, (I_s^+,I_s^-)\}$ in $R^1\subset R^2$, each one labeled with
pair of numbers $(i_k^+,i_k^-)$, $1 \leq i_k^+,i_k^- \leq n$. To simplify the notation
we will often denote such object simply by the ordered sequence of labels. A morphism
in $\K_n$ is given by a generalized Kirby tangle each component of which is either
closed or intersects at most one of the planes $R^2 \times \{0\}$ and $R^2 \times
\{1\}$ in a pair of intervals $(I^+_k,I^-_k)$, in such a way that the label of each
interval coincides with the label of the framed component to which it belongs. 
Moreover, any open component has a half-integer framing. The composition of two
morphisms is obtained by identifying the target of the first one with the sourse of
the second and rescaling. Observe that our notion of admissible tangles is more
restrictive than the one in \cite{MaP92} and \cite{KL01}, but this is all we need. 
$\K_n$ is a braided monoidal category with respect to juxtaposition of diagrams. In
particular, $\one$ is the empty set, the identity on $(i,j)$ and the braiding
morphisms $\gamma_{(i,j),(i',j')}$ and $\gamma_{(i',j'),(i,j)}^{-1}$ are presented
in Figure \ref{morphK/fig}, while $\iota_A,\,{}_A\iota$ and the associativity morphisms
are the identities.

\begin{Figure}[htb]{morphK/fig}{}{Identity and the braiding morphisms in $\K_n$.}
\centerline{\fig{}{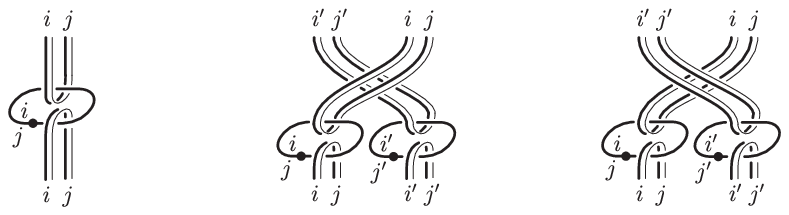}}
\end{Figure}

We recall from the Introduction that a morphism $K$ in $\K_n$ is closed if it has the
empty object as sourse and target. Moreover, $K$ is complete if the labels of the
dotted components (taken as ordered pairs of indices) together with the identities
$(i,i)$, $i = 1 , \dots, n$, generate all the groupoid $\G_n$. Equivalently, $K$ is
complete if and only if the graph having $n$ ordered vertices and one edge connecting
the $i$-th and the $j$-th vertex for any dotted component of label $(i,j)$, is
connected. To see that the notion of completeness is well-posed with respect to
equivalences of morphisms, we observe that the only move in which a dotted component
may appear or disappear is the bottom move in Figure \ref{diag4/fig} where if there $i
\neq j$, neccesarily the cancelling framed component passes through other 1-handes,
the corresponding edges of which connect the $i$-th and the $j$-th vertex of the
graph. In particular, removing or adding the edge corresponding to the cancelling
dotted component does not change the connectedness of the graph.

Let $\hat\K_n^c$ denote the set of closed complete morphisms in $\K_n$. Of course, a
closed morphism  is complete if and only if the corresponding handlebody is connected.
Then the stabilization with an 1-handle of label $(n+1, n)$ shown in the upper part
of Figure \ref{diag4/fig} defines a bijective map $\up_n^{n+1}:\hat\K_n^c \to
\hat\K_{n+1}^c$, while the procedure of reducing the number of 0-handles given in
\ref{red-kirby/par} represents its inverse $\down_n^{n+1}: \hat\K_{n+1}^c \to
\hat\K_n^c$. Observe that ${\up_n^{n+1}}K \in \hat\K_{n+1}^c$ and $K \in \hat\K_n^c$
describe the same 4-dimensional 2-handlebody up to 2-deformation. 
\end{block}

\begin{block}
Even if it is irrelevant for the present work, we would like to point out that the
category $\K_n$ is equivalent to the category of $n$-labeled $3+1$ cobordisms $\cob_n$
(for $n=1$ see Section 1 in \cite{Ke99}). The objects in $\cob_n$ are oriented
3-dimensional 1-handlebodies with $n$ 0-handles. Given two such handlebodies $M_s$ and
$M_r$ respectively with $s$ and $r$ 1-handles, let $N(M_s, M_r)$ denote the
3-dimensional 1-handlebody obtained from $M_s \sqcup M_r$ by attaching for any $i \leq
n$ a single 1-handle connecting the $i$-th 0-handles of $M_s$ and $M_r$. We can cancel
the new 1-handles against some of the 0-handles thinking of $N(M_s, M_r)$ as a
3-dimensional 1-handlebody with $n$  0-handles as well. Then a morphism $W: M_s \to
M_r$ is a relative (4-dimensional) 2-handlebody build on $N(M_s, M_r)$. Note that the
term relatitve handlebody build on $M$, is usually limited to the case when $M$ is a
closed 3-manifold, but the generalization is straightforward: we attach 1- and
2-handles on $N(M_s, M_r) \times \{1\} \subset N(M_s, M_r) \times [0,1]$. Then $W_1,
W_2: M_s \to M_r$ are called 2-equivalent if they are obtained from each other by a
2-deformation, i.e. changing the attaching maps of the handles through isotopy or
adding/deleting a canceling pair of 1/2-handles. Obviously these handle moves are
limited to $N(M_s, M_r) \times \{1\}$. The composition of two relative cobordims $W:
M_s \to M_r$ and $W': M_r \to M_t$ is obtained by gluing through an orientation
reversing homeomorphism $M_r \times \{0\} \subset N(M_s, M_r) \times [0,1]$ and $M_r
\times \{0\} \subset N(M_r, M_t) \times [0,1]$ and canonically identifying the
resulting manifold as $N(M_s, M_t) \times [0,1]$ on which we have attached $r$
4-dimensional 1-handles. 

\begin{Figure}[htb]{diag7/fig}{}
 {Identifying a morphism in $\K_n$ as relative cobordism.}
\centerline{\fig{}{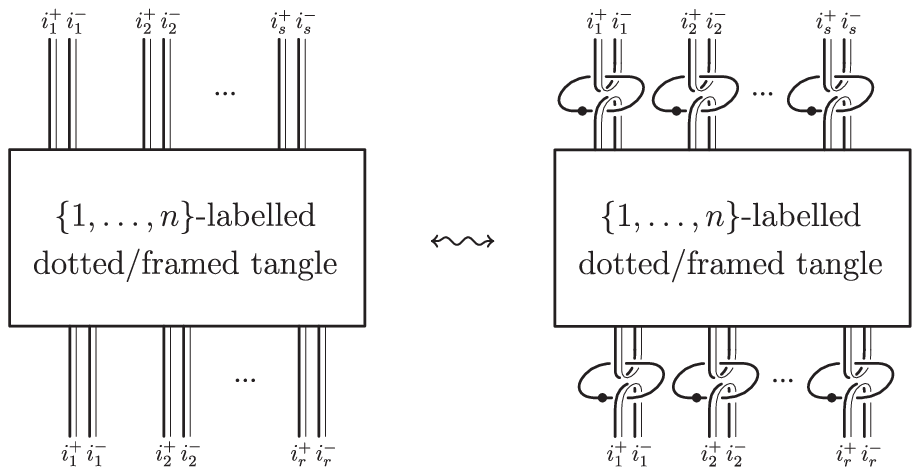}}\vskip-3pt
\end{Figure}

Any morphism $K: \{(i_1^+, i_1^-), \dots, (i_s^+, i_s^-)\} \to \{(j_1^+, j_1^-), \dots,
(j_r^+, j_r^-)\}$ can be easily seen as describing a cobordism $W(K): M_s \to M_r$ by
composing it (if necessary) with the identity morphisms as shown in Figure
\ref{diag7/fig}.
\end{block}

\begin{block}\label{kirby-3man/par}
 The main theorem of Kirby calculus \cite{Ki89} asserts that two orientable
4-dimensional 2-handlebodies have diffeomorphic boundaries and the same signature if
and only if they are related by 2-deformations and 1/2-handle trading, while two
orientable 4-dimensional 2-handlebodies have diffeomorphic boundaries if and only if
they are related by 2-deformations, positive/negative blowing up/down and 1/2-handle
trading.

\begin{Figure}[b]{diag6/fig}{}{Blowing up/down and 1/2-handle trading.}
\centerline{\fig{}{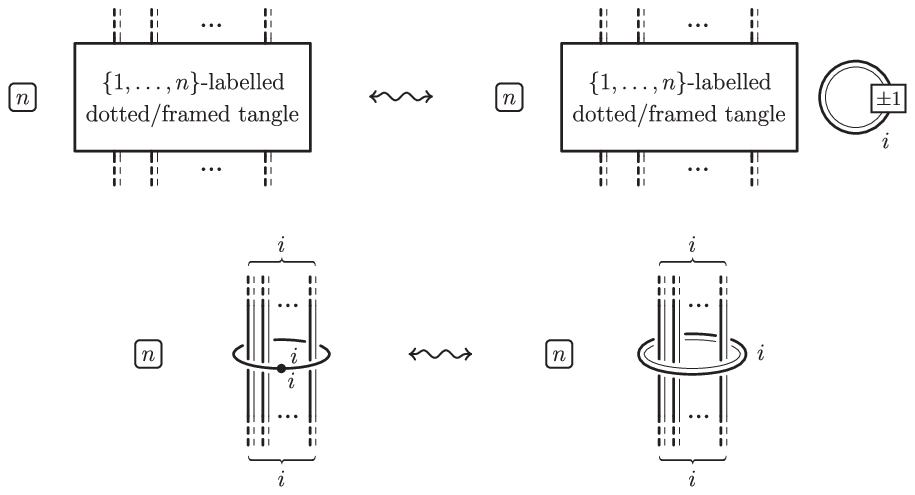}}
\end{Figure}

\begin{Figure}[htb]{defn-boundK1/fig}{}
 {Reducing blowing up/down and 1/2-handle trading to special cases.}
\vskip4pt\centerline{\fig{}{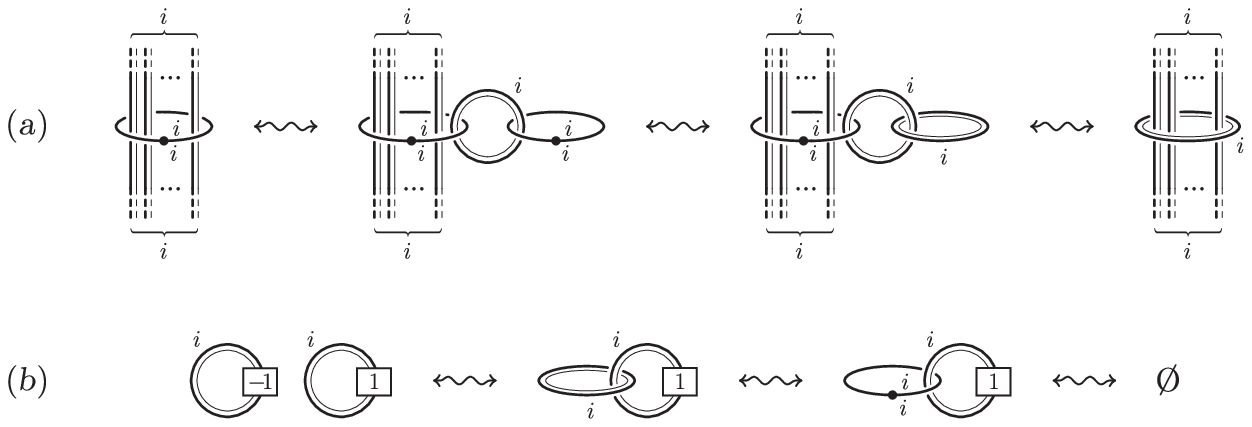}}
\end{Figure}

In terms of generalized Kirby diagrams these last three modifications can be realized
by the moves in Figure \ref{diag6/fig}. These moves essentially coincide with the
corresponding ones for ordinary Kirby diagrams (with $i = n = 1$), being the involved
labels all the same.

As it is shown in Figure \ref{defn-boundK1/fig} \(a), modulo 2-deformations, 1/2-handle
trading is equivalent to 1/2-handle trading for a canceling pair, which is the move
presented 
\mypagebreak
in Figure \ref{defn-boundK/fig} \(a). Moreover, the relation in Figure
\ref{defn-boundK1/fig} \(b) implies that in $\partial^{\star}K_n$ a negative blow up
can be obtained through 2-deformations and 1/2-handle trading from the positive one
(Figure \ref{defn-boundK/fig} \(b)). Therefore we define the category
$\partial^\star\K_n$ to be the $\K_n$ modulo the relation in Figure
\ref{defn-boundK/fig} \(a), while $\partial \K_n$ to be the $\K_n$ modulo the relations
in Figure \ref{defn-boundK/fig} \(a) and \(b). From \cite{Ki89} and the discussion
above it follows that the closed morphisms $\partial^\star \hat\K_n$ in $\partial^\star
\K_n$ describe framed 3-manifolds, while the closed morphisms $\partial\hat\K_n$ in
$\partial \K_n$ describe 3-manifolds. 
\begin{Figure}[htb]{defn-boundK/fig}{}
 {Additional relations in $\partial^\star \K_n$ and $\partial \K_n$.}
\centerline{\fig{}{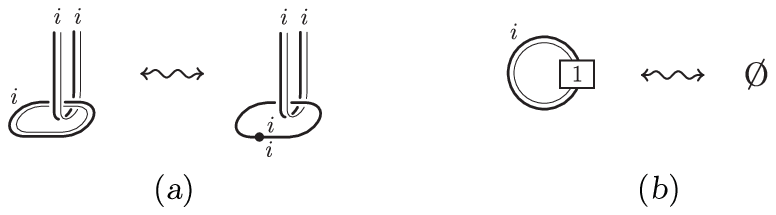}}\vskip-2pt
\end{Figure}
\end{block}

\vskip-12pt\vskip0pt

\section{The category $\S_n$ of labeled \rs-tangles%
\label{surfaces/sec}}

We review the notion of ribbon surface and 1-isotopy of such surfaces from
\cite{BP}. A smooth compact surface $F \subset B^4$ with $\Bd F \subset S^3$ is
called a {\sl ribbon surface} if the Euclidean norm restricts to a Morse function on
$F$ with no local maxima in $\Int F$. Assuming $F \subset R^4_- \subset R^4_- \cup
\{\infty\} \simeq B^4$, this property is topologically\break equivalent to the fact
that the fourth Cartesian coordinate restricts to a Morse height function on $F$
with no local maxima in $\Int F$. Such a surface $F \subset R^4_-$ can be
horizontally (preserving the height function) isotoped to make its orthogonal
projection into $R^3$ a self-transversal immersed surface, whose double points form
disjoint arcs as in Figure \ref{ribbsurf1/fig} \(a).

\begin{Figure}[htb]{ribbsurf1/fig}{}{}\smallskip
\vskip4pt\centerline{\fig{}{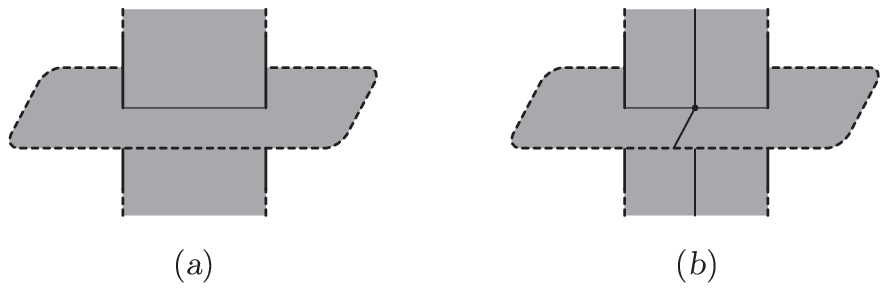}}\vskip-4pt
\end{Figure}

We refer to such a projection as a {\sl 3-dimensional diagram} of $F$. Actually, any
immersed compact surface $F \subset R^3$ with no closed components and all
self-intersections of which are as above, is the diagram of a ribbon surface
uniquely determined up to vertical isotopy. This can be obtained by pushing $\Int F$
down inside $\Int R^4_-$ in such a way that all self-intersections disappear. 

\begin{Figure}[htb]{ribbsurf5/fig}{}{1-isotopy moves.}
\centerline{\fig{1-isotopy moves}{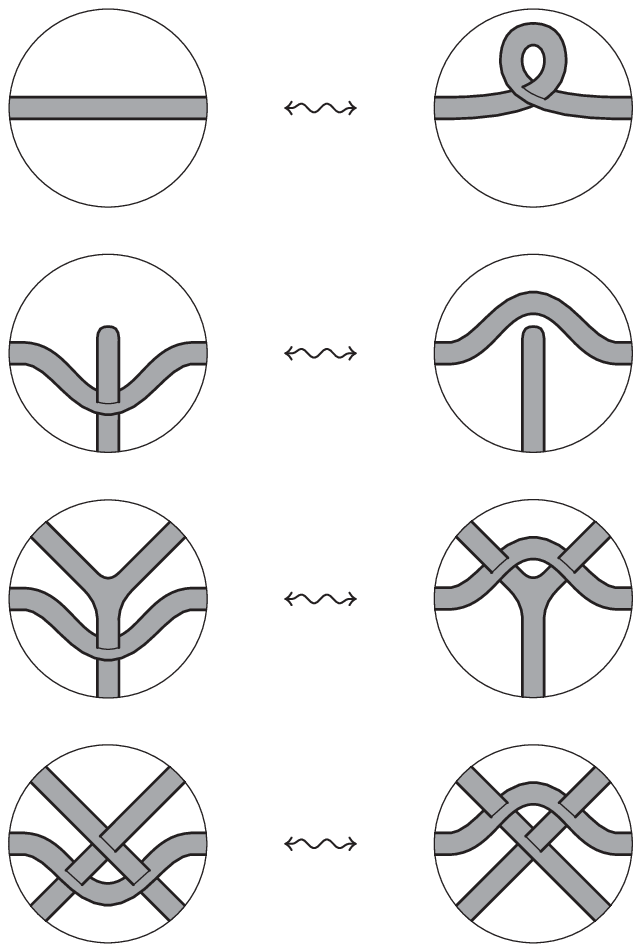}}
\end{Figure}

In, [BP] we introduced {\sl 1-isotopy} between ribbon surfaces, to be generated by
3-dimensional diagram isotopy in $R^3$ (that means isotopy preserving ribbon
intersections) and the 1-isotopy moves depicted in Figure \ref{ribbsurf5/fig}.

\medskip

Since a ribbon surface $F$ has no closed components, any 3-dimensional diagram of
it, considered as a 2-dimensional complex in $R^3$, collapses to a graph $G$. We choose
\mypagebreak
$G$ to be the projection in $R^3$ of a smooth simple spine $\tilde G$ of
$F$ (simple means that all the vertices have valency one or three), which meets at
exactly one point each arc projecting to a ribbon intersection of the 3-dimensional
diagram of $F$ as shown in Figure \ref{ribbsurf1/fig} \(b). The inverse image in
$\tilde G$ of such a point consists of a single uni-valent vertex and a point in the
interior of an edge, while the restriction of the projection to $\tilde G$ with such
uni-valent vertices removed, is injective.

Therefore, $G$ turns out to have only vertices of valency 1 and 3. We call {\sl
singular vertices} the tri-valent vertices located at the ribbon intersections,
and {\sl flat vertices} all the other vertices. Moreover, we assume that $G$
has distinct tangent lines at each flat tri-valent vertex while the tangent lines
to two of the edges at a singular vertex coincide (since those two edges form a unique
edge of $\tilde G$) and differ from the tangent line to the third edge at such
vertex.

Up to a further horizontal isotopy, we can contract the 3-dimensional diagram of
$F$ to a narrow regular neighborhood of the graph $G$. Moreover, by considering a
planar diagram of $G$, we can assume that the diagram of $F$ is contained in the
projection plane, except for a finite number of positive/negative ribbon half-twists,
ribbon intersections and ribbon crossings, as the ones depicted in Figure
\ref{ribbsurf2/fig}.

\begin{Figure}[htb]{ribbsurf2/fig}{}{}
\centerline{\fig{}{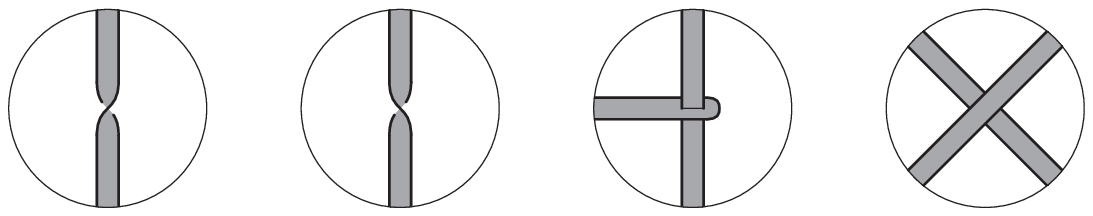}}
\end{Figure}

In this way we get a new diagram of $F$, consisting of a certain number of regions as
the ones presented in Figures \ref{ribbsurf2/fig} and \ref{ribbsurf3/fig}, suitably
connected by flat bands contained in the projection plane. We call such diagram a
{\sl planar diagram} of $F$.

\begin{Figure}[htb]{ribbsurf3/fig}{}{}
\vskip6pt\centerline{\fig{}{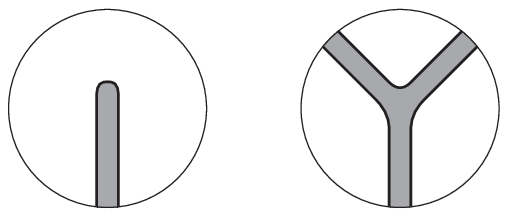}}
\end{Figure}

Actually, a planar diagram of $F$ arises as a diagram of the pair $(F,G)$ and this is
the right way to think of it. However, we omit to draw the graph $G$ in the pictures
of a planar diagram, since it can be trivially recovered, up to diagram isotopy, as
the core of the diagram itself. In particular, the singular vertices and the
projection crossings are located at the centers of the third circle in Figure
\ref{ribbsurf2/fig}, while the flat vertices of $G$ are located at the centers of the
circles in Figure \ref{ribbsurf3/fig}. To be precise, there are two choices in
recovering the graph $G$ at a singular vertex, as shown in Figure \ref{ribbsurf4/fig}.
They give the same graph diagram of each other, but differ for the way the graph is
embedded in the 3-dimensional diagram of the surface. We consider the move depicted in
Figure \ref{ribbsurf4/fig} as an equivalence move for the pair $(F,G)$. Up to this
move, $G$ is uniquely determined also as a graph in the 3-dimensional diagram of $F$.

\begin{Figure}[htb]{ribbsurf4/fig}{}{}
\vskip6pt\centerline{\fig{}{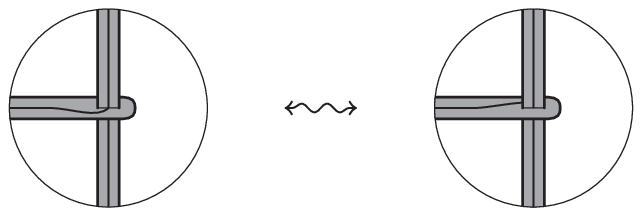}}
\end{Figure}

Like the 3-dimensional diagrams, also the planar ones uniquely determine the ribbon
surface $F$ up to vertical isotopy. Here, by vertical isotopy we mean an isotopy
which preserves the first 2 coordinates. In other words, the 3-dimensional height
functions (as well as the 4-dimensional one) is left undetermined for a planar
diagram. Of course, such height function is required to be consistent with the
restrictions deriving from the local configurations in Figure \ref{ribbsurf2/fig}.

\medskip

In the following, {\sl ribbon surfaces will be always represented by planar diagrams
and considered up to vertical isotopy} (preserving the first 2 coordinates).

\medskip

 By a {\sl ribbon surface tangle} (or simply {\sl \rs-tangle}) we mean the intersection
of a pair $(F,G)$, consisting  of a ribbon surface $F$ and its core graph $G$, with
$R^2 \times [0,1]\subset R^3$, where $F$ intersects transversely $R^2 \times \{0,1\}$
in some sets of arcs, called boundary arcs, and $G$ intersects transversely $R^2 \times
\{0,1\}$ in some set of points, one for each boundary arc, different from its vertices.

Two \rs-tangles are {\sl 1-isotopic} if they are related by 1-isotopy inside $R^2
\times [0,1]$ leaving fixed a small regular neighborhood of the boundary arcs.
The \rs-tangles form a category $\S$ whose objects are sets of intervals in
$R^2$, and whose morphisms are \rs-tangles in which source and target are respectively
given by the intersection of the \rs-tangle with $R^2 \times \{0\}$ and
$R^2 \times \{1\}$. Then the composition of two morphisms is obtained by identifying
the target of the first one with the source of the second and rescaling. We also
define the product of two \rs-tangles by horizontal juxtaposition.

\begin{Figure}[b]{ribbsurf13/fig}{}{Elementary diagrams in $\S$}
\centerline{\ \fig{}{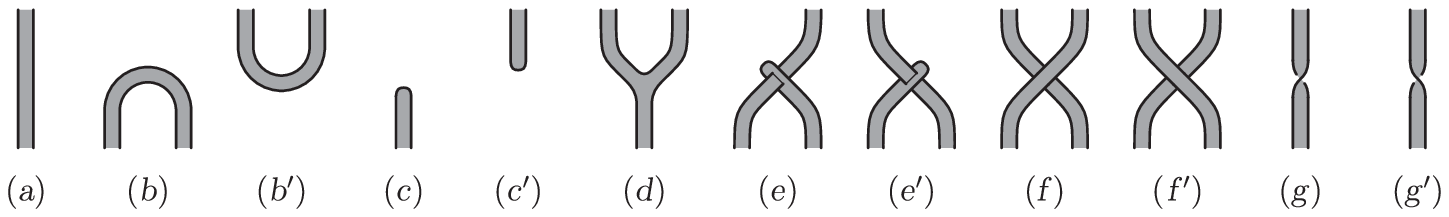}}\vskip-4pt
\end{Figure}

\begin{Figure}[b]{ribbsurf14/fig}{}{Planar isotopy relations in $\S$}
\medskip\centerline{\fig{Moves 2}{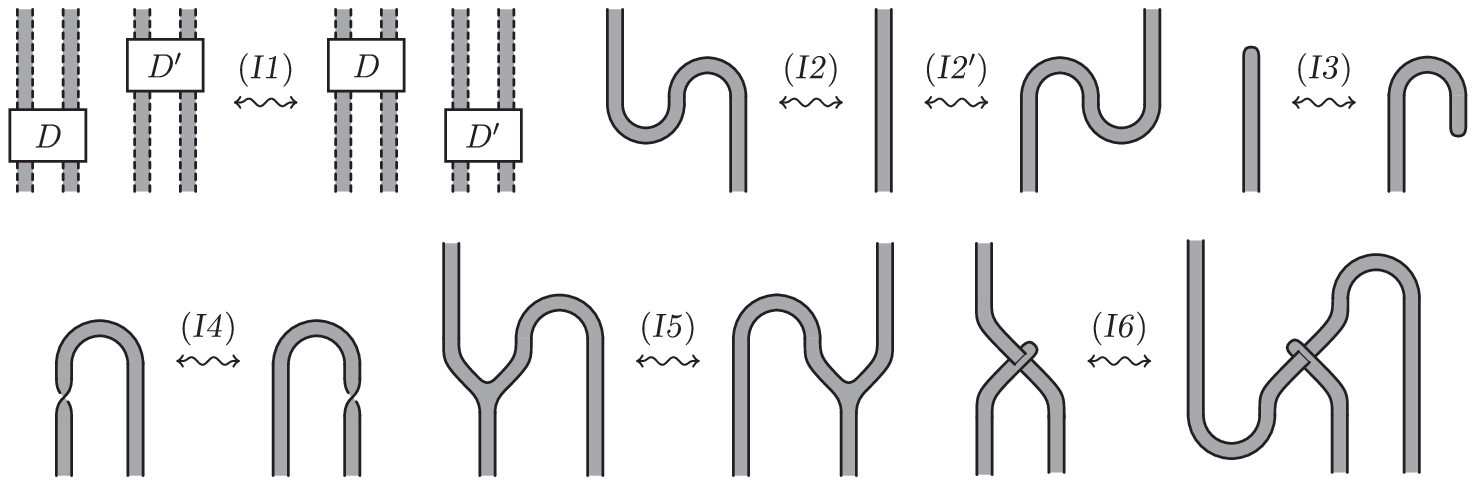}}
\end{Figure}

\begin{Figure}[b]{ribbsurf15/fig}{}{3-dimensional isotopy relations in $\S$}
\medskip\centerline{\fig{Moves 2}{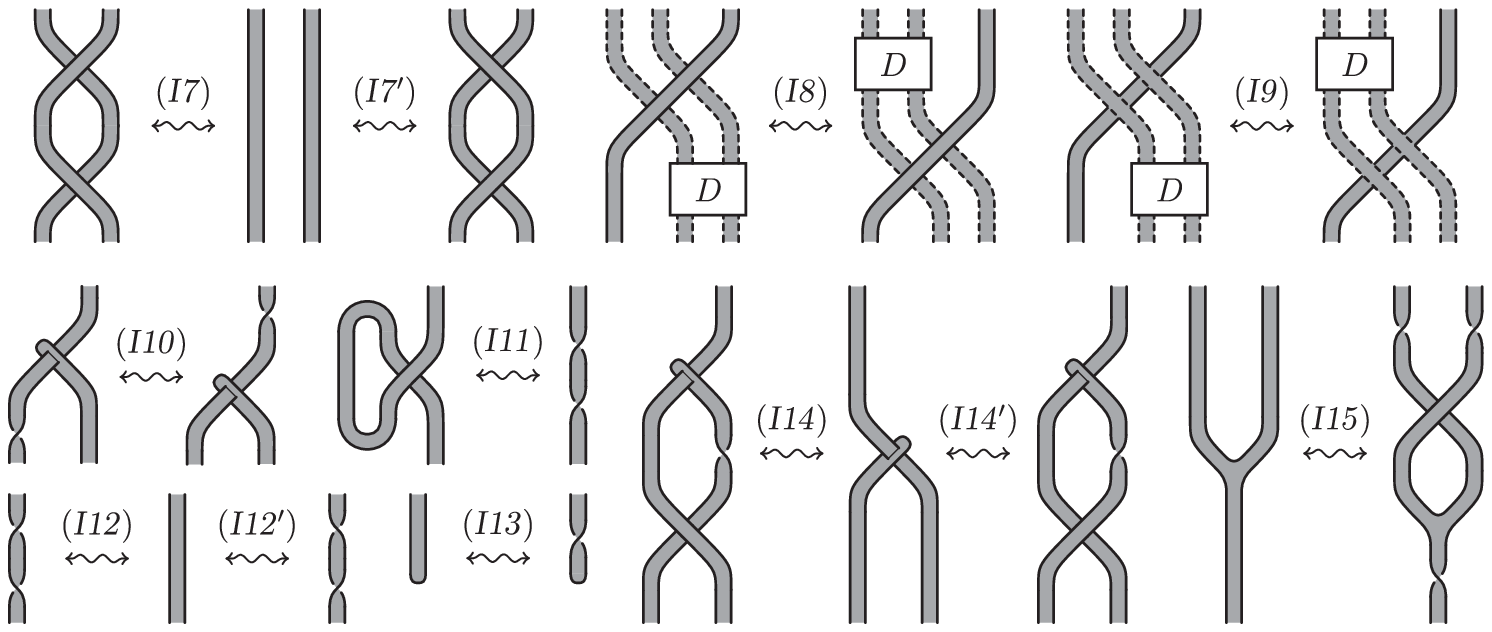}}
\end{Figure}

\begin{proposition}\label{surfcat/theo}
 $\S$ is equivalent to the category of planar diagrams, whose objects are finite
sequences of disjoint intervals in $R$, and whose morphisms are iterated products
and compositions of the elementary planar diagrams presented in Figure
\ref{ribbsurf13/fig}, modulo plane isotopies preserving the $y$-coordinate and the
moves presented in Figures \ref{ribbsurf14/fig}--\ref{ribbsurf17/fig}, where $D$ and
$D'$ in moves \(I1), \(I8) and \(I9) correspond to any of the elementary diagrams
in Figure \ref{ribbsurf13/fig}.
\end{proposition}

\begin{Figure}[htb]{ribbsurf16/fig}{}{Graph changing relations in $\S$}
\centerline{\fig{Moves 3}{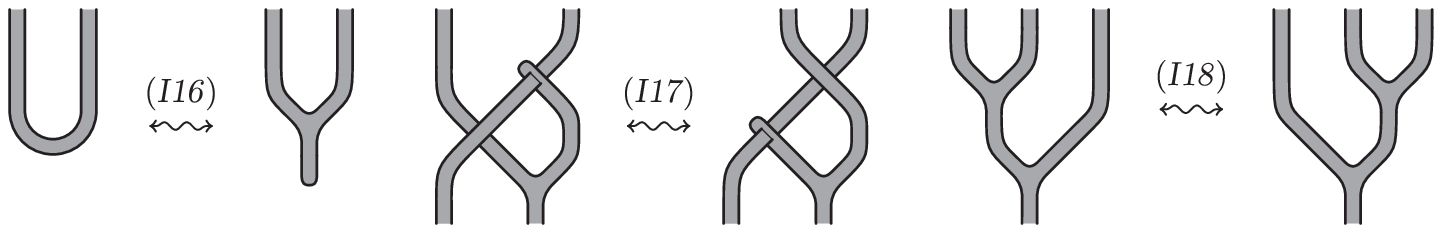}}\vskip-2pt
\end{Figure}

\begin{Figure}[htb]{ribbsurf17/fig}{}{1-isotopy relations in $\S$}
\centerline{\fig{Moves 4}{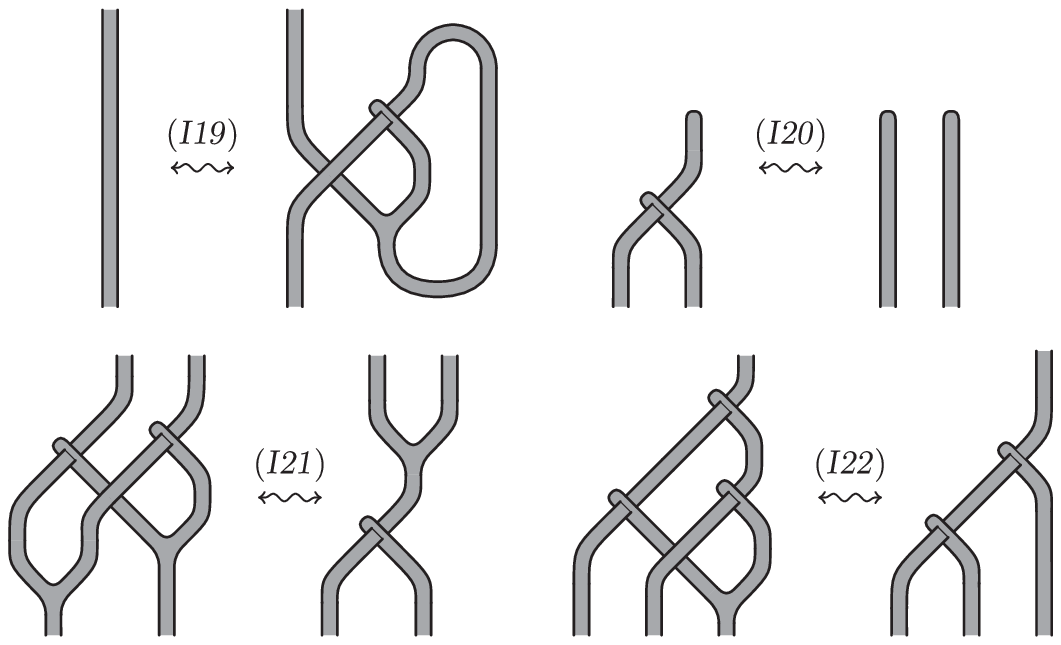}}\vskip-2pt
\end{Figure}

The proof of Proposition \ref{surfcat/theo} is technical but quite standard and the
details are presented in Appendix I (p. \pageref{appendix1/sec}). Here, we only make
few observations: moves \(I1) to \(I6) come from isotopy of (oriented with respect to
the $y$-axis) planar diagrams in the projection plane (but suffice for that only up to
the subsequent moves); moves \(I7) to \(I15) describe the changes under 3-dimensional
isotopies preserving the ribbon intersections; moves \(I16) to \(I18) describe the
changes which occur in the planar diagram when one changes the core graph $G$ of $F$;
moves \(I19) to \(I22) represent in terms of standard planar diagrams the 1-isotopy
moves in Figure \ref{ribbsurf5/fig}.

\begin{block}\label{autonomous}
Proposition \ref{surfcat/theo} implies that $\S$ is a braided monoidal category with
respect to juxtaposition of diagrams: $\one$ is again the empty set, the braiding
morphisms $\gamma_{I,I'}$ and $\gamma_{I,I'}^{-1}$ are presented in Figure
\ref{ribbsurf3/fig} \(f) and \(f'), while $\iota_A,\,{}_A\iota$ and the associativity
morphisms are the identities. In fact, as we will show, $\S$ carries also autonomous
and tortile structures. Recall from \cite{Sh94} that an {\sl autonomous} braided
monoidal category $(\C,\,\diam,\, \one,\, \iota_A,\,{}_A\iota,\, \alpha,\gamma,
\Lambda, \lambda)$ is a braided monoidal category $\C$ in which every object $A$ has a
right dual $(A^\ast, \Lambda_A, \lambda_A)$, where the morphisms
$$
\begin{array}{rl}
 \Lambda_A: &\one \to A^\ast \diam A \quad\text{({\sl coform}\/)}\\[2pt]
 \lambda_A: &A \diam A^\ast \to \one \quad\text{({\sl form}\/)}
\end{array}
$$
are such that the compositions
$$
\begin{array}{ccccccccccc}
 A &\mapright{{}_A\iota^{-1}} &A \diam \one &\mapright{\id \diam \Lambda_A}
 &A \diam (A^\ast \diam A) &\mapright{\alpha^{-1}} &(A \diam A^\ast)\diam A
 &\mapright{\lambda_A \diam \id} & \one\diam A &\mapright{\iota_A} &A,\\[4pt]
 A^\ast &\mapright{\iota^{-1}_A} &\one\diam A^\ast&\mapright{\Lambda_A \diam \id}
 &(A^\ast \diam A)\diam A^\ast &\mapright{\alpha} &A^\ast \diam (A\diam A^\ast)
 &\mapright{\id \diam \lambda_A} &A^\ast\diam \one &\mapright{{}_A\iota} &A^\ast,
\end{array}
$$
 are the identities. Then, given any morphism $F: A \to B$ in $\C$, we define its {\sl
dual} $F^\ast: B^\ast \to A^\ast$ as follows:
$$
F^\ast={}_{A}\iota^{-1}\circ(\id_{A^\ast}\diam\lambda_B)\circ
\alpha_{A^\ast,B,B^\ast}\circ((\id_{A^\ast}\diam F)\diamond \id_{B^\ast})
\circ (\Lambda_A\diam\id_{B^\ast})\circ\iota^{-1}_{B^\ast}.
$$

A {\sl twist} for a braided monoidal category is a family of natural isomorphisms
$\theta_A: A \to A$ such that $\theta_{\one} = \id_{\one}$ and
$\theta_{A\diam B} = \gamma_{B,A} \circ (\theta_B\diam\theta_A) \circ
\gamma_{A,B}$.

An autonomous braided monoidal category is called {\sl tortile} (the terminology
is from \cite{Sh94}) if it is equipped with a distinguished twist such that 
$\theta_{A^\ast} = (\theta_A)^\ast$ for any object $A$.

\begin{Figure}[htb]{dualtwist/fig}{}{}
\vskip-6pt\centerline{\fig{}{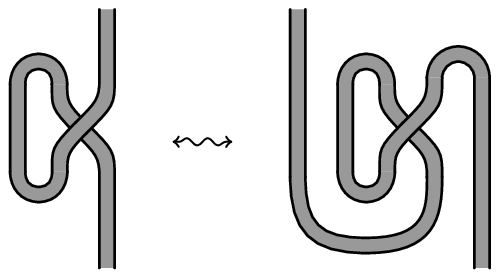}}
\end{Figure}

A classical example of tortile category is the category $\T$ of framed tangles.
We remind that  $\T$ has the same objects as $\S$, but the morphisms are compositions
of products of the elementary diagrams \(a), \(b), \(b'), \(f) and \(f') in Figure
\ref{ribbsurf3/fig} modulo the relations \(I1), \(I2-2'), \(I7-7'), \(I8), \(I9)
and the relation in Figure \ref{dualtwist/fig}. The autonomous structure on $\T$ is
determined by $A^\ast=A$ for any object $A$, while $\lambda_I$ and $\Lambda_I$ for
a single interval $I$ are represented by the diagrams \(b) and \(b') in Figure 
\ref{ribbsurf3/fig} and then extended to any other object by requiring that 
$\Lambda_{A \diam B} = (\id_{B}\diam\Lambda_A\diam\id_B) \circ \Lambda_B$ and
$\lambda_{A \diam B} = \lambda_A \circ (\id_A \diam \lambda_B \diam \id_{A})$.
Then the conditions on the form and coform in \ref{autonomous} reduce to relations
\(I2-2'). The twist is defined by $\theta_A =(\lambda_A \diam\id_A) \circ (\id_{A}
\diam \gamma_{A,A}) \circ (\Lambda_{A} \diam \id_A )$. In particular, $\theta_I$ is
presented on the left in Figure  \ref{dualtwist/fig} while the relation in the same
figure represents the condition $\theta_{I^\ast} = (\theta_I)^\ast$. The fact that
this condition is satisfied for any object $A$ and that $\theta_{A\diam B} =
\gamma_{B,A} \circ (\theta_B\diam\theta_A) \circ \gamma_{A,B}$, is a straightforward
application of the relations \(i7-9). 
\end{block}

\begin{proposition}
 $\S$ is a tortile category and the map which sends any tangle in $\T$ to the
corresponding rs-tangle in $\S$ defines a functor $\T\mapsto\S$ which is a strict
tortile functor.
\end{proposition}

\begin{proof} In order to show that the map is a functor we only need to observe
that the relation in Figure \ref{dualtwist/fig} is satisfied in $\S$. Indeed, this
follows easily by applying \(I11), \(I4) and \(I2). Now we define the form,
coform and the twist in $\S$ to be equal to the images of the correspoinding
morphisms in $\T$. Those would make $\S$ into a tortile category if the twist is
natural in $\S$, i.e if for any morphism $F: A \to B$ in $\S$, $\theta_B\circ
F = F \circ \theta_A$. This follows from the naturality of the twist in $\T$ and
relations \(I13-15).
\end{proof}

\begin{block}\label{labeledS/defn}
An \rs-tangle labeled in the symmetric group $\Sigma_n$, is an \rs-tangle in which
each boundary arc and each connected component of its 3-dimensional diagrams with
double points removed, is labeled by a transposition in $\Sigma_n$, according to the
following rules: the label of any boundary arc coincides with the label of the
component to which it belongs; at any ribbon intersection, the label of the ribbon
entering into a disk, gets conjugated by the label of this disk. 

Then the category $\S_n$ of {\sl labeled \rs-tangles}, is defined to have as objects
the finite sets $\{I_i,(k_i\;l_i)\}_i$ of intervals $I_i \subset R^2$ labeled by
transpositions $(k_i\;l_i) \in \Sigma_n$, and as morphisms the equivalence classes of
labeled \rs-tangles modulo the labeled version of the defining relations of $\S$ and
the relations \(R1) and \(R2) presented in Figure \ref{ribbon-m/fig}, where $i$, $j$,
$k$ and $l$ are all distinct.

A labeling of a planar diagram of an \rs-tangle induces a labeling of the planar
diagram of its core graph. We call a singular vertex of such graph uni-, bi- or
tri-colored according to how many different transpositions appear as labels of the
edges attached to it. Then move \(R1) in Figure \ref{ribbon-m/fig} states that a
tri-colored singular vertex carries a cyclic symmetry, while move \(R2) allows to
create or remove in $\S_n$ bi-colored singular vertices. 
 
$\S_n$ is also a monoidal tortile category, where the product corresponds to disjoint
union, the identity corresponds to the empty set and the tortile structure is induced
by that of $\S$.

We recall that a morphism in $\S_n$ is called complete if the transpositions occuring
as labels of any planar diagram of it generate all $\Sigma_n$. Notice that this
property does not depend on the choice of the planar diagram representation of the
morphism, being invariant under the defining relations in Figures
\ref{ribbsurf14/fig}--\ref{ribbsurf17/fig}.

As already mentioned in the Introduction, one can define a stabilization map
$\up_n^{n+1}: \hat\S_n^c \to \hat\S_{n+1}^c$ by sending a ribbon surface to its
disjoint union with a disk labeled $(n\;n+1)$. Moreover, according to Proposition 4.2
in \cite{BP}, $\up_n^{n+1}$ is invertible for $n > m \geq 3$ and its inverse
$\down_n^{n+1}: \hat\S^c_{n+1} \to \hat\S^c_n$ is called a reduction map.
\end{block}

 Finally, we define $\partial^\star \S_n$ to be quotient of $\S_n$ modulo the relation
\(T) in Figure \ref{defn-boundS/fig}, and $\partial \S_n$ to be quotient of
$\partial^\star \S_n$ modulo the relations \(P) and \(P') presented in the same figure.

\vskip-\lastskip
\vskip0pt

\section{From labeled ribbon surfaces to Kirby diagrams%
\label{surf-kirby/sec}}

In this section we recall from Section 2 of \cite{BP}, the definition of the bijective
map $F \mapsto K_F$ sending any labeled ribbon surface $F$ to a generalized Kirby
diagram $K_F$ in $\K_n$.

The main idea is that an 1-handlebody structure on the labeled ribbon surface $F$,
representing a $n$-fold simple branched covering $p:M \to B^4$, naturally induces a
2-handlebody structure on $M$ defined up to 2-deformations. In particular, the
construction of $K_F$ requires to make some choices: 
\begin{itemize}\itemsep\smallskipamount
\item[{\sl a}\/)]\vskip-\lastskip\smallskip
a choice of an adapted 1-handlebody structure on $F$, i.e. a decomposition 
$$F = D_1 \cup \dots \cup D_p \cup B_1 \cup \dots \cup B_q,$$ 
where the $D_h$'s are disjoint flat disks (the 0-handles of $F$) and the $B_h$'s are
disjoint bands attached to $F_0 = D_1 \cup \dots \cup D_p$ (the 1-handles of $F$);
in particular, any ribbon intersection occurs between a band and a disk;
\item[{\sl b}\/)] a choice of orientations of the disks $D_h$;
\item[{\sl c}\/)] for every ribbon intersection contained in a disk $D_h$, a
choice of an arc $\alpha$ in $D_h$, joining the ribbon intersection with $\Bd
D_h\cap \Bd F$; such arcs are taken disjoint from each other.
\end{itemize}\vskip-\lastskip\medskip

Then the following steps summarize the procedure defined in \cite{BP} for constructing
$K_F$ in the case when $F$ is presented by a planar diagram. In particular, in this
case a choice of an orientation of the projection plane induces an orientation of all
disks (since they lie in this plane):
\begin{itemize}\itemsep\smallskipamount
\item[{\sl a}\/)] \vskip-\lastskip\smallskip
replace any disk $D_h$ by a dotted unknot coinciding with $\Bd D_h$; if $D_h$ is
labeled $(i\;j) \in \Sigma_n$ with $i < j$, then assign to its upper side the label $i$
and to its lower side, the label $j$;
\item[{\sl b}\/)] cut any band $B_h$ at the ribbon intersections by removing a small
neighborhood of the intersection arc; then take two small positive and negative
vertical shifts of the resulting pieces of $B_h$, disjoint from the original surface;
if a piece is labeled $(i\;j) \in \Sigma_n$ with $i < j$, then label by $i$ its upper
shift and by $j$ its lower shift; Figure \ref{defnKF5/fig} ({\sl a}\/) shows how the
two shifts look like in a neighborhood of a positive half-twist (in this figure, as
well as in the next one, the projection plane is depicted from a perspective point of
view);
\item[{\sl c}\/)] connect the two shifts of any attaching arc of a band $B_h$ by 
a band which forms a ribbon intersection with the attaching disk, as shown in Figure
\ref{defnKF5/fig} ({\sl b}\/);
\begin{Figure}[htb]{defnKF5/fig}{}{($i < j$)}
\centerline{\fig{}{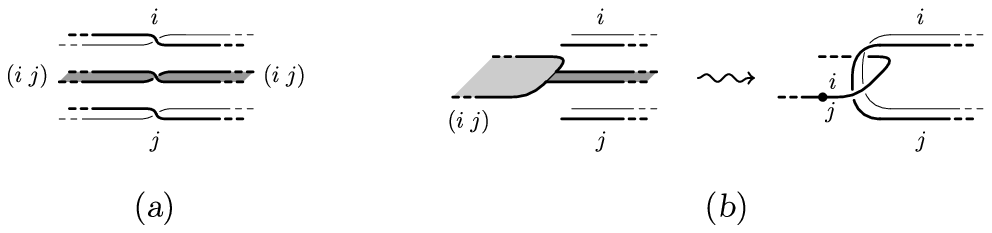}}\vskip-4pt
\end{Figure}
\begin{Figure}[b]{defnKF6/fig}{}{($i < j < k$)}
\vskip-4pt\centerline{\fig{}{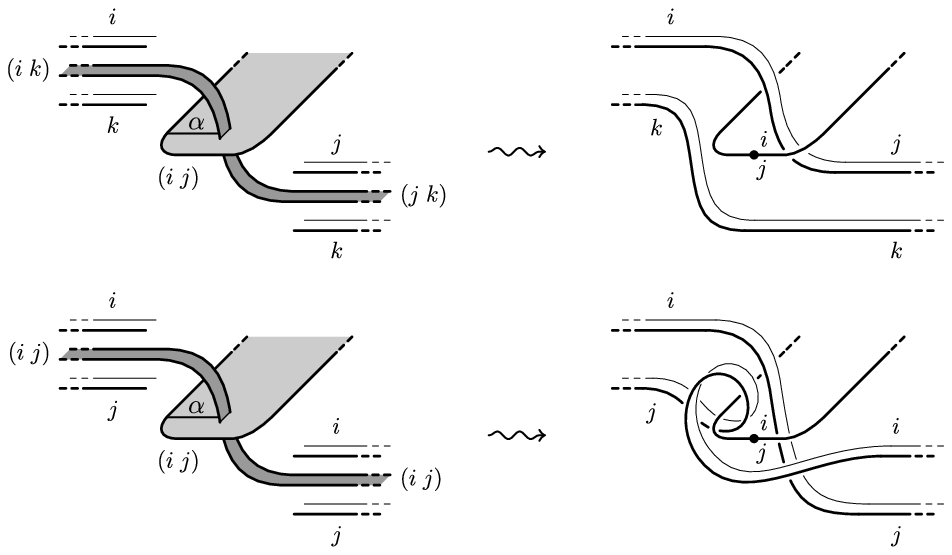}}\vskip-4pt
\end{Figure}
\item[{\sl d}\/)] for each ribbon intersection of $F$ contained in a disk $D_h$,
connect in pairs the shifts approaching $D_h$ from opposite sides, according to the
following rule: if two shifts have the same label and this is different from the ones
of $D_h$, then connect them with a flat band disjoint from $D_h$; if two shifts have
different labels coinciding with the ones of $D_h$, then connect them by a band which
forms a ribbon intersection with $D_h$ and is flat except for a twist around $\Bd D_h$
in the case when the two shifts approach $D_h$ from the wrong sides with respect to
the labeling; the whole construction is carried out in a regular neighborhood of the
arc $\alpha$ and Figure \ref{defnKF6/fig} shows two special cases.
\end{itemize}\vskip-\lastskip\smallskip

Lemma 2.3 and Proposition 2.4 in \cite{BP} state that $K_F$ is well defined up to
2-deformation, i.e. making different choices of 1-handlebody structure,
orientations or arcs, as well as changing $F$ by 1-isotopy or ribbon moves, changes
$K_F$ only by 2-deformation. For more details on this point, as well as for the
proof that the restriction of the map to the complete morphisms is bijective, we
refer the interested reader to \cite{BP}, observing that our final result relies on
the bijectivity of the map, but does not need the explicit form of the inverse, which
can be found in \cite{BP}.

\begin{Figure}[b]{defnKF1/fig}{}{Definition of $K_F$ ($i < j$ and $i' < j'$).}
\centerline{\fig{}{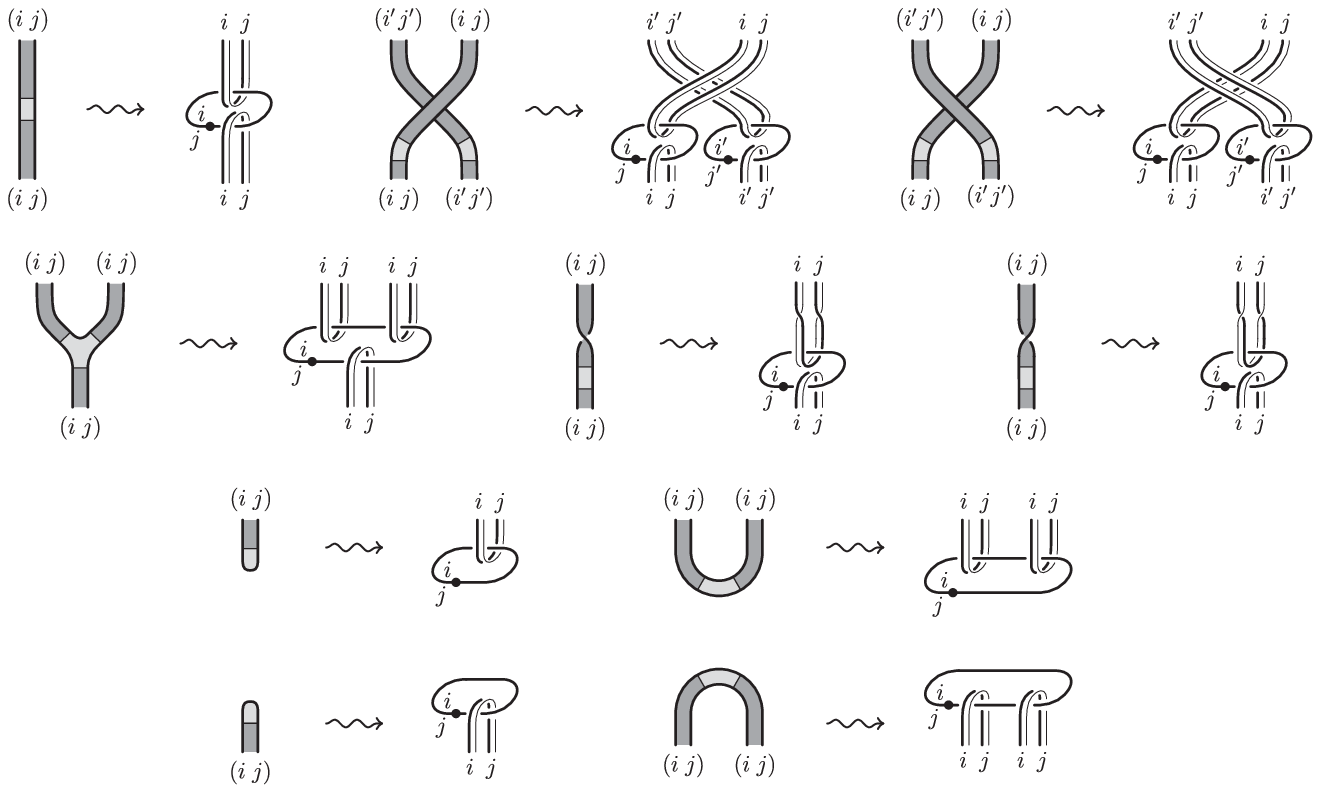}}\vskip-3pt
\end{Figure}

\begin{Figure}[t]{defnKF2/fig}{}
 {Definition of $K_F$ ($i < j < k$, $h < l$ and $\{i,j\}\cap\{h,l\} = \emptyset$).}
\vskip6pt\centerline{\fig{}{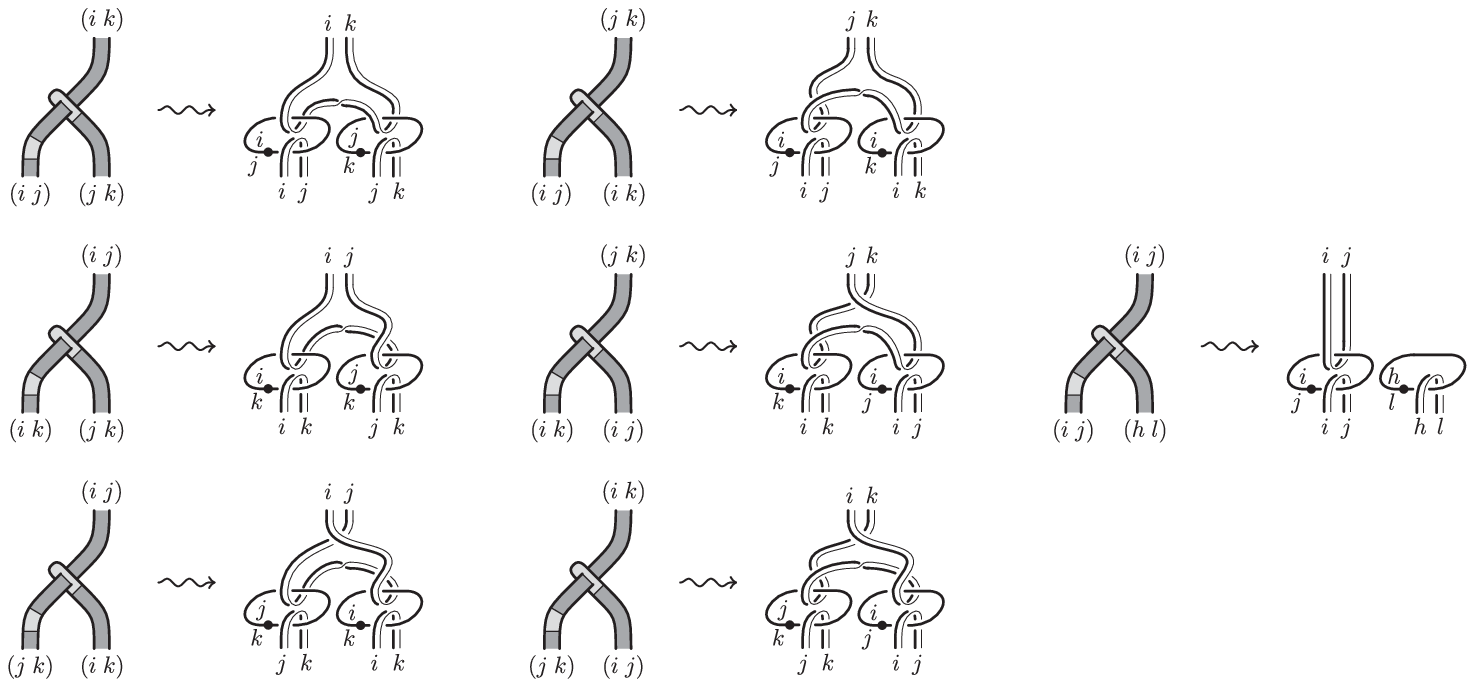}}\vskip-3pt
\end{Figure}

\begin{Figure}[htb]{defnKF3/fig}{}{Definition of $ K_F$ ($i < j$).}
\centerline{\fig{}{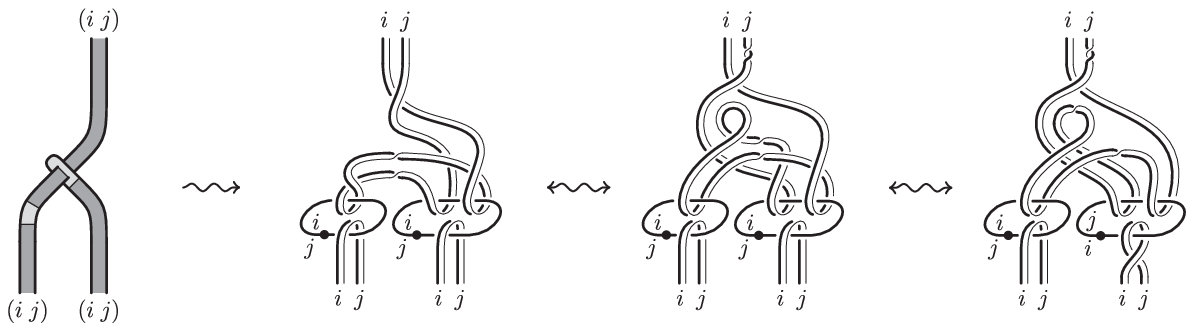}}\vskip-8pt
\end{Figure}

Instead, given any labeled ribbon surface $F\in\hat\S_n^c$, we will need an explicit
form for $K_F$, which requires a specific choice of adapted 1-handlebody structure for
any such surface. Since the morphisms in $\hat\S_n^c$ are compositions of the
elementary morphisms presented by the planar diagrams in Figure \ref{ribbsurf3/fig},
we can define an 1-handlebody structure on each morphism by indicating its
intersection with the elementary morphisms. This is done in Figures \ref{defnKF1/fig},
\ref{defnKF2/fig} and \ref{defnKF3/fig}, where the disks are denoted with lighter gray
and the bands with heavier gray color. Observe that there are two types of disks: the
ones which deformation retract onto neighborhoods of the vertices of the core graph,
and the others which divide the ribbons in such a way that none of them contains two
boundary arcs. Moreover, since any band forms at most one ribbon intersection with any
disk $D_h$, and in this case $D_h$ has a single band attached to it, the choice of the 
arcs $\alpha$ is unique up to isotopy, so we omit it.

The application of the construction of $K_F$ to the particular 1-handlebody structure
of $F$ chosen above is shown in Figures \ref{defnKF1/fig}, \ref{defnKF2/fig} and
\ref{defnKF3/fig}, where for later use the image  of the uni-colored singular vertex in
Figure \ref{defnKF3/fig} has been transformed through the isotopy move presented in
Figure \ref{defnKF4/fig}. In particular, given $F\in\hat\S_n^c$ as a composition of
elementary diagrams, $K_F$ is the formal composition of the corresponding generalized
Kirby tangles on the right in those figures. Observe that the results in \cite{BP}
assure us that such composition is well defined morphism in $\hat{K}_n^c$. 
 By proving Theorem \ref{eq-alg-kirby/theo}, we will prove that the map $\S_n\to\K_n$
defined in Figures \ref{defnKF1/fig}, \ref{defnKF2/fig} and \ref{defnKF3/fig} is a
composition of two braided monoidal functors: $\Psi_n: \S_n \to \H_n$ and $\Phi_n: \H_n
\to \K_n$, and therefore it is a braided monoidal functor itself.

\begin{Figure}[htb]{defnKF4/fig}{}{Reversing a dotted unknot in $\K_n$ ($i \neq j$)}
\vskip4pt\centerline{\fig{}{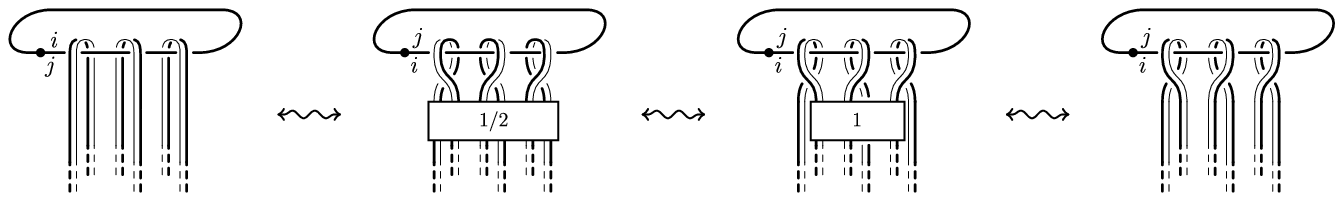}}
\end{Figure}

\vskip-12pt\vskip0pt

\section{The universal braided Hopf algebras $\H(\G)$ and $\H^u(\G)$%
\label{hopf-alg/section}}

\begin{block}\label{groupoid}
 We remind that a groupoid $\G$ is a small category in which every morphism is
invertible. Hopefully without causing a confusion, we will indicate by $\G$ also the
set of the morphisms of $\G$, endowed with the partial binary operation given by the
composition, for which we adopt the multiplicative notation from left to right. The
identity of $i \in \Obj\G$ will be denoted by $1_i \in \G$, while the inverse of $g
\in \G$ will be denoted by $\bar g \in \G$. Moreover, given two objects $i,j \in
\Obj\G$, we denote by $\G(i,j) \subset \G$ the subset of morphisms from $i$ to $j$.
Consequently, if $g \in \G(i,j)$ and $h \in \G(j,k)$ then $gh \in \G(i,k)$. In
particular, $g \bar g = 1_i$ and $\bar g g = 1_j$, and sometimes the identity morphisms
will be represented in this way. A groupoid is called {\sl connected} if $\G(i,j)$  is
non-empty for any $i,j \in \Obj\G$. Given two groupoids $\G \subset \G'$, we say that
$\G$ is {\sl full} in $\G'$ if $\G(i,j) = \G'(i,j)$ for all $i,j \in \Obj \G$.

\medskip
We have the following basic examples of (connected) groupoids:
\begin{itemize}\itemsep\smallskipamount
\item[{\sl a}\/)] \vskip-\lastskip\smallskip
any group $G$, considered as a groupoid $\G$ with a single object $\one$, i.e.
$\Obj\G = \{\one\}$ and $\G(\one,\one) = G$;
\item[{\sl b}\/)] 
the groupoid $\G_n$ with $n \in {\Bbb N}$, such that $\Obj\G_n = \{1,2,\dots,n\}$ and
$\G_n(i,j)$ consists of a unique arrow $(i,j):i \to j$ for any $i,j \in \Obj \G$;
in particular, the composition is given by $(i,j)(j,k) = (i,k)$ and ${\bar{(i,j)\!}\,}
= (j,i)$.
\end{itemize}\vskip-\lastskip\medskip

The $\G_n$'s are exactly the groupoids which we will need later.
\end{block}\bigskip

\begin{block}\label{defnh}
 For a groupoid $\G$ and $x \in \G(i,k)$, define the functor $\_^x: \G \to \G$ as
follows: 
\begin{itemize}\itemsep\smallskipamount
\item[{\sl a}\/)] \vskip-\lastskip\smallskip
$k^x = i$ and $j^x = j$ if $j \in \Obj\G - \{k\}$;
\item[{\sl b}\/)]\strut\vskip-20pt
$$g^x = \left\{
\arraycolsep0pt
\begin{array}{ll}
	xg \bar x &\quad\text{if } g \in \G(k,k),\\
 xg &\quad\text{if } g \in \G(k,l) \text{ with } l \neq k,\\
	g \bar x &\quad\text{if } g \in \G(j,k) \text{ with } j \neq k,\\
 g &\quad\text{if } g \in \G(j,l) \text{ with } j,l \neq k.
\end{array}
\right.
$$
\end{itemize}

Then the following statements hold:
\begin{itemize}\itemsep\smallskipamount
\item[{\sl c}\/)] \vskip-\lastskip\smallskip
 $\bar{g^x} = \bar g^x$ for any $g \in \G$;
\item[{\sl d}\/)] if $k\neq i$ then the image of the functor $\_^x$ is the full
subgroupoid $\G^{\bs k}$ of $\G$ with $\Obj \G^{\bs k} = \Obj\G - \{k\}$ (we use the
notation $\G^{\bs k}$ instead of $\G^x$, to emphasize that this groupoid depends only
on the target of $x$); 
\item[{\sl e}\/)] $\_^x$ restricts to the identity on $\G^{\bs k}$ and to an
equivalence of categories $\G^{\bs i} \to \G^{\bs k}$, whose inverse $\G^{\bs k} \to
\G^{\bs i}$ is the corresponding restriction of $\_^{\bar x}$;
\item[{\sl f}\/)] for any $x \in \G(i,k)$ and $y \in \G(j,k)$, there exists a natural
transformation\break $N:\_^x \to \_^y$ such that $N(k) = x \bar y$ and $N(l) = 1_l$
if $l \neq k$.
\end{itemize}\vskip-20pt
\end{block}

\begin{block}\label{hcat}
 Given a groupoid $\G$ and a braided monoidal category $\C$, a {\sl Hopf $\G$-algebra}
in $\C$ is a family of objects $H = \{H_g\}_{g \in \G}$ in $\C$, equipped with the
following families of morphisms in $\C$ (here and in the sequel we will often write 
$g$ instead of $H_g$. In particular we will use the notations $\id_g = \id_{H_g}$,
$\gamma_{g,h} = \gamma_{H_g,H_h}$, ${}_g\iota = {}_{H_g}\iota$ and $\iota_g =
\iota_{H_g}$; moreover, based on the MacLane's coherence result for monoidal categories
(p. 161 in \cite{McL}), we will omit the associativity morphisms since they can be
filled in a unique way):

\medskip\noindent
 a {\sl comultiplication} $\Delta = \{\Delta_g: H_g \to H_g \diam H_g\}_{g \in \G}$,
such that for any $g \in \G$
\vskip-4pt
$$ (\Delta_g\diam\id_g) \circ \Delta_g = (\id_g\diam\Delta_g) \circ \Delta_g;
\eqno{\(a1)}$$
\smallskip\noindent
 a {\sl counit} $\epsilon = \{\epsilon_g: H_g \to \one\}_{g \in \G}$, such that
 for any $g \in \G$
\vskip-4pt
$$(\epsilon_g\diam\id_g)\circ\Delta_g = \id_g = (\id_g\diam\epsilon_g)\circ\Delta_g;
\eqno{\(a2-2')}$$
\smallskip\noindent
 a {\sl multiplication} $m = \{m_{g,h}: H_g \diam H_h \to H_{gh}\}_{g,h,gh \in
\G}$ ($m_{g,h}$ is defined when $g$ and $h$ are composable in $\G$), such that if
$f,g,h,fgh \in \G$ then
\vskip-3pt
$$m_{fg,h} \circ (m_{f,g} \diam \id_h) = m_{f,gh} \circ (\id_f \diam m_{g,h}),
\eqno{\(a3)}$$
$$(m_{g,h} \diam m_{g,h}) \circ (\id_g \diam \gamma_{g,h}\diam \id_h) \circ
(\Delta_g \diam \Delta_h) = \Delta_{gh} \circ m_{g,h},
\eqno{\(a5)}$$
$$\epsilon_{gh} \circ m_{g,h} = \epsilon_g \diam \epsilon_h;                  
\eqno{\(a6)}$$
\smallskip\noindent
 a {\sl unit} $\eta = \{\eta_i: \one \to H_{1_i}\}_{i \in \Obj \G}$, such
that if $g \in \G(i,j)$ then
\vskip-3pt
$$m_{g,{1_j}}\circ(\id_g\diam\eta_j) = \id_g = m_{{1_i},g}\circ(\eta_i\diam\id_g),
\eqno{\(a4-4')}$$
$$\Delta_{1_i} \circ \eta_i = \eta_i \diam \eta_i,
\eqno{\(a7)}$$
$$\epsilon_{1_i} \circ \eta_i = \id_{\one};
\eqno{\(a8)}$$
\smallskip\noindent
 an {\sl antipode} $S = \{S_g: H_g \to H_{\bar g}\}_{g \in \G}$ and its inverse $\bar 
S = \{\bar S_g: H_g \to H_{\bar g}\}_{g \in \G}$, such that if $g\in \G(i,j)$ then
\vskip-9pt
$$m_{\bar g,g}\circ(S_g\diam\id_g)\circ\Delta_g = \eta_{1_j}\circ\epsilon_g,
\eqno{\(s1)}$$
$$m_{g,\bar g}\circ(\id_g\diam S_g)\circ\Delta_g = \eta_{1_i} \circ \epsilon_g,
\eqno{\(s1')}$$
$$S_{\bar g}\circ\bar S_g =\bar S_{\bar g}\circ S_g=\id_g.
\eqno{\(s2-2')}$$

\smallskip
The definition above is a straightforward generalization of a categorical group
Hopf algebra given in \cite{V01}. Observe that an ordinary braided Hopf algebra in
$\C$ is a Hopf $\G = \G_1$-algebra, where $\G_1$ is the trivial groupoid with a
single object and a single morphism. In particular, $H_{1_i}$ is a braided Hopf
algebra in $\C$ for any $i \in \Obj \G$.
\end{block}

\vskip-12pt\vskip0pt

\begin{block}\label{integral}
 Let $\C$ be a braided monoidal category and $H = \{H_g\}_{g \in \G}$ be a Hopf
$\G$-algebra in $\C$. By a categorical {\sl left} (resp. {\sl right}\/) 
{\sl cointegral} of $H$ we mean a family $l = \{l_i: H_{1_i} \to \one \}_{i \in \Obj
\G}$ of morphisms in $\C$, such that for any $i \in \Obj \G$
$$
\begin{array}{c}
 (\id_{1_i} \diam l_i) \circ \Delta_{1_i} = \eta_i \circ l_i: 
 H_{1_i} \to H_{1_i}\\
 (\text{resp. } (l_i \diam \id_{1_i}) \circ \Delta_{1_i} = \eta_i \circ l_i: 
 H_{1_i} \to H_{1_i}).
\end{array} 
\eqno{\(i1-1')}
$$
\smallskip\noindent
 On the other hand, by categorical {\sl right} (resp. {\sl left}\/) {\sl integral} of
$H$ we mean a family $L = \{L_g: \one \to H_g\}_{g \in \G}$ of morphisms in $\C$, such
that if $g,h,gh \in \G$ then
$$
\begin{array}{c}
 m_{g,h} \circ (L_g \diam \id_h) = L_{gh} \circ \epsilon_h: 
 H_h \to H_{gh}\\
 (\text{resp. } m_{g,h} \circ (\id_g \diam L_h) = L_{gh} \circ \epsilon_h: 
 H_g \to H_{gh}).
\end{array}
\eqno{\(i2-2')}
$$
\smallskip\noindent
 If $l$ (resp. $L$) is both right and left categorical cointegral (integral) of $H$, 
 we call it simply a cointegral (integral) of $H$.
\end{block}

\vskip3pt\vskip0pt

\begin{block}
 Given a groupoid $\G$, let $\H(\G)$ be the free braided monoidal category with a
Hopf $\G$-algebra $H = \{H_g\}_{g \in \G}$ in it and with a left cointegral $l$ and a
right integral $L$ of $H$, such that for any $i \in\Obj\G$ 
$$l_i \circ L_{1_i} = \id_{\one} = l_i \circ S_{1_i} \circ L_{1_i}.
\eqno{\(i3-3')}$$

We refer to $\H(\G)$ as the {\sl universal Hopf $\G$-algebra}. In particular, if
$\C$ is any braided monoidal category with a Hopf $\G$-algebra $A = \{A_g\}_{g \in
\G}$ in it which has a left cointegral and a right integral satisfying the condition
above, then there exists a braided monoidal functor $\H(\G)\to\C$ sending $H_g$ to
$A_g$.

\begin{Figure}[b]{elemdiag/fig}{}{Elementary diagrams in $\H(\G)$.}
\vskip-6pt\centerline{\fig{}{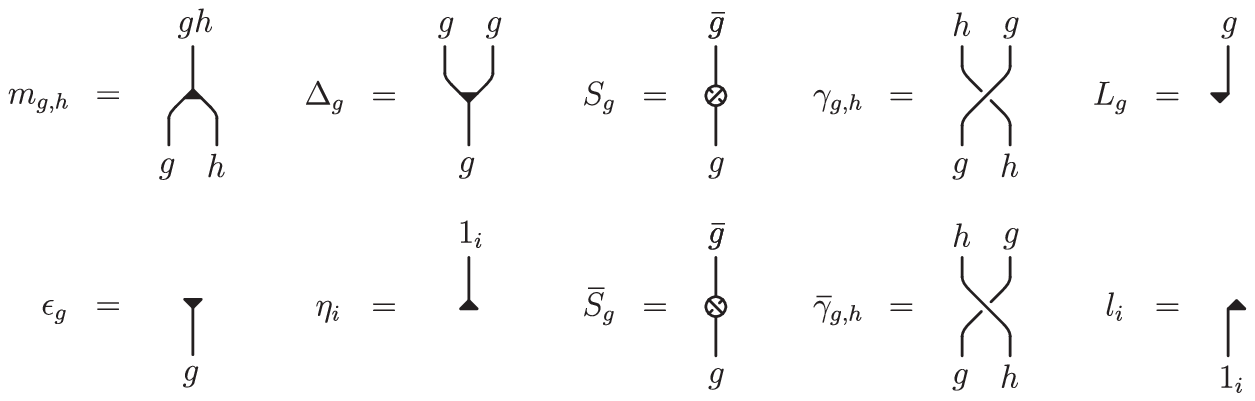}}\vskip-6pt
\end{Figure}

\begin{Figure}[htb]{isot/fig}{}{Braid axioms for $\H(\G)$.}
\vskip4pt\centerline{\fig{}{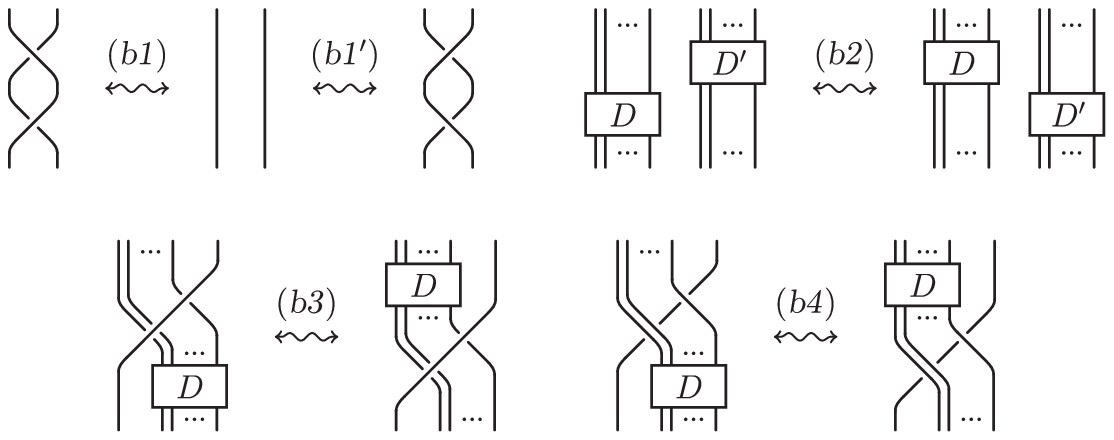}}\vskip-3pt
\end{Figure}

\begin{Figure}[htb]{algebra/fig}{}{Bi-algebra axioms for $\H(\G)$.}
\vskip4pt\centerline{\fig{}{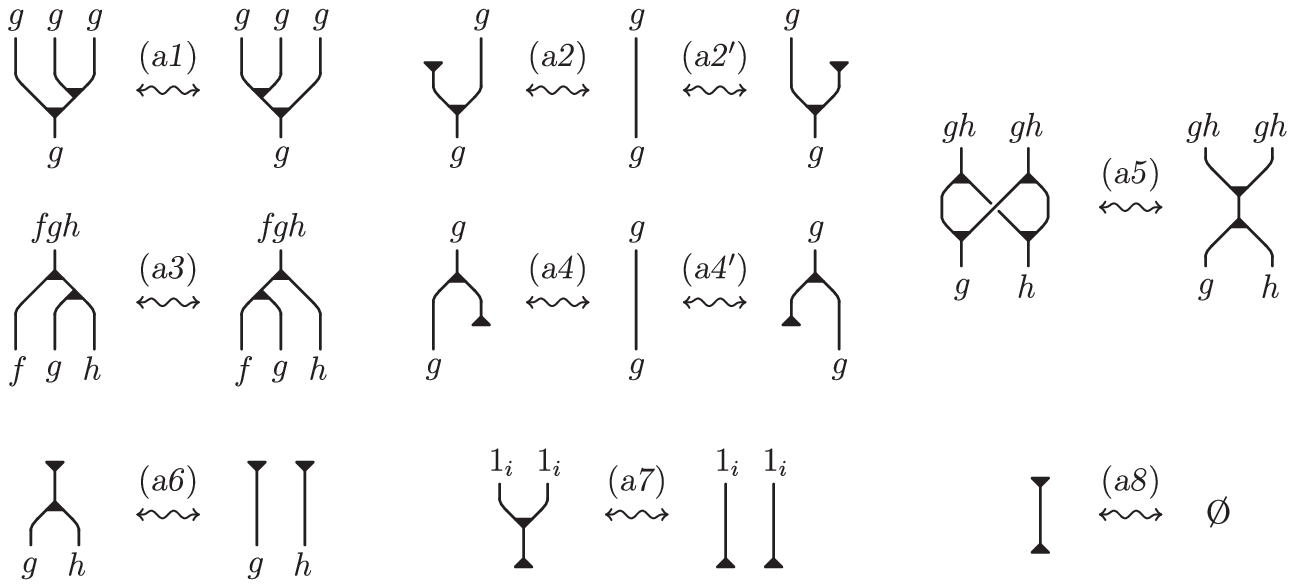}}\vskip-3pt
\end{Figure}

\begin{Figure}[t]{antipode/fig}{}{Antipode axioms for $\H(\G)$ ($g \in \G(i,j)$).}
\centerline{\fig{}{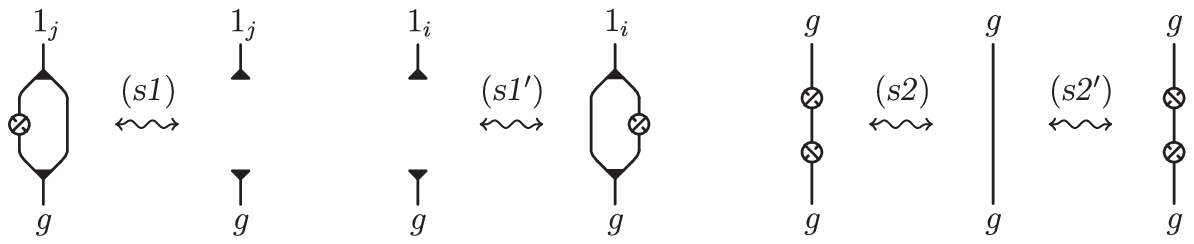}}\vskip-6pt
\end{Figure}

\begin{Figure}[htb]{integrals/fig}{}{Integral axioms for $\H(\G)$.}
\vskip4pt\centerline{\fig{}{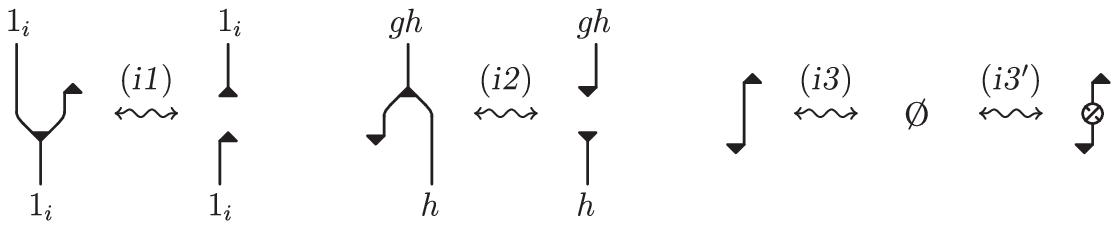}}\vskip-6pt
\end{Figure}

Analogously to \cite{Ke02}, $\H(\G)$ can be described as a category of planar
diagrams in $R\times [0,1]$. The objects of $\H(\G)$ are sequences of points in $R$
labeled by elements in $\G$, and the morphisms are iterated products and
compositions of the elementary diagrams presented in Figure \ref{elemdiag/fig},
modulo the relations presented in Figures \ref{isot/fig}, \ref{algebra/fig},
\ref{antipode/fig} and \ref{integrals/fig} and plane isotopies which preserve the
$y$-coordinate. We remind that the composition of diagrams $D_2 \circ D_1$
is obtained by stacking $D_2$ on the  top of $D_1$ and then rescaling, while the
product $D_1 \diam D_2$ is given by the horizontal juxtaposition of $D_1$ and $D_2$.

The plane diagrams and the relations between them listed above are going to be our main
tool, so it would be useful to stop and reflect, trying to catch some general rules
which would allow faster reading of the somewhat complicated compositions which will
appear later. The first observation is that the plane diagrams used represent
projections of embedded graphs in $R^3$ with uni-, bi- and tri-valent  vertices and the
projections of the edges do not contain local extremis. The vertices correspond to the
defining morphisms in the algebra, and often we will call them with the name of the
corresponding morphism. For example the bi-valent vertices (which have one incoming
and one outgoing edge) will also be called antipode vertices. The uni-valent vertices
are divided in unit vertices (corresponding to $\eta$ and $\epsilon$) and integral
vertices (corresponding to $l$ and $\Lambda$). 
	
The notation for the vertices reflects their interaction in the relations and later
their interpretation in terms of Kirby tangles. In particular, the uni- and tri-valent
vertices are represented by small triangles which may point up (positively polarized)
or down (negatively polarized). Then relations \(a8) and \(i3) can be summarized by
saying that two uni-valent vertices of opposite  polarizations, connected by an edge,
annihilate each other. Analogously, relations \(a6), \(a7), \(i1) and \(i2) can be
summarized by saying that a uni-valent vertex connected to a tri-valent vertex of
opposite polarization, cancels this last one, creating two uni-valent vertices of the
same polarization as its own. 
	
Relations \(a2) and \(a4) can be put together as well, by saying that if a uni-valent
vertex is connected to a tri-valent vertex of the same polarization, we can delete both
vertices and the edge between them, fusing the remaining two edges into a single one.
But we point out that there are two other possibilities to connect a uni- and
tri-valent vertices of the same polarization, which do not appear in the relations
above. Later, in Lemma \ref{lemma/auto} (cf. Figure \ref{form/fig}), we will see that
the diagrammatic language can be generalized to extend the statement to those two cases
as well, but a ``correction'' appears.
	
The remaining three relations \(a1), \(a3) and \(a5) in Figure \ref{algebra/fig},
concern the diagrams in which two tri-valent vertices are connected by a single edge.
Observe that if the vertices have the same polarization, we can slide one through the
other, while if they have opposite polarization, the diagram splits as shown in \(a5).
\end{block}	

The proofs of most of our theorems require to show that some morphisms in the universal
algebra are equivalent, meaning that the graph diagram of one of them can be obtained
from the graph diagram of the other by applying a sequence of the defining relations
(moves) in the algebra.  We will outline the main steps in this procedure by drawing
in sequence the intermedate diagrams, and for each step we will indicate in the
corresponding order, the main moves needed to transform the diagram on the left into
the one on the right. Actually, some steps can be understood more easily by starting
from the diagram on the right and reading the moves in the reverse order. Notice, that
the moves represent equivalences of diagrams and we use the same notation for them and
their inverses. In the captions of the figures the reader will find (in square
brackets) the reference to the pages where those moves are defined. For example, in
the first step in Figure \ref{pr-antipode-pf/fig} we obtain the diagram on the left
from the one on the right by first applying moves \(a1-3) in Figure \ref{algebra/fig}
on p. \pageref{algebra/fig}, then move \(s1) in Figure \ref{antipode/fig} on p.
\pageref{antipode/fig} and finally moves \(a2-4') in Figure \ref{algebra/fig} on p.
\pageref{algebra/fig}. To be precise, before applying \(s1) and \(a4'), we use the
braid axioms in Figure \ref{isot/fig}, but this in general will not be indicated.

\begin{Figure}[htb]{pr-antipode-pf/fig}{}
 {Proof of Lemma \ref{lemma/auto} 
 [{\sl a}/\pageref{algebra/fig}, 
 {\sl s}/\pageref{antipode/fig}].}
\centerline{\fig{}{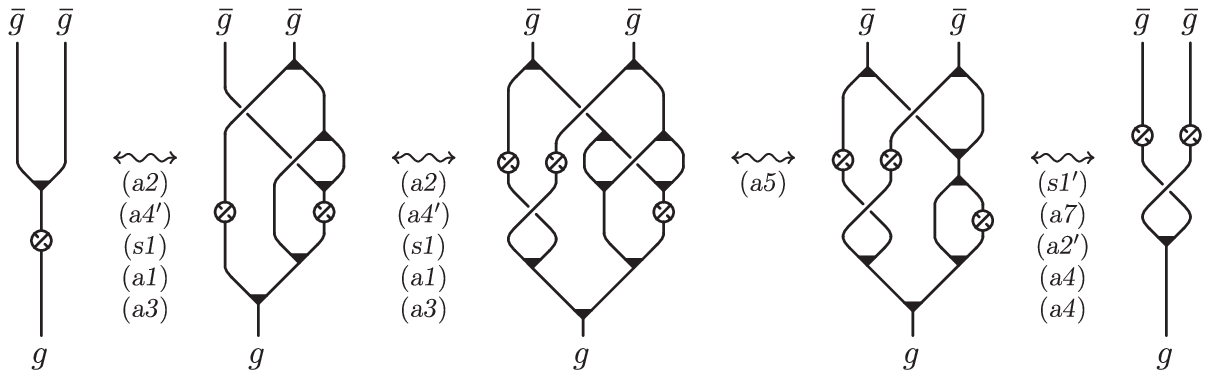}}\vskip-4pt
\end{Figure}

 In Figure \ref{pr-antipode-pf/fig} it is proved the first of the following properties
of the antipode in $\H(\G)$, which hold for any $g,h \in \G$ such that $gh$ is defined
and for any $i \in \G$:
\vskip-3pt
{$$\Delta_{\bar g} \circ S_g = (S_g \diam S_g) \circ \gamma_{g,g} \circ
\Delta_g: H_g
\to H_{\bar g} \diam H_{\bar g},
\eqno{\(s3)}$$
$$S_{gh} \circ m_{g,h} = m_{\bar h,\bar g} \circ (S_h \diam S_g) \circ \gamma_{g,h}:
H_g \diam H_h \to H_{\bar h \bar g},
\eqno{\(s4)}$$
$$S_{1_i} \circ \eta_i = \eta_i,
\eqno{\(s5)}$$
$$\epsilon_{\bar g} \circ S_g = \epsilon_g.
\eqno{\(s6)}$$}
\vskip-9pt\noindent

\begin{Figure}[b]{pr-antipode/fig}{}{Properties of the antipode in $\H(\G)$.}
\centerline{\fig{}{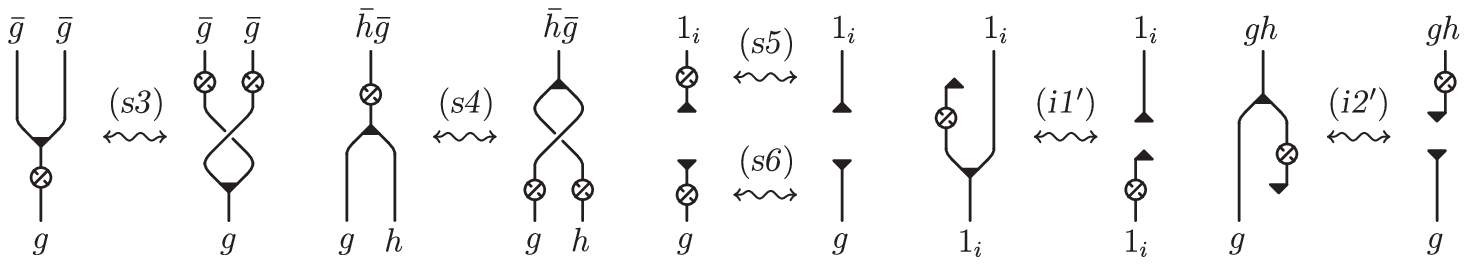}}\vskip-6pt
\end{Figure}

In terms of diagrams these relations are presented in Figure \ref{pr-antipode/fig}.
They are analogous to the case of braided and group Hopf algebras (see \cite{Ke02} and
\cite{V01}) and the rest of the proofs are left as an exercise. 

Using these properties one immediately sees that if $l = \{l_j: H_{1_j} \to \one \}_{j
\in \Obj \G}$ is a left cointegral of $H$, then (cf. Figure
\ref{pr-antipode/fig}): 
$$
l \circ S = \{l_j \circ S_{1_j}: H_{1_j} \to \one \}_{j \in \Obj \G}
 \text{ is a right cointegral of $H$.} \eqno{\(i1')}
$$
Analogously, if $L = \{L_g: \one \to
H_g\}_{g \in \G}$ is a right integral of $H$, then (cf. Figure
\ref{pr-antipode/fig} again):
$$
S \circ L =
\{S_g \circ L_g: \one \to H_{\bar g}\}_{g \in \G}\text{ is a left
integral of $H$.} \eqno{\(i2')}
$$

\smallskip

The next lemma extends Lemma 7 of \cite{Ke02} to possibly non-unimodular categories.

\begin{lemma} \label{lemma/auto}
$\H(\G)$ is an autonomous category such that for every $g \in \G$
(cf. Figure \ref{form/fig}):
\vskip-9pt
$$H_g^\ast = H_g\text{,}$$
$$\Lambda_{H_g} = \Lambda_g = \Delta_g \circ L_g\text{,}
\eqno{\(f1)}$$
$$\lambda_{H_g} = \lambda_g = l_{g \bar g} \circ m_{g,\bar g}\circ(\id_g \diam S_g).
\eqno{\(f2)}$$\vskip-\lastskip
\end{lemma}

\vskip-\lastskip
\vskip0pt

\begin{proof}[Lemma \ref{lemma/auto}]
 Figure \ref{pf-auto/fig} shows that $\Lambda_{H_g}$ and $\lambda_{H_g}$ satisfy the
relations in Paragraph \ref{autonomous}. Then it suffices to observe that such
relations propagate to $\Lambda_A$ and $\lambda_A$ for any object $A$ in $\H(\G)$,
where $\Lambda_A$ and $\lambda_A$ are defined inductively by the following identities
(observe that the definition is well-posed, giving equivalent results for different
decompositions $A \diam B = A' \diam B'$):
\vskip-6pt
$$(A \diam B)^\ast = B^\ast \diam A^\ast,$$
$$\Lambda_{A \diam B} = (\id_{B^\ast}\diam\Lambda_A\diam\id_B) \circ \Lambda_B,$$
$$\lambda_{A \diam B} = \lambda_A \circ (\id_A \diam \lambda_B \diam \id_{A^\ast}).
\hbox to 0pt{$\quad\square$\hss}$$\phantom{\proved}
\end{proof}

\begin{Figure}[htb]{pf-auto/fig}{}
 {Proof of Lemma \ref{lemma/auto} 
 [{\sl a}/\pageref{algebra/fig}, 
 {\sl i}/\pageref{integrals/fig}-\pageref{pr-antipode/fig},
 {\sl s}/\pageref{antipode/fig}-\pageref{pr-antipode/fig}].}
\vskip-20pt\centerline{\fig{}{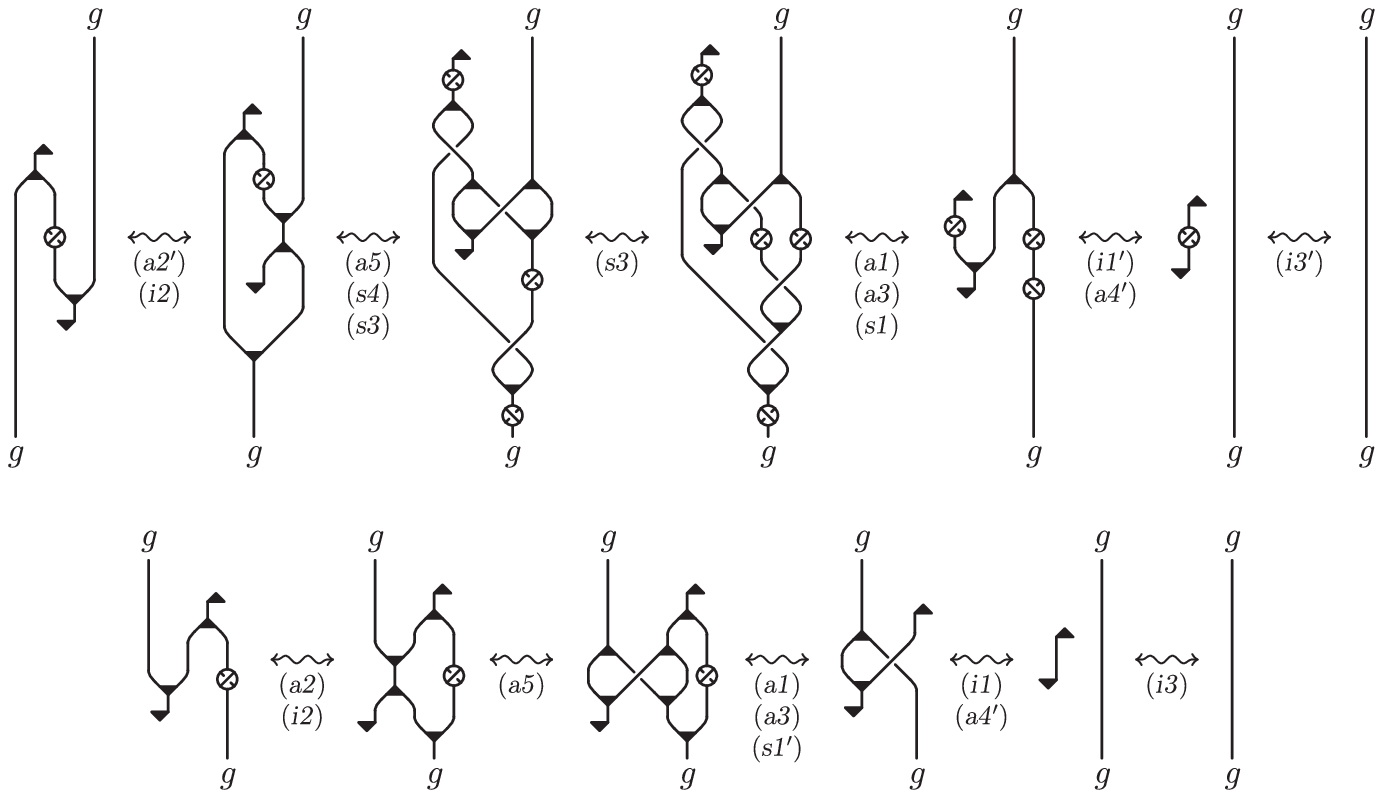}}\vskip-4pt
\end{Figure}

\vskip-4pt

To simplify notations we will often write $\lambda_g$ and
$\Lambda_g$ instead of $\lambda_{H_g}$ and $\Lambda_{H_g}$.

Lemma \ref{lemma/auto} implies that in the diagrams representing the morphisms of
$\H(\G)$, it is appropriate to use the notations \(f1) and \(f2) in Figure
\ref{form/fig} for the coform and the form. In fact, the identities in Paragraph
\ref{autonomous} reduce to the standard ``pulling the string'' (relation \(f3-3') in
Figure \ref{form/fig}), which together with the braid axioms in Figure \ref{isot/fig}
realize regular isotopy of strings.

\begin{Figure}[htb]{form/fig}{}{Coform and form in $\H(\G)$.}
\centerline{\fig{}{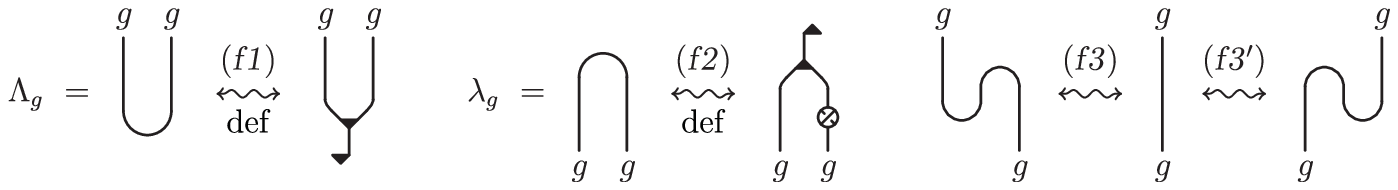}}\vskip-4pt
\end{Figure}

\medskip

Let $\Mor(A,B)$ denote the set of morphisms $A \to B$ in $\H(\G)$. Then we define
the ``rotation'' map 
$$\rot:\Mor(H_g \diam A,B\diam H_h) \to \Mor(A \diam  H_h,H_g \diam B)$$
by the equation
$$\rot(F) = (\id_{H_g \diam B} \diam \lambda_{h}) \circ (\id_{g} \diam F \diam \id_{h})
\circ (\Lambda_{g} \diam \id_{A \diam H_h}).$$

The following lemma shows that the negatively polarized tri-valent vertices are
invariant under such rotation, while for the positively polarized tri-valent vertices
this is equivalent to applying the antipode on the outgoing edge and its inverse on the
right incoming edge. 

\begin{Figure}[b]{cycl0/fig}{}{Properties of coform and form in $\H(\G)$.}
\centerline{\fig{}{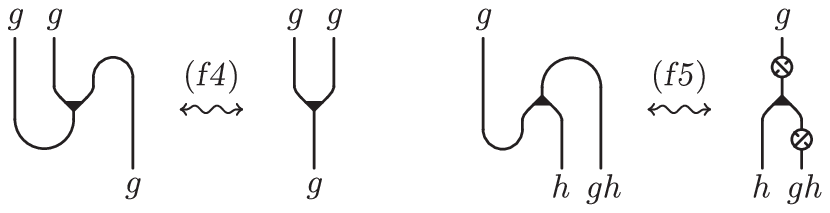}}\vskip-6pt
\end{Figure}

\begin{Figure}[htb]{cycl/fig}{}
 {Proof of Lemma \ref{lemma/cycl} 
 [{\sl a}/\pageref{algebra/fig}, {\sl f}/\pageref{form/fig}].}
\centerline{\fig{}{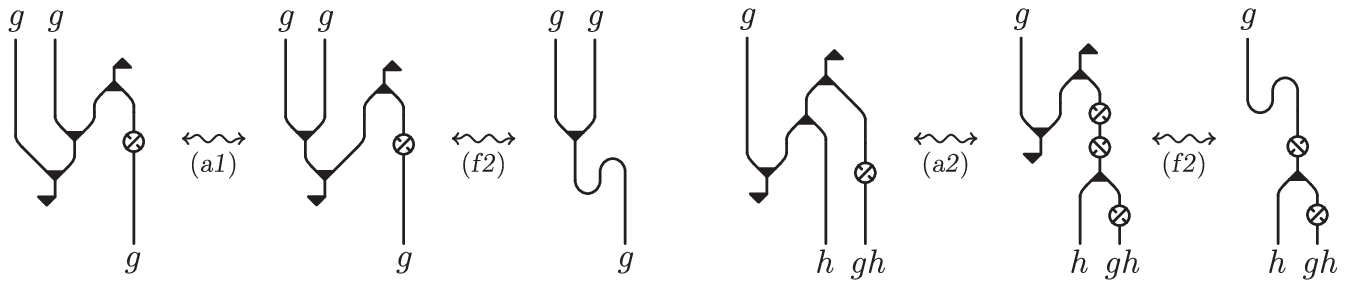}}\vskip-6pt
\end{Figure}

\begin{lemma} \label{lemma/cycl}
 If $g,h,gh\in\G$ ($g$ and $h$ are composable) then (cf. Figure \ref{cycl0/fig}):
\vskip-6pt
$$\rot(\Delta_g) = \Delta_g: H_g \to H_g \diam H_g\text{,}
\eqno{\(f4)}$$
$$\rot(m_{g,h}) = \bar S_{\bar g} \circ m_{h,\bar h \bar g} \circ (\id_g \diam
S_{gh}) :  H_h \diam H_{gh} \to H_g.
\eqno{\(f5)}$$\vskip-12pt
\end{lemma}

\begin{proof}
 See Figure \ref{cycl/fig}.
\end{proof}

\begin{block}\label{unimodular}
A Hopf $\G$-algebra $H$ in a braided monoidal category $\C$ is called {\sl unimodular}
if $H$ has $S$-invariant integral and cointegral, i.e.
\vskip-3pt
$$S_g \circ L_g = L_{\bar g},\eqno{\(i4)}$$ 
$$l_{\bar g} \circ S_g = l_g.\eqno{\(i5)}$$
\smallskip

\begin{Figure}[htb]{unimodular/fig}{}{Integral axioms for $\H^u(\G)$.}
\vskip-3pt\centerline{\fig{}{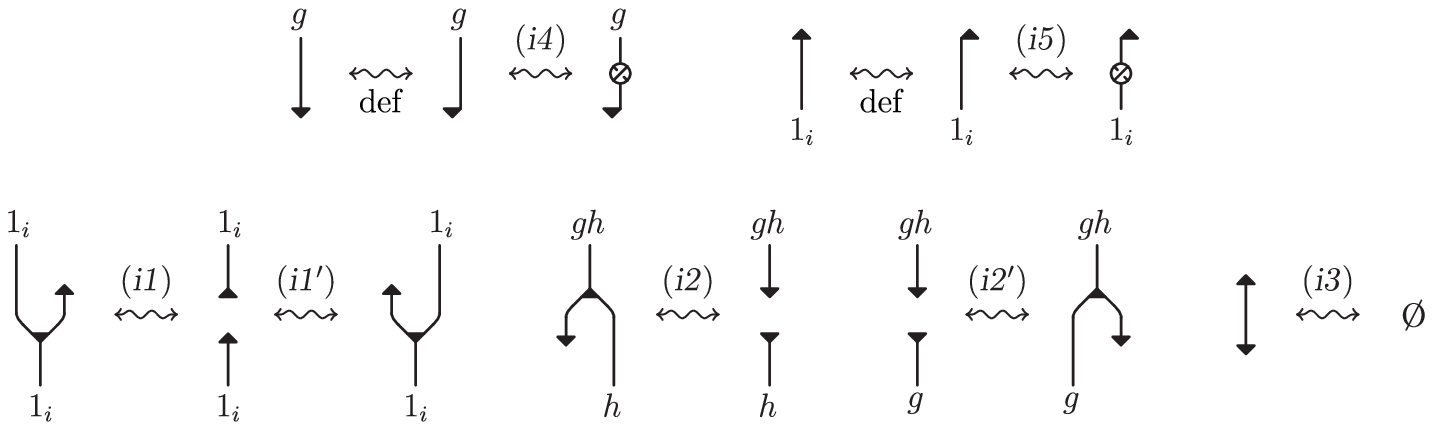}}\vskip-6pt
\end{Figure}

Given a groupoid $\G$, let $\H^u(\G)$ be the quotient of $\H(\G)$ modulo the
relations \(i4) and \(i5) presented in Figure \ref{unimodular/fig}. We refer to
$\H^u(\G)$ as the {\sl universal unimodular Hopf $\G$-algebra}. Moreover, in $\H^u(\G)$
we change the notation for the integral vertices by connecting the edge to the middle
point of the base of the triangle, to reflect that the corresponding integral is
two-sided (cf. the bottom line in Figure \ref{unimodular/fig}). 

Using the integral axioms \(i1)-\(i5) and \(a8) in Figure \ref{algebra/fig}, it is easy
to see that the uni-valent vertices of the same polarization are dual to each other
with respect to the form/coform as shown in Figure \ref{form-uni/fig}.

\begin{Figure}[htb]{form-uni/fig}{}
 {Duality of uni-valent vertices with the same polarization in $\H^u(\G)$.}
\centerline{\fig{}{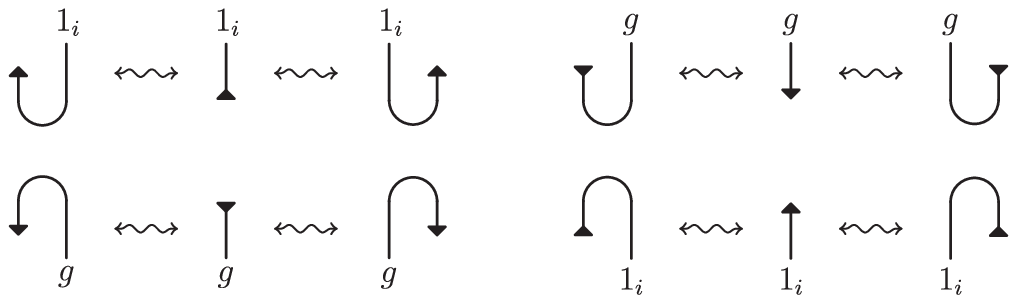}}\vskip-6pt
\end{Figure}

The next lemma is a version of Lemma 8 in \cite{Ke02}, but for completeness we
present the proof.

\begin{lemma} \label{h-tortile}  
$\H^u(\G)$ is a tortile category where for any $g \in \G$ (cf. Figure
\ref{h-tortile1/fig}): 
\vskip-6pt
$$\theta_{H_g} = S_{\bar g} \circ S_g\text{,} 
\eqno{\(f6-6')}$$
$$(\id_g \diam S_g) \circ\Lambda_{\bar g} = (S_{\bar g} \diam \id_g) \circ \Lambda_g
\text{,}
\eqno{\(f7)}$$
$$\lambda_g \circ(\id_g \diam S_{\bar g}) = \lambda_{\bar g} \circ (S_g \diam
\id_{\bar g}).
\eqno{\(f8)}$$\vskip-12pt
\end{lemma}

\begin{proof}
Define $\theta_A  =(\lambda_{A^\ast}\!\diam\id_A) \circ (\id_{A^\ast}\!\diam
\gamma_{A,A}) \circ (\Lambda_{A}\diam\id_A )$ for any object $A$ in $\H^u(\G)$ (cf.
Paragraph \ref{autonomous}). 
In particular, $\theta_g$ is represented by the diagram on the left in move \(f6) in
Figure \ref{h-tortile1/fig}. This definition guarantees that $\theta_1 = \id_1$ and
that
\mypagebreak
the identity $\theta_{A \diam B} = \gamma_{B,A} \circ (\theta_B\diam\theta_A) \circ
\gamma_{A,B}$ holds, up to isotopy moves, for any two objects $A,B$ in $\H^u(\G)$.
Therefore, in order to see that $\theta$ makes $\H^u(\G)$ into a tortile category, it
is enough to show its naturality and that $\theta_{A^\ast} = (\theta_A)^\ast$ for any
$A$. 

\begin{Figure}[t]{h-tortile1/fig}{}
 {Additional properties of coform and form in $\H^u(\G)$.}
\centerline{\fig{}{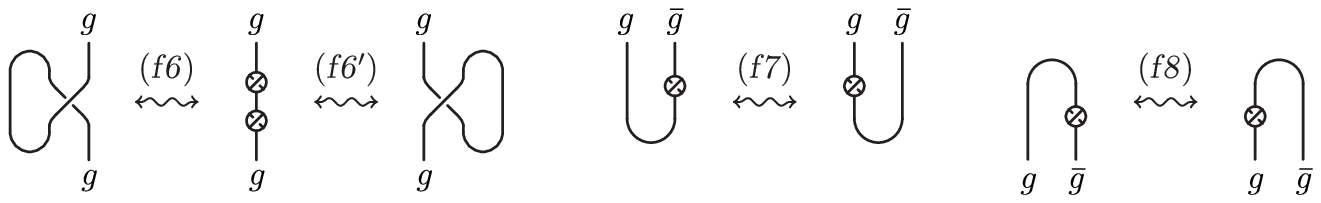}}\vskip-6pt
\end{Figure}

\begin{Figure}[htb]{h-tortile2/fig}{}{}
\vskip-6pt\centerline{\fig{}{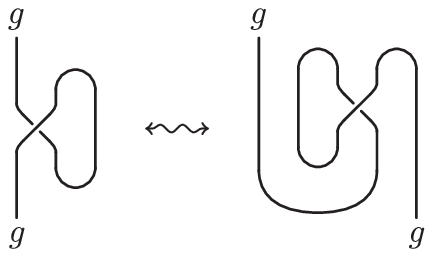}}\vskip-6pt
\end{Figure}

We will first prove the last identity. Actually, through an induction argument, one can
easily see that if this identity is true for $A=H_g$, then it is true for any other
object $A$ in $\H(\G)$. As it is shown in Figure \ref{h-tortile2/fig}, the diagram on
the right of move \(f6') in Figure \ref{h-tortile1/fig} coincides with
$(\theta_g)^\ast$  up to isotopy moves. In particular,  the identity for $A=H_g$ would
follow if we show \(f6) and \(f6'). Now, in Figure \ref{h-tortile/fig}  we prove moves
\(f6) and \(f7). Then move \(f6') follows by applying \(f6) and \(f7) to right diagram
in Figure \ref{h-tortile2/fig}.

\begin{Figure}[htb]{h-tortile/fig}{}
 {Proof of Lemma \ref{h-tortile} 
 [{\sl f}/\pageref{form/fig}, 
 {\sl s}/\pageref{pr-antipode/fig}-\pageref{antipode/fig}].}
\vskip-2pt\centerline{\fig{}{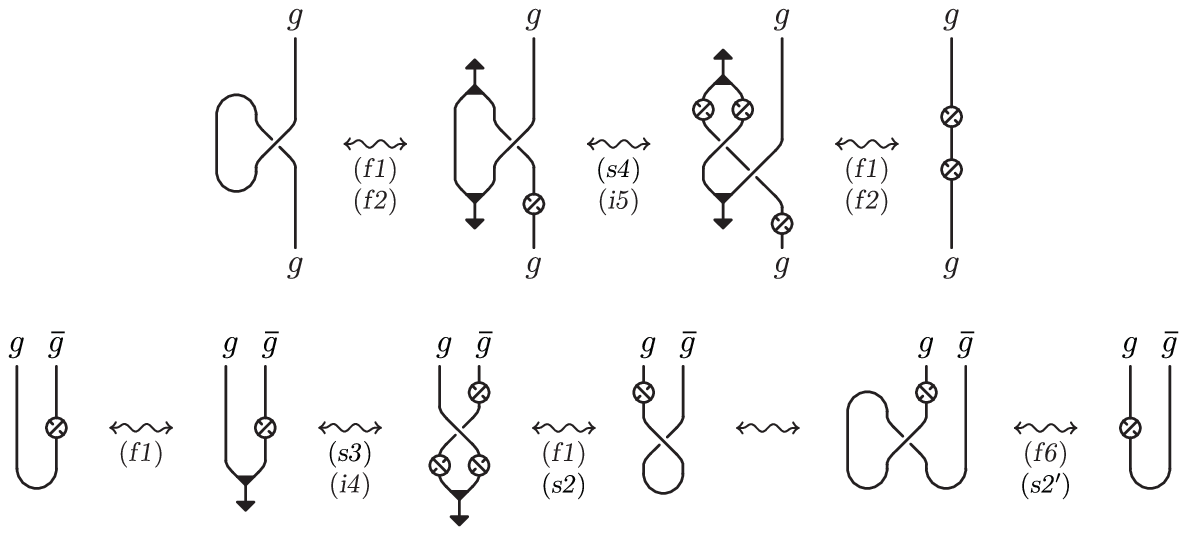}}\vskip-4pt
\end{Figure}

It is left to show the naturality of $\theta$, i.e. that for any morphism $F:A \to B$
in $\H^u(\G)$ one has $\theta_B \circ F = F \circ \theta_A$. Actually, using the fact
that $\theta_{A\diam B} = \gamma_{B,A} \circ (\theta_B\diam\theta_A)\circ\gamma_{A,B}$
and isotopy moves, one sees that it is enough to show this identity when $F$ is any
elementary morphism and in this case it follows from \(f6), \(s3-6) in Figure
\ref{pr-antipode/fig} and \(i4-5) in Figure \ref{unimodular/fig}.

To complete the proof of the theorem, observe that move \(f8) in Figure
\ref{h-tortile1/fig} follows from \(f7) and the properties \(f3-3') of the form and
coform presented in Figure \ref{form/fig}.
\end{proof}
\end{block}

\vglue-30pt\vskip0pt
\section{The universal ribbon Hopf algebra $\H^r(\G)$%
\label{ribbon-alg/section}}

\begin{block}\label{copairing}
 Let $\C$ be a braided monoidal category and $(A, m^A,\eta^A, \Delta^A, \epsilon^A,
S^A)$ and $(B, m^B,\eta^B, \Delta^B, \epsilon^B, S^B)$ be Hopf algebras in it (over
the trivial groupoid). Then a morphism $\sigma_{A,B}: \one \to A \diam B$ is called a
{\sl Hopf copairing} if the following conditions are satisfied:
\vskip-3pt
$$(\Delta^A \diam \id_B) \circ \sigma_{A,B} = (\id_{A \diam A} \diam m^B) \circ
(\id_A \diam \sigma_{A,B} \diam \id_B) \circ \sigma_{A,B},$$
$$(\id_A \diam \Delta^B) \circ \sigma_{A,B} = (m^A \diam \id_{B \diam B})
\circ (\id_A \diam \sigma_{A,B} \diam \id_B) \circ \sigma_{A,B},$$
$$(\epsilon^A \diam \id_B) \circ \sigma_{A,B} = \eta^B \quad \text{and} \quad 
(\id_A \diam \epsilon^B) \circ \sigma_{A,B} = \eta^A.$$
\vskip3pt\noindent
The Hopf copairing $\sigma_{A,B}$ is called {\sl trivial} if $\sigma_{A,B} = \eta^A
\diam \eta^B$.\smallskip
\end{block}

\begin{block}\label{ribbon}
A unimodular Hopf $\G$-algebra $H$ in a braided monoidal category $\C$ is called {\sl
ribbon} if there exists a family $v =\{v_g: H_g\to H_g\}_{g \in \G}$ of invertible
morphisms in $\C$, called {\sl ribbon morphisms}, such that:
\vskip-6pt 
{$$\epsilon_g \circ v_g = \epsilon_g \quad (\epsilon\text{-invariance}),
\eqno{\(r2)}$$
$$v_g \circ L_g = L_g \quad (L\text{-invariance}),
\eqno{\(r3)}$$
$$S_g \circ v_g = v_{\bar g} \circ S_g \quad (S\text{-invariance}),
\eqno{\(r4)}$$
$$m_{g,h} \circ (v_g \diam \id_h) = v_{gh} \circ m_{g,h} = m_{g,h} \circ (\id_g \diam
v_h) \quad (\text{centrality}),
\eqno{\(r5-5')}$$
$$\gamma_{g,h} \circ (\id_g \diam v_h ) = (v_h \diam \id_h) \circ \gamma_{g,h},
\eqno{\(b3)}$$
$$\bar\gamma_{g,h} \circ (\id_g\diam v_h ) = (v_h \diam \id_h) \circ \bar\gamma_{g,h};
\eqno{\(b4)}$$
\vskip3pt}
\smallskip\noindent
moreover, the family $\sigma = \{\sigma_{i,j}: \one \to H_{1_i} \diam H_{1_j}\}_{i,j
\in \Obj \G}$ of morphisms, defined by the identity
$$\arraycolsep0pt
\sigma_{i,j} = \left\{
\begin{array}{ll} 
 (v_{1_j}^{-1} \diam (v_{1_j}^{-1} \circ S_{1_j})) \circ \Delta_{1_j} \circ v_{1_j}
 \circ \eta_j &\quad \text{if } i = j,\\
 \eta_i \diam \eta_j &\quad \text{if } i\neq j,
\end{array}\right.
\eqno{\(r6-7)}$$
satisfies the following properties:
\vskip-6pt
{$$\hbox to 0pt{%
\hss$\sigma_{i,j}\text{ is a Hopf copairing for any } i,j \in \Obj \G$ \hss}
\eqno{\(r8-8'-9-9')}$$
$$\Delta_g \circ v_g^{-1} = \mu_{g,g} \circ (v_g^{-1} \diam
v_g^{-1}) \circ \bar \gamma_{g,g} \circ \Delta_g: H_g \to H_g \diam H_g,
\eqno{\(r10)}$$
$$\arraycolsep0pt
\begin{array}{rl}
(m_{1_j,h} \diam m_{g,1_i}) \circ {}& (S_{1_j} \diam (\mu_{h,g} \circ \bar
\gamma_{g,h} \circ \mu_{g,h}) \diam S_{1_i}) \circ (\rho^l_{g,j} \diam
\rho^r_{h,i}) = {} \\[4pt] & {} = \gamma_{g,h}: H_g \diam H_h \to H_h \diam H_g,
\end{array}
\eqno{\(r11)}$$}%
 where $\bar g g = 1_i$, $h \bar h = 1_j$ and:
$$
\begin{array}{rl}
 \rho^r_{g,j} &= (m_{g,1_i} \diam \id_{1_j}) \circ 
 (\id_g \diam \sigma_{i,j}): H_g \to H_g \diam H_{1_j},\\[4pt]
 \rho^l_{h,i} &= (\id_{1_i} \diam m_{1_j,h}) \circ 
 (\sigma_{i,j} \diam \id_h): H_h \to H_{1_i} \diam H_h,\\[4pt]
 \mu_{g,h} &= (m_{g,1_i} \diam m_{1_j,h}) \circ (\id_g \diam \sigma_{i,j} \diam
 \id_h) : H_g \diam H_h \to H_g \diam H_h.
\end{array}
$$
 It is easy to check that $\rho^r_{g,j}$ (resp. $\rho^l_{g,j}$) make $H_g$ into a right
(resp. left) $H_{1_j}\mkern-1mu$-comodule. We will refer to the morphisms
$\sigma_{i,j}$'s as {\sl copairing morphisms} or simply {\sl copairings}.\vskip-12pt
\end{block}

\begin{block}\label{defnHr/par}
 Given a groupoid $\G$, let $\H^r(\G)$ be the free braided monoidal category with a
ribbon Hopf $\G$-algebra $H$. We refer to $\H^r(\G)$ as the {\sl universal ribbon
Hopf $\,\G$-algebra}. $\H^r(\G)$ has the same objects as $\H(\G)$ and its morphisms
are iterated compositions and products of the elementary morphisms presented in
Figures \ref{elemdiag/fig} and \ref{diag-v/fig}, modulo the relations presented in
Figures \ref{isot/fig} (where now $D$ can be also a ribbon morphism),
\ref{algebra/fig}, \ref{antipode/fig}, \ref{unimodular/fig} and the additional
relations presented in Figures \ref{ribbon1/fig} and \ref{ribbon2/fig}, which express
the ribbon axioms \(r6)--\(r11) in \ref{ribbon}. Look at Figure \ref{diag-rho/fig}
for the diagramatic presentations of the morphisms $\rho^r_{g,j}$, $\rho^l_{h,i}$,
$\mu_{g,h}$ (here we also present the morphism $\mu^{-1}_{g,h}$ defined and shown to
be the inverse of $\mu_{g,h}$ in Proposition \ref{coform-s/theo}). 

\begin{Figure}[htb]{diag-v/fig}{}{Additional elementary diagrams in $\H^r(\G)$.}
\vskip-4pt\centerline{\fig{}{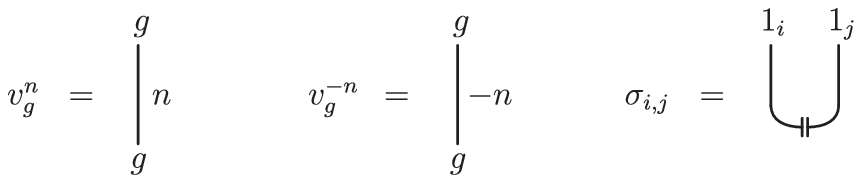}}\vskip-6pt
\end{Figure}
\begin{Figure}[htb]{ribbon1/fig}{}{Additional axioms for $\H^r(\G)$ -- I.}
\vskip-4pt\centerline{\fig{}{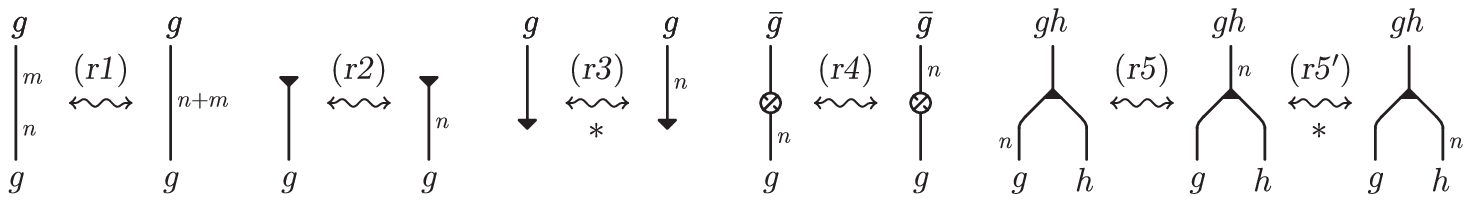}}\vskip-6pt
\end{Figure}
\begin{Figure}[htb]{ribbon2/fig}{}{Additional axioms for $\H^r(\G)$ -- II 
 ($i \neq j$).}
\centerline{\fig{}{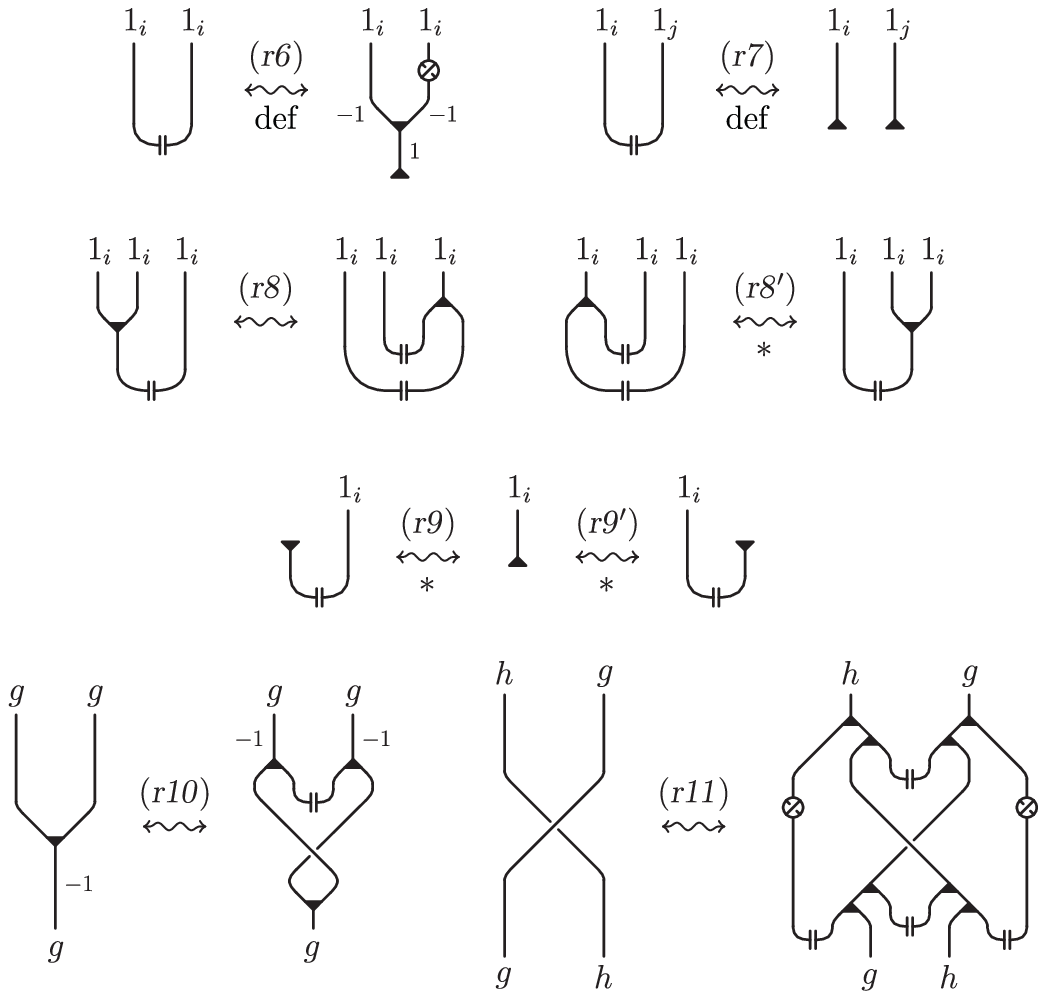}}\vskip-6pt
\end{Figure}
\begin{Figure}[htb]{diag-rho/fig}{}{}
\centerline{\fig{}{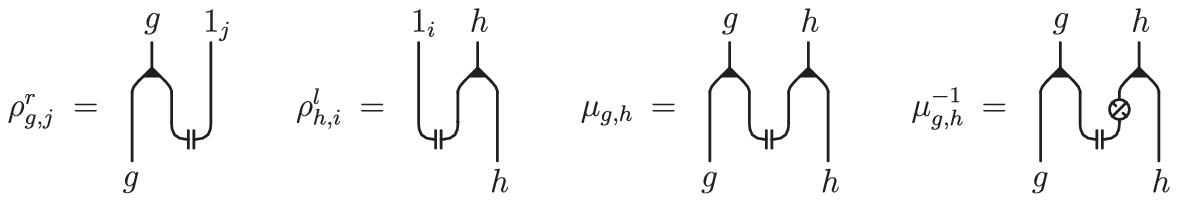}}\vskip-6pt
\end{Figure}

We make here few useful observations about the axioms of the ribbon algebra and their
diagramatic interpretation:
\begin{itemize}\itemsep\smallskipamount
\item[{\sl a}\/)] \vskip-\lastskip\smallskip
 In the diagrams  $v_g^n$ is interpreted as weight of an edge. Then relations \(r2) and
\(r3) say that the weight of an edge attached to a negatively polarized uni-valent
vertex can be changed arbitrarily, while \(r5-5') imply that the weight of one of the
edges attached to a positevely polarized tri-valent vertex can be transferred to any
other edge attached at that vertex.
\item[{\sl b}\/)] In Figure \ref{ribbon2/fig}, the relations  \(r8-8') and \(r9-9') 
expressing the fact that $\sigma_{i,j}$ is a Hopf copairing, are shown only for $i =
j$, since for $i\neq j$ the copairing is trivial  and they follow from the rest of the
axioms.
\item[{\sl c}\/)] For some combinations of the labels, moves \(r10) and \(r11) simplify
because trivial copairings appear. In particular, move \(r11) reduces to a crossing
change when $g \in \G(i,j)$ and $h\in\G(k,l)$ with $\{i,j\} \cap \{k,l\} =
\emptyset$. \label{r11/note}
\item[{\sl d}\/)] An interesting and important question is which of the axioms for $v$
and $\sigma$ in \ref{ribbon} are independent. In particular, the relations \(r3),
\(r5'), \(r8') and \(r9-9') (marked with a star below the arrow), are consequences of
the rest and have been listed as axioms only for convenience. In fact: \(r3) can be
obtained from \(r2) by using \(r5) and the duality of the negative uni-valent vertices
(Figure \ref{form-uni/fig}); \(r9-9') immediately follow from the definition \(r6-7) of
the copairing, and the properties \(a2-2'), \(s5-6) and \(r2) (Figures
\ref{algebra/fig} and \ref{pr-antipode/fig}); \(r5') can be shown to be equivalent to
\(r5) by using property \(s4) of the antipode (Figure \ref{pr-antipode/fig}).
The proof that \(r8') derives from the rest of the axioms is presented in Appendix
(Section \ref{appendix/sec}). On the other hand, it seems unlikely that the relations
\(r10) and \(r11), the new ones with respect to the relations in \cite{Ke02}, are
consequences of the rest, even if, as we will see later, they can be presented in
different equivalent forms (cf. Corollary \ref{eqaxioms/theo}). 
\end{itemize}\vskip-20pt
\end{block}

\begin{proposition}\label{formalext/theo} 
 Any functor $\phi: \G \to \G'$ between groupoids which is injective on the set of
objects can be extended to a  functor  $\Upsilon_\phi:\H^r(\G) \to
\H^r(\G')$. Moreover, if $\phi$ is faithful (an embedding) then $\Upsilon_\phi$ is
also faithful.
\end{proposition}
	
 In particular, when $\iota: \G\subset \G'$ is an inculsion of groupoids, we can
identify $\H(\G)$ and $\Upsilon_\iota(\H(\G))$ through the isomorphism of categories
$\Upsilon_\iota: \H(\G) \to \Upsilon_\iota(\H(\G))$ and write $\H(\G) \subset \H(\G')$.

\begin{proof}
 The extension $\Upsilon_\phi$ of $\phi$ is formally defined as $\Upsilon_\phi(m_{g,h})
= m_{\phi(g),\phi(h)}$, $\Upsilon_\phi(\eta_i) = \eta_{\phi(i)}$, etc. 
To see that $\Upsilon_\phi$ is well defined, we need to check that all relations for
$\H^r(\G)$ are satisfied in the image. The only problem we might have would be with
relation \(r7) in Figure \ref{ribbon2/fig}, if $i \neq j$ and $\phi(i) = \phi(j)$. But
this cannot happen since $\phi$ is injective on objects. This concludes the first part
of the proposition since the functoriality of $\Upsilon_\phi$ is obvious.

At this point, it is left to show that when $\phi$ is injective on the set of
morphisms, then $\Upsilon_\phi$ is injective on  morphisms as well. First of all, in
this case $\phi$ induces an isomorphism of categories $\hat{\phi}:\G\mapsto \G^{\phi}$,
where $\G^{\phi}=\phi(\G)$. Then by definition, $\Upsilon_{\hat\phi} \circ
\Upsilon_{\hat\phi^{-1}}$ and $\Upsilon_{\hat\phi^{-1}}\circ\Upsilon_{\hat\phi}$ are
the identities; therefore $\Upsilon _{\hat\phi}: \H(\G)\mapsto \H(\G^{\phi})$ is an
isomorphism of categories as well. Moreover, $\Upsilon_\phi=\Upsilon_\iota\circ
\Upsilon_{\hat\phi}$ where $\iota:\G^{\phi}\subset \G'$ is the corresponding inclusion.
Hence, the statement would follow if we show that the functor $\Upsilon_\iota$ is
injective on the set of morphisms.

Now let $F,F':A \to B$ be two morphsims in $\H(\G^{\phi})$ such that $\Upsilon_\iota(F)
= \Upsilon_\iota(F')$ in $\H(\G')$; in particular, $F$ and $F'$ are represented by two
diagrams labeled in $\G^{\phi}$, which are related by a sequence of moves in $\H(\G')$.
Observe that when we apply a relation move to a diagram representing a morphism of
$\H(\G')$, the only new labels that can appear are identities of $\G'$ and products of
labels already occurring in it or their inverses. Since $\G^{\phi}$ is a subcategory of
$\G'$, this implies that the only labels not belonging to $\G^{\phi}$ that can occur
in the intermediate diagrams of the sequence are identities $1_i$ with $i \in \Obj \G'
- \Obj \G^{\phi}$. The parts of the graph diagram carrying such labels interact with
the rest of the diagram only through move \(r7) in Figure \ref{ribbon2/fig}; in
particular, if to any intermediate diagram we apply move \(r7) to change back any
trivial copairing into two units, the part of the graph labeled by $1_i$'s with $i \in
\Obj \G' - \Obj \phi(\G)$, forms a component, disjoint from the rest. By deleting such
component, we obtain a new sequence of diagrams between $F$ and $F'$ related by moves
in $\H(\G^{\phi})$ which proves that $F = F'$ (in $\H(\G^{\phi})$).
\end{proof}

\smallskip

\begin{proposition}\label{coform-s/theo} 
For every $i,j \in \Obj \G$ and $g,h \in \G$, we have (cf. Figure \ref{coform-s/fig}):
\vskip-6pt
{$$(S_{1_i}\diam\id_j) \circ \sigma_{i,j} = (\id_i \diam S_{1_j}) \circ \sigma_{i,j}
\text{,}
\eqno{\(p1)}$$
$$\mu_{g,h} \circ \mu^{-1}_{g,h}  = \id_{g \diam h} = \mu^{-1}_{g,h} \circ \mu_{g,h}
\text{,}
\eqno{\(p2-2')}$$}
\vskip-9pt\noindent
with $\mu^{-1}_{g,h}$ defined by the following identity, where we are assuming $\bar g
g = 1_i$ and $h \bar h = 1_j$ (cf. Figure \ref{diag-rho/fig}):
\vskip-18pt
$$\mu^{-1}_{g,h} = (m_{g,1_i} \diam m_{1_j,h}) \circ (\id_g \diam ((\id_{1_i} \diam
S_{1_j}) \circ \sigma_{i,j}) \diam \id_h): H_g \diam H_h \to H_g \diam H_h.$$ 
\end{proposition}

\begin{Figure}[htb]{coform-s/fig}{}{Some more relations in $\H^r(\G)$ -- I.}
\vskip-24pt\centerline{\fig{}{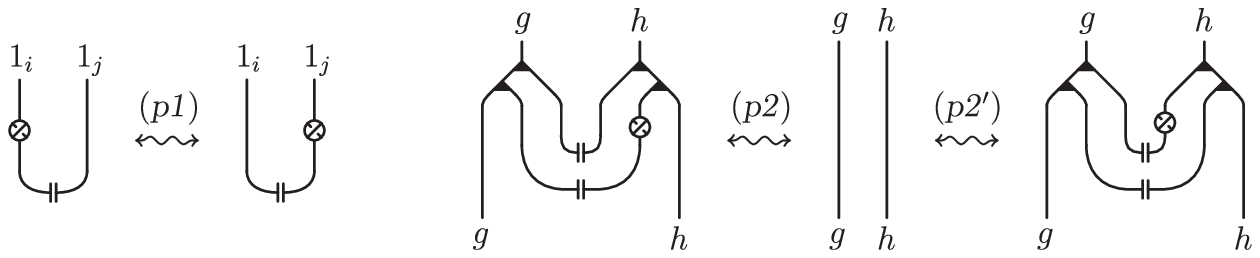}}\vskip-6pt
\end{Figure}

\begin{proof}
\(p2) is proved in Figure \ref{coform-spf/fig}. The proof of \(p2') is analogous. 
To prove \(p1) let us consider the morphism $\bar\mu^{-1}_{g,h}$, which is the same as
$\mu^{-1}_{g,h}$ but with the antipode moved to the left side of $\sigma_{i,j}$:
$$\bar\mu^{-1}_{g,h} = (m_{g,1_i} \diam m_{1_j,h}) \circ (\id_g \diam ((S_{1_i} \diam
\id_{1_j}) \circ \sigma_{i,j}) \diam \id_h): H_g \diam H_h \to H_g \diam H_h.$$ 
Then, by replacing \(r8'), \(r9') and \(s1') in Figure \ref{coform-spf/fig}
respectively with \(r8), \(r9) and \(s1), we obtain that \(p2) and \(p2') are still
valid if $\mu^{-1}_{g,h}$ is replaced by $\bar\mu^{-1}_{g,h}$. Hence
$\bar\mu^{-1}_{g,h}$ and $\mu^{-1}_{g,h}$ are equal, being both two-sided inverses of
$\mu_{g,h}$. This gives \(p1).
\end{proof}

\begin{Figure}[htb]{coform-spf/fig}{}
 {Proof of relation \(p2) in Figure \ref{coform-s/fig}
 [{\sl a}/\pageref{algebra/fig}, {\sl r}/\pageref{ribbon2/fig},
 {\sl s}/\pageref{antipode/fig}].}
\centerline{\fig{}{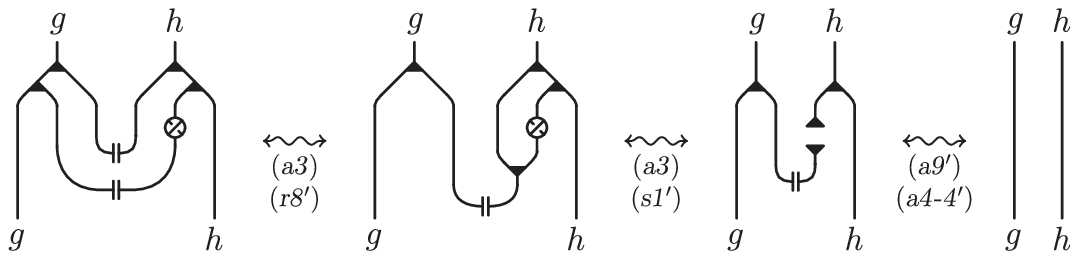}}\vskip-6pt
\end{Figure}

\begin{Figure}[b]{ribbon-isot/fig}{}{Additional isotopy relations in $\H^r(\G)$.}
\centerline{\fig{}{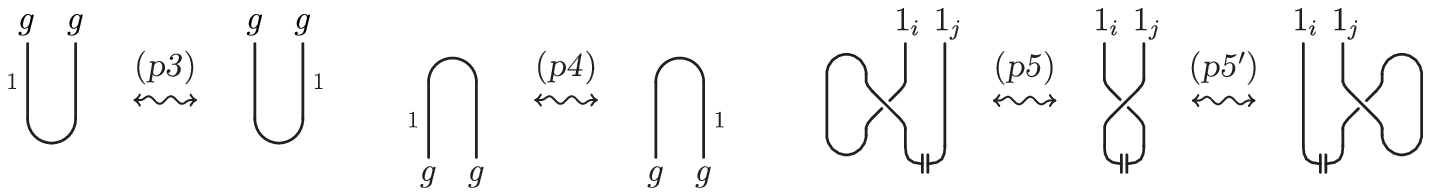}}\vskip-4pt
\end{Figure}

\begin{Figure}[htb]{ribbon-tot/fig}{}{Some more relations in $\H^r(\G)$ -- II.}
\vskip4pt\centerline{\fig{}{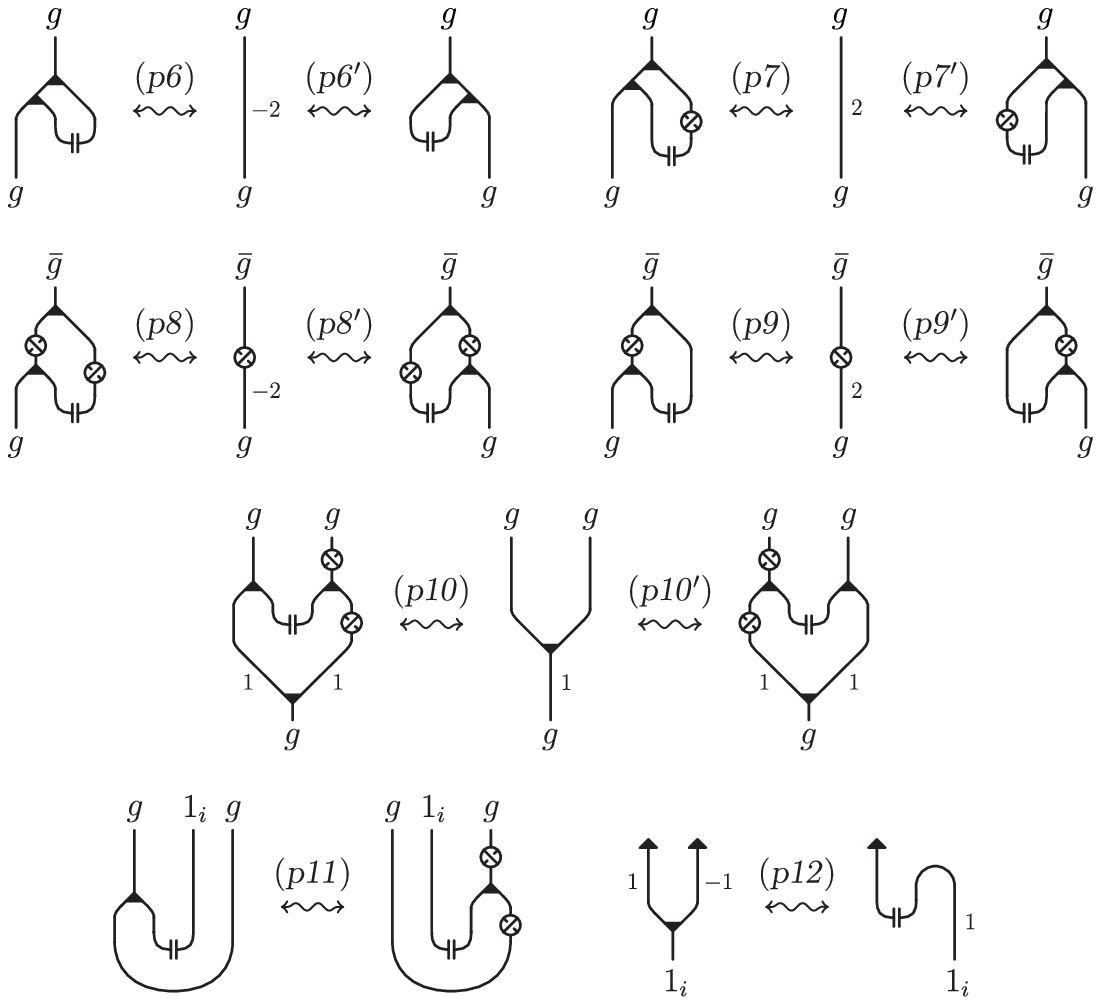}}\vskip-4pt
\end{Figure}

 In Figures \ref{ribbon-isot/fig} and \ref{ribbon-tot/fig} we list some more useful
relations satisfied in $\H^r(\G)$, but in order not to make this section excessively
technical, we collect their proofs in Appendix (Section \ref{appendix/sec}). We only
point out that \(p3-3'), \(p5-5') and \(p10-10') -- \(p12) follow directly from the
definition of $\sigma_{i,j}$ and axioms \(r1)--\(r5'), but do not depend on the extra
axioms in Figure \ref{ribbon2/fig}, in particular on the fact that $\sigma_{i,j}$ is a
Hopf copairing. Observe also that \(p1) in Figure \ref{coform-s/fig} and the moves in
Figure \ref{ribbon-isot/fig} allow us to extend the notion of isotopy to $\H^r(\G)$. 
As it is shown in Appendix \ref{appendix/sec}, they imply that if a string connecting
two polarized vertices, is divided by a copairing morphism, this morphisms can be
shifted anywhere between those vertices. Therefore, we will say that two graph
diagrams in $\H^r(\G)$ are {\sl isotopic} or that one of them is obtained from the
other through {\sl isotopy} if  it is obtained by applying a sequence of moves
\(b1-4)  in Figure \ref{isot/fig}, \(f3-3') in Figure \ref{form/fig},
\(f6)${}\circ{}$\(f6') and \(f7-8) in Figure \ref{h-tortile1/fig}, \(p1) in Figure
\ref{coform-s/fig} and \(p3-5') in Figure \ref{ribbon-isot/fig}.

\medskip

An immediate consequence of Proposition \ref{coform-s/theo} is the following.

\begin{corollary}\label{eqaxioms/theo} 
 The axiom \(r10) can be replaced by either \(r10'), \(r12) or \(r12') in Figure
\ref{ribbon4/fig}, while axiom \(r11) can be replaced by \(r11') in Figure
\ref{ribbon5/fig}.
\end{corollary}\vskip-18pt
 
\begin{Figure}[htb]{ribbon4/fig}{}{Equivalent presentations of \(r10).}
\vskip6pt\centerline{\fig{}{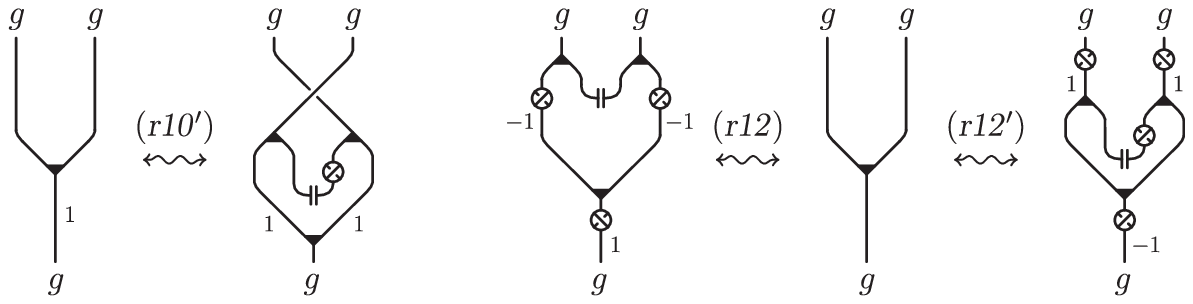}}\vskip-4pt
\end{Figure}

\begin{Figure}[htb]{ribbon5/fig}{}{Equivalent presentation of \(r11).}
\vskip-12pt\centerline{\fig{}{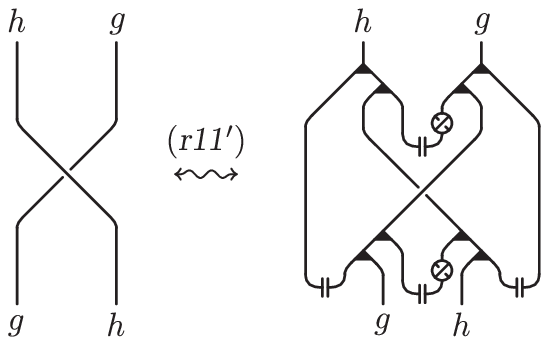}}\vskip-4pt
\end{Figure}

\begin{proof}
The equivalence between \(r10) and \(r12) derives from \(s3) in Figure
\ref{pr-antipode/fig} after having composed both sides of \(r10) on the bottom with
$v_g$. \(r10') is obtained from \(r10) by composing both sides on the bottom with
$v_g$ and on the top with the invertible morphism $\gamma_{g,g} \circ \mu^{-1}_{g,g}
\circ (v_g \diam v_g)$. Analogously, \(r12') is obtained from \(r12) by replacing
$g$ with $\bar g$ and composing both sides on the bottom and on the top respectively
with the invertible morphisms $S_g \circ v_g^{-1}$ and $(\bar S_g \diam
\bar S_g) \circ \mu^{-1}_{g,g} \circ (v_g \diam v_g)$.

To see that \(r11) and \(r11') are equivalent, it suffices to observe that the right
sides of the two relations are inverse of each other by using \(p2-2').
\end{proof}

\begin{Figure}[b]{repT/fig}{}
 {Additional relations for $T_g$ -- I ($g \in \G(i,j)$, $i \neq j$).}
\centerline{\fig{}{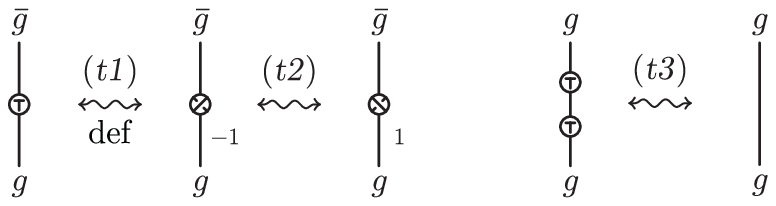}}\vskip-6pt
\end{Figure}

\begin{block}\label{defnT/par}
As we have seen at page \pageref{r11/note} for move \(r11), some relations involving
the copairing simplify significantly for suitable choices of the labeling, due to the
fact that $\sigma_{i,j}$ is trivial if $i \neq j$.

In particular, for $g \in \G(i,j)$ with $i \neq j$ relation \(p8) in Figure
\ref{ribbon-tot/fig} implies that $S_g \circ v_g^{-1} = \bar S_g \circ v_g$ and in
this case we define (cf. Figure \ref{repT/fig})
$$T_g = S_g \circ v_g^{-1} = \bar S_g \circ v_g
\quad \text{if } g \in \G(i,j) \text{ with } i \neq j.
\eqno{\(t1-2)}$$
 Then we have
$$T_{\bar g} \circ T_g = \id_g.
\eqno{\(t3)}$$ 

For any object $A$ in $\H^r(\G)$, let $V_A: A \to A^\ast$ be the morphism in
$\H^r(\G)$, defined inductively by the following identities (the definition is
well-posed, giving equivalent results for different decompositions $A \diam B = A'
\diam B'$):
$$V_{\one}=\id_{\one}\,,\!\quad V_{H_g} = v_g\quad \text{and} \quad V_{A \diam B} =
(V_B \diam V_A) \circ \gamma_{A,B}.$$
 This morphism is obviously invertible and for any $F: A \to B$ we define: 
$$\rev(F) = V_B^{-1} \circ F \circ V_A.$$\vskip-18pt
\end{block}

\begin{proposition} \label{T/theo} 
 For any $g \in \G(i,j)$ and $h \in \G(j,k)$ with $i\neq j\neq k \neq i$, we
have (cf. Figure \ref{T-theo/fig}):
\vskip-6pt
$$(T_{\bar g} \circ T_{\bar g}) \circ \Delta_{\bar g} \circ T_g = \Delta_g =
\rev(\Delta_g)\text{;}
\eqno{\(t4-5)}$$
$$\rot(m_{g,h}) = T_{\bar g} \circ m_{h,\bar{gh}} \circ (\id_h \diam T_{gh})\text{;}
\eqno{\(t6)}$$ 
$$\rev(m_{g,h}) = T_{\bar{gh}} \circ m_{\bar h,\bar g} \circ (T_h \diam T_g).
\eqno{\(t7)}$$\vskip-18pt
\end{proposition}

\begin{Figure}[htb]{T-theo/fig}{}{Additional relations for $T_g$ -- II
 ($g \in \G(i,j)$, $h \in \G(j,k)$, $i \neq j\neq k \neq i$).\kern-5pt}
\vskip-4pt\centerline{\fig{}{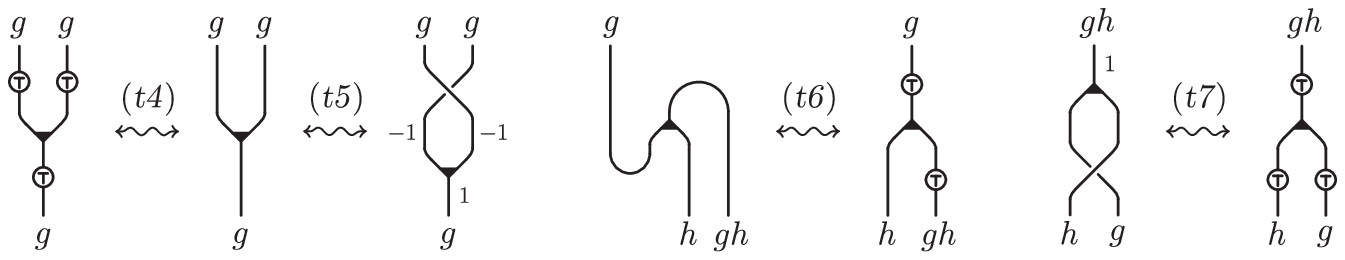}}\vskip-6pt
\end{Figure}

\begin{proof}
\(t4) and \(t5) rewrite \(r12) and \(r10) (in Figures \ref{ribbon4/fig} and
\ref{ribbon2/fig}) when $g \in \G(i,j)$ with $i\neq j$, i.e. when $\sigma_{i,j}$ is
trivial. \(t6) and \(t7) rewrite \(f5) and \(s4) (in Figures \ref{cycl0/fig} and
\ref{pr-antipode/fig}), under the further assumption that $h \in \G(j,k)$ and $i 
\neq k \neq j$.
\end{proof}

\medskip

\begin{block}\label{defn-boundH/par}
 A ribbon Hopf $\G$-algebra $H$ in a braided monoidal category $\C$ is called {\sl
selfdual} if 
$$(l_i \diam \id_{1_i}) \circ \sigma_{i,i} = L_{1_i} = (\id_{1_i} \diam l_i) \circ
\sigma_{i,i}.
\eqno{\(d1-1')}$$
\medskip\noindent
 A selfdual ribbon Hopf $\G$-algebra $H$ is called {\sl boundary} if
$$l_i \circ v_{1_i} \circ \eta_i = \id_{1_i} = l_i \circ v^{-1}_{1_i} \circ \eta_i.
\eqno{\(d2-2')}$$

\begin{Figure}[htb]{defn-boundH/fig}{}
 {Axioms for $\partial^\star \H^r(\G)$ and $\partial \H^r(\G)$.}
\centerline{\fig{}{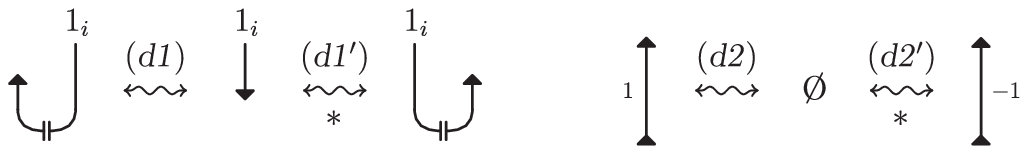}}\vskip-4pt
\end{Figure}
	
We define $\partial^\star\H^r(\G)$ (resp. $\partial \H^r(\G)$) to be the quotient
category of $\H^r(\G)$ modulo the relation \(d1) presented in Figure
\ref{defn-boundH/fig} (resp. the relations \(d1) and \(d2) in the same figure).
We call $\partial^\star\H^r(\G)$ (resp. $\partial \H^r(\G)$) the {\sl universal
selfdual} (resp. {\sl boundary}) ribbon Hopf algebra.

\break

Actually \(d1') and \(d2') (marked with a star below the arrow) have been listed as
axioms only for convenience. Indeed, \(d1') follows from \(d1) and \(p5-5') in Figure
\ref{ribbon-isot/fig}, while \(d2') follows from \(d2) and the relation shown in
Figure \ref{boundH-pf/fig}.
\end{block}

\begin{Figure}[htb]{boundH-pf/fig}{}
 {A relation in $\partial^\star\H^r(\G)$
 [{\sl a}/\pageref{algebra/fig},
 {\sl d}/\pageref{defn-boundH/fig}, {\sl p}/\pageref{ribbon-tot/fig}].}
\centerline{\fig{}{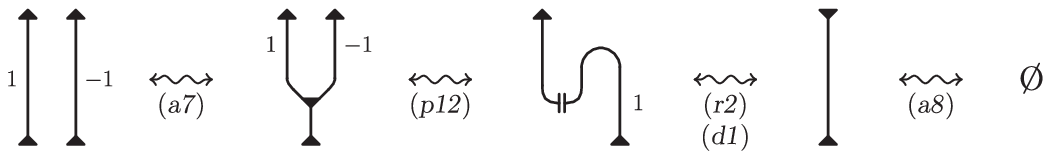}}\vskip-4pt
\end{Figure}

As it was observed by Kerler in \cite{Ke02}, the relation \(d1) in Figure
\ref{defn-boundH/fig} implies the nondegeneracy of the copairing and the duality of 
the multiplication and comultiplication morphisms in $\partial^{\star}\H^r(\G)$. 
In particular, if one defines pairing morphisms $H_{1_i} \diam H_{1_j} \to \one$ as
\(d3) in Figure \ref{selfdual/fig}, then the relations \(d4-4') and \(d5) presented 
in the same figure can be easily derived by using isotopy moves and move \(r8') in
Figure \ref{ribbon2/fig}.

\begin{Figure}[htb]{selfdual/fig}{}
 {Some more relations in $\partial^{\star}\H^r(\G)$.}
\centerline{\fig{}{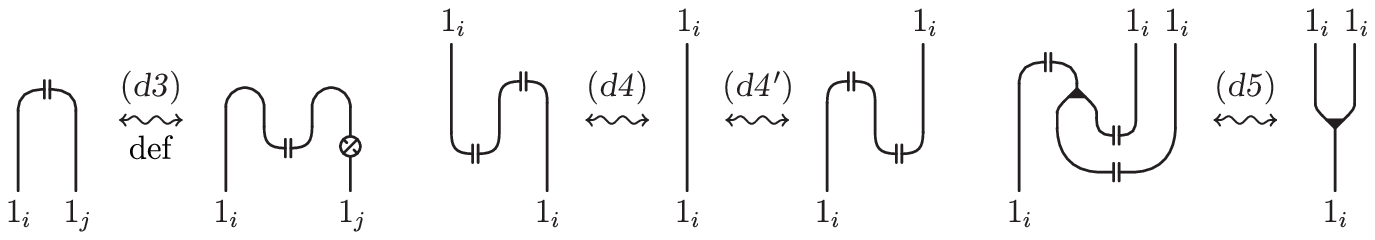}}\vskip-4pt
\end{Figure}

\vskip-9pt\vskip0pt

\section{From the algebra to the generalized Kirby diagrams%
\label{alg-kirby/sec}}\vskip-3pt\vskip0pt
\begin{block}\label{alg-kirby/par}
{\sl Proof of Theorem \ref{alg-kirby/theo}.}\kern1ex
 The braided monoidal functor $\Phi_n: \H_n^r \to \K_n$ is defined by sending
$H_{(i,j)}$ to the ordered pair $(i,j)$ of labeled intervals, while the images of the
elementary morphisms in $\H_n^r$ are presented in Figure \ref{defnPhi/fig}.

\begin{Figure}[htb]{defnPhi/fig}{}
 {Definition of $\Phi_n: \H_n^r \to \K_n$.}
\centerline{\fig{}{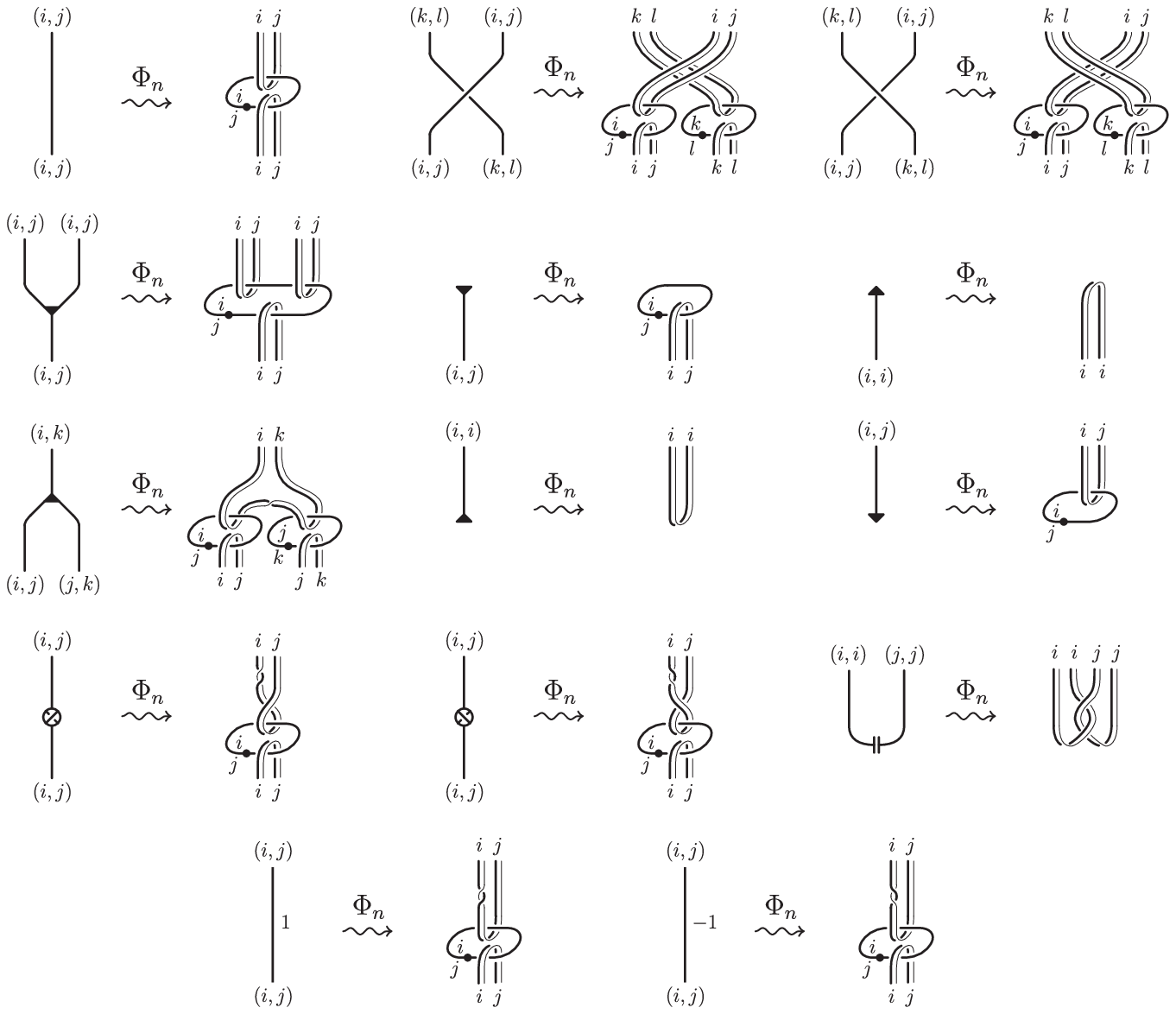}}\vskip-4pt
\end{Figure}

The reader can check that the images of the form and the coform $\lambda_{(i,j)}$ and
$\Lambda_{(i,j)}$ are equivalent in $\K_n$ to the ones presented in Figure
\ref{defnPhi7/fig}.

\begin{Figure}[htb]{defnPhi7/fig}{}
 {$\Phi_n(\lambda_{(i,j)})$ and $\Phi_n(\Lambda_{(i,j)})$.}
\centerline{\fig{}{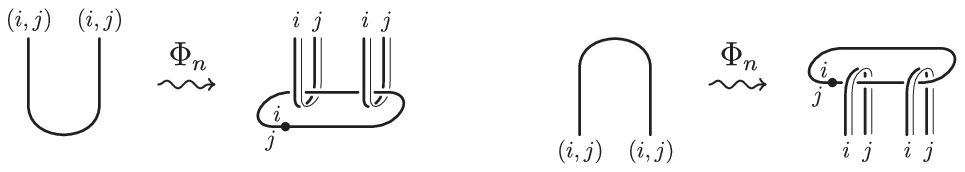}}\vskip-4pt
\end{Figure}

The proof of the theorem is an extension of the well known fact that the category of
admissible tangles contains a braided Hopf algebra object (see \cite{Ke02, H00}). In
particular, here we work with the category of generalized  tangles, i.e. tangles with
dotted components and different labels, and correspondingly a groupoid Hopf algebra.
Moreover we need to check the extra ribbon axioms in \ref{ribbon}.

Most of the Hopf algebra axioms are very easy to check. For example \(a2-2'), \(a6)
and \(a8) in Figure \ref{algebra/fig} and \(i2), \(i3) in Figure
\ref{unimodular/fig} follow directly from the deletion of 1/2-canceling pairs
(the bottom move in Figure \ref{diag4/fig}). The same is true for \(a1), \(a3),
\(a4-4'), \(a7) and \(i1-1'), but one needs to make one or two handle slides before
deleting. \(s2-2') and \(i4) reduce to an isotopy. \(a5) in Figure
\ref{algebra/fig} and \(s1) in Figure \ref{antipode/fig} are shown in Figures
\ref{defnPhi1/fig} and \ref{defnPhi2/fig} respectively (\(s1') is analogous).

\begin{Figure}[b]{defnPhi1/fig}{}
 {Proof of \(a5) in Figure \ref{algebra/fig}.}
\centerline{\fig{}{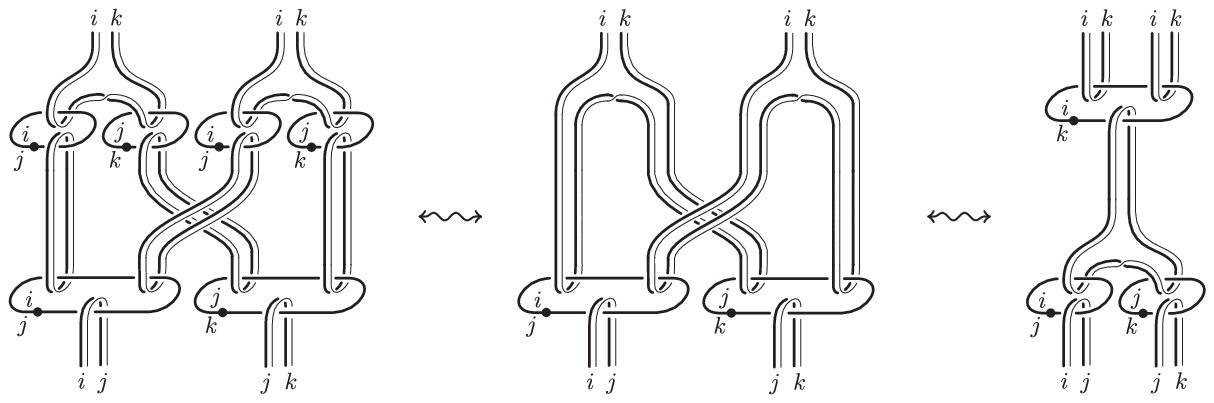}}\vskip-6pt
\end{Figure}

\begin{Figure}[htb]{defnPhi2/fig}{}
 {Proof of \(s1) in Figure \ref{antipode/fig}.}
\centerline{\fig{}{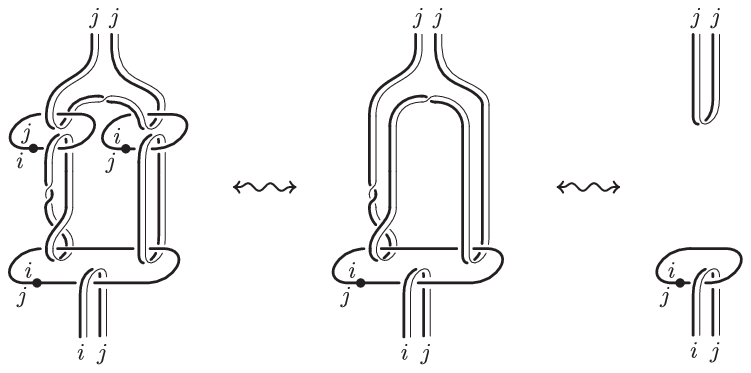}}\vskip-6pt
\end{Figure}

\begin{Figure}[b]{defnPhi8/fig}{}
 {Proof of \(r5-5') in Figure \ref{ribbon1/fig}.}
\vskip-2pt\centerline{\fig{}{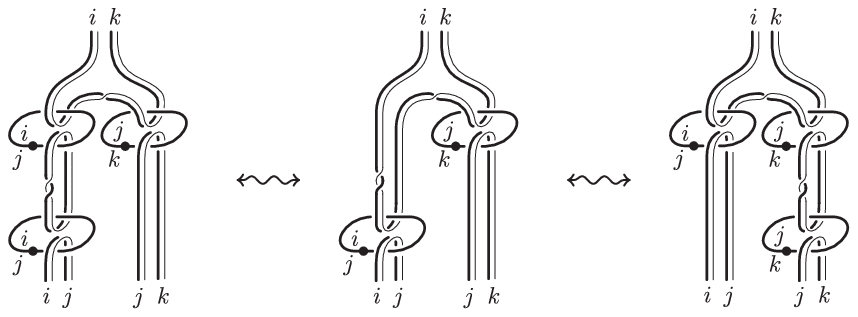}}\vskip-6pt 
\end{Figure}

\begin{Figure}[htb]{defnPhi4/fig}{}
 {Proof of \(r6) in Figure \ref{ribbon2/fig}.}
\centerline{\fig{}{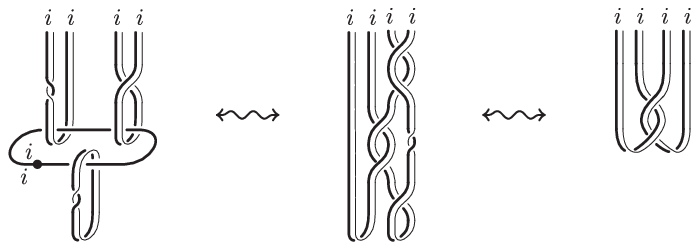}}\vskip-4pt
\end{Figure}

\begin{Figure}[htb]{defnPhi3/fig}{}
 {Proof of \(r8) in Figure \ref{ribbon2/fig}.}
\centerline{\fig{}{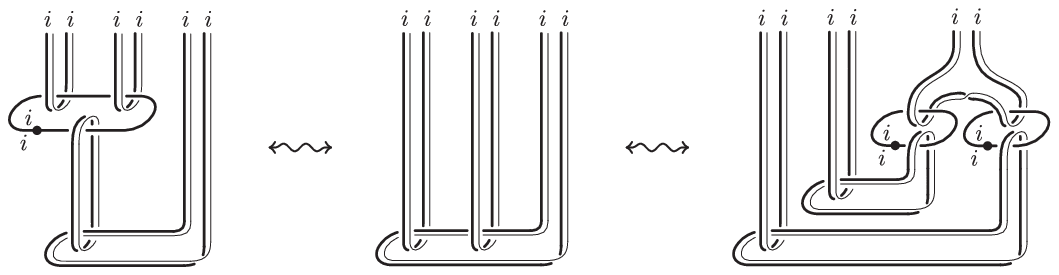}}\vskip-2pt
\end{Figure}

\begin{Figure}[htb]{defnPhi5/fig}{}
 {Proof of \(r10) in Figure \ref{ribbon2/fig}.}
\centerline{\fig{}{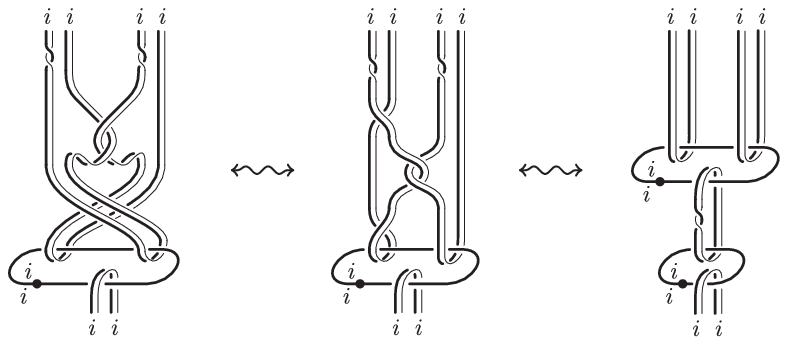}}\vskip-6pt
\end{Figure}

\begin{Figure}[htb]{defnPhi6/fig}{}
 {Proof of \(r11) in Figure \ref{ribbon2/fig}.}
\centerline{\fig{}{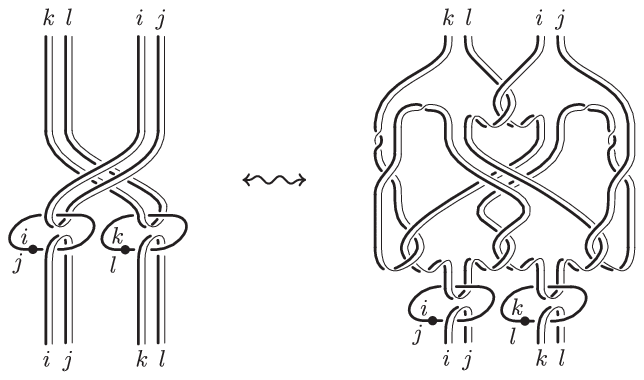}}\vskip-6pt
\end{Figure}

The ribbon axiom \(r7) in Figure \ref{ribbon2/fig} follows from the last move in Figure
\ref{diag3/fig}, which allows to change crossings between components with different
labels. The ribbon axioms \(r2), \(r3), \(r9) and \(r9') follow from the deletion
of 1/2-canceling pairs, while in showing \(r1) and \(r4) one needs to make a handle
slide before deleting. The rest of the ribbon axioms are shown in Figures
\ref{defnPhi8/fig} to \ref{defnPhi6/fig}. $\;\square$
\end{block}

Observe  that under the functor $\Phi_n$, the relations \(d1) and \(d2) in Figure
\ref{defn-boundH/fig} translate directly in \(a) and \(b) in Figure
\ref{defn-boundK/fig}. This leads to the following proposition.

\begin{proposition} \label{boundK/theo} 
$\Phi_n$ induces functors	on the quotient categories
\vskip-6pt
$$\partial^\star\Phi_n: \partial^\star\H^r_n \to \partial^\star\K_n
\quad \text{and} \quad
\partial\Phi_n: \partial\H^r_n \to \partial\K_n.$$
\vskip-12pt
\end{proposition}

\begin{block}\label{Phi-inv/par}
The restriction of $\Phi_n$ to the set of complete closed morphisms in $\H^r_n$ is
surjective, and this will follow from Theorem \ref{eq-alg-kirby/theo}, once it is
proved. Nevertheless, we sketch here the procedure which allows to find the
preimage under $\Phi_1$ of any Kirby diagram $K\in \hat\K_1$. Obviously, this has a
great practical value, since it allows to calculate invariants of 4-dimensional
2-handlebodies directly from their surgery presentation.
 
Let $K = \cup_j L_j$, where some of the $L_j$'s form a trivial link of dotted
unknots and the others are framed knots. Consider a planar diagram of $K$, where the
dotted components project trivially, in such a way that the crossings are either
between two framed components or between one framed and one dotted component. 
By changing a finite number of such crossings, say $C_1, \dots, C_l$, one can obtain
a new Kirby diagram $K' = \cup_j L_j'$ which is trivial as a link. 
Without loss of generality, $K'$ can be assumed to coincide with $K$ outside $E_1\cup
\dots \cup E_l$, where each $E_i$ is a cylinder projecting onto a small circular
neighborhood of $C_i$. One of such cylinders $E_i$, together with the relative portions
of the diagrams $K$ and $K'$, is depicted in Figure \ref{crossing1/fig} \(a) and \(b),
where $j$ and $k$ may or may not be distinct.

\begin{Figure}[htb]{crossing1/fig}{}{}
\centerline{\fig{}{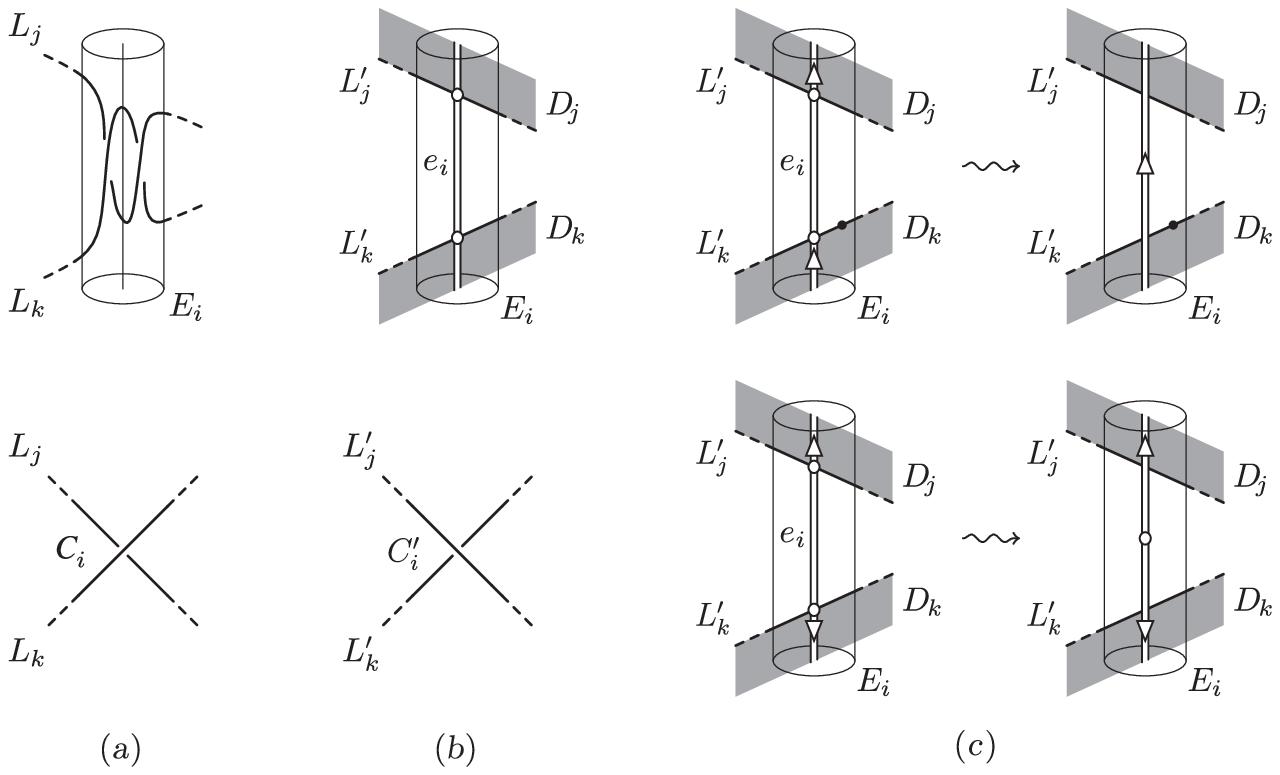}}
\end{Figure}

Now let $D = \cup_j D_j$ be a disjoint union of disks embedded in $R^3$ and bounded by
$\cup_j L_j'$. Embed in each $D_j$ a connected rooted uni/tri-valent tree whose
uni-valent vertices (different from the root) are the preimages of the $C_i'$'s in
$L_i'$, and the rest of the tree lies in the interior of the disk. Then connect through
a vertical (with respect to the projection plane) edge $e_i$ the two preimages of each
crossing $C_i'$ (cf. Figure \ref{crossing1/fig} \(b)). Orient all edges of the trees
in such a way that if $L_j'$ is a dotted component, each tri-valent vertex has one
incoming and two outgoing edges, and if $L_j'$ is a framed component, each tri-valent
vertex has two incoming and one outgoing edge. In this way one obtains a graph embedded
in $R^3$ whose bi-valent vertices correspond to the preimages of $C_i'$. Remove these
vertices from the graph. When the edges of the trees, attached to these preimages, are
oriented in a consistent way as shown on the top in Figure \ref{crossing1/fig} \(c),
this induces an orientation of the resulting edge. Otherwise divide this edge by a
single bi-valent vertex as shown on the bottom in Figure \ref{crossing1/fig} \(c). In
this way one obtains a graph $G$ embedded in $R^3$ whose bi-valent vertices stay above
those crossings $C_i'$ which involve two framed components.

\begin{Figure}[b]{Phi-inv/fig}{}
 {Surjectivity of $\Phi_1$ on the set of closed morphisms.}
\vskip-6pt\centerline{\fig{}{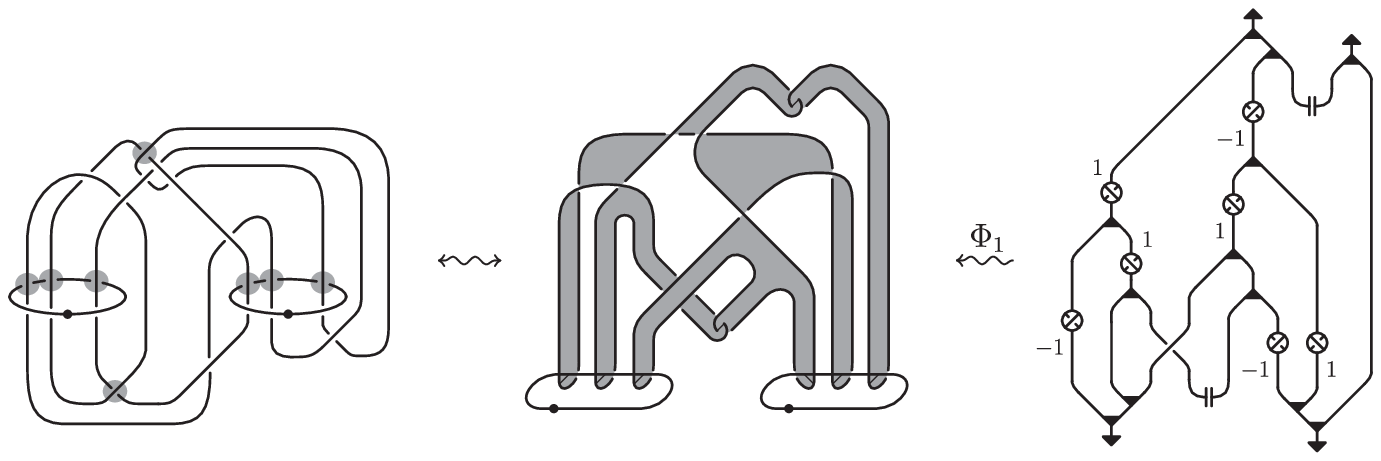}}\vskip-6pt
\end{Figure}		

At this point, consider a closed regular neighborhood $N_G$ of $G$ in $R^3$ extending
the regular neighborhood $E=\cup_iE_i$ of the edges $e_i$ and such that $J = N_G\cap D$
is a regular neighborhood of the trees in $D$. Let $I = \partial D \cap J$ be the set
of intervals in which  $J$ intersects the boundary of the  disks $D_j$. Then the link
$$K'' = (K\cap E)\cup \left(\partial J - \mbox{Int}\,I \right),$$
 is isotopic to $K$ through an isotopy which restricts to a deformation retraction on
each $D_j$. Isotope further $N_G$ (inducing an isotopy of G and $K''$) in such a way
that G is in regular position with respect to the projection plane, and in a
neighborhood of each root and tri-valent vertex the projection of the edges on the
$y$-axis is increasing. Then isotope $K''$ in $N_G$, to put it in regular position
with respect to the projection plane and transform its framing into blackboard
framing. Finally, observe that, up to isotopy and creations of canceling pairs of
1/2-handles, $K''$ is composed by the elementary diagrams on the right in Figures
\ref{defnPhi/fig} and \ref{defnPhi4/fig}, where the bi-valent vertices correspond to
the images of the copairing morphism. An example (with blackboard framing) is
presented in Figure \ref{Phi-inv/fig} where the crossings $C_i$ are encircled by a
small gray disk.
\end{block}

\section{The reduction map%
\label{reduction/sec}}
 
This section is dedicated to the proof of Theorem \ref{reduction/theo}, i.e. given
$m < n$, we construct a bijective map $\down_m^n: \hat\H^{r,c}_n \to \hat\H^{r,c}_m$
such that $\down_m^n \circ \Phi_n = \Phi_m \circ \down_m^n$. In other words, we
realize in $\hat\H^{r,c}_n$ (actually in $\hat\H^{r,c}(\G)$ for any $\G$) the
algorithm for reducing the labels (0-handles) of a generalized Kirby diagram
described in \ref{red-kirby/par}.

According to the algorithm in \ref{red-kirby/par}, in order to cancel the $n$-th
0-handle in a diagram $K$, one first separates the part of the diagram contained in
this 0-handle by ``pulling it up'', and then slides it over an 1-handle $U$ of label
$x = (i_0,n)$ with $i_0\neq n$. As a result the labels in the diagram change
from $j$ to $j^x$, where $\_^x: \G_n \to \G_{n-1}$ denotes the functor constructed in
Paragraph \ref{defnh}; in particular, $n^x = i_0$ and $j^x = j$ for any $j \neq n$. 
Eventually, one cancels the $n$-th 0-handle and $U$. The resulting diagram $K^U \in
\hat\K_{n-1}$ does not contain the label $n$. 

Observe that $K^U$ can be thought as obtained from $K$ through a sequence of the
following local changes:
\begin{itemize}\itemsep\smallskipamount
\item[{\sl a}\/)]\vskip-\lastskip\smallskip
 (if necessary) flip over in the projection plane any dotted components of label
$(i_0, n)$ in such a way that the upper face is labeled $n$;
\item[{\sl b}\/)] change the sign of all crossings at which a component labeled $n$
runs under a component labeled $i_0$; 
\item[{\sl c}\/)] change all labels from $j$ to $j^x$;
\item[{\sl d}\/)] cancel the component $U$. 
\end{itemize}\vskip-\lastskip\smallskip

The modifications {\sl a}\/) and {\sl b}\/) produce equivalent diagrams in $\hat\K_n$,
while {\sl c}\/) replaces any label $n$ with $i_0$, producing this way a diagram in
$\Re_\K^x(K) \in \hat\K_{n-1}$.\break Since $\Re_\K^x(K)$ is obtained from $K$ through
local changes which depend on $x$, but not on $U$, we may think of it as defining a
map $\Re_\K^x: \K_n \to \K_{n-1}$. Figure \ref{rh-kirby1/fig} shows
$\Re_\K^x(\Phi_n(\Delta_{(i_0, n)}))$, $\Re_\K^x(\Phi_n(S_{(i_0, n)}))$ and
$\Re_\K^x(\Phi_n(\bar S_{(n,i_0)}))$ (here we apply the move described in Figure
\ref{defnKF4/fig} to flip the dotted components). Observe that $\Re_\K^x$ respects the
monoidal structures of the two categories, a structure which is violated by {\sl d}\/).

\begin{Figure}[hbt]{rh-kirby1/fig}{}
 {($x \in \G(i_0,n)$ with $i_0 \neq n$)}
\centerline{\fig{}{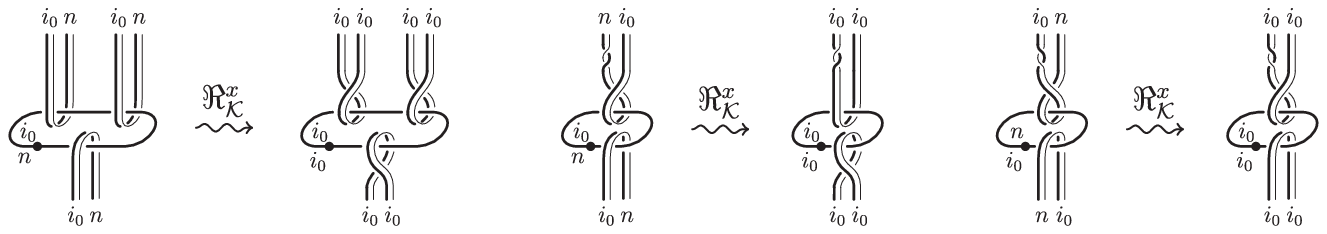}}\vskip-6pt
\end{Figure}

On the other hand, given a different morphism $y = (j_0,n)$ with $i_0\neq j_0\neq n$,
$\Re^x_\K(K)$ and $\Re^y_\K(K)$ may not be equivalent as diagrams in $\K_{n-1}$
even if $K$ is closed. Indeed, if $K$ is made out of a single 1-handle of label
$(n, n)$, then $\Re^x_\K(K)$ consists in a single 1-handle of index $(i_0, i_0)$,
while $\Re^y_\K(K)$ consists in a single 1-handle of index $(j_0, j_0)$. Those two
cannot be related unless it is present another 1-handle of label $y \bar x = (j_0,
i_0)$, when $\Re^x_\K(K)$ can be slided over this 1-handle obtaining $\Re^y_\K(K)$.
In particular, if $K$ is a complete morphism in $\hat\K_n^c$ and $V$ is any dotted
component of $K$ with label $y = (j_0,n)$, then $K^U$ can be transformed into $K^V$
through 1-handle slides (cf. Figure \ref{example-red1/fig}). Therefore the reduction
map ${\down_{n-1}^n}K = K^U$ is well defined on the subset of closed complete
morphisms.

\medskip

Now we want to mimic the above reduction procedure in $\H^r(\G)$ where $\G$ an
arbitrary groupoid. Notice that in this case there is no natural order on $\Obj \G$, so
the indices $i_0$ and $n$ above will be replaced by two (possibly coinciding) objects
$i_0,k_0 \in \Obj\G$. 

Given $x \in \G(i_0,k_0)$, we extend the functor $\_^x :\G \to \G$ defined in Paragraph
\ref{defnh} to a map $\_^x : \Obj\H^r(\G) \to \Obj\H^r(\G)$, by requiring that $(A
\diam B)^x = A^x \diam B^x$ for any $A,B \in \Obj\H^r(\G)$.

\smallskip

Our first goal is to prove the following theorem, where $\G^{\bs i}$ denotes the full
subgroupoid of $\G$ with $\Obj \G^{\bs i} = \Obj \G - \{i\}$ (cf. Paragraph
\ref{defnh}) and $\H^r(\G^{\bs i}) \subset \H^r(\G)$ are the universal ribbon Hopf
algebras constructed respectively on $\G^{\bs i}$ and $\G$, with the inclusion given
by Proposition \ref{formalext/theo}.

\begin{theorem} \label{functR/theo} 
 Let $\G$ be a groupoid. For any $x \in \G(i_0,k_0)$ there exists a functor
$$\Re^x: \H^r(\G) \to \H^r(\G)$$
 which coincides with $\_^x: \Obj \H^r(\G) \to \Obj \H^r(\G)$ on the objects, such
that:
\begin{itemize}\itemsep\smallskipamount
\item[{\sl a}\/)]\vskip-\lastskip\smallskip
given any other $y \in \G(j_0,k_0)$ and any object $A$ in
$\H^r(\G)$ there exists an invertible morphism
$$\xi_A^{x,y}: H_{y\bar x} \diam A^x \to H_{y\bar x} \diam A^y,$$\vskip6pt\noindent
	such that for any morphism $F: A \to B$ in $\H^r(\G)$ (cf. Figure \ref{natuR/fig}):
$$\xi_B^{x,y} \circ (\id_{y\bar x} \diam \Re^x(F)) = (\id_{y\bar x} \diam 
\Re^y(F)) \circ \xi_A^{x,y}\text{;}
\eqno{\(q1)}$$
\item[{\sl b}\/)] $\Re^x$ restricts to the identity on $\H^r(\G^{\bs k_0})$
and to an equivalence of categories $\Re^x_| : \H^r(\G^{\bs i_0}) \to \H^r(\G^{\bs
k_0})$ whose inverse is $\Re^{\bar x}_|: \H^r(\G^{\bs k_0}) \to \H^r(\G^{\bs i_0})$;
\item[{\sl c}\/)] if $i_0\neq k_0$ then $\Re^x(\H^r(\G)) \subset \H^r(\G^{\bs k_0})$,
hence we have a retraction functor
$$\Re^x:\H^r(\G) \to \H^r(\G^{\bs k_0}).$$
\end{itemize}
Moreover, $\Re^x$ induces functors on the quotient categories
$$\partial^\star\Re^x:\partial^\star\H^r(\G) \to \partial^\star\H^r(\G)
\quad \text{and} \quad \partial\Re^x:\partial\H^r(\G) \to \partial\H^r(\G).$$
\end{theorem}\vskip-18pt

We observe that in general $\Re^x$ and $\Re^y$ are not naturally equivalent, but 
{\sl a}\/) can be interpreted as a weaker version of that.

In the next corollary we isolate two particular cases of \(q1) which will be
relevant for our purposes.

\begin{corollary} \label{commutingxi/theo}
 For any $x \in \G(i_0,k_0)$ and any morphism $F: A \to B$ in $\H^r(\G)$, we have
 (cf. Figure \ref{natuR/fig}):
\vskip-6pt
$$\xi_B^{1_{k_0\!},x} \circ (\id_x \diam F) = (\id_x \diam \Re^x(F))
\circ \xi_A^{1_{k_0\!},x},
\eqno{\(q2)}$$
$$\xi_B^{x,x} \circ (\id_{1_{i_0}\!} \diam \Re^x(F)) = (\id_{1_{i_0}\!}
\diam \Re^x(F)) \circ \xi_A^{x,x}.
\eqno{\(q3)}$$\vskip-12pt
\end{corollary}

{\sc Definition of $\xi$.} As it is clear from the discussion above, $\xi$ should
be the algebraic analog of sliding over 1-handle or pushing through a dotted
component. The definition itself is somewhat heavy and consists in the following
steps. Given
\mypagebreak
$x \in \G(i_0,k_0)$ and $y \in \G(j_0,k_0)$, we define $\zeta_A^{x,y}:
H_{y\bar x} \diam A^x \to A^y$ for any object $A$ in $\H^r(\G)$ inductively by the
following identities (cf. Figure \ref{defnzeta/fig}) where $g \in \G(i,j)$:

$$\arraycolsep0pt\kern-6pt
\zeta_{H_g}^{x,y} = \zeta_g^{x,y} = \left\{
\begin{array}{ll} 
 \epsilon_{y\bar x} \diam \id_{g^x} &\ \ \text{if } i \neq k_0 \neq j,\\[2pt]
 m_{y\bar x,g^x} &\ \ \text{if } i = k_0 \neq j,\\[2pt]
 m_{g^x,x\bar y} \circ (\id_{g^x} \diam S_{y\bar x}) \circ 
 \gamma_{y\bar x,g^x} &\ \ \text{if } i \neq k_0 = j,\\[2pt]
 m_{y\bar x,g^x x\bar y} \circ (\id_{y\bar x} \diam m_{g^x,x\bar y}) \circ {}\\[2pt]
 \quad {} \circ (\id_{y\bar x} \diam ((\id_{g^x} \diam S_{y\bar x}) \circ
 \gamma_{y\bar x,g^x})) \circ (\Delta_{y\bar x} \diam \id_{g^x}) 
 &\ \ \text{if } i = k_0 = j;
\end{array}\right.
\eqno{\(q4)}
$$
$$
 \zeta_{A  \diam H_g}^{x,y} = (\zeta_A^{x,y} \diam \zeta_g^{x,y})
 \circ (\id_{y\bar x} \diam \gamma_{y\bar x,A^x} \diam 
 \id_{g^x}) \circ (\Delta_{y\bar x} \diam \id_{A^x \diam g^x}).
 \eqno{\(q5)}
$$

\begin{Figure}[htb]{defnzeta/fig}{}
 {($x \in \G(i_0,k_0)$, $y \in \G(j_0,k_0)$, $g \in \G(i,j)$ and $A = H_{g_1}\! \diam
  H_{g_2}\! \diam \dots \diam H_{g_r}$).\kern-10pt}
\vskip-3pt\centerline{\ \fig{}{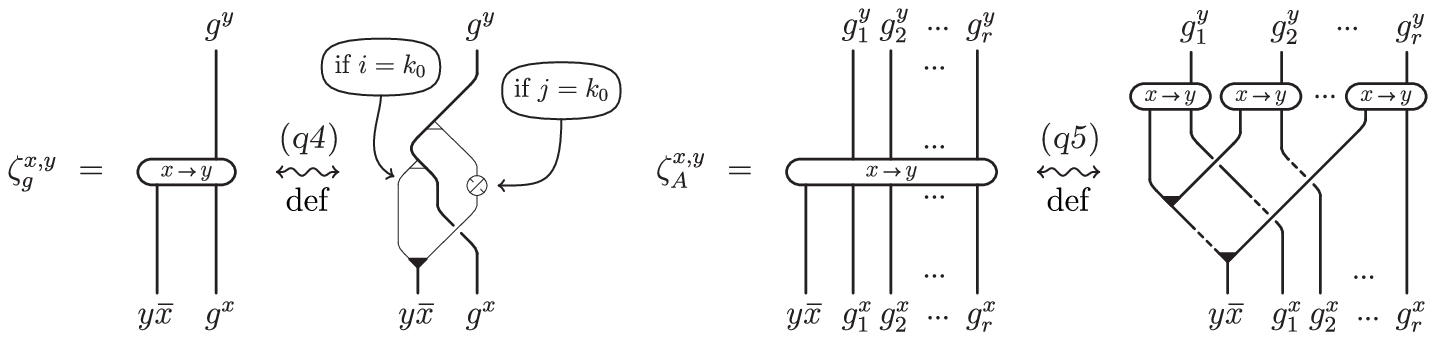}}\vskip-6pt
\end{Figure}

\begin{proposition} \label{actzeta/theo}
	If $x \in \G(i_0,k_0)$, $y \in \G(j_0,k_0)$ and $z \in \G(l_0,k_0)$, then for any
 object $A$ in $\H^r(\G)$ we have (cf. Figure \ref{actzeta/fig}):
$$\zeta^{y,z}_A \circ (\id_{z \bar y} \diam \zeta^{x,y}_A) = \zeta^{x,z}_A \circ
(m_{z \bar y, y \bar x} \diam \id_{A^x}): H_{z \bar y} \diam H_{y \bar x} \diam H_{A^x}
\to H_{A^z}.\kern-2pt
\eqno{\(q6)}$$
 In particular, $\zeta^{1_i,1_i}_{1_i}\!:\! H_{1_i} \!\diam\! H_{1_i}\! \to H_{1_i}$
makes $H_{1_i}$ into a $H_{1_i\!}$-module for any $i \in \Obj \G$.
\end{proposition}\vskip-6pt

\begin{Figure}[htb]{actzeta/fig}{}
 {($x \in \G(i_0,k_0)$, $y \in \G(j_0,k_0)$ and $z \in \G(l_0,k_0)$)}
\vskip-3pt\centerline{\fig{}{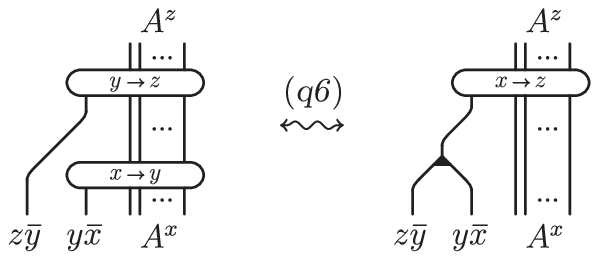}}\vskip-6pt
\end{Figure}

\begin{Figure}[htb]{actzeta-pf/fig}{}
 {Proof of \(q6) when $A = H_g$ with $g \in \G(k_0,k_0)$
 [{\sl a}/\pageref{algebra/fig}, {\sl s}/\pageref{pr-antipode/fig}].}
\vskip-6pt\centerline{\fig{}{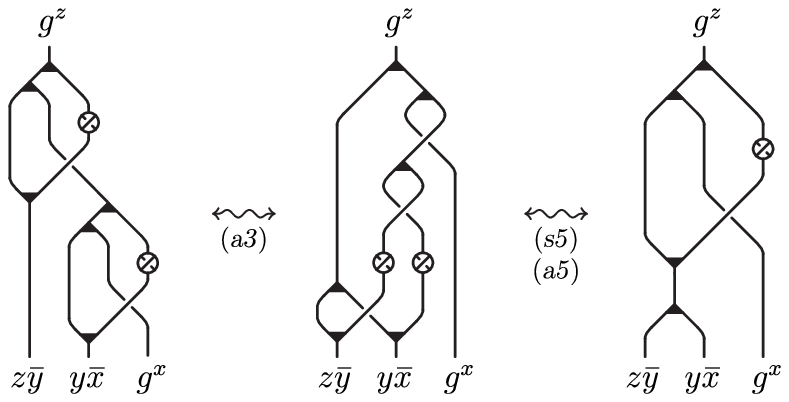}}\vskip-6pt
\end{Figure}

\begin{Figure}[htb]{actzeta-pfind/fig}{}
 {Proof of \(q6) -- the inductive step 
 [{\sl a}/\pageref{algebra/fig}, 
  {\sl q}/\pageref{defnzeta/fig}-\pageref{actzeta/fig}]}
\centerline{\fig{}{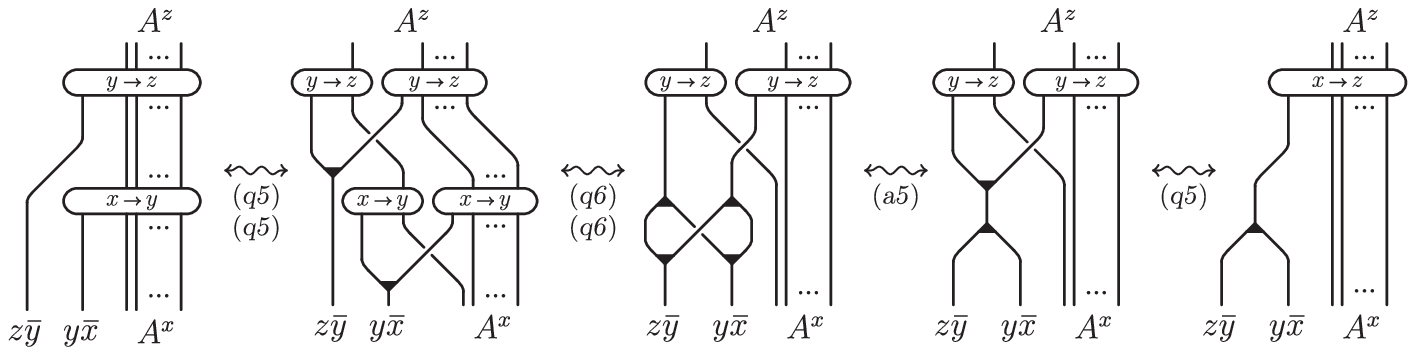}}\vskip-6pt
\end{Figure}

\begin{proof} 
 First of all, let us consider the special case when $A = H_g$ with $g \in \G(i,j)$. If
$j \neq k_0$, the statement is equivalent to \(a3) or \(a6) in Figure
\ref{algebra/fig}. If $j = k_0$ and $i \neq i_0$, it reduces to the antipode property
\(s3) in Figure \ref{pr-antipode/fig}. Finally, the case when $i = j = k_0$ is shown in
Figure \ref{actzeta-pf/fig}. At this point, the general case follows by the inductive
argument shown in Figure \ref{actzeta-pfind/fig}.
\end{proof}

Now we define the morphism $\xi_A^{x,y}$ and its inverse (cf. Figure
\ref{defnxi/fig}):
\vskip-6pt
$$\xi_A^{x,y} = (\id_{y\bar x} \diam \zeta_A^{x,y}) \circ (\Delta_{y\bar x} \diam
\id_{A^x}): H_{y\bar x} \diam A^x \to H_{y\bar x} \diam A^y,$$
$$(\xi_A^{x,y})^{-1} = (\id_{y\bar x} \diam \zeta_A^{y,x}) \circ (((\id_{y\bar x}
\diam S_{y\bar x}) \circ \Delta_{y\bar x}) \diam \id_{A^y}).$$

\begin{Figure}[htb]{defnxi/fig}{}{}
\vskip-6pt\centerline{\fig{}{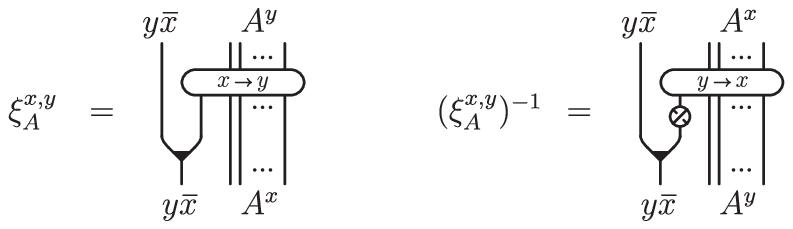}}\vskip-6pt
\end{Figure}

The proof that $(\xi_A^{x,y})^{-1} \!\circ \xi_A^{x,y} = \id_{y \bar x} \diam
\id_{A^x}$ is given in Figure \ref{xiinv-pf/fig}. In a similar way, using \(s1')
instead of \(s1), one can prove that $\xi_A^{x,y} \circ (\xi_A^{x,y})^{-1} =
\id_{y \bar x} \diam \id_{A^x}$.

\begin{Figure}[htb]{xiinv-pf/fig}{}
 {Invertibility of $\xi^{x,y}_A$
 [{\sl a}/\pageref{algebra/fig},
  {\sl q}/\pageref{actzeta/fig},
  {\sl s}/\pageref{antipode/fig}-\pageref{pr-antipode/fig}].}
\centerline{\fig{}{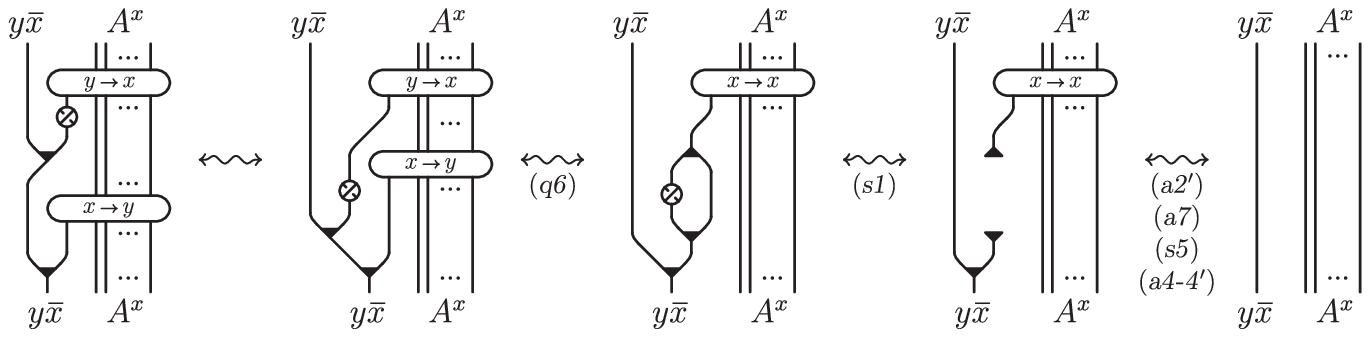}}\vskip-6pt
\end{Figure}

With the notation introduced in Figure \ref{defnxi/fig}, the relations of Theorem
\ref{functR/theo} and Corollary \ref{commutingxi/theo} look like in Figure
\ref{natuR/fig}.

\begin{Figure}[htb]{natuR/fig}{}
 {($x \in \G(i_0,k_0)$, $y \in \G(j_0,k_0)$, $F: A \to B$)}
\centerline{\fig{}{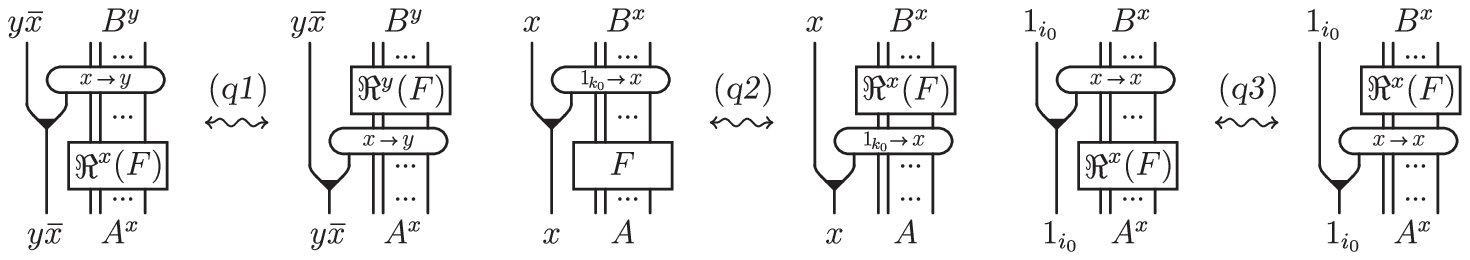}}\vskip-6pt
\end{Figure}

\medskip

We will need also the morphisms
\vskip-6pt
$$\bar{\zeta}_A^{x,y} = \zeta_A^{x,y} \circ \bar \gamma_{A^x,y\bar x}: A^x \diam 
H_{y\bar x} \to A^y,$$
$$\bar{\xi}_A^{x,y} = (\bar{\zeta}_A^{x,y} \diam \id_{y\bar x}) \circ (\id_{A^x}
\diam \Delta_{y\bar x}): A^x \diam H_{y\bar x} \to A^y \diam H_{y\bar x},$$
which are obtained from $\zeta_A^{x,y}$ and $\xi_A^{x,y}$ by simply ``pulling''
the string labeled $y\bar x$ to the right as shown in Figure \ref{defnxi-bar/fig}.
$\bar{\xi}_A^{x,y}$ is invertible with
$$(\bar{\xi}_A^{x,y})^{-1} =
(\bar{\zeta}_A^{y,x} \diam \id_{y\bar x}) \circ (\id_{A^y} \diam ((\bar S_{y\bar x}
\diam \id_{y\bar x}) \circ \Delta_{y\bar x})),$$
	as shown in Figure \ref{defnxi-bar/fig}.

\begin{Figure}[htb]{defnxi-bar/fig}{}{}
\vskip-4pt\centerline{\fig{}{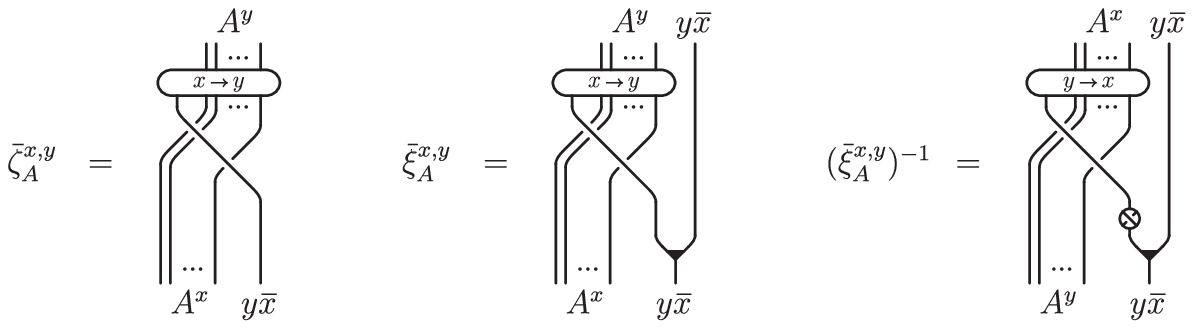}}\vskip-4pt
\end{Figure}

\begin{Figure}[b]{xi-kirby/fig}{}
 {$\Phi_n(\zeta_{g \diam h}^{x,y})$ for $x = (i_0,k_0)$, $y = (j_0,k_0)$, 
 $g = (k_0,k_0)$, $h = (k_0,i)$, $i \neq k_0$.}
\vskip-4pt\centerline{\fig{}{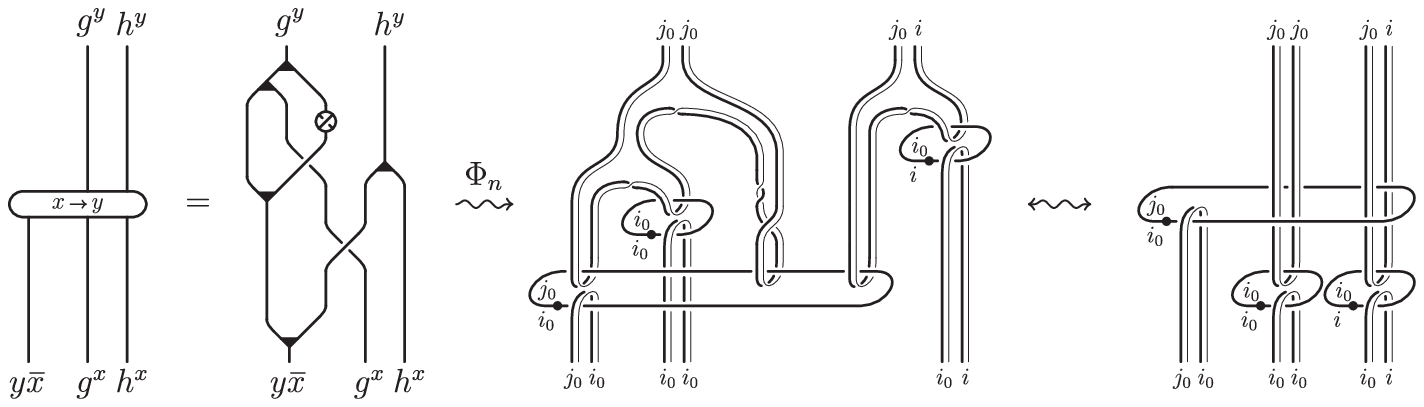}}\vskip-6pt
\end{Figure}

In order to see that $\zeta_A^{x,y}$ and $ \xi_A^{x,y}$ are the algebraic analog of
sliding over 1-handles, in Figure \ref{xi-kirby/fig} we represent the image 
$\Phi_n(\zeta_{g \diam h}^{x,y})$ in $\K_n$ when $x = (i_0,k_0)$, $y = (j_0,k_0)$, 
$g = (k_0,k_0)$ and $h = (k_0,i)$ with $i\neq k_0$ are morphisms in $\G_n$. It should
be obvious now that in terms of Kirby diagrams, $\zeta_A^{x,y}$ corresponds to a dotted
component $U$ of label $y \bar x = (j_0,i_0)$ embracing all framed strings of label
$k_0^x = i_0$ and transforming these labels into $k_0^y = j_0$. The reader can see
that $\Phi_n(\xi_A^{x,y})$ is the same, but there is an extra framed arc passing
through $U$, which is a reflected image of the left most one of labels $(i_0,j_0)$
with respect to a horizontal diameter of $U$.

Now we can interpret \(q2) and \(q3) in Figure \ref{natuR/fig} in terms of Kirby
tangles. In particular, \(q2) implies that $\Re^x_\K(F)$ is obtained from $F$ by
pushing up the part contained in the $k_0$'th 0-handle through a dotted component of
label $x = (i_0,k_0)$, changing this way the index $k_0$ to $i_0$. On the other hand,
\(q3) states that if we have 
\mypagebreak
already done this, the same part can be pushed through a
dotted component of label $(i_0,i_0)$ without further changes.

\medskip

{\sc Definition of $\Re^x$.} Let $x \in \G(i_0,k_0)$ as above. For $i_0 = k_0$ we
define $\Re^x$ to be the formal extension of $\_^x$ to $\H^r(\G)$, in the sense that
$\Re^x(\Delta_g) = \Delta_{g^x}$, $\Re^x(\gamma_{g,h}) = \gamma_{g^x,h^x}$,
$\Re^x(S_g) = S_{g^x}$ and so on. 

However, when $i_0 \neq k_0$ such formal extension would run into problems. 
In fact the copairing $\sigma_{i,j}$ is defined to be trivial if $i \neq j$, so the
axiom \(r7) in Figure \ref{ribbon2/fig} would not be satisfied in the image if
$i^x = j^x$. Thus we need to make some corrections in order to have a functor.
Being $\Re^x$ the algebraic analog of $\Re^x_\K$, the nature of these corrections can
be understood by reviewing the definition of $\Re^x_\K$ and Figure
\ref{rh-kirby1/fig}. 

Namely, for $i_0 \neq k_0$ we define $\Re^x:\H^r(\G) \to \H^r(\G)$ by putting, for any
$i \in \Obj \G$ and any $g,h \in \G$ (assumed to be composable in $\G$ when dealing
with $m_{g,h}$):
\vskip-6pt
$$\arraycolsep0pt
\begin{array}{l}
 \Re^x(\Delta_g) = \left\{%
 \begin{array}{ll} 
  \mu^{-1}_{g^x,g^x} \circ \Delta_{g^x} &\quad\text{if } g \in \G(i_0,k_0),\\[2pt]
 	\Delta_{g^x} &\quad\text{otherwise};
 \end{array}	\right.\\[16pt]
 \Re^x(\gamma_{g,h}) = \left\{%
 \begin{array}{l}
  \begin{array}{ll}
   \gamma_{g^x\!,h^x} &\quad\text{if } g \in \G(k_0,k_0),\\[4pt]
   (\bar{\zeta}^{x,x}_h \diam \id_{g^x}) \circ
    (\id_{h^x} \diam \rho^l_{g^x\!,i_0}) \circ \gamma_{g^x\!,h^x}
    &\quad\text{if } g \in \G(i,k_0) \text{, } i \neq k_0,\\[3pt]
   \gamma_{g^x\!,h^x} \circ (\id_{g^x} \diam \zeta^{x,x}_h) \circ {}\\[2pt]
    \kern4mm {} \circ (\id_{g^x} \diam S_{1_{i_0}} \diam \id_{h^x}) \circ
    (\rho^r_{g^x\!,i_0} \diam \id_{h^x})
    &\quad\text{if } g \in \G(k_0,j) \text{, } j \neq k_0,
  \end{array}\\[31pt]
  (\bar{\zeta}^{x,x}_h \diam \id_{g^x}) \circ
   (\id_{h^x} \diam \rho^l_{g^x\!,i_0}) \circ \gamma_{g^x\!,h^x} \circ {}\\[2pt]
   \kern4mm {} \circ (\id_{g^x} \diam \zeta^{x,x}_h) \circ (\id_{g^x} \diam 
   S_{1_{i_0}} \diam \id_{h^x}) \circ (\rho^r_{g^x\!,i_0} \diam \id_{h^x})
   \quad\text{otherwise;}
 \end{array}\right.\\[54pt]
 \Re^x(\bar \gamma_{g,h}) = \left\{%
 \begin{array}{l}
  \begin{array}{ll}
   \bar \gamma_{g^x\!,h^x} &\quad\text{if } h \in \G(k_0,k_0),\\[4pt]
   (\id_{h^x} \diam \zeta^{x,x}_g) \circ
    (\rho^r_{h^x\!,i_0} \diam \id_{g^x}) \circ \bar \gamma_{g^x\!,h^x}
    &\quad\text{if } h \in \G(k_0,l) \text{, } l \neq k_0,\\[3pt]
   \bar \gamma_{g^x\!,h^x} \circ (\bar{\zeta}^{x,x}_g \diam \id_{h^x}) \circ {}\\[2pt]
    \kern4mm {} \circ (\id_{g^x} \diam \bar S_{1_{i_0}} \diam \id_{h^x}) 
    \circ (\id_{g^x} \diam \rho^l_{h^x\!,i_0})
    &\quad\text{if } h \in \G(k,k_0) \text{, } k \neq k_0,
   \end{array}\\[31pt]
  (\id_{h^x} \diam \zeta^{x,x}_g) \circ
   (\rho^r_{h^x\!,i_0} \diam \id_{g^x}) \circ \bar \gamma_{g^x\!,h^x} \circ {}\\[2pt]
   \kern4mm {} \circ (\bar{\zeta}^{x,x}_g \diam \id_{h^x}) \circ (\id_{g^x}
   \diam \bar S_{1_{i_0}} \diam \id_{h^x}) \circ (\id_{g^x} \diam
   \rho^l_{h^x\!,i_0})
   \quad\text{otherwise;}
 \end{array}\right.\\[54pt]
 \Re^x(S_g) = \left\{%
 \begin{array}{ll} 
  \bar S_{g^x} v_{g^x}^2 &\quad\text{if } g \in \G(i_0,k_0),\\[2pt]
  S_{g^x} &\quad\text{otherwise};
 \end{array}	\right.\\[16pt]
 \Re^x(\bar S_g) = \left\{%
 \begin{array}{ll} 
  S_{g^x} v_{g^x}^{-2} &\quad\text{if } g \in \G(k_0,i_0),\\[2pt]
  \bar S_{g^x} &\quad\text{otherwise};
 \end{array}\right.\\[14pt]
 \begin{array}{ll}
  \hbox to 40mm{$\Re^x(\eta_i) = \eta_{i^x};$\hss}
  &\Re^x(\epsilon_g) = \epsilon_{g^x};\\[6pt]
  \hbox to 40mm{$\Re^x(l_i) = l_{i^x};$\hss}
  &\Re^x(L_g) = L_{g^x};\\[6pt]
  \Re^x(m_{g,h}) = m_{g^x\!,h^x};
  &\Re^x(v_g) = v_{g^x};
 \end{array}\\[28pt]
 \Re^x(\sigma_{i,j}) = \left\{%
 \begin{array}{ll} 
  \sigma_{i^x\!,j^x} &\quad\text{if $\{i, j\}\neq \{i_0, k_0\}$},\\[2pt]
  \eta_{i_0} \diam \eta_{i_0} &\quad\text{otherwise}.
 \end{array}\right. \label{imsigma/par}
\end{array}$$

\smallskip

$\Re^x(\Delta_g)$, $\Re^x(\gamma_{g,h})$ and $\Re^x(\bar \gamma_{g,h})$ are presented 
in Figures \ref{delta/fig} and \ref{cross/fig}. Moreover, up to equivalence in
$\H^r(\G)$, one can see that $\Re^x(\gamma_{g,A})$ and $\Re^x(\bar \gamma_{A,h})$ are
given by the same expressions for any object $A$ in $\H^r(\G)$ (cf. Figure
\ref{cross11/fig}). The proof of this statement for $\Re^x(\gamma_{g,A})$ with $A =
h_1 \diam h_2$ in the case when $g \in \G(i,j)$ with $i \neq k_0$ and $j \neq k_0$ is
illustrated in Figure
\ref{cross12/fig}. The other cases are simpler and the generalization to an arbitrary
$A$ is straightforward. The proof for $\Re^x(\bar \gamma_{A,h})$ is similar.

\begin{Figure}[htb]{delta/fig}{}
 {$\Re^x(\Delta_g)$ when $x \in \G(i_0,k_0)$ with $k_0 \neq i_0$ ($g \in \G(i,j)$) 
 [{\sl r}/\pageref{ribbon4/fig}].}
\centerline{\fig{}{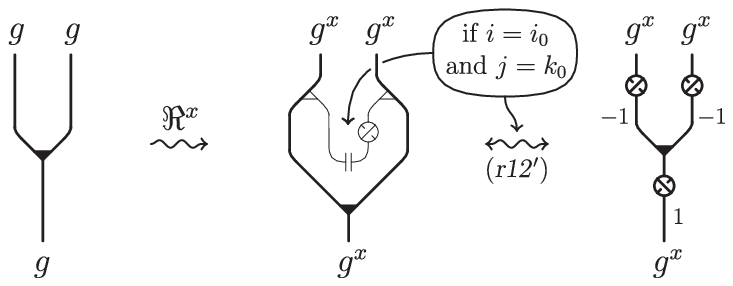}}\vskip-6pt
\end{Figure}

\begin{Figure}[htb]{cross/fig}{}
 {$\Re^x(\gamma_{g,h})$ and $\Re^x(\bar\gamma_{g,h})$ when $x \in
 \G(i_0,k_0)$ with $k_0 \neq i_0$ ($g \in \G(i,j)$ and $h \in \G(k,l)$) 
 [{\sl p}/\pageref{ribbon-isot/fig},
  {\sl r}/\pageref{ribbon2/fig}, 
  {\sl s}/\pageref{pr-antipode/fig}].}
\vskip-12pt\centerline{\kern5mm\fig{}{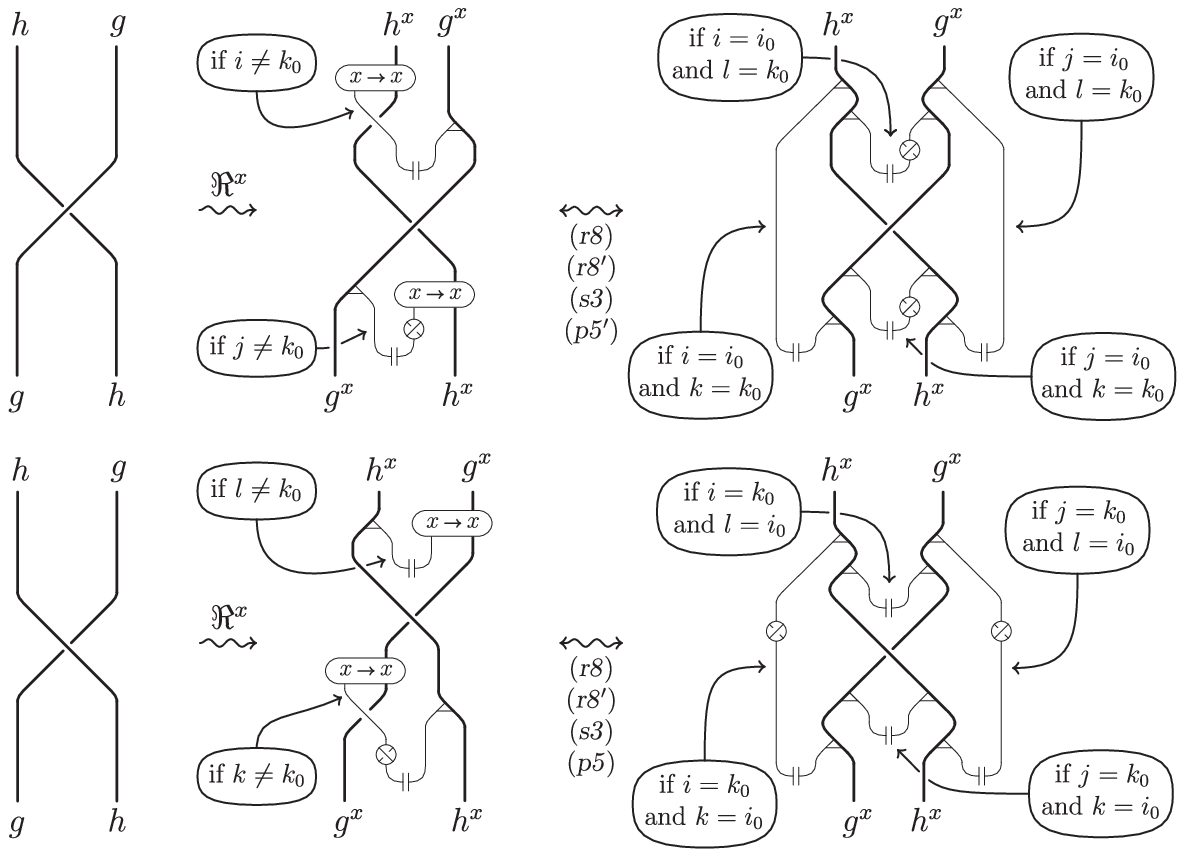}}\vskip-9pt
\end{Figure}

\begin{Figure}[htb]{cross11/fig}{}
 {$\Re^x(\gamma_{g,A})$ and $\Re^x(\bar\gamma_{A,h})$ when $x \in
 \G(i_0,k_0)$ with $k_0 \neq i_0$ ($g \in \G(i,j)$ and $h \in \G(k,l)$).}
\vskip-12pt\centerline{\ \fig{}{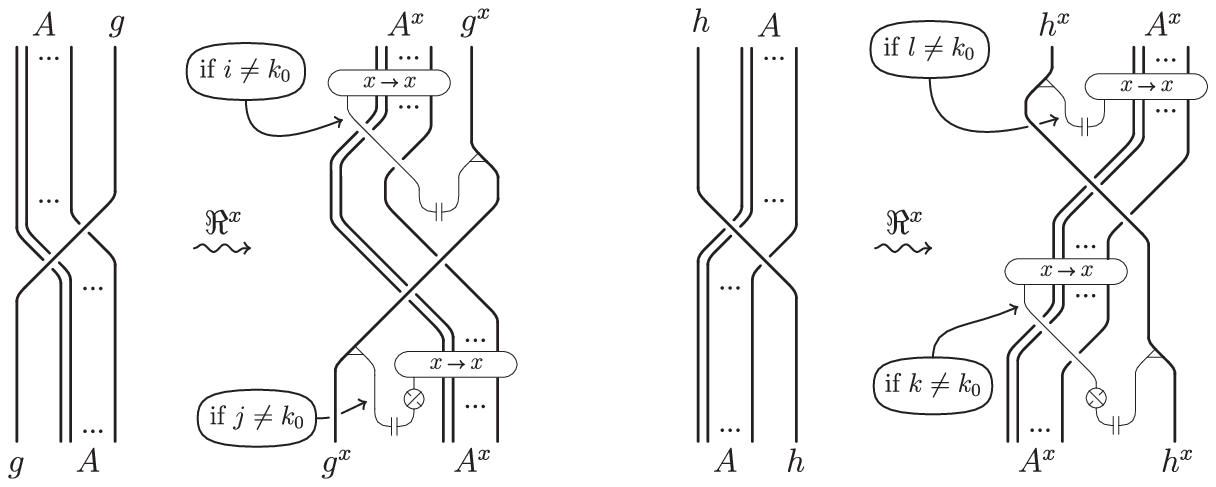}}\vskip-6pt
\end{Figure}

\begin{Figure}[htb]{cross12/fig}{}
 {($x \in \G(i_0,k_0)$, $k_0 \neq i_0$, $g \in \G(i,j)$, $i \neq k_0$, $j\neq k_0$)
 [{\sl q}/\pageref{defnzeta/fig},
  {\sl r}/\pageref{ribbon2/fig}, 
  {\sl s}/\pageref{pr-antipode/fig}].}
\centerline{\fig{}{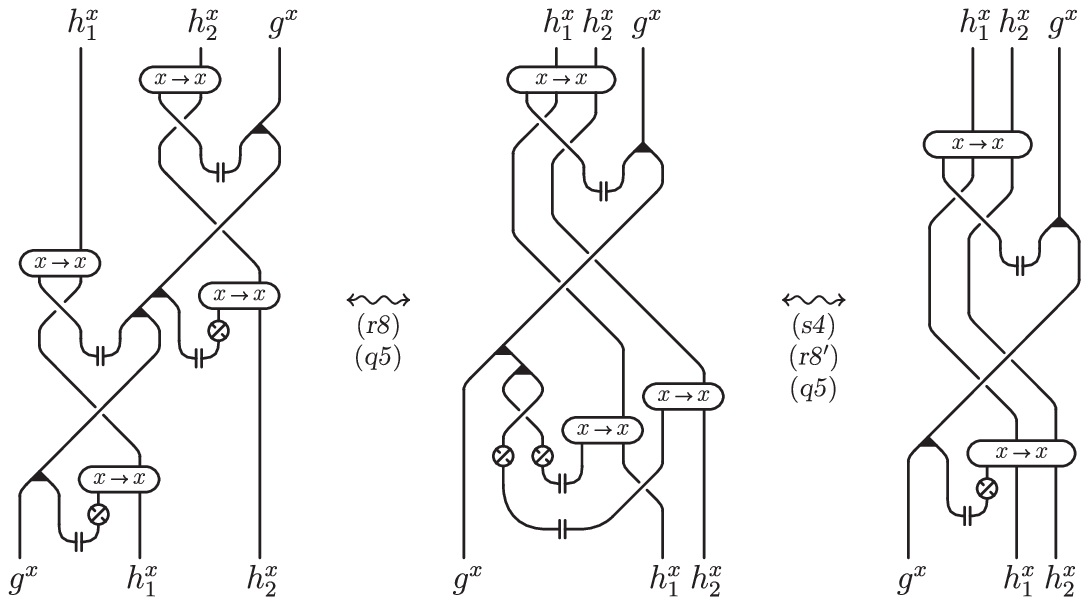}}\vskip-6pt
\end{Figure}

\begin{proof}[Theorem \ref{functR/theo}]
First of all we prove that \(q1) in {\sl a}\/) holds for any $x \in \G(i_0,k_0)$,
$y \in \G(j_0,k_0)$ and any elementary morphism $F: A \to B$, with $\xi^{x,y}_A$ and
$\xi^{x,y}_B$ defined as above. More precisely, we will prove the identity
$$\zeta_B^{x,y} \circ (\id_{y\bar x} \diam \Re^x(F)) =\Re^y(F)\circ\zeta_A^{x,y}, 
\eqno{\(q1')}$$
 which, by composing with $\Delta_{y\bar x} \diam \id_{A^x}$, implies \(q1).

Notice that if $F$ is invertible then the relation \(q1') for $F^{-1}$ follows from the
one for $F$, once we know that $\Re^x(F^{-1})$ is the inverse of $\Re^x(F)$. It can
be easily verified that this is true when $F$ is $\gamma_{g,h}$ and $S_g$, by using
\(p2-2') in Figure \ref{coform-s/fig} and \(q6) in Figure \ref{actzeta/fig}. 
Moreover, since the relations \(r6) and \(r7) in Figure \ref{ribbon2/fig} are
trivially satisfied in the image, the identity \(q1') for $\sigma_{i,j}$ follows from
the ones for the other elementary morphisms. Therefore, it suffices to prove that
\(q1') holds for $F = \eta_{1_i}$, $l_{1_i}$, $\Delta_g$, $S_g$, $\epsilon_g$, $L_g$,
$v_g$, $m_{g,h},\gamma_{g,h}$.

\begin{description}\itemsep\smallskipamount
\item[\ms{\sl $F = \eta_i,l_i$}.]
If $i \neq k_0$ there is nothing to prove, being $\xi^{x,y}_{1_i}$ trivial and
$\Re^x(F) = \Re^y(F) = F$. The case when $i = k_0$ is shown in Figure \ref{invxi1/fig}.

\begin{Figure}[htb]{invxi1/fig}{}
 {Proof of \(q1') for $F = \eta_{k_0}$ and $F = l_{k_0}$
 ($x \in \G(i_0,k_0)$ and $y \in \G(j_0,k_0)$)
 [{\sl a}/\pageref{algebra/fig},
  {\sl i}/\pageref{unimodular/fig},
  {\sl f}/\pageref{h-tortile1/fig},
  {\sl s}/\pageref{antipode/fig}-\pageref{pr-antipode/fig}].}
\vskip-3pt\centerline{\fig{}{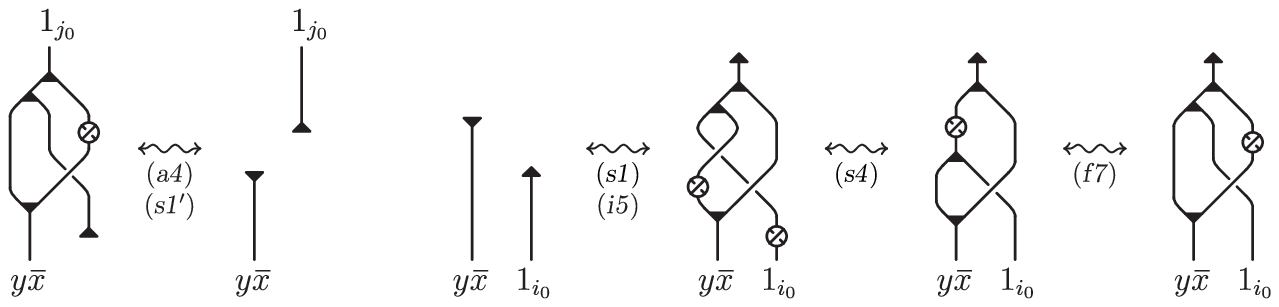}}\vskip-6pt
\end{Figure}

\begin{Figure}[htb]{invxi32/fig}{}
 {Proof of \(q1') for $F = \Delta_g$ with $g \in \G(i,k_0)$ 
 and $i \neq k_0$ ($x \in \G(i_0,k_0)$ and $y \in \G(j_0,k_0)$)
 [{\sl a}/\pageref{algebra/fig},
  {\sl r}/\pageref{ribbon1/fig},
  {\sl s}/\pageref{pr-antipode/fig}].}
\centerline{\fig{}{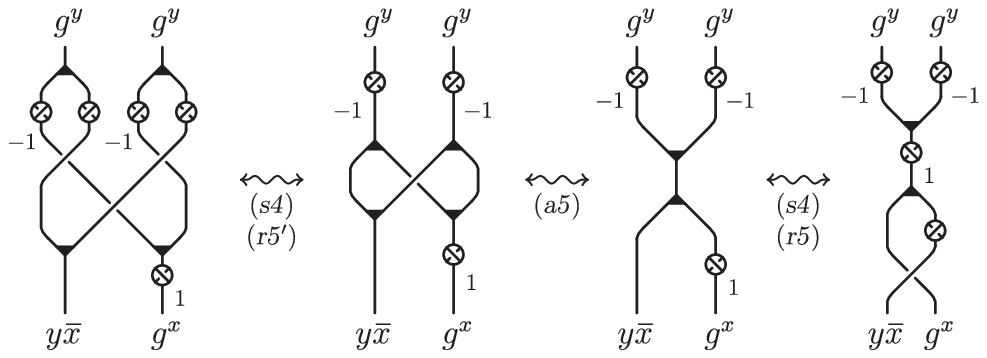}}\vskip-6pt
\end{Figure}

\begin{Figure}[htb]{invxi3/fig}{}
 {Proof of \(q1') for $F = \Delta_g$ with $g \in \G(k_0,k_0)$
 ($x \in \G(i_0,k_0)$ and $y \in \G(j_0,k_0)$)
 [{\sl a}/\pageref{algebra/fig}, 
  {\sl p}/\pageref{ribbon-tot/fig},
  {\sl r}/\pageref{ribbon2/fig},
  {\sl s}/\pageref{pr-antipode/fig},
 step \($\ast$) in Figure \ref{invxi4/fig}].}
\centerline{\fig{}{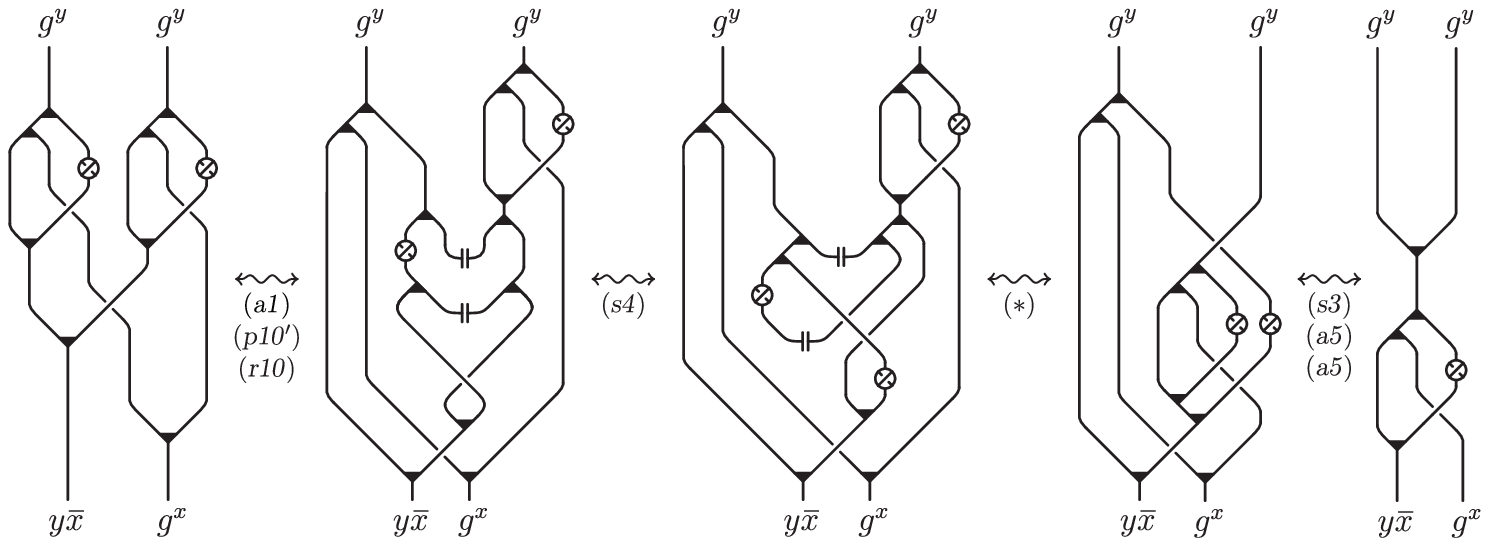}}\vskip-6pt
\end{Figure}

\begin{Figure}[htb]{invxi4/fig}{}
 {Proof of step $(\ast)$ in Figure \ref{invxi3/fig}
  [{\sl p}/\pageref{coform-s/fig}-\pageref{ribbon-tot/fig}, 
   {\sl q}/\pageref{actzeta/fig},
   {\sl r}/\pageref{ribbon2/fig},
   {\sl s}/\pageref{pr-antipode/fig}].}
\centerline{\fig{}{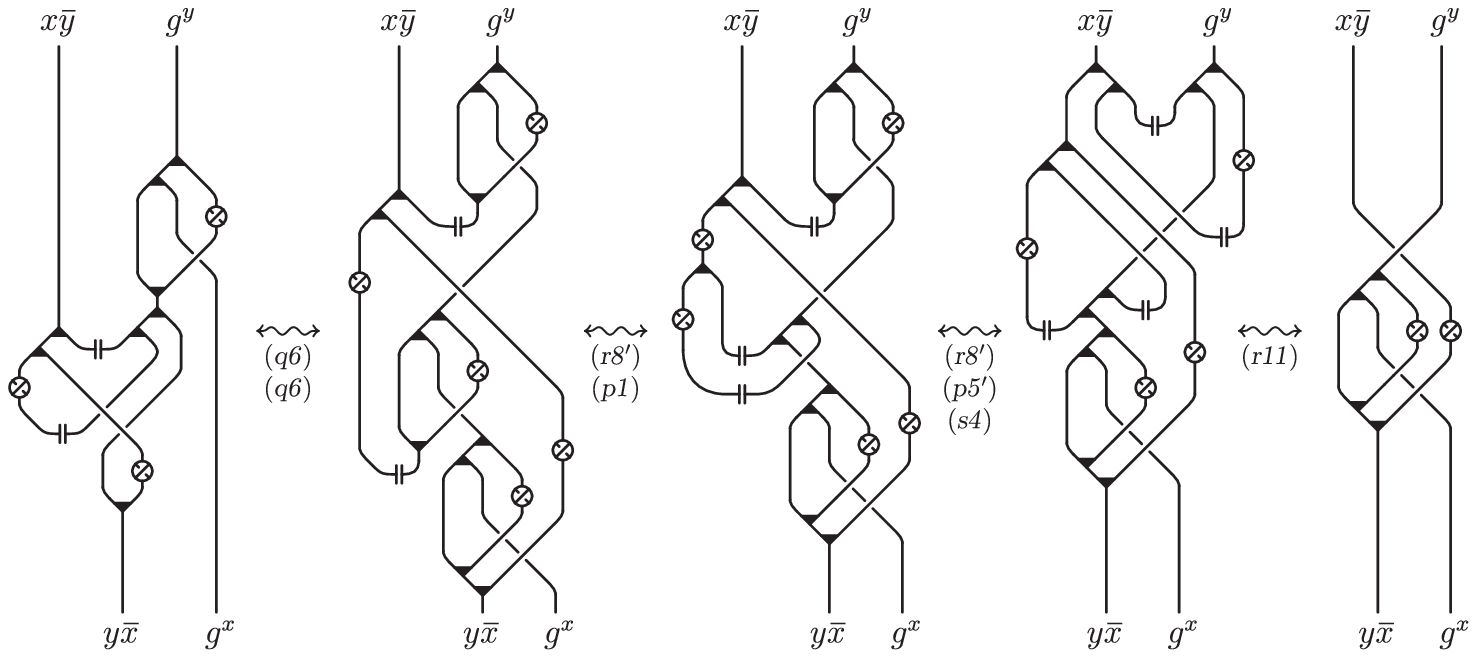}}\vskip-6pt
\end{Figure}

\item[\ms{\sl $F = \Delta_g$}.]
If $g \in \G(i,j)$ with $i,j \neq k_0$ there is nothing to prove, being $\xi^{x,y}_g$
trivial. For $g \in \G(k_0,j)$ with $j \neq k_0$, the statement essentially reduces to
the relation \(a5) in Figure \ref{algebra/fig}. The proof for $g \in \G(i,k_0)$ with 
$i \neq k_0$ is presented in Figure \ref{invxi32/fig}, where one can see that the last
diagram is exactly $\Re^y(\Delta_g)\circ\zeta_g^{x,y}$, by confronting it with Figure
\ref{delta/fig} when $i= i_0$, or applying \(t4) in Figure \ref{T-theo/fig} when
$i \neq i_0$. The case $g \in \G(k_0,k_0)$ is shown in Figures \ref{invxi3/fig} and
\ref{invxi4/fig}.

\item[\ms{\sl $F = S_g$}.]\vglue-15pt
As above, there is nothing to prove for $g \in \G(i,j)$ with $i,j \neq k_0$. When $g
\in \G(i,k_0)$ or $g \in \G(k_0,j)$ with $i,j \neq k_0$, the statement is equivalent
to the property \(s4) of the antipode in Figure \ref{pr-antipode/fig} (see Figure
\ref{invxi33/fig} for the former case). Figure~\ref{invxi2/fig} addresses the case of
$g \in \G(k_0,k_0)$.

\begin{Figure}[htb]{invxi33/fig}{}
 {Proof of \(q1') for $F = S_g$ with $g \in \G(i,k_0)$ and $i \neq k_0$
 ($x \in \G(i_0,k_0)$ and $y \in \G(j_0,k_0)$)
 [{\sl r}/\pageref{ribbon1/fig}, 
  {\sl s}/\pageref{pr-antipode/fig},
  {\sl t}/\pageref{repT/fig}].}
\centerline{\fig{}{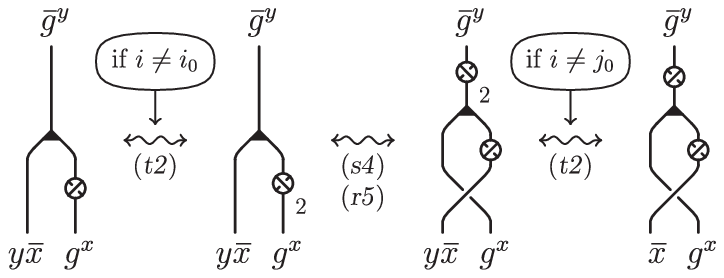}}\vskip-7pt
\end{Figure}

\begin{Figure}[htb]{invxi2/fig}{}
 {Proof of \(q1') for $F = S_g$ with $g \in \G(k_0,k_0)$
 ($x \in \G(i_0,k_0)$ and $y \in \G(j_0,k_0)$)
 [{\sl a}/\pageref{algebra/fig},
  {\sl p}/\pageref{coform-s/fig}-\pageref{ribbon-tot/fig}, 
  {\sl r}/\pageref{ribbon1/fig}-\pageref{ribbon4/fig}-\pageref{ribbon5/fig},
  {\sl s}/\pageref{pr-antipode/fig}].}
\vskip-6pt\centerline{\fig{}{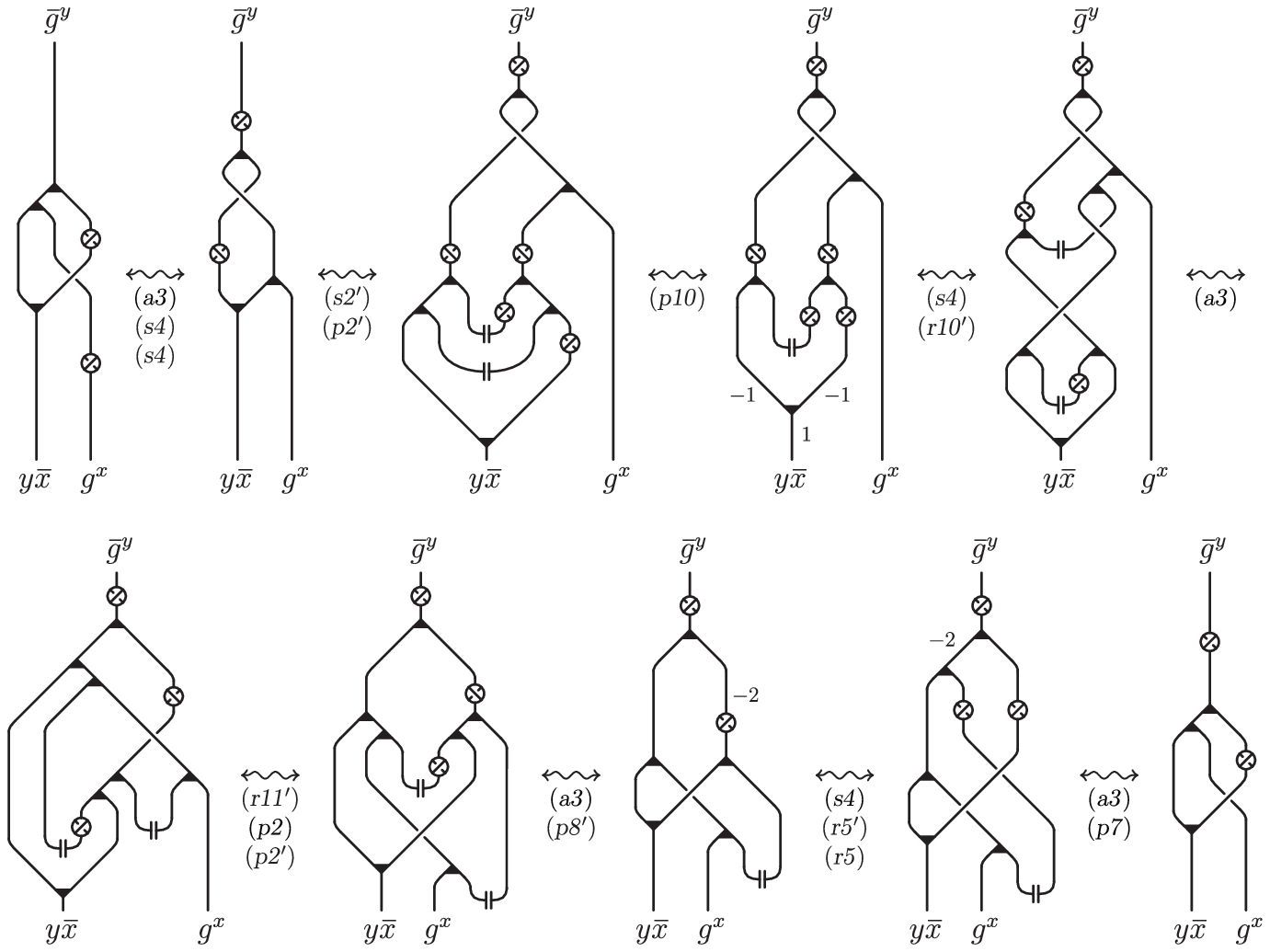}}\vskip-7pt
\end{Figure}

\item[\ms{\sl $F = \epsilon_g, L_g, v_g$}.]
The statements follow respectively from \(a6) in Figure \ref{algebra/fig}, \(i2-2') in
Figure \ref{unimodular/fig} and \(r5-5') in Figure \ref{ribbon1/fig}.

\item[\ms{\sl $F = m_{g,h}$}.]
See Figure \ref{invxi31/fig} for the case in which $g,h \in \G(k_0,k_0)$. The other
cases are simpler and anologous.

\begin{Figure}[htb]{invxi31/fig}{}
 {Proof of \(q1') for $F = m_{g,h}$ with $g,h \in \G(k_0,k_0)$
 ($x \in \G(i_0,k_0)$ and $y \in \G(j_0,k_0)$) 
 [{\sl a}/\pageref{algebra/fig}, {\sl s}/\pageref{antipode/fig},].}
\centerline{\fig{}{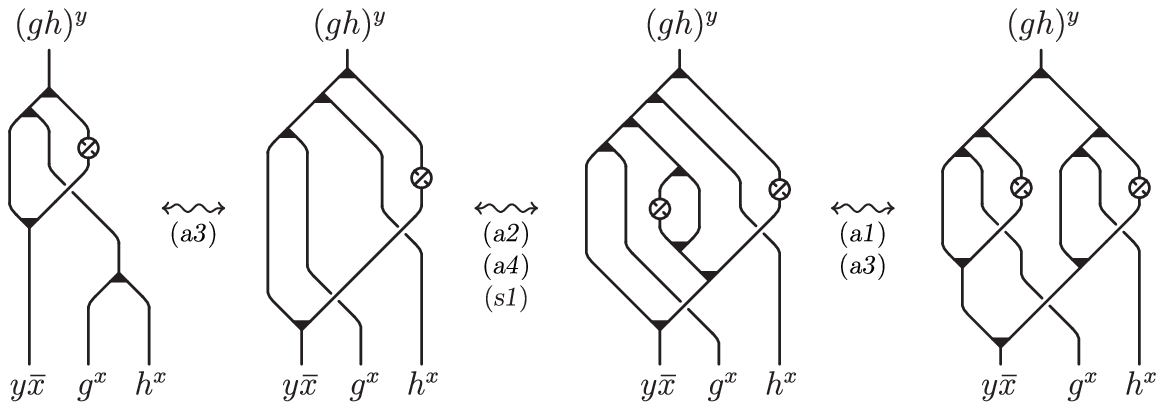}}\vskip-6pt
\end{Figure}

\begin{Figure}[b]{invxi7/fig}{}
 {Proof of \(q1') for $F = \gamma_{g,h}$ -- I
 ($g \in \G(i,j)$, $h \in \G(k,l)$, $x \in \G(i_0,k_0)$ and $y \in \G(j_0,k_0)$)
 [{\sl q}/\pageref{actzeta/fig}].}
\vskip-6pt\centerline{\fig{}{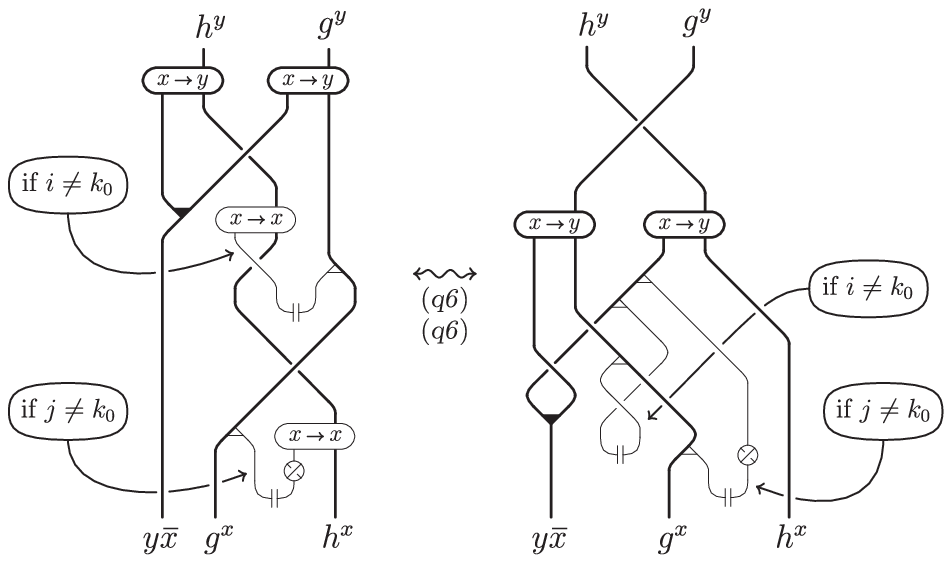}}\vskip-6pt
\end{Figure}

\begin{Figure}[htb]{invxi11/fig}{}
 {Proof of \(q1') for $F = \gamma_{g,h}$  with $g \in \G(i,j)$, $i \neq k_0$ and $j 
 \neq k_0$ -- II ($h \in \G(k,l)$, $x \in \G(i_0,k_0)$ and $y \in \G(j_0,k_0)$)
 [{\sl f}/\pageref{h-tortile1/fig},
  {\sl p}/\pageref{coform-s/fig}-\pageref{ribbon-isot/fig},
  {\sl q}/\pageref{actzeta/fig},
  {\sl r}/\pageref{ribbon5/fig},
  {\sl s}/\pageref{pr-antipode/fig},
  step \($\ast\ast$) follows from \(q1') for $F = S$ and $x,y = 1_{j_0}$ (cf. Figure
  \ref{invxi12/fig})].}
\centerline{\fig{}{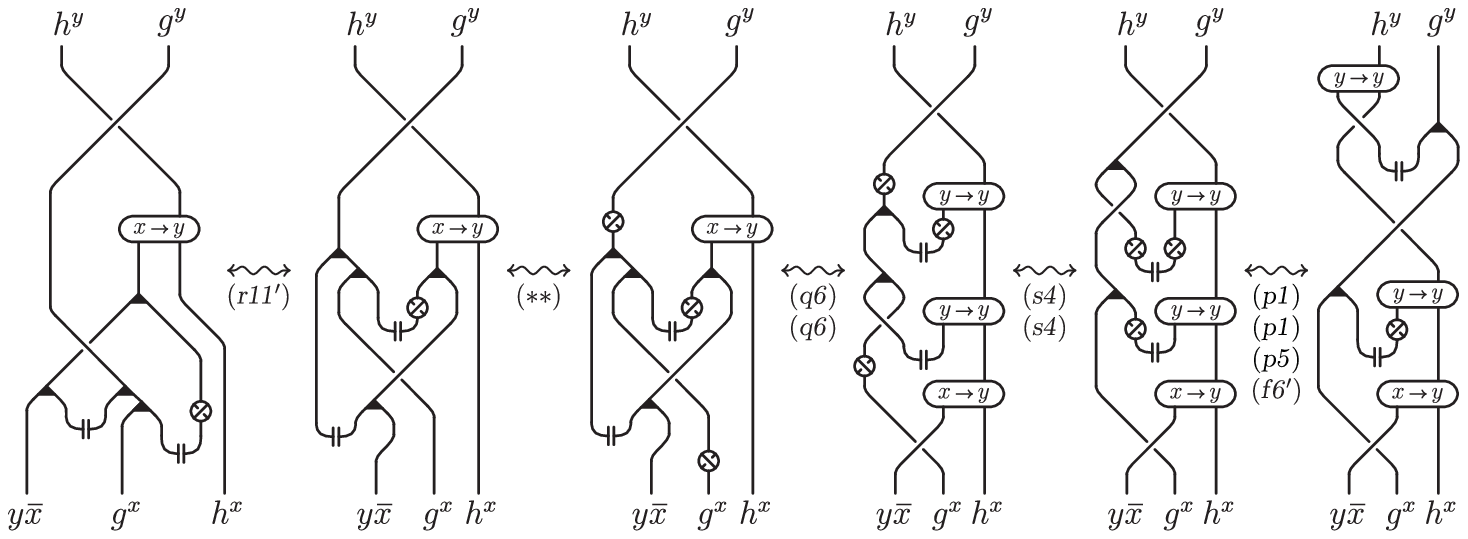}}\vskip-6pt
\end{Figure}

\begin{Figure}[htb]{invxi12/fig}{}
 {($g \in \G(i,j)$, $x \in \G(i_0,k_0)$, $y \in \G(j_0,k_0)$, $i \neq k_0$ and $j \neq
 k_0$)
 [{\sl a}/\pageref{algebra/fig},
  {\sl p}/\pageref{coform-s/fig},
  {\sl r}/\pageref{ribbon2/fig}].}
\centerline{\fig{}{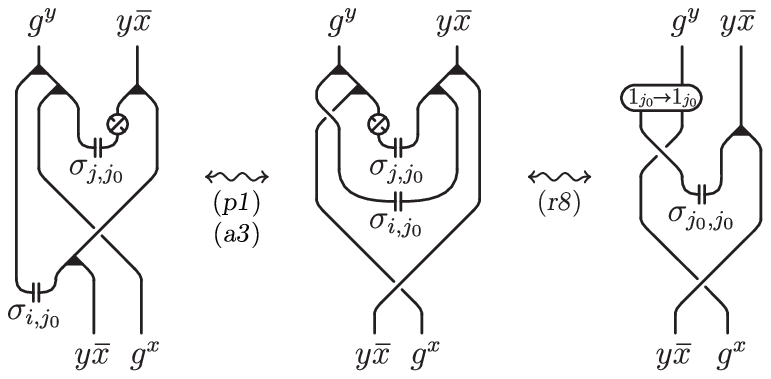}}\vskip-6pt
\end{Figure}

\begin{Figure}[htb]{invxi8/fig}{}
 {Proof of \(q1') for $F = \gamma_{g,h}$  with $g \in \G(k_0,j)$, $j \neq k_0$ -- II
 ($h \in \G(k,l)$, $x \in \G(i_0,k_0)$ and $y \in \G(j_0,k_0)$)
 [{\sl f}/\pageref{h-tortile1/fig},
  {\sl p}/\pageref{ribbon-tot/fig},
  {\sl q}/\pageref{actzeta/fig},
  {\sl r}/\pageref{ribbon2/fig}-\pageref{ribbon5/fig},
  {\sl s}/\pageref{pr-antipode/fig}].}
\centerline{\fig{}{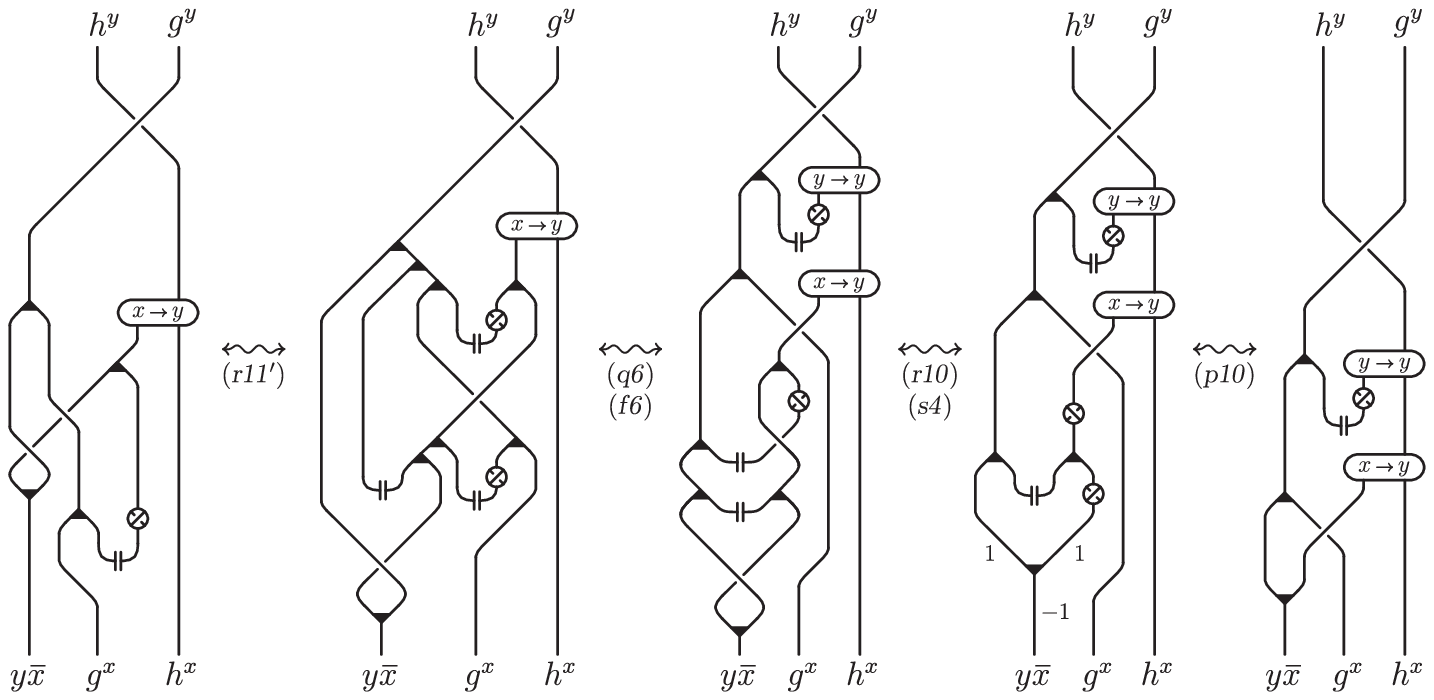}}\vskip-6pt
\end{Figure}

\begin{Figure}[htb]{invxi9/fig}{}
 {Proof of \(q1') for $F = \gamma_{g,h}$  with $g \in \G(i, k_0)$, $i \neq k_0$ -- II
 ($h \in \G(k,l)$, $x \in \G(i_0,k_0)$ and $y \in \G(j_0,k_0)$)
 [{\sl f}/\pageref{h-tortile1/fig},
 {\sl s}/\pageref{pr-antipode/fig},
 step \({$\ast{\ast}\ast$}) as in Figure \ref{invxi8/fig}].}
\centerline{\fig{}{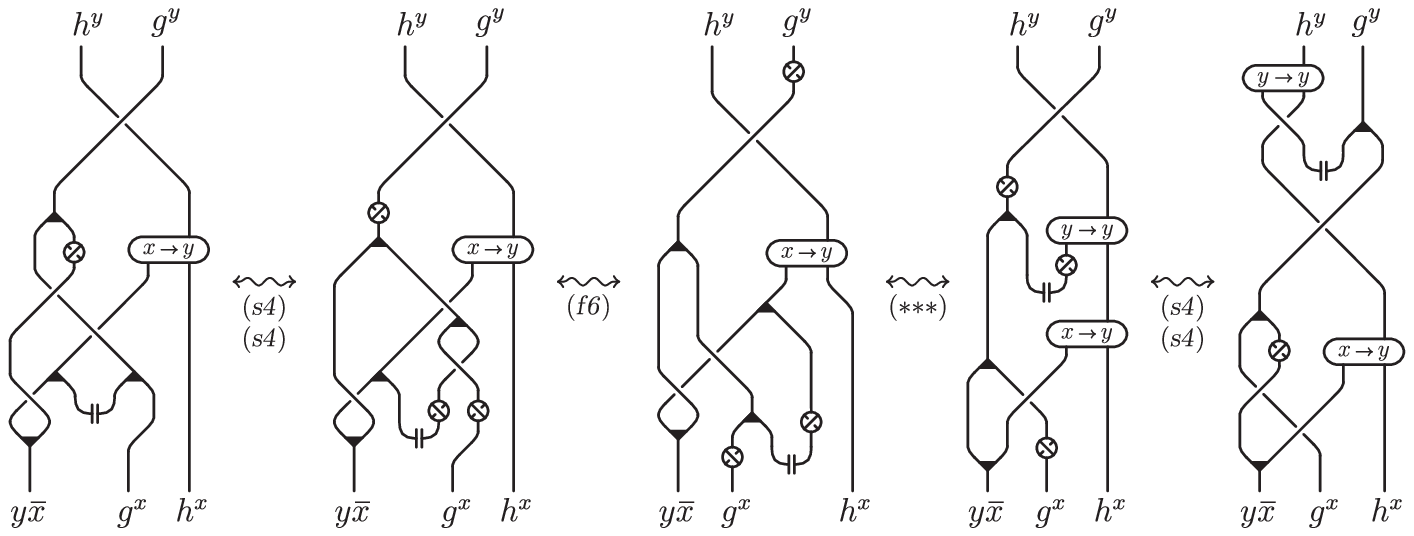}}\vskip-6pt
\end{Figure}

\begin{Figure}[b]{invxi10/fig}{}
 {Proof of \(q1') for $F = \gamma_{g,h}$  with $g \in \G(k_0,k_0)$ -- II
 ($h \in \G(k,l)$, $x \in \G(i_0,k_0)$ and $y \in \G(j_0,k_0)$)
 [{\sl f}/\pageref{h-tortile1/fig},
  {\sl p}/\pageref{ribbon-tot/fig},
  {\sl q}/\pageref{actzeta/fig},
  {\sl r}/\pageref{ribbon2/fig}-\pageref{ribbon5/fig}].}
\vskip-24pt\centerline{\fig{}{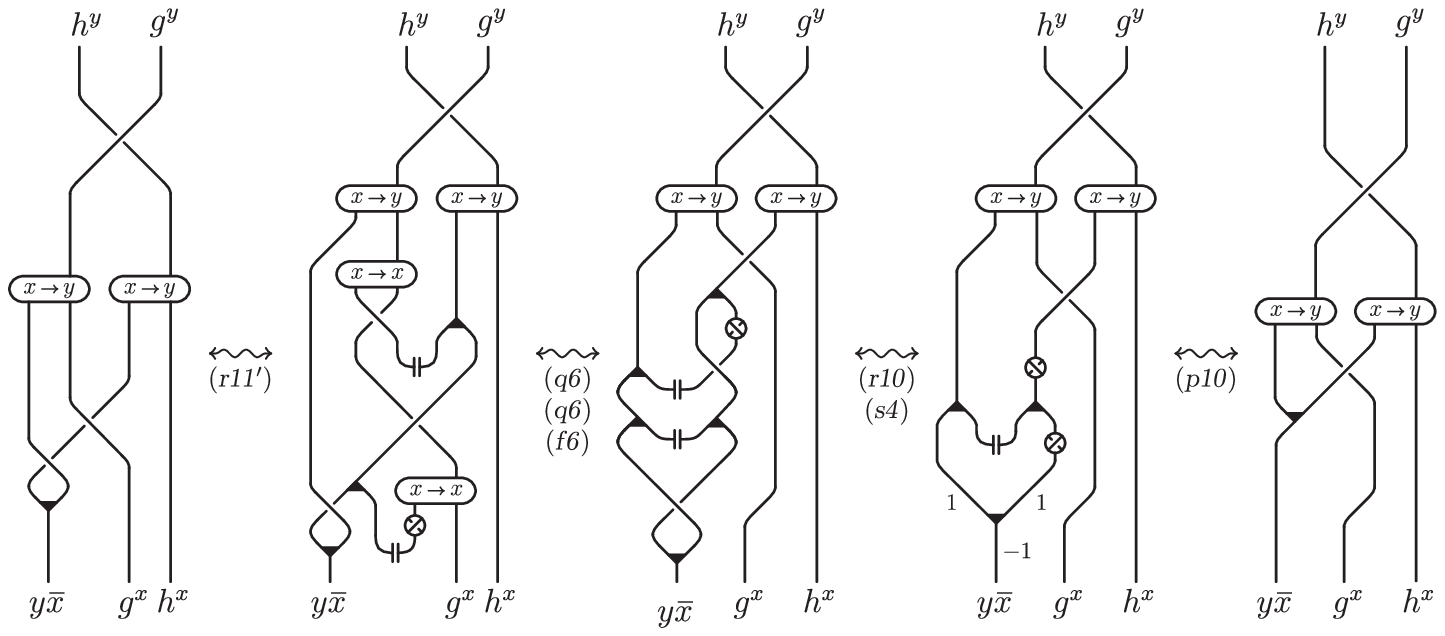}}\vskip-6pt
\end{Figure}

\item[\ms{\sl $F = \gamma_{g,h}$}.]
Let $g \in \G(i,j)$ and $h \in \G(k,l)$. We can transform the left side of \(q1') as
shown in Figure \ref{invxi7/fig}. Now the idea is to change the ``wrong'' crossing
between the edges labeled $g^x$ and $y\overline x$, through move \(r13') in Figure
\ref{ribbon5/fig}, and then use the copairings that appear to transform
$\gamma_{y\overline x,y\overline x} \circ \Delta_{y\overline x}$ into
$\Delta_{y\overline x}$. This is done separately for the four cases: $i \neq k_0$ and
$j \neq k_0$ in Figures \ref{invxi11/fig} and \ref{invxi12/fig}; $i = k_0$ and $j \neq
k_0$ in Figure \ref{invxi8/fig}; $i \neq k_0$ and $j = k_0$ in Figure \ref{invxi9/fig};
$i = k_0$ and $j = k_0$ in Figure \ref{invxi10/fig}.
\end{description}

At this point, having established property {\sl a}\/) for $F$ being an elementary
morphism, we are ready to show that $\Re^x$ is a functor, i.e. it preserves all the
relations between the elementary morphisms in $\H^r(\G)$ presented in Figures
\ref{isot/fig}, \ref{algebra/fig}, \ref{antipode/fig}, \ref{unimodular/fig},
\ref{ribbon1/fig} and \ref{ribbon2/fig}. This is trivially true when $x \in
\G(k_0,k_0)$. Hence, we focus on the case when $x \in \G(i_0,k_0)$, $i_0 \neq k_0$.

We start with the isotopy moves in Figure \ref{isot/fig}. As it was already observed,
moves \(b1-1') follow from \(p2-2') in Figure \ref{coform-s/fig} and \(q6) in Figure
\ref{actzeta/fig}.

Moves \(b3) and \(b4) in Figure \ref{isot/fig} need to be shown only when $D$ is an
elementary morphism $F:A \to B$. In this case, according to Figure \ref{cross11/fig},
the statement follows from the identity $\zeta_B^{x,x} \circ (\id_{x \bar x} \diam
\Re^x(F)) = \Re^x(F) \circ \zeta_A^{x,x}$, which is equivalent to the specialization of
the relation \(q1') proved above for $y=x$.

We continue with the bi-algebra axioms presented in Figure \ref{algebra/fig}. The
only non-trivial	ones are \(a1), \(a2) and \(a5) when $g$ or $h$ are in
$\G(i_0,k_0)$. In this case, \(a1) and \(a2) follow directly from the fact that
$\Re^x(\Delta_g) = (T_{\bar g^x} \diam T_{\bar g^x}) \circ \Delta_{\bar g^x} \circ
T_{g^x}^{-1}$ (cf. the definition of $T$ in Paragraph \ref{defnT/par} and the
presentation of $\Re^x(\Delta_g)$ in Figure \ref{delta/fig}). The proof of \(a5)
for $g \in \G(i_0,k_0)$ and $h \in \G(k_0,k_0)$ is presented in Figure
\ref{funct-bia/fig}. The proofs for different choices of $g$ and $h$ are simpler or
analogous.

\begin{Figure}[htb]{funct-bia/fig}{}
 {$\Re^x$ preserves \(a5) when $g \in \G(i_0,k_0)$ and $h \in \G(k_0,k_0)$
 \hbox{($x \in \G(i_0,k_0)$)\kern-6pt}
 [{\sl a}/\pageref{algebra/fig},{\sl p}/\pageref{coform-s/fig}].}
\centerline{\fig{}{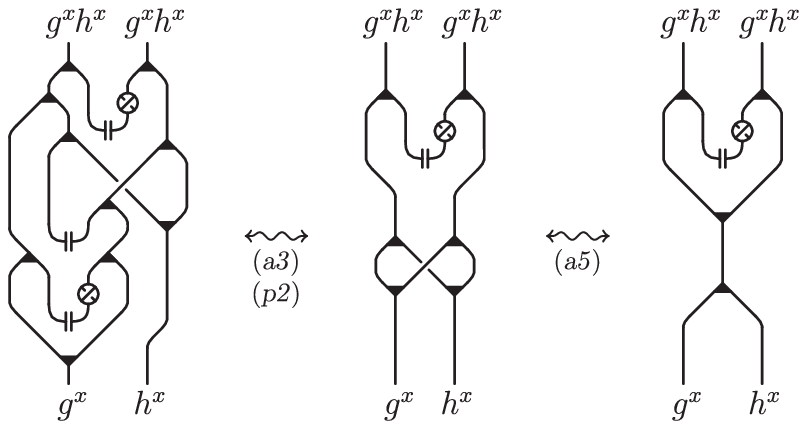}}\vskip-6pt
\end{Figure}

From the antipode axioms in Figure \ref{antipode/fig} the only non-trivial ones are
\(s1) and \(s1') for $g \in \G(i_0,k_0)$. Figure \ref{funct-anti/fig} proves \(s1) in
this case, while the proof of \(s1') is obtained by reflecting all the diagrams in the
same figure with respect to the $y$-axis.

\begin{Figure}[htb]{funct-anti/fig}{}
 {$\Re^x$ preserves \(s1) when $g \in \G(i_0,k_0)$
 ($x \in \G(i_0,k_0)$)
 [{\sl p}/\pageref{coform-s/fig}-\pageref{ribbon-tot/fig},
  {\sl s}/\pageref{antipode/fig}].}
\centerline{\fig{}{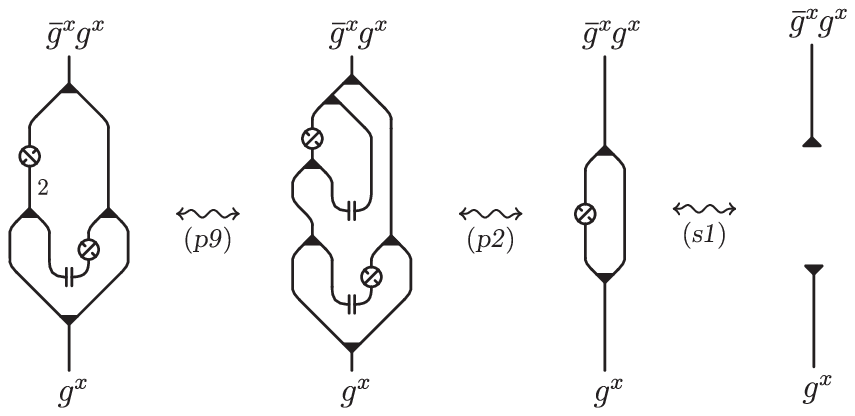}}\vskip-6pt
\end{Figure}

The integral axioms in Figure \ref{unimodular/fig} are trivially preserved because
$\Re^x(m_{g,h}) = m_{g^x,h^x}$ and $\Re^x(\Delta_{1_i}) = \Delta_{1_i^x}$.

\begin{Figure}[b]{funct-delta/fig}{}
 {$\Re^x$ preserves \(r10) when $g \in \G(i_0,k_0)$
 ($x \in \G(i_0,k_0)$)
 [{\sl p}/\pageref{coform-s/fig},
  {\sl r}/\pageref{ribbon2/fig}].}
\vskip-6pt\centerline{\fig{}{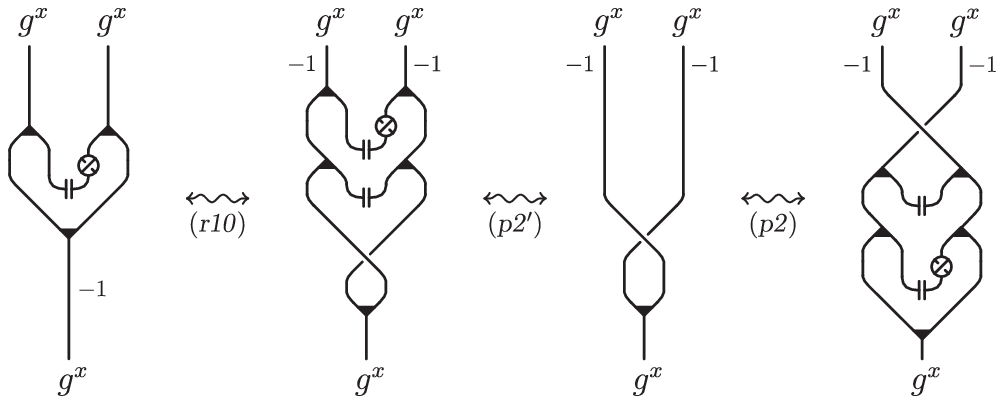}}\vskip-6pt
\end{Figure}

\begin{Figure}[htb]{funct-cross/fig}{}
 {$\Re^x$ preserves \(r11)
 ($g \in \G(i,j)$, $h \in \G(k,l)$, $x \in \G(i_0,k_0)$)
 [{\sl p}/\pageref{coform-s/fig},
  {\sl r}/\pageref{ribbon2/fig}-\pageref{ribbon5/fig}].}
\vskip-3pt\centerline{\fig{}{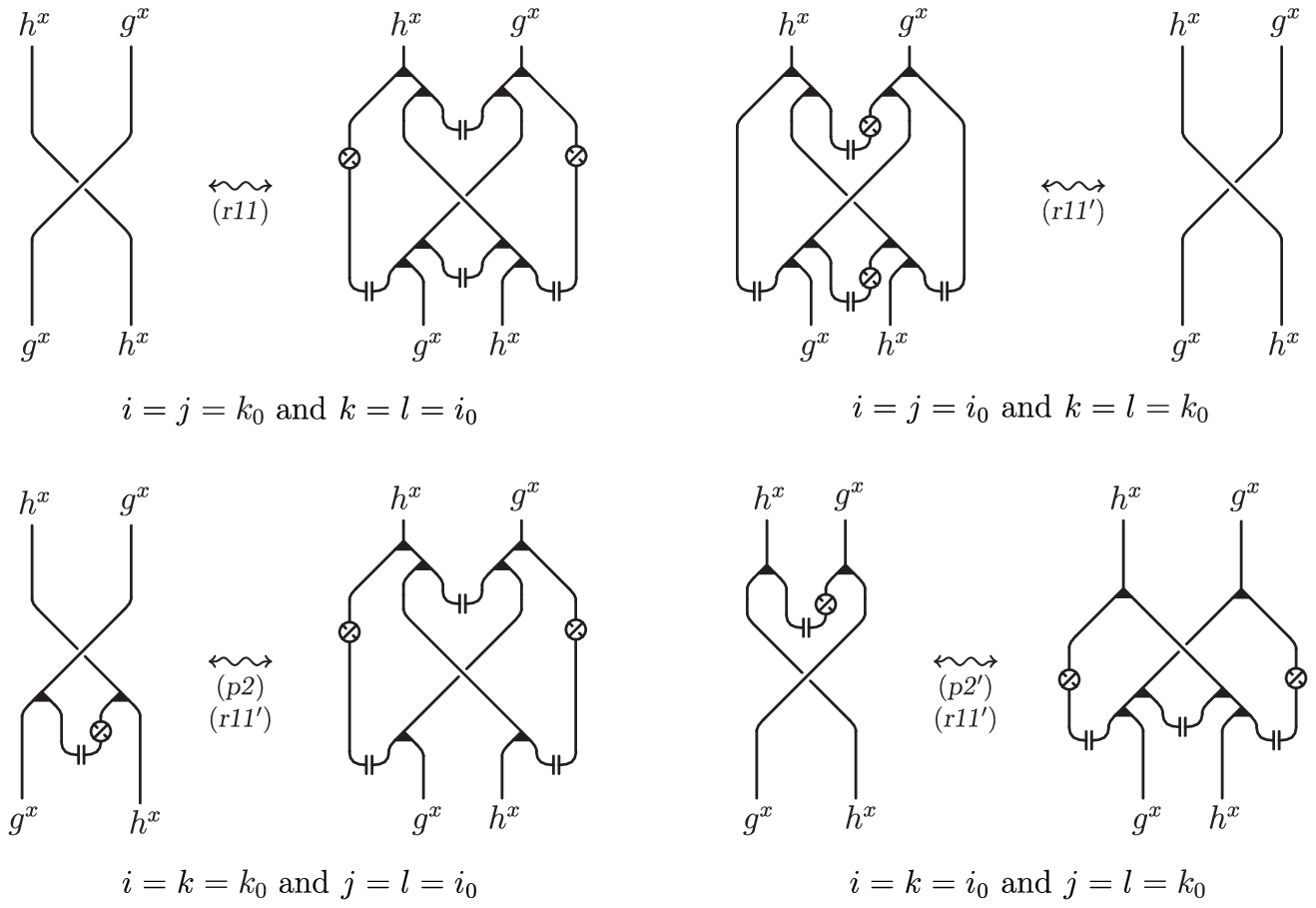}}\vskip-3pt
\end{Figure}

\break

Also the axioms of the ribbon structure in Figure \ref{ribbon1/fig} are trivially
preserved, while the only non-trivial ones from Figure \ref{ribbon2/fig} are \(r10)
when $g \in \G(i_0,k_0)$ and \(r11) with $g \in \G(i,j)$ and $h \in \G(k,l)$ such that
some of $i,j$ are equal to $k_0$ and some of $k,l$ are equal to $i_0$ or vice versa.
Figure \ref{funct-delta/fig} deals with the above metioned case of \(r10). Some of the
cases of \(r11) are presented in Figure \ref{funct-cross/fig} (cf. the expression for
$\Re^x(\sigma_{i,j})$ on page \pageref{imsigma/par}) and the rest are analogous. This
concludes the proof of the functoriality of $\Re^x$. Moreover, the functoriality
together with the fact that {\sl a}\/) holds for $F$ being an elementary morphism,
imply {\sl a}\/) for any other morphism $F$.

Properties {\sl b}\/) and {\sl c}\/) derive from the definition of $\Re^x$; in
particular, from the fact that $\Re^x$ coincides with the formal extension of $\_^x$ 
on the elementary morphisms of $\H^r(\G)$ which do not involve $i_0$ or $k_0$ (see
right side of Figure \ref{cross/fig} for the crossings), taking into account
\ref{defnh} {\sl e}\/).

Finally, the last part of the theorem follows from the fact that the additional
relations defining the quotient categories $\partial^\star \H^r(\G)$ and
$\partial \H^r(\G)$ involve only morphisms in $H_{1_i}$ for $i \in \Obj \G$.
\end{proof}

Notice that, given any $x \in \G(i_0,k_0)$ and $y \in \G(j_0,l_0)$ with $i_0 \neq k_0
\neq l_0 \neq j_0$, both compositions of functors $\Re^{y^x} \!\circ \Re^x$ and
$\Re^{x^y} \!\circ \Re^y$ retract $\H^r(\G)$ to the subcategory $\H^r(\G^{\bs k_0 \bs
l_0}) = \H^r(\G^{\bs l_0 \bs k_0})$, provided that $i_0 \neq l_0$ or $j_0 \neq k_0$
(if $i_0 = l_0$ and $j_0 = k_0$ one has $y^x \in \G(l_0,l_0)$ and $x^y \in
\G(k_0,k_0)$).

The next proposition asserts that these two retractions are naturally equivalent.
Since we will focus later on closed morphisms, this fact will not be used in the
present work. Nevertheless, we include it for its significance, and also in view 
of a possible extension of what follows to non-closed morphisms.

\begin{proposition} \label{natuR/theo} 
 Given $i_0 \neq k_0 \neq l_0 \neq j_0$ in $\Obj \G$, such that either $i_0 \neq l_0$
or $j_0 \neq k_0$, and morphisms $x \in \G(i_0,k_0)$ and $y \in \G(j_0,l_0)$, there
exists a natural transformation $\nu^{x,y}: \Re^{x^y} \!\circ \Re^y \to \Re^{y^x}
\!\circ
\Re^x$. In particular, for any morphism $F: A \to B$ in $\H^r(\G)$, we have the
following commutative diagram:
\medskip\centerline{\fig{}{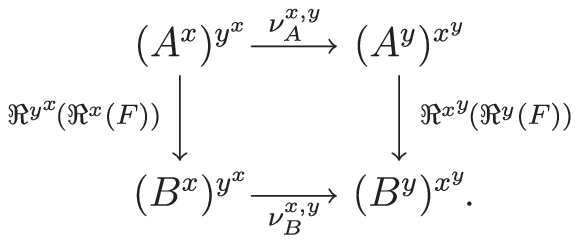}}\vskip-12pt
\end{proposition}

\begin{proof}
Let $x$ and $y$ be as in the statement. To simplify the notation we put $\Re^{y,x} =
\Re^{y^x} \circ \Re^x$ and $\Re^{x,y} = \Re^{x^y} \circ \Re^y$.
Given any object $A$ in $\H^r(\G)$, it is easy to check that $\Re^{y,x}(A) =
(A^x)^{y^x} = (A^y)^{x^y} = \Re^{x,y}(A)$, and we denote this image by $A'$. In
particular, we have $H_x' = H_{1_{i_0}^y}$ and $H_y' = H_{1_{j_0}^x}$.

\smallskip

We define the morphism $\nu^{x,y}_A: A' \to A'$ as follows:
$$\nu^{x,y}_A = (\epsilon_{1_{i_0}^y} \!\diam \epsilon_{1_{j_0}^x} \!\diam \id_{A'})
\,\circ\, \Re^{y,x}(\Theta_A \circ\, \Xi_A^{-1}) \,\circ\, (\eta^{1_{i_0}^y}
\!\diam \eta^{1_{j_0}^x} \!\diam \id_{A'}),$$
where $\Xi_A: H_x \diam H_y \diam A \to H_x \diam H_y \diam A'$
and $\Theta_A: H_x \diam H_y \diam A \to H_x \diam H_y \diam A'$ 
are the invertible morphisms given by 
$$\Xi_A = ((\xi^{1_{i_0},x}_{y^x})^{-1} \!\diam \id_{A'}) \,\circ\, (\id_x \diam
\xi^{1_{j_0}^x,y^x}_{A^x}) \,\circ\, \xi^{1_{i_0},x}_{H_y \diam A},$$
$$\Theta_A = ((\gamma_{x,y} \,\circ\, (\xi^{1_{j_0},y}_{x^y})^{-1}) \diam \id_{A'}) 
\,\circ\, (\id_y \diam \xi^{1_{i_0}^y,x^y}_{A^y}) \,\circ\, 
\xi^{1_{j_0},y}_{H_x \diam A} \,\circ (\bar\gamma_{x,y} \diam \id_A).$$
 Then, for any morphism $F: A \to B$ in $\H^r(\G)$, the relation \(q2) implies that
$$\Xi_B \,\circ\, (\id_x \diam \id_y \diam F) = (\id_x \diam \id_y \diam
\Re^{y,x}(F)) \,\circ\, \Xi_A,$$ 
$$\Theta_B \circ (\id_x \diam \id_y \diam F) =\break (\id_x \diam \id_y \diam
\Re^{x,y}(F)) \circ \Theta_A.$$ 
 From the former equality, it also follows that 
$$\Xi_B^{-1} \circ\, (\id_x \diam \id_y \diam \Re^{y,x}(F)) = (\id_x
\diam \id_y \diam F) \,\circ\, \Xi_A^{-1}.$$ 
 Summing up, we have 
$$\Theta_B \,\circ\, \Xi_B^{-1} \circ\, (\id_x \diam \id_y \diam
\Re^{y,x}(F)) = (\id_x \diam \id_y \diam \Re^{x,y}(F)) \,\circ\, \Theta_A
\,\circ\, \Xi_A^{-1}.$$
 By applying $\Re^{y,x}$ to both the members of this identity in
$\H^r(\G)$, we get the following identity in
$\H^r(\G^{\bs k_0
\bs l_0})$: 
$$\Re^{y,x}(\Theta_B \circ \Xi_B^{-1}) \,\circ\,
(\id_{1_{i_0}^y} \diam \id_{1_{j_0}^x} \diam \Re^{y,x}(F)) = 
(\id_{1_{i_0}^y} \diam \id_{1_{j_0}^x} \diam \Re^{x,y}(F))
\,\circ\, \Re^{y,x}(\Theta_A \circ \Xi_A^{-1}),$$
where we have used  the fact that $\Re^{y,x}\circ \Re^{y,x}=\Re^{y,x}$ and
$\Re^{y,x}\circ \Re^{x,y}=\Re^{x,y}$ (cf. Theorem \ref{functR/theo} {\sl b}\/) and
{\sl c}\/)).

Finally, we compose with $\epsilon_{1_{i_0}^y} \diam \epsilon_{1_{j_0}^x} \diam
\id_{B'}$ on the left and with $\eta_{1_{i_0}^y} \diam \eta_{1_{j_0}^x} \diam \id_{A'}$
on the right, to obtain the desired identity: $\nu^{x,y}_B \,\circ\, \Re^{y,x}(F) =
\Re^{x,y}(F)\,\circ\,\nu^{x,y}_A$.
\end{proof}

Now we proceed with the main result of this section, namely Theorem
\ref{gen-reduction/theo}, which is a reformulation of the first part of Theorem
\ref{reduction/theo}, concerning reduction and stabilization maps, in the general
context of the ribbon Hopf algebras $\H^r(\G)$, where $\G$ is an arbitrary (finitely
generated) connected groupoid $\G$ (cf. Paragraph \ref{groupoid}). But before stating
the theorem we need some preliminaries.

Recall that, given a connected groupoid $\G$, a morphism $F$ in $\H^r(\G)$ is called
complete if the elements of $\G$ occurring as labels in a graph diagram of $F$ together
with the identities of $\G$ generate all $\G$, i.e. any element of $\G$ which is not
an identity can be obtained as a product of such labels or their inverses. Observe that
the definition is independent on the choice of the particular diagram representation of
$F$, since the only new labels that can appear or disappear in the equivalence moves of
$\H^r(\G)$ are identities of $\G$ or products and inverses of labels already occuring.
 
On the other hand, given any $g \in \G$, we can use equivalence moves to change any
diagram  representing a complete morphism $F$, into a new one where $g$ occurs as a
label. Namely, move \(s2) allows us to create a new edge whose label is the inverse of
a preexisting one, while applying the move in Figure \ref{create-edg/fig}, possibly
after isotopy,  allows us to create a new edge whose label is the product of two
preexisting labels. Since by the completeness of $F$, $g$ is a product of labels which
appear in the original diagram or their inverses, we can eventually create an edge
labeled by $g$.

\begin{Figure}[htb]{create-edg/fig}{}
{($g \in \G(i,j)$, $h \in \G(j,k)$)
 [{\sl a}/\pageref{algebra/fig}].}
\centerline{\fig{}{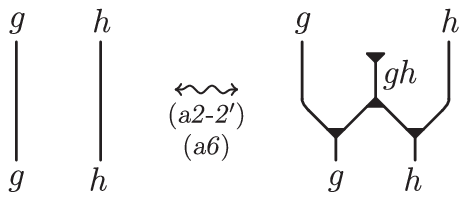}}\vskip-6pt
\end{Figure}

We denote by $\hat\H^{r,c}(\G)$ the set of closed complete morphisms in $\H^r(\G)$.
Notice that the set $\hat\H^{r,c}(\G)$ of closed complete morphisms in $\H^r(\G)$ is
non-empty if and only if $\G$ is finitely generated, which in turn implies that $\Obj
\G$ is finite.

\begin{Figure}[b]{creadisk/fig}{}{Isotoping $F$ to $F' \circ L_g$.}
\centerline{\fig{}{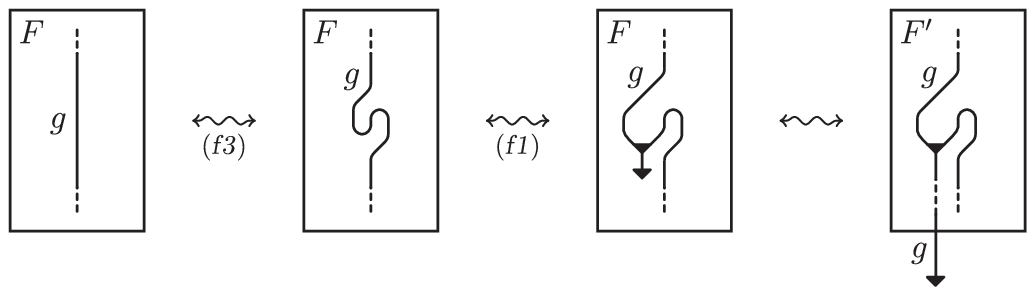}}\vskip-12pt
\end{Figure}

\begin{lemma} \label{completeness/theo} 
 Let $\G$ be a connected groupoid and $F \in \hat\H^{r,c}(\G)$. Then for any $g \in \G$
there exists a morphism $F'$ in $\H^r(\G)$ such that $F = F' \circ L_g$.
\end{lemma}
 
\begin{proof}
 According to what was said above, we can assume that $F$ contains an edge labeled $g$.
Then $F$ can be isotoped to the form $F' \circ L_g$, for a suitable morphism $F'$ in
$\H^r(\G)$, as shown in Figure \ref{creadisk/fig}.
\end{proof}

\begin{theorem} \label{gen-reduction/theo}
 Any full inclusion $\G \subset \G'$ of finitely generated connected groupoids uniquely
determines two bijective maps 
$$\up_\G^{\G'}: \hat\H^{r,c}(\G) \to \hat\H^{r,c}(\G') \quad \text{and} \quad
\down_\G^{\G'}: \hat\H^{r,c}(\G') \to \hat\H^{r,c}(\G),$$
which we call respectively stabilization map and reduction map, such that:
\begin{itemize}\itemsep\smallskipamount
\item[{\sl a}\/)]\vskip-\lastskip\medskip
 $\up_\G^{\G'}$ and $\down_\G^{\G'}$ are inverse of each other;
\item[{\sl b}\/)] 
$\up_\G^\G = \down_\G^\G = \id_{\hat\H^{r,c}(\G)}$ and moreover
$$\up_{\G'}^{\G''} \circ\, \up_\G^{\G'} = \up_\G^{\G''} \quad \text{and} \quad
\down_\G^{\G'} \circ\, \down_{\G'}^{\G''} = \down_\G^{\G''}$$
for any other full inclusion $\G' \subset \G''$ of groupoids as above;
\item[{\sl c}\/)]
 if $x \in \G(i_0,k_0)$ with $i_0 \neq k_0$, then
\vskip-6pt
$$
\begin{array}{c}
 \up_{\G^{\bs k_0}}^\G(F) = F \diam (\epsilon_x \circ L_x) \ \ \text{for any
  $F \in \hat\H^{r,c}(\G^{\bs k_0})$,}\\[6pt]
 \down_{\G^{\bs k_0}}^\G(F) = \Re^x(F') \circ \eta_{1_{i_0}} \ \ \text{for any
  $F = F' \circ L_x \in \hat\H^{r,c}(\G)$.}
\end{array}
$$
\end{itemize}\smallskip
Moreover, $\up_\G^{\G'}$ and $\down_\G^{\G'}$ induce bijective maps between the sets
of closed complete morphisms of the quotient categories
\vskip-12pt
$$
\begin{array}{c}
 \partial^\star{\up_\G^{\G'}}: \partial^\star\hat\H^{r,c}(\G) \to
  \partial^\star\hat\H^{r,c}(\G') \quad \text{and} \quad
  \partial^\star{\down_\G^{\G'}}: \partial^\star\hat\H^{r,c}(\G') \to
  \partial^\star\hat\H^{r,c}(\G)\text{,}\\[6pt]
 \partial{\up_\G^{\G'}}: \partial\hat\H^{r,c}(\G) \to \partial\hat\H^{r,c}(\G')
  \quad \text{and} \quad
  \partial{\down_\G^{\G'}}: \partial\hat\H^{r,c}(\G') \to \partial\hat\H^{r,c}(\G).
\end{array}
$$
\vskip-15pt
\end{theorem}

\begin{proof}
 We begin by verifying that {\sl c}\/) gives a good definition for the two maps
$\up_{\G^{\bs k_0}}^\G : \hat\H^{r,c}(\G^{\bs k_0}) \to \hat\H^{r,c}(\G)$ and
$\down_{\G^{\bs k_0}}^\G : \hat\H^{r,c}(\G) \to \hat\H^{r,c}(\G^{\bs k_0})$.

Concerning $\up_{\G^{\bs k_0}}^\G$, we have to prove that for any $F \in \hat\H^{r,c}
(\G^{\bs k_0})$ the morphism $F \diam (\epsilon_x \circ L_x)$ belongs to
$\hat\H^{r,c}(\G)$ and does not depend on the choice of $x$. The only non-trivial
point in the first assertion is the completeness of $F \diam (\epsilon_x \circ L_x)$
in $\H^r(\G) $. This follows from the completeness of $F$ and from the fact that
$\G^{\bs k_0}$ together with $x$ generate all $\G$ (cf. Paragraph \ref{defnh}). The
independence on $x$ is proved in Figure~\ref{indep-x/fig}, where $y \in \G(j_0,k_0)$
with $j_0\neq k_0$ and $F''$ is a morphism in $\H^r(\G^{\bs k_0})$ such that $F = F''
\circ y \bar x$, whose existence is guaranteed by Proposition
\ref{completeness/theo}.

\begin{Figure}[htb]{indep-x/fig}{}
 {$F \diam (\epsilon_x \circ L_x) = F \diam (\epsilon_y \circ L_y)$ in $\H^{r}(\G)$
 [{\sl a}/\pageref{algebra/fig}, 
  {\sl i}/\pageref{unimodular/fig}].}
\vskip3pt\centerline{\fig{}{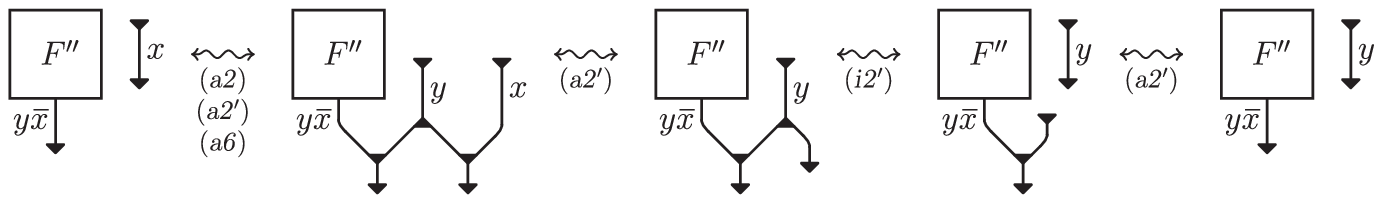}}\vskip-3pt
\end{Figure}

Passing to $\down_{\G^{\bs k_0}}^\G$, we have to prove that for any $F \in
\hat\H^{r,c}(\G)$ and any decomposition $F = F' \circ L_x$ the morphism $\Re^x(F) \circ
\eta_{1_{i_0}}$ belongs to $\hat\H^{r,c}(\G^{\bs k_0})$ and does not\break depend on
the choice of the decomposition. As above, for the first assertion it is enough to
observe that $\Re^x(F) \circ\, \eta_{1_{i_0}}$ is complete in $\H^r(\G^{\bs k_0})$,
since the labels occurring in it include the images under the functor $\_^x$ of those
occurring in $F$, so they generate all $\G^{\bs k_0}$.  To see that  two different
decompositions $F = F' \circ L_x$ with $x\in\G(i_0,k_0)$ and $F = F'' \circ L_y$ with
$y \in \G(j_0,k_0)$ give rise to the same reduction of $F$,
\mypagebreak 
without loss of
generality, we may assume that those decompositions are obtained by isotoping edges,
as in Figure \ref{creadisk/fig}, of two  diagrams representing $F$. Assume first, that
those two diagrams are the same, i.e that the decompositions are obtained by starting
with the same diagram and  making two different choices of an edge and an isotopy. By
isotoping both  edges simultaneously we obtain a third decomposition $F = F'''\circ
(L_x \diam L_y)$ (the center diagram in Figure \ref{indepdisk/fig}) such that, up to
isotopy, $F' = F''' \circ L_y$ and $F'' = F'''\circ L_x$. Then the functoriality of
$\Re^x$ and $\Re^y$ implies that $\Re^x(F''') \circ (\eta_{1_{i_0}} \diam L_{y \bar
x})$ is equivalent to $\Re^x(F') \circ \eta_{1_{i_0}}$ and $\Re^y(F''') \circ
(L_{x\bar y} \diam \eta_{1_{j_0}})$ is equivalent to $\Re^y(F'') \circ \eta_{1_{j_0}}$
in $\H^{r,c}(\G^{\bs k_0})$ (the  equivalences on the left and on the right in Figure
\ref{indepdisk/fig}). Therefore it is enough to show that $\Re^x(F''') \circ
(\eta_{1_{i_0}}\! \diam L_{y\bar x})$ and $\Re^y(F''') \circ (L_{x\bar y} \diam
\eta_{1_{j_0}})$ are equivalent in $\H^{r,c}(\G^{\bs k_0})$, which is done in Figure
\ref{pfindepdisk/fig}. Now the independence on the particular diagram  representing
$F$, follows from the fact that all  elementary relations are local. So, by the
argument above, when applying a relation to the diagram of $F$, we can always choose
the isotoped edge producing the decomposition $F = F' \circ L_x$, away from the
support of the relation and think of the  relation as if it were applied to the
diagram of $F'$. Now the statement follows from the functoriality of $\Re^x$. This
concludes the proof that the stabilization and the reduction maps in c) are well
defined.

\begin{Figure}[htb]{indepdisk/fig}{}
 {Two different reductions of $F = F''' \circ (L_x \diam L_y)$.}
\vglue3pt\centerline{\fig{}{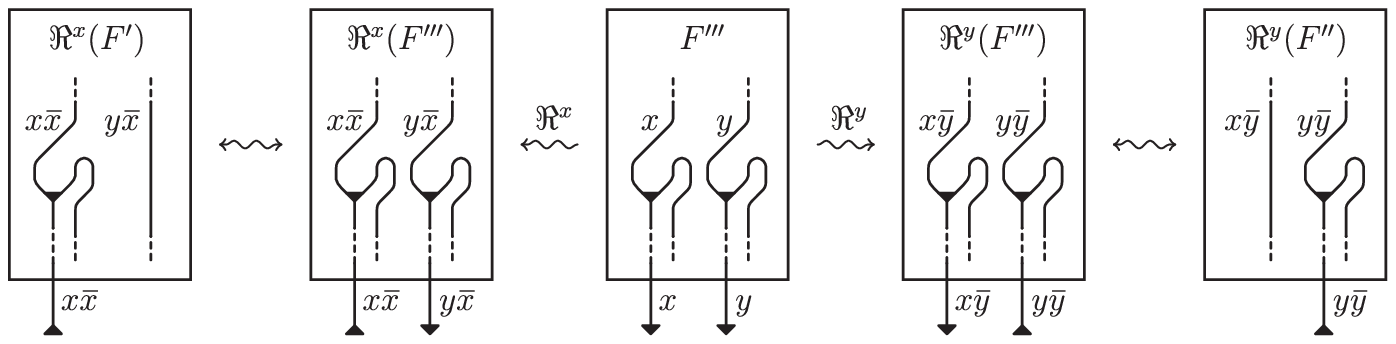}}\vskip-3pt
\end{Figure}			

\begin{Figure}[htb]{pfindepdisk/fig}{}
 {$\Re^x(F''') \circ (\eta_{x\bar x} \diam L_{y\bar x}) = \Re^y(F''') \circ (L_{x\bar
 y} \diam \eta_{y\bar y})$ 
 [{\sl a}/\pageref{algebra/fig}, 
  {\sl f}/\pageref{h-tortile1/fig}, 
  {\sl i}/\pageref{unimodular/fig},\break
  {\sl q}/\pageref{natuR/fig}, 
  {\sl s}/\pageref{pr-antipode/fig}].}
\vskip-18pt\centerline{\fig{}{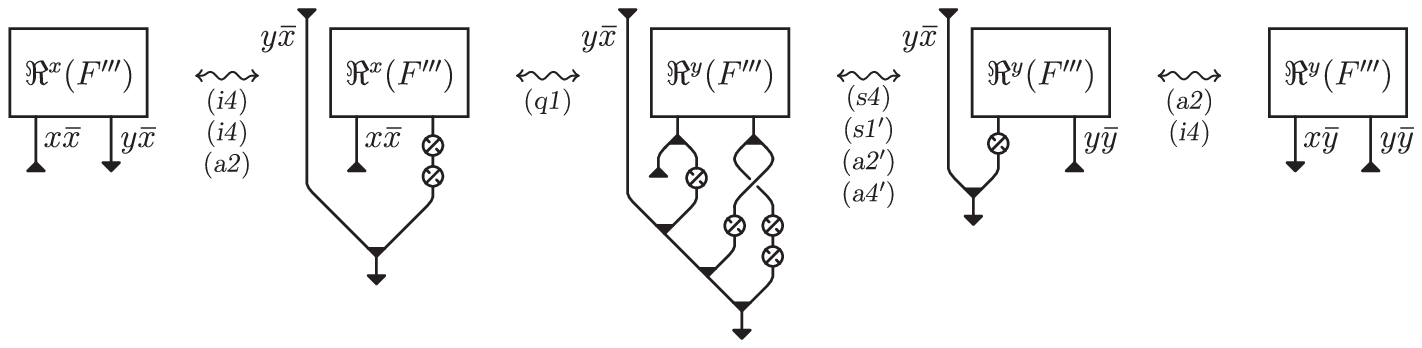}}\vskip-6pt
\end{Figure}

Now we want to show that the maps $\up_{\G^{\bs k_0}}^\G$ and $\down_{\G^{\bs k_0}}^\G$
defined in {\sl c}\/) are inverse of each other. Indeed, given $F \in \hat\H^{r,c}
(\G^{\bs k_0})$ and $x \in \G(i_0,k_0)$ with $i_0 \neq k_0$, Theorem \ref{functR/theo}
{\sl b}\/) implies that $\Re^x(F)=F$. Then $\down_{\G^{\bs k_0}}^{\G}
({\up^{\G}_{\G^{\bs k_0}}}(F)) = \down_{\G^{\bs k_0}}^{\G}(F \diam (\epsilon_x \circ
L_x)) = \Re^x(F \diam \epsilon_x) \circ \eta_{1_{i_0}} = F \diam (\epsilon_{1_{i_0}}\!
\circ \eta_{1_{i_0}}) = F$. On the other hand, Figure \ref{updown/fig} shows that
$F=\up_{\G^{\bs k_0}}^{\G}({\down^{\G}_{\G^{\bs k_0}}}(F))$ for any $F \in
\hat\H^{r,c}(\G)$.
	
At this point, we define the maps $\up_\G^{\G'}$ and $\down_\G^{\G'}$ by iteration.
Namely, given the full inclusion of groupoids $\G \subset \G'$ with $\Obj \G' - \Obj 
\G = \{k_1, \dots, k_n\}$, we consider the sequence $\G'=G_n \supset G_n^{\bs
k_n}=G_{n-1} \supset \dots \supset G_2^{\bs k_2}=G_1 \supset G_1^{\bs k_1}=\G$, and
define

\begin{Figure}[htb]{updown/fig}{}
 {$F =\up_{\G^{\bs k_0}}^{\G}({\down^{\G}_{\G^{\bs k_0}}}(F))$
[{\sl a}/\pageref{algebra/fig},
 {\sl f}/\pageref{h-tortile1/fig},
 {\sl i}/\pageref{unimodular/fig},
 {\sl q}/\pageref{natuR/fig}, 
 {\sl s}/\pageref{pr-antipode/fig}].}
\centerline{\fig{}{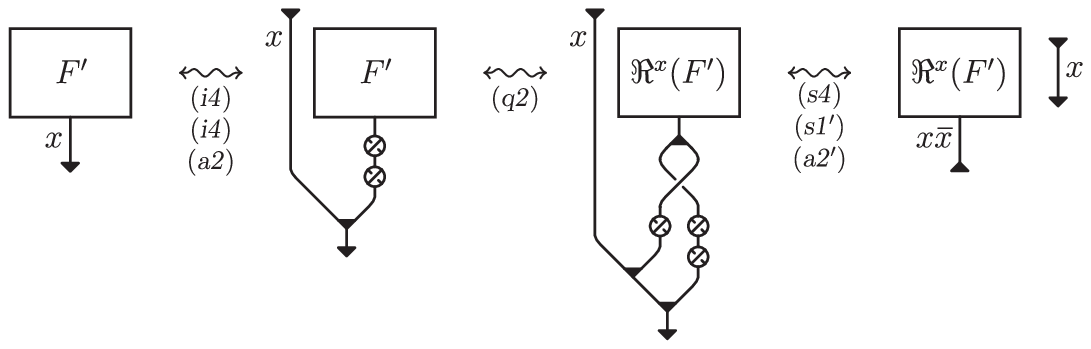}}\vskip-6pt
\end{Figure}

\break\vglue-30pt
$$\up_\G^{\G'} = \up_{G_n^{\bs k_n}}^{G_n} \!\circ \dots \circ \up_{G_2^{\bs
k_2}}^{G_2} \!\circ \up_{G_1^{\bs k_1}}^{G_1} \quad \text{and} \quad
\down_\G^{\G'} = \down_{G_1^{\bs k_1}}^{G_1} \!\circ \down_{G_2^{\bs k_2}}^{G_2} 
\!\circ \dots \circ \down_{G_n^{\bs k_n}}^{G_n}.$$

Observe that the maps $\up_\G^{\G'}$ and $\down_\G^{\G'}$ are inverses of each other
and since the definition of first map is obviously independent on the choice of
ordering of $\Obj\G' - \Obj \G$, the same must be true for the second one. Then
properties {\sl a}\/) and {\sl b}\/) are trivially satisfied.

Finally, the induced maps between the sets of closed complete morphisms of the quotient
categories are defined similarly to $\up_\G^{\G'}$ and $\down_\G^{\G'}$, by replacing
$\Re^x$ with the induced functors $\partial^\star \Re^x$ and $\partial\Re^x$ (cf.
Theorem \ref{functR/theo}). 
\end{proof}

\smallskip
\begin{block}\label{reduction/par}
{\sl Proof of Theorem \ref{reduction/theo}.}\kern1ex
 We remind that $\H^r_n = \H^r(\G_n)$, where $\G_n$ is the groupoid with $\Obj \G_n =
\{1,2,\dots,n\}$ and $\G_n(i,j) = \{(i,j)\}$ for any $1 \leq i,j \leq n$.
Then Theorem \ref{reduction/theo} is essentially the specialization of Theorem
\ref{gen-reduction/theo} to the case when $\G = \G_m$, $\G' = \G_n$ and the full
inclusion $\G_m \subset \G_n$ is the canonical one; moreover, the notations $\up_m^n
=\up_{\G_m}^{\G_n}$ and $\down_m^n = \down_{\G_m}^{\G_n}$ have been used. Then the
equality $\Phi_n \circ \up_m^n = \up_m^n \circ \Phi_m$ follows directly from the
definitions of the functors $\Phi_k$ and of the stabilization maps. $\;\square$
\end{block}

\section{From surfaces to the algebra%
\label{surf-alg/sec}}

This section is dedicated to the proof of Theorem \ref{defnpsi/theo}. First of all,
given a strict total order $<$ on $\Obj\G_n$, we define the functor $\Psi_n^<: \S_n
\to \H^r_n$ in the following way: on the objects $\Psi_n^<$ is uniquely determined by
the identities
$$\Psi_n^<((i\;j)) = H_{(i, j)} \ \ \text{if $i < j$},$$
$$\Psi_n^<(A \diam B) = \Psi_n^<(A) \diam \Psi_n^<(B),$$
while Figures \ref{surfalg1/fig} and \ref{surfalg2/fig} describe the images under
$\Psi_n^<$ of the elementary morphisms of $\S_n$ (the image of any labeling of the
morphism presented Figure \ref{ribbsurf13/fig} \(e'), is defined through relation
\(I6) in Figure \ref{ribbsurf14/fig}).

\begin{Figure}[htb]{surfalg1/fig}{}
 {Definition of $\Psi_n^<$ -- I ($i < j$, $i' < j'$).}
\centerline{\ \fig{}{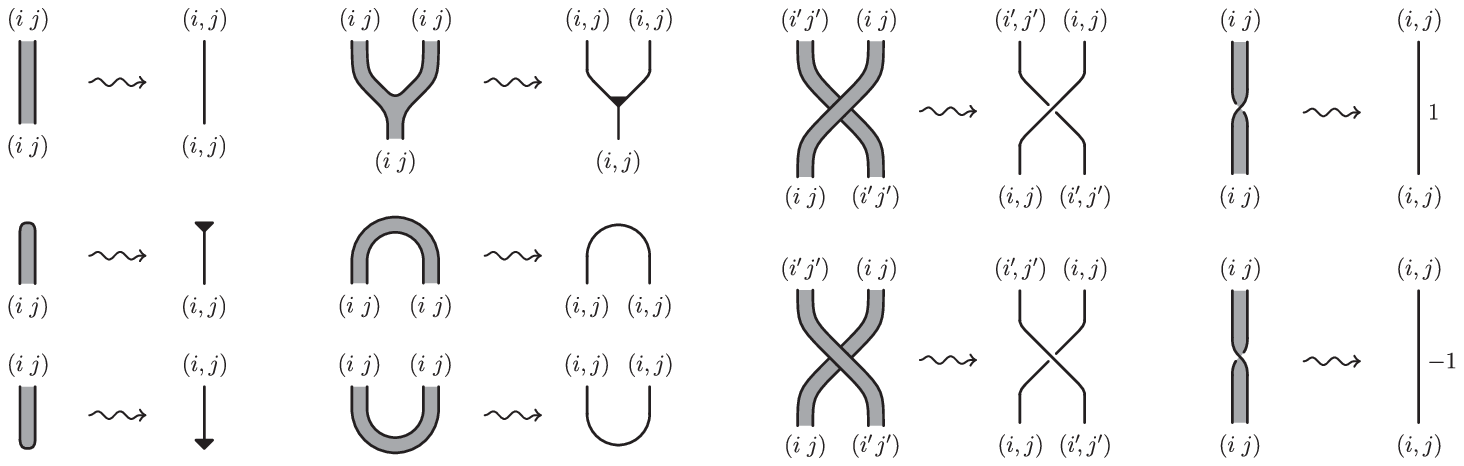}}\vskip-3pt
\end{Figure}

\begin{Figure}[htb]{surfalg2/fig}{}
 {Definition of $\Psi_n^<$ -- II ($i < j < k$, $h < l$ and
  $\{i,j\} \cap \{h,l\} = \emptyset$).}
\vskip-9pt\centerline{\fig{}{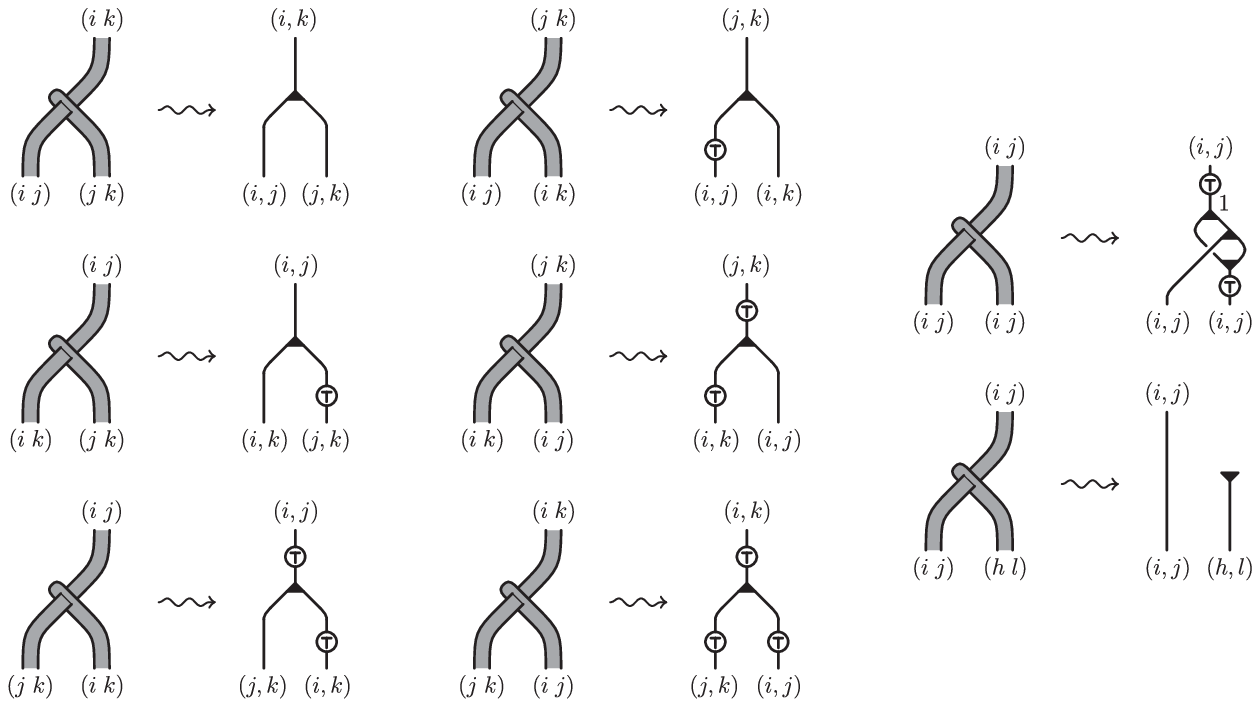}}\vskip-3pt
\end{Figure}

Then Theorem \ref{defnpsi/theo} states that $\Psi_n^<$ is a braided monoidal functor
and if $<'$ is another strict total order, there is a natural equivalence $\tau:
\Psi_n^< \to \Psi_n^{<'}$ which is identity on the empty set.

The proof of the theorem will make use of the relations \(t8) and \(t9), which are
introduced in the next proposition for an arbitrary groupoid $\G$. Given $g \in
\G(i,j)$ with $i \neq j$, let $T_g: H_g \to H_{\bar g}$ be defined as in Paragraph
\ref{defnT/par}. We remind that $T_{\bar g} \circ T_g = \id_g$ (cf. \(t3) in Figure
\ref{repT/fig}).

\begin{proposition} \label{cusp/theo}
	Let $\H^r(\G)$ be the universal ribbon Hopf $\G$-algebra. For any $g,h\in \G(i,j)$
with $i\neq j$, define the morphisms $c^\pm_{g, h}: H_g \diam H_{h} \to H_{h g
h}$ as follows:
\vskip-6pt
$$c^+_{g, h} = m_{h, g h} \circ (\id_{h} \diam m_{g, h}) \circ
(\gamma_{g,h} \diam \id_{h}) \circ (\id_g \diam \Delta_{h})\text{,}$$
$$c^-_{g, h} = m_{h, g h} \circ (\id_{h} \diam m_{g, h}) \circ
(\bar \gamma_{g,h} \diam \id_{h}) \circ (\id_g \diam \Delta_{h}).$$
\vskip6pt\noindent
 Then (cf. Figure \ref{cusp-0/fig})
\vskip-9pt
$$c^+_{g,h} \circ (\id_g \diam v_{h}) = c^-_{g,h}\text{,}
\eqno{\(t8)}$$
$$c^+_{g, h} = T_{\bar h \bar g \bar h} \circ c^+_{\bar g,\bar h} \circ (T_g \diam
T_{h}).
\eqno{\(t9)}$$\vskip-\lastskip
\end{proposition}\vskip-12pt

\begin{Figure}[htb]{cusp-0/fig}{}
 {Relations of Proposition \ref{cusp/theo} ($g,h \in \G(i,j)$, $i \neq j$).}
\centerline{\fig{}{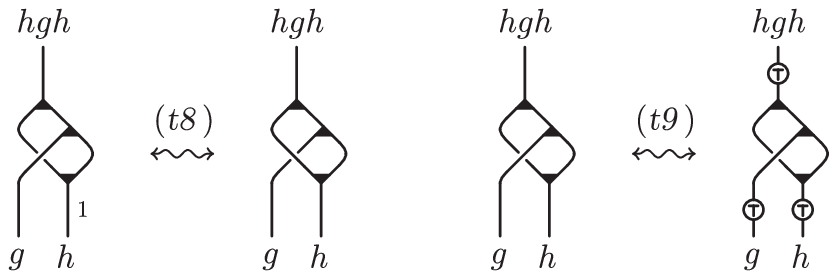}}\vskip-6pt
\end{Figure}

\begin{proof} 
 See Figure \ref{cusp-new/fig}.
\end{proof}

\begin{Figure}[htb]{cusp-new/fig}{}
 {Proof of Proposition \ref{cusp/theo} 
 [{\sl a}/\pageref{algebra/fig},
  {\sl p}/\pageref{ribbon-tot/fig},
  {\sl r}/\pageref{ribbon1/fig},
  {\sl s}/\pageref{pr-antipode/fig},
  {\sl t}/\pageref{repT/fig}-\pageref{cusp-0/fig}].}
\centerline{\fig{}{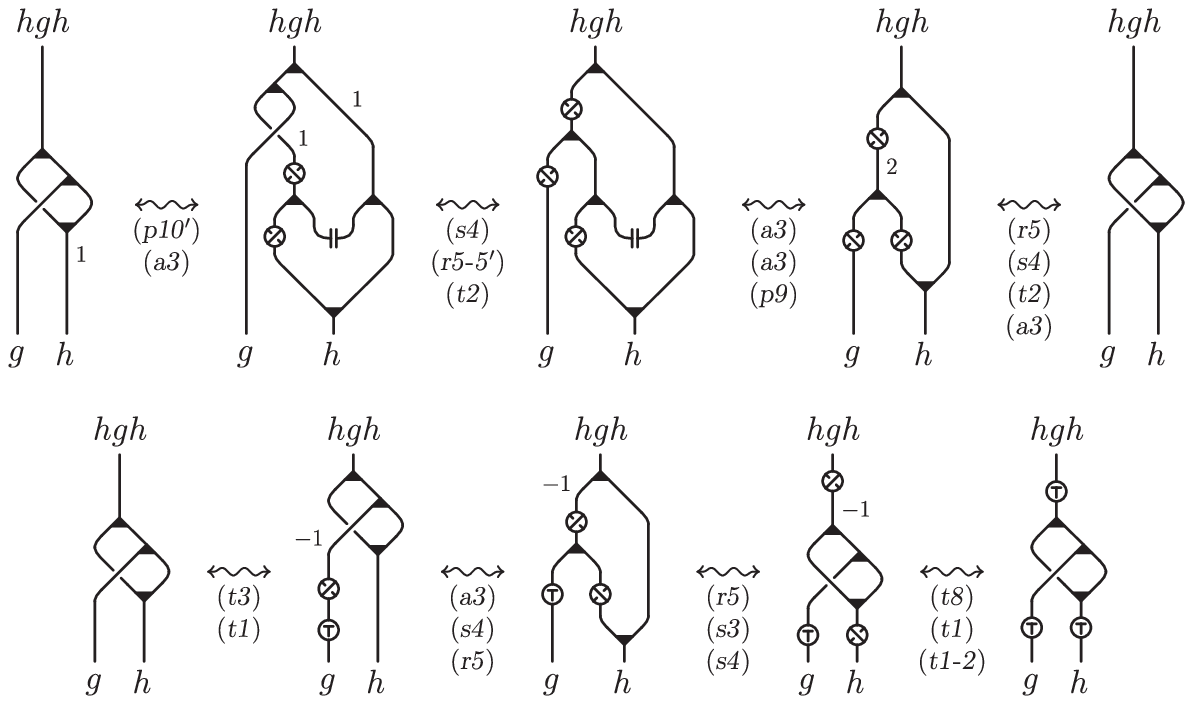}}\vskip-6pt
\end{Figure}

Observe that, when $\G = \G_n$ and $g = h = (i, j)$, $S_{\bar h g \bar h} \circ
c^+_{g,\bar h} \circ (\id_g \diam T_h)$ is exactly the image of the uni-colored
singular vertex under $\Psi_n$ (cf. Figure \ref{surfalg2/fig}). 

\begin{block}\label{connect-surf/par}
{\sl Proof of Theorem \ref{defnpsi/theo}.}\kern1ex
Given any two strict total orders $<$ and $<'$ on $\Obj\G_n$, for any object $A$ in
$\S_n$ we define the invertible morphism $\tau_A: \Psi_n^<(A) \to \Psi_n^{<'}(A)$ by
induction, according to the relations

$$
\tau_{(i,j)} = \left\{%
\begin{array}{ll}
 \id_{(i,j)}: H_{(i,j)} \to H_{(i,j)} &\ \ \text{if } i<j \text{ and } i<'j,\\[2pt]
 T_{(i,j)}: H_{(i,j)} \to H_{(j,i)} &\ \ \text{if } i<j \text{ and } j<'i;
\end{array}\right.$$
\vskip-4pt
$$\tau_{\one}=\id_{\one},\quad \tau_{A \diam B} = \tau_A \diam \tau_B.$$

\vskip3pt
Then for any morphism $F: A \to B$ in $\S_n$, represented as a fixed composition
of products of elementary morphisms, we have the following commutative diagram:

\vskip9pt\centerline{\fig{}{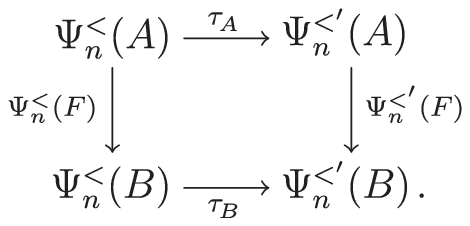}}\vskip6pt

\noindent
 Indeed, the commutativity in the case when $F$ is an elementary morphism can be easily
derived from the definition of $\Psi^{<}_n$, by using the relations \(t3) in Figure
\ref{repT/fig}, \(t4) in Figure \ref{T-theo/fig} and \(t9) in Figure \ref{cusp-0/fig},
while the extension to an arbitrary composition is allowed by the relation \(t3). 

Now we have to show that $\Psi_n^<$ is well-defined on the level of morphisms. In other
words, that it sends labeled rs-tangles related by any of the defining moves of $\S_n$
(the labeled versions of the moves in Figures \ref{ribbsurf14/fig},
\ref{ribbsurf15/fig}, \ref{ribbsurf16/fig} and \ref{ribbsurf17/fig} together with the
moves in Figure \ref{ribbon-m/fig}) to labeled graph diagrams representing the same 
morphism in $\H^r_n$. In showing that, by the commutativity of the above diagram, we
can choose the most convenient order $<$ for each single move.

We remind that move \(I6) in Figures \ref{ribbsurf14/fig} was used to define the
images under $\Psi_n^<$ of the morphisms on the left hand side, so there is nothing to
prove for it. Analogously, for \(R2) in Figure \ref{ribbon-m/fig} it suffices to
look at the definition of $\Psi_n^<$.

Here below we indicate how the verification goes for the remaining moves (actually,
most of them essentially rewrite the algebra axioms and the relations in $\H^r_n$
proved in Sections \ref{hopf-alg/section} and \ref{ribbon-alg/section}).

\begin{description}\itemsep\smallskipamount 
\item[\ms{\(I1), \(I7-7') and \(I8-9)}]
 follow from the braid axioms in Figure \ref{isot/fig} (p. \pageref{isot/fig}).
\item[\ms{\(I2-2')}] correspond to axioms \(f3-3') in Figure \ref{form/fig} 
 (p. \pageref{form/fig}).
\item[\ms{\(I3)}] follows from one of the botton-left duality moves in Figure
 \ref{form-uni/fig} (p. \pageref{form-uni/fig}).
\item[\ms{\(I4)}] follows from relation \(p4) in Figure \ref{ribbon-isot/fig} (p.
 \pageref{ribbon-isot/fig}).
\item[\ms{\(I5)}] follows from move \(f4) in Figure \ref{cycl0/fig} (p.
 \pageref{cycl0/fig}).
\item[\ms{\(I10), \(I12-12') and \(I13)}] follow from the ribbon axioms in Figure
\ref{ribbon1/fig} (p. \pageref{ribbon1/fig}).
\item[\ms{\(I11)}] follows from relations \(f6') in Figure \ref{h-tortile/fig}
 (p. \pageref{h-tortile1/fig}) and \(t2) in Figure \ref{repT/fig}
 (p. \pageref{repT/fig}).
 
\item[\ms{\sl\(I14-14')}] for a bi-colored singular vertex are trivial, while for
a tri-colored singular vertex follow directly from \(t6) and \(t7) in Figure
\ref{T-theo/fig} (p. \pageref{T-theo/fig}). The proof for a uni-colored singular vertex
is presented in Figures \ref{proof-i14/fig} and \ref{proof-i14a/fig}, where Figure
\ref{proof-i14/fig} shows that the image under $\Psi_n^<$ of the labeled rs-tangle in
the middle of \(I14-14') is equivalent to the third graph diagram in Figure
\ref{proof-i14a/fig}.

\begin{Figure}[htb]{proof-i14/fig}{}
 {Proof of \(I14-14') in the uni-colored case -- I ($i < j$) 
 [{\sl a}/\pageref{algebra/fig}, 
  {\sl f}/\pageref{form/fig}-\pageref{cycl0/fig}-\pageref{h-tortile1/fig}, 
  {\sl p}/\pageref{ribbon-isot/fig},
  {\sl r}/\pageref{ribbon1/fig}, 
  {\sl s}/\pageref{pr-antipode/fig}, 
  {\sl t}/\pageref{repT/fig}].}
\centerline{\fig{}{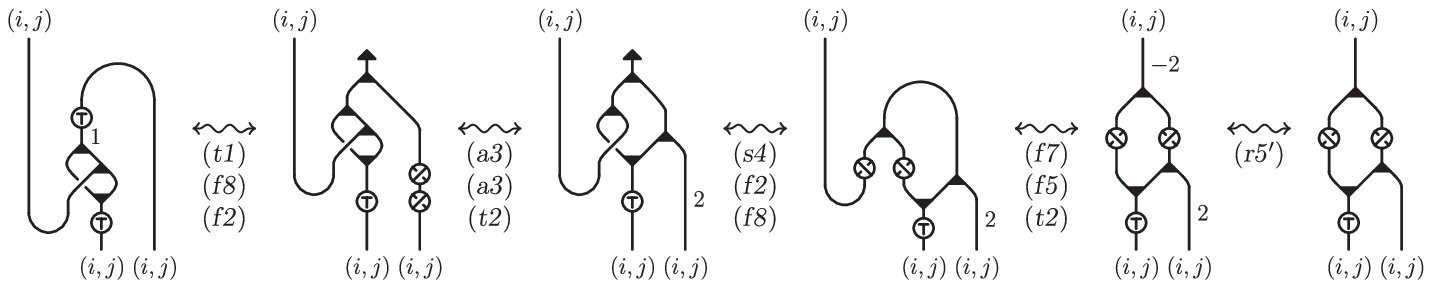}}\vskip-6pt
\end{Figure}

\begin{Figure}[htb]{proof-i14a/fig}{}
 {Proof of \(I14-14') in the uni-colored case -- II ($i < j$) 
 [{\sl a}/\pageref{algebra/fig},
  {\sl r}/\pageref{ribbon1/fig}, 
  {\sl s}/\pageref{pr-antipode/fig},
  {\sl t}/\pageref{repT/fig}].}
\centerline{\fig{}{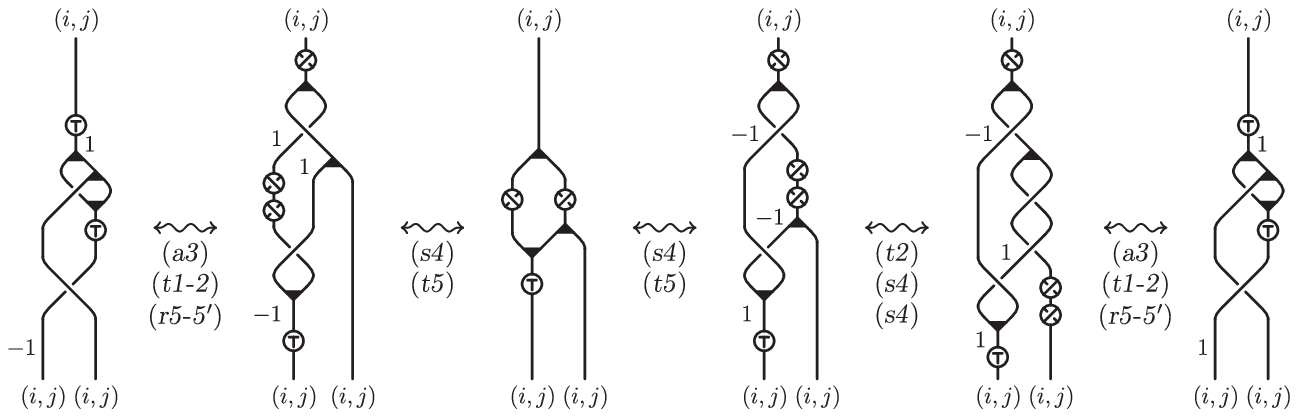}}\vskip-6pt
\end{Figure}

\item[\ms{\sl\(I15)}] follows from relation \(t5) in Figure \ref{T-theo/fig} (p.
\pageref{T-theo/fig}).

\item[\ms{\sl\(I16)}] follows from the definition of the coform \(f1) in Figure
 \ref{form/fig} (p. \pageref{form/fig}).
 
\item[\ms{\sl\(I17)}]
 for a bi-colored singular vertex reduces to a crossing change, so it follows from
axioms \(r11) in Figure \ref{ribbon2/fig} (p. \pageref{ribbon2/fig}). For a tri-colored
singular vertex the proof is presented in Figure \ref{proof-i17/fig}, while for a
uni-colored singular vertex it is presented in Figure \ref{proof-i17a/fig}.
 
\begin{Figure}[b]{proof-i17/fig}{}
 {Proof of \(I17) in the tri-colored case ($i < j < k$)
 [{\sl f}/\pageref{T-theo/fig},
  {\sl p}/\pageref{coform-s/fig}-\pageref{ribbon-tot/fig},
  {\sl r}/\pageref{ribbon1/fig}-\pageref{ribbon5/fig},
  {\sl s}/\pageref{pr-antipode/fig},
  {\sl t}/\pageref{repT/fig}].}
\vskip-6pt\centerline{\fig{}{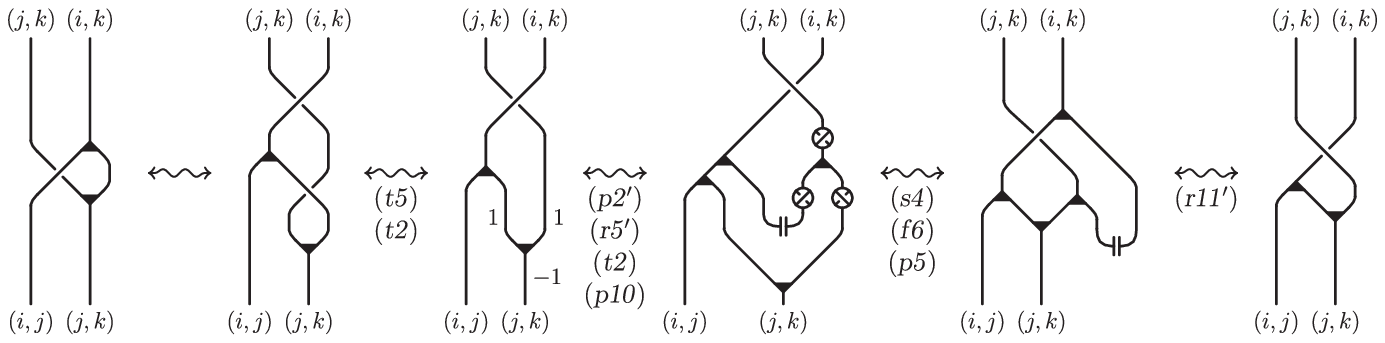}}\vskip-6pt
\end{Figure}

\begin{Figure}[htb]{proof-i17a/fig}{}
 {Proof of \(I17) in the uni-colored case ($i < j$)
 [{\sl a}/\pageref{algebra/fig},
  {\sl p}/\pageref{coform-s/fig}-\pageref{ribbon-tot/fig},
  {\sl r}/\pageref{ribbon1/fig}-\pageref{ribbon2/fig}-\pageref{ribbon5/fig},
  {\sl s}/\pageref{pr-antipode/fig},
  {\sl t}/\pageref{repT/fig}].}
\vskip-6pt\centerline{\fig{}{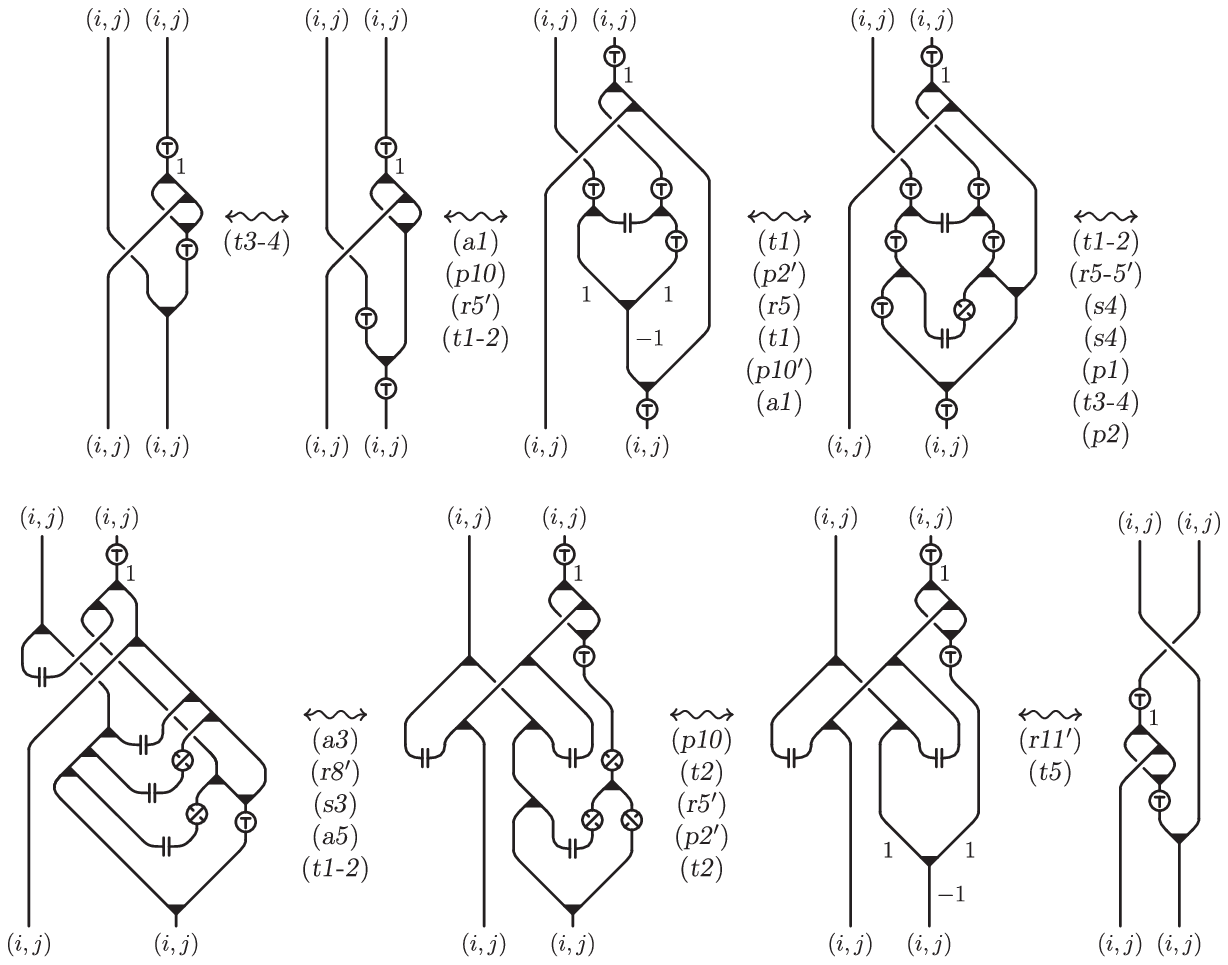}}\vskip-6pt
\end{Figure}

\item[\ms{\sl\(I18)}] corresponds to the bi-algebra axiom \(a1) in
 Figure \ref{algebra/fig} (p. \pageref{algebra/fig}).
\item[\ms{\sl\(I19)}] is presented in Figure \ref{proof-i19/fig}.

\begin{Figure}[htb]{proof-i19/fig}{}
 {Proof of \(I19) ($i < j$) 
 [{\sl a}/\pageref{algebra/fig}, 
  {\sl f}/\pageref{form/fig}-\pageref{h-tortile1/fig},
  {\sl r}/\pageref{ribbon1/fig},
  {\sl s}/\pageref{antipode/fig},
  {\sl t}/\pageref{repT/fig}].}
\centerline{\fig{}{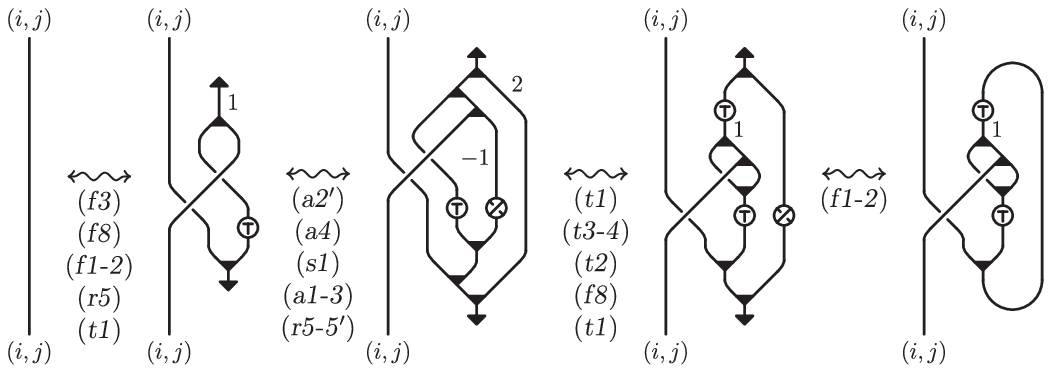}}\vskip-6pt
\end{Figure}

\item[\ms{\sl\(I20)}]
 follows from \(a6), \(s6) and \(r2) respectively in Figures \ref{algebra/fig} (p.
\pageref{algebra/fig}), \ref{pr-antipode/fig} (p. \pageref{pr-antipode/fig}) and
\ref{ribbon1/fig} (p. \pageref{ribbon1/fig}).
\item[\ms{\sl\(I21)}] for bi-colored singular vertices is trivial, while for
tri-colored singular vertices, with the proper choice of the order $<$, it corresponds
to the bi-algebra axiom \(a5) in Figure \ref{algebra/fig} (p. \pageref{algebra/fig}).
The proof for uni-colored singular vertices is presented in Figure \ref{proof-i21/fig}.

\begin{Figure}[htb]{proof-i21/fig}{}
 {Proof \(I21) in the uni-colored case ($i < j$) 
 [{\sl a}/\pageref{algebra/fig},
  {\sl p}/\pageref{ribbon-isot/fig},
  {\sl r}/\pageref{ribbon1/fig}-\pageref{ribbon2/fig},
  {\sl s}/\pageref{pr-antipode/fig},
  {\sl t}/\pageref{T-theo/fig}].}
\centerline{\ \fig{}{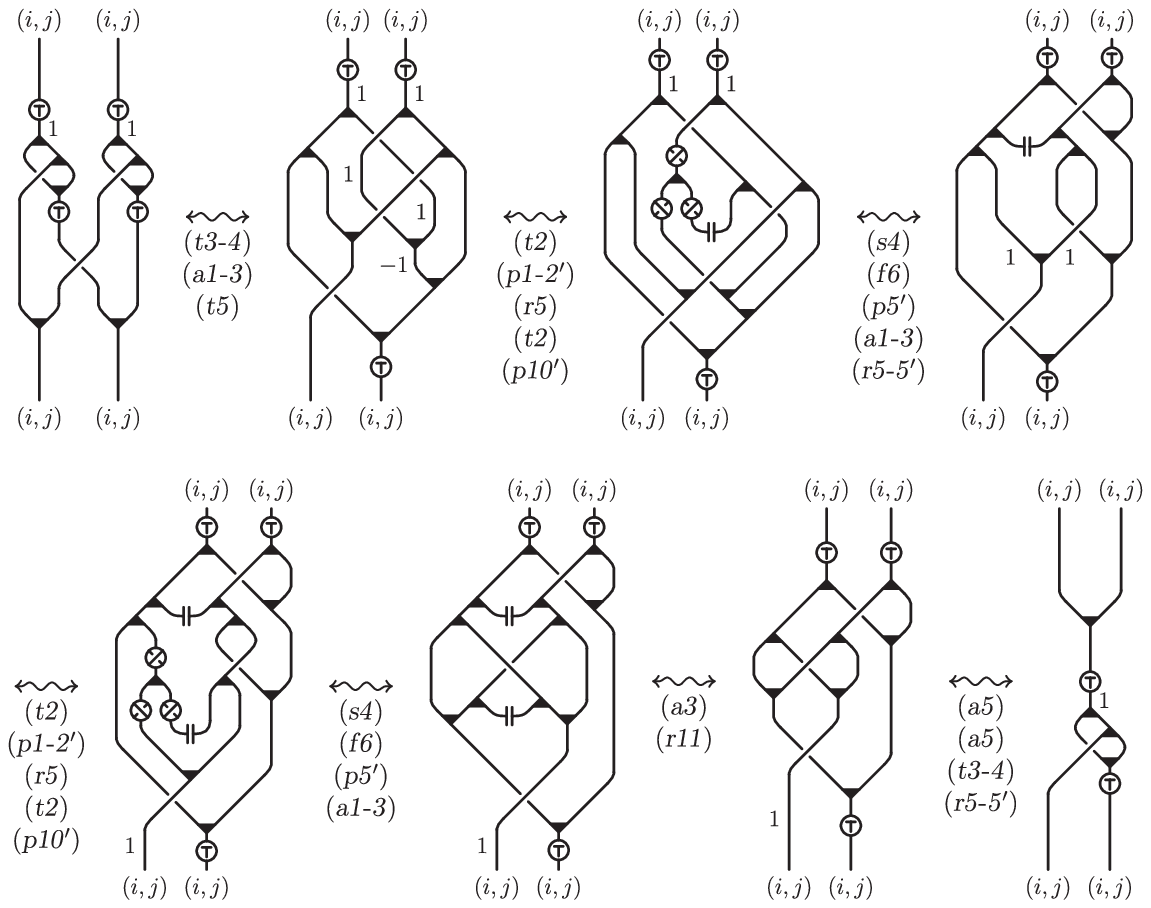}}\vskip-6pt
\end{Figure}

\begin{Figure}[t]{proof-i22/fig}{}
 {Proof of \(I22) in the non trivial cases when a bi-colored singular vertex occurs, 
  ($i < j < k < l$)
 [{\sl a}/\pageref{algebra/fig},
  {\sl r}/\pageref{ribbon1/fig}-\pageref{ribbon2/fig},
  {\sl s}/\pageref{antipode/fig},
  {\sl t}/\pageref{repT/fig}].}
\centerline{\fig{}{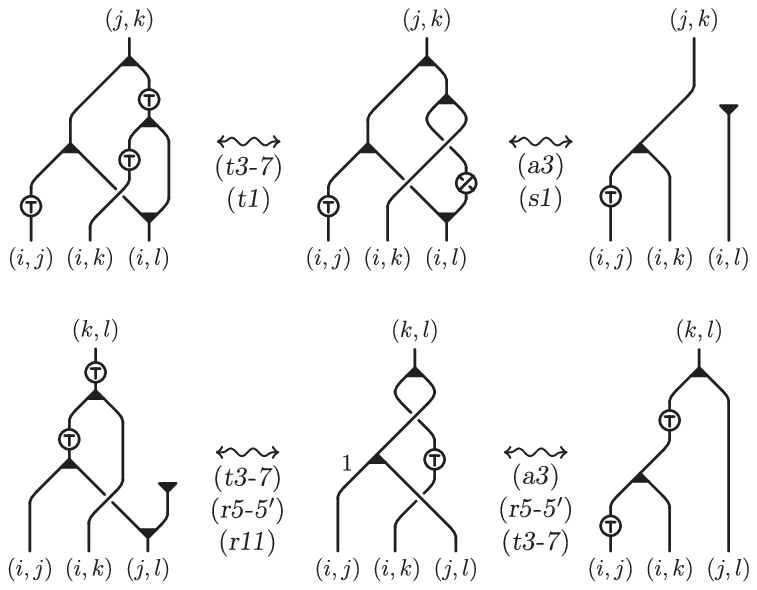}}\vskip-6pt
\end{Figure}

\begin{Figure}[b]{proof-i22a/fig}{}
 {Proof of \(I22) in the cases when one uni-colored singular vertex occurs 
 ($i < j < k$)
 [{\sl a}/\pageref{algebra/fig},
  {\sl r}/\pageref{ribbon1/fig},
  {\sl s}/\pageref{antipode/fig},
  {\sl t}/\pageref{repT/fig}-\pageref{cusp-0/fig}].}
\vskip-24pt\centerline{\fig{}{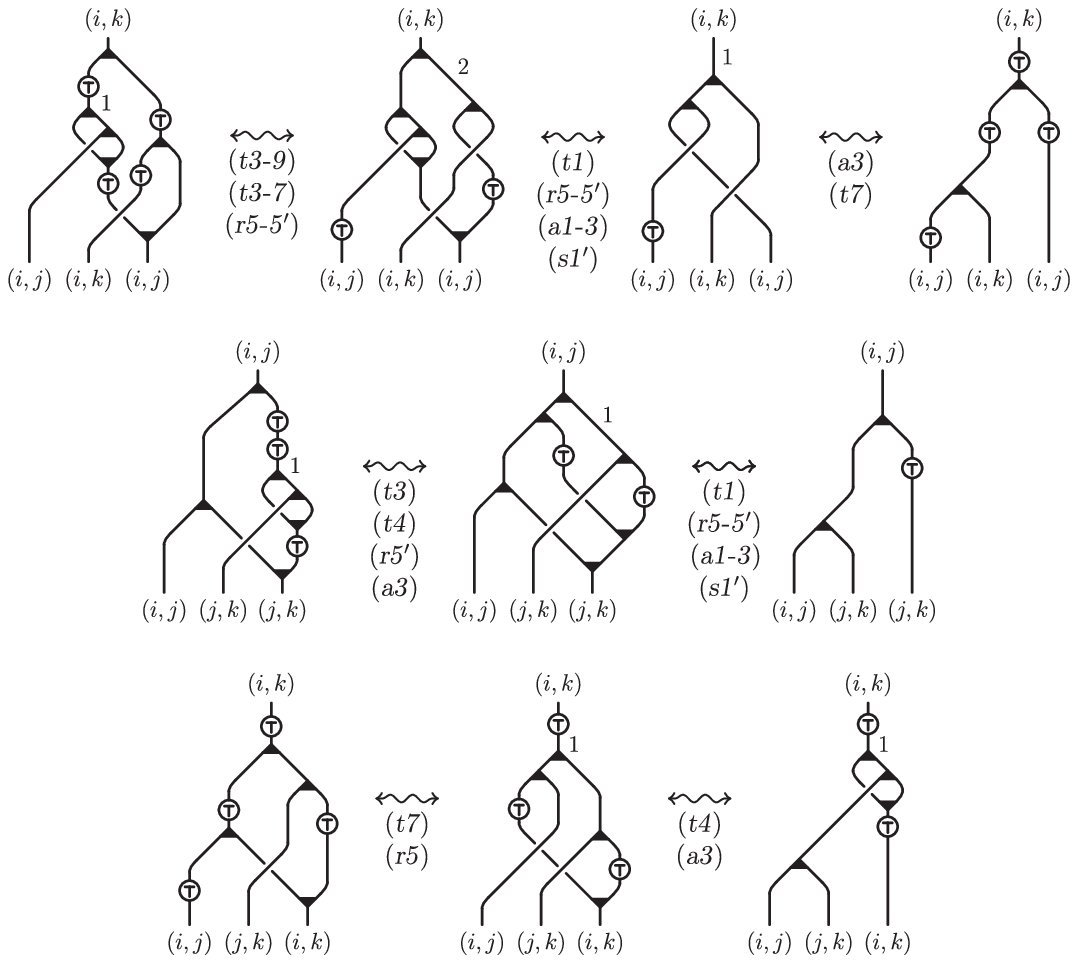}}\vskip-6pt
\end{Figure}

\begin{Figure}[b]{proof-i22b/fig}{}
 {Proof of \(I22) in the cases when two uni-colored singular vertices occur
  ($i < j < k$)
  [{\sl a}/\pageref{algebra/fig},
   {\sl r}/\pageref{ribbon1/fig},
   {\sl s}/\pageref{antipode/fig},
   {\sl t}/\pageref{repT/fig}-\pageref{cusp-0/fig}].}
\vskip6pt\centerline{\fig{}{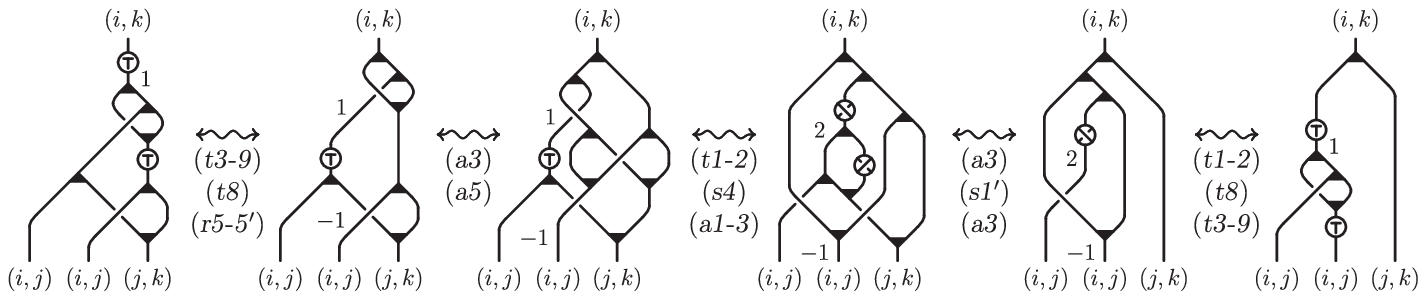}}\vskip-6pt
\end{Figure}

\begin{Figure}[t]{proof-i22c/fig}{}
 {Proof of \(I22) in the uni-colored case ($i < j$)
 [{\sl a}/\pageref{algebra/fig},
  {\sl p}/\pageref{coform-s/fig}-\pageref{ribbon-isot/fig},
  {\sl r}/\pageref{ribbon1/fig}-\pageref{ribbon2/fig},
  {\sl s}/\pageref{antipode/fig}-\pageref{pr-antipode/fig},
  {\sl t}/\pageref{repT/fig}-\pageref{cusp-0/fig}].}
\centerline{\fig{}{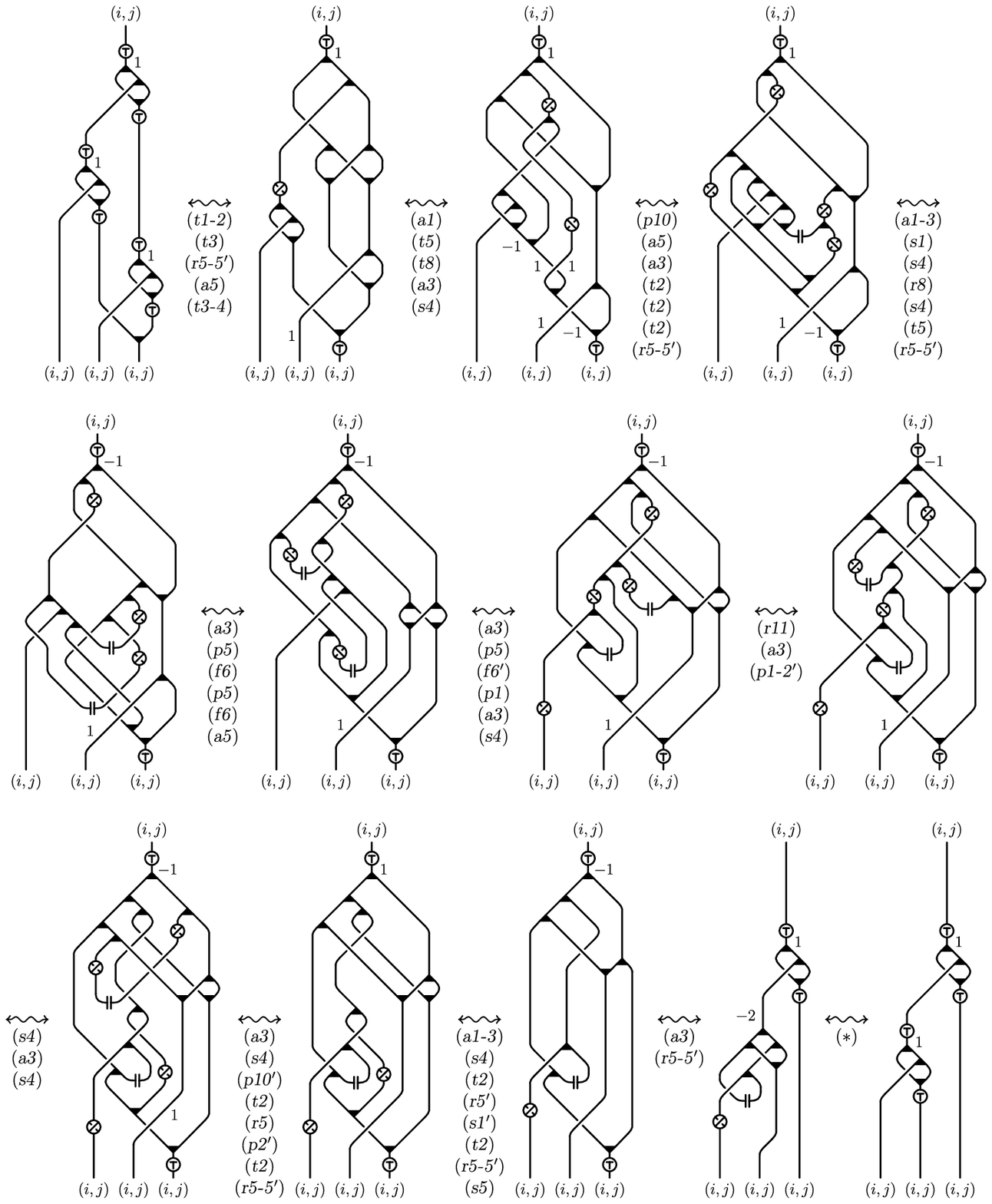}}\vskip-4pt
\end{Figure}

\begin{Figure}[b]{proof-i22d/fig}{}
 {Proof of step \($\ast$) in Figure \ref{proof-i22c/fig}
 [{\sl a}/\pageref{algebra/fig},
  {\sl p}/\pageref{ribbon-tot/fig},
  {\sl r}/\pageref{ribbon1/fig}-\pageref{ribbon2/fig},
  {\sl s}/\pageref{pr-antipode/fig},
  {\sl t}/\pageref{repT/fig}-\pageref{cusp-0/fig}].}
\vskip-25pt\centerline{\fig{}{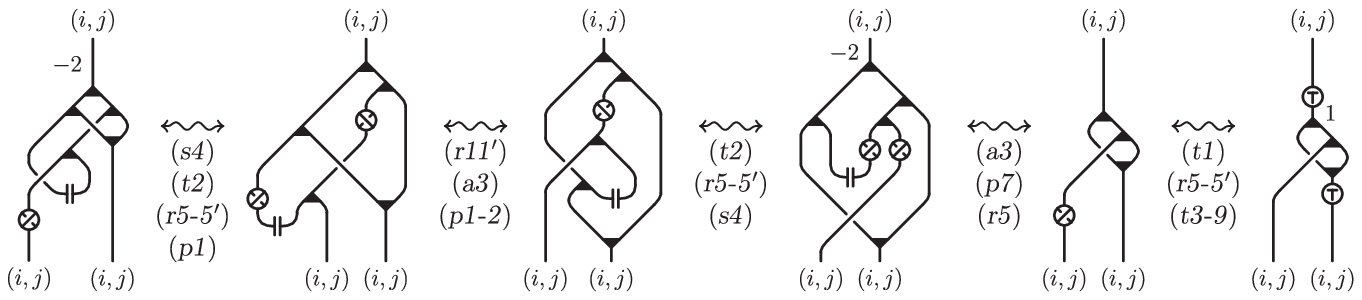}}\vskip-4pt
\end{Figure}

\item[\ms{\sl\(I22)}]\vskip-6pt
is the most complicated relation to deal with, since the source of the involved
morphisms consists of three intervals which can be labeled independently from each
other, so there are many different cases. First of all, we observe that the presence of
disjoint labels allows us to simplify the relation by using move \(R2) to remove the
bi-colored singular vertices. In particular, such simplification makes \(I22) trivial
when the label of the rightmost interval is disjoint from the other two, while it
makes \(I22) easily reducible to other moves (namely \(I3), \(I5), \(I16) and \(I20))
when the leftmost and the rightmost intervals have disjoint labels. Up to conjugation,
the only remaining labelings of the three intervals which include a pair of disjoint
labels are given by the sequences $(i\;j)\,,(j\;k)\,,(k\;l)$ and
$(i\;j)\,,(i\;k)\,,(i\;l)$ where $i,j,k,l$ are all distinct. The first case
corresponds, modulo move \(R2), to the bi-algebra axiom \(a3) in Figure
\ref{algebra/fig} (p. \pageref{algebra/fig}), while the second case is treated in the
top line of Figure \ref{proof-i22/fig}. The bottom line of the same figure concerns the
unique remaining case when a bi-colored singular vertex occurs (even if there is no
pair of dijoint labels in the source). The rest of the cases are grouped depending on
the number of uni-colored singular vertices: the three cases with only one such vertex
are presented in Figure \ref{proof-i22a/fig}; the single case with two uni-colored
vertices is considered in Figure \ref{proof-i22b/fig}; finally, the last case in which
all singular vertices are uni-colored is shown in Figures \ref{proof-i22c/fig} and
\ref{proof-i22d/fig}.

\item[\ms{\sl\(R1)}] follows from the relations \(t3) in Figure \ref{repT/fig} (p.
\pageref{repT/fig}) and \(t6) in Figure \ref{T-theo/fig} (p. \pageref{T-theo/fig}),
taking into account that $T$ propagates through the form and the coform, due to
\(f7-8) in Figure \ref{h-tortile1/fig} (p. \pageref{h-tortile1/fig}) and \(p3-4) in
Figure \ref{ribbon-isot/fig} (p. \pageref{ribbon-isot/fig}). $\;\square$
\end{description}\vskip-\lastskip
\end{block}\vskip-12pt
\vskip0pt

\begin{block}\label{4-mani/par}
{\sl Proof of Theorem \ref{eq-alg-kirby/theo}.}\kern1ex
The fact that $\Phi_n(\Psi_n(F)) = K_F$ for any $F \in \hat\S^c_n$ can be seen by
comparing the definitions of $\Phi_n$ (cf. Figures \ref{defnPhi/fig} and
\ref{defnPhi7/fig}) and $\Psi_n$ (cf. Figures \ref{surfalg1/fig} and
\ref{surfalg2/fig}) with the description of $K_F$ given in Section \ref{surf-kirby/sec}
(cf. Figures \ref{defnKF1/fig} and \ref{defnKF2/fig}). Hence, as discussed after the
statement of the theorem (at p. \pageref{eq-alg-kirby/theo}), it suffices to verify
that the map $\Psi_n: \hat\S^c_n \to \hat\H^{r,c}_n$ is surjective for $n \geq 4$.
Actually, we will do this for any $n \geq 3$.

Let $F$ be an arbitrary morphism in $\hat\H^{r,c}_n$ with $n \geq 3$, represented by a
given diagram (without using copairings and form/coform notation). An edge of the
diagram will be called an $i$-{\sl edge}, $1 \leq i \leq n$, if it is labeled $(i,i)$
and it is not attached to one positive tri-valent and one positive integral vertex.
Moreover, a vertex will be called an $i$-{\sl vertex} if it is not a positive integral
vertex and all edges attached to it are labeled $(i,i)$. 

As a preliminary step, we will show how to transform the diagram representing $F$ into
an equivalent one, where no $i$-edges appear for any $1 \leq i \leq n$. Actually, the
figures below deal only with edges of zero weight and not containing the antipode
morphism, but the generalization to other weights or to the presence of the antipode
is straightforward. Observe also that since $F$ is complete and $n \geq 2$, we may
assume that near a given $i$-edge there is a counit $\epsilon_{(i,j)}$ with $i\neq j$
(such a counit can be obtained by move \(a2) from an edge labeled $(i,j)$ and then
isotoped everywhere). 

We start by eliminating all uni-valent $i$-vertices as described in the first three
moves in Figure \ref{elim-mono1/fig}. Then by applying, if necessary, the edge breaking
shown in the fourth move in the same figure, we obtain a diagram  where all $i$-edges
connect two tri-valent vertices such that at most one of them is an $i$-vertex. 

\begin{Figure}[htb]{elim-mono1/fig}{}
 {Eliminating  uni-valent $i$-vertices ($i \neq j$)
 [{\sl a}/\pageref{algebra/fig},
  {\sl s}/\pageref{antipode/fig}].}
\centerline{\fig{}{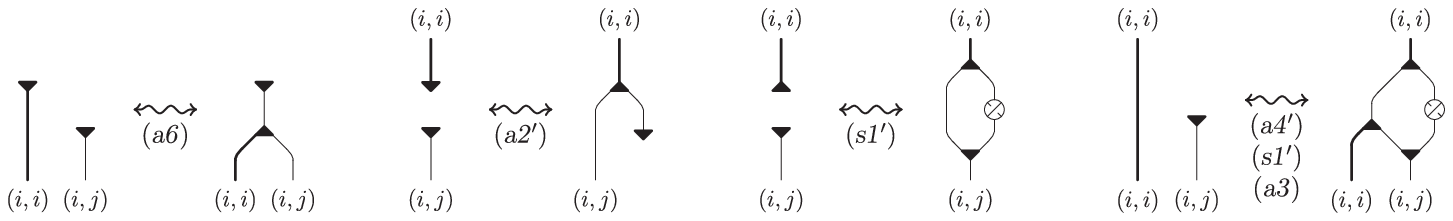}}\vskip-4pt
\end{Figure}

\begin{Figure}[b]{elim-mono2/fig}{}
 {Eliminating  tri-valent $i$-vertices ($i \neq j$)
 [{\sl a}/\pageref{algebra/fig}].}
\centerline{\fig{}{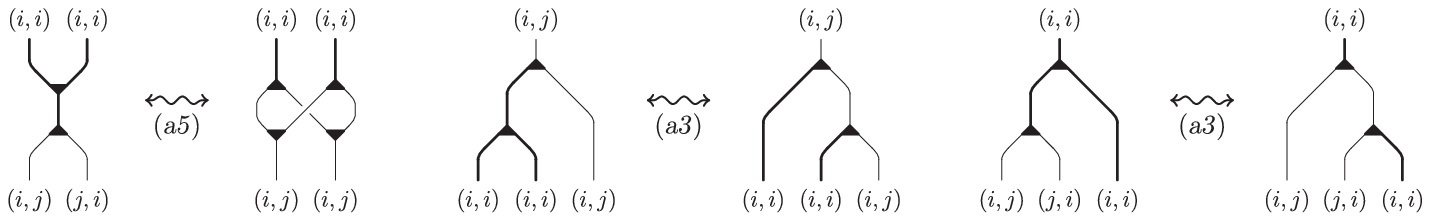}}\vskip-4pt
\end{Figure}

We proceed by eliminating the tri-valent $i$-vertices, starting with the negative ones
as shown in the leftmost move in Figure \ref{elim-mono2/fig} and then using the two
other moves in the same figure (or their vertical reflections) to eliminate the
positive $i$-vertices as well.

At this point, the only remaining $i$-edges, connect two positive tri-valent vertices
none of which is an $i$-vertex. Such edges are eliminated through the moves shown in
Figure \ref{elim-mono3/fig} (or their vertical reflections), where, since $n \geq 3$
and $F$ is complete, without loss of generality we have assumed that there exists 
close by a counit of label $(j, k)$ with $i \neq j \neq k \neq i$.

\begin{Figure}[hbt]{elim-mono3/fig}{}
 {Eliminating $i$-edges connecting two tri-valent vertices none of which is an
  $i$-vertex ($i \neq j \neq k \neq i$)
 [{\sl a}/\pageref{algebra/fig}, 
  {\sl s}/\pageref{antipode/fig}].}
\centerline{\fig{}{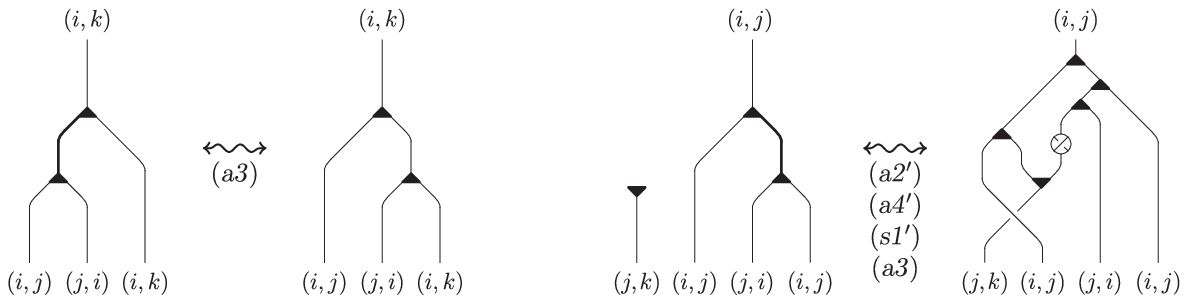}}\vskip-4pt
\end{Figure}

In this way we have represented the morphism $F$ by a diagram in which any edge labeled
$(i,i)$ is attached to one positive tri-valent and one positive integral vertex.
Then, by using the form notation introduced in Figure \ref{form/fig}, we can eliminate
those exceptional edges as well, obtaining a diagram of $F$ whose labels are all of
the type $(i,j)$ with $i \neq j$.

Now, given any order $<$ on $\Obj \G_n$, we modify this last diagram of $F$ in such a
way that the labels $(i,j)$ with $j < i$ are concentrated at short arcs of zero weight
delimited on one end by a uni- or tri-valent vertex and on the other end by an
antipode. This can be done by using relation \(t2) in Figure \ref{repT/fig} to convert
all the inverses of antipodes into antipodes and  to create/eliminate pairs of
antipodes along the same edge (at the cost of adding some twists) and relations \(f7-8)
in Figure \ref{h-tortile1/fig}, \(p1-4) in Figure \ref{ribbon-isot/fig} and \(r4) in
Figure \ref{ribbon1/fig} to slide antipodes and twists along edges.

Then, we remove those of the short arcs above, which are attached to uni-valent
vertices and to negative tri-valent vertices by applying moves \(i4-5) in Figure
\ref{unimodular/fig}, \(s5-6) in Figure \ref{pr-antipode/fig} and \(t4) in Figure
\ref{T-theo/fig}. The antipodes still present in the diagram, bound short arcs
attached to positive tri-valent vertices and using relation \(t1) in Figure
\ref{repT/fig}, we express them in terms of the morphisms $T_{(i,j)}$.

We finish the proof of Theorem \ref{eq-alg-kirby/theo} by observing that the resulting
diagram of $F$ is a composition of the diagrams on the left side of Figures
\ref{surfalg1/fig} and \ref{surfalg2/fig}, hence $F$ is in the image of $\Psi_n$.
$\;\square$
\end{block}

\begin{block}\label{3-mani/par}
{\sl Proof of Theorem \ref{3-mani/theo}.}\kern1ex
According to Proposition \ref{boundK/theo}, $\Phi_n$ induces functors
$\partial^\star\Phi_n: \partial^\star\H^r_n \to \partial^\star\K_n$ and 
$\partial\Phi_n: \partial\H^r_n \to \partial\K_n$ between the corresponding quotient
categories. Moreover, by the last part of Theorem \ref{gen-reduction/theo}, $\down_m^n$
and ${\up_m^n}$ induce well defined bijective maps between the closed complete
morphisms in these quotient categories for any $m < n$. 
In order to complete the proof of Theorem \ref{3-mani/theo} we only need to show that
the functor $\Psi_n$ induces functors
\vskip-6pt
$$\partial^\star\Psi_n: \partial^\star\S_n \to \partial^\star\H^r_n
\quad \text{and} \quad
\partial\Psi_n: \partial\S_n \to \partial\H^r_n,$$
in other words that the images under $\Psi_n$ of the two sides of the relation \(T)
(resp. relations \(T), \(P) and \(P')) in Figure \ref{defn-boundS/fig} are equivalent
in $\partial^\star\H_n^r$ (resp. $\partial\H_n^r$). This is done in Figures
\ref{threeman1/fig} (resp. \ref{threeman2/fig}). $\;\square$
\end{block}

\begin{Figure}[htb]{threeman1/fig}{}
 {Proof of \(T) in $\partial^\star\H^r_n$ ($i < j$)
 [{\sl a}/\pageref{algebra/fig},
  {\sl f}/\pageref{form/fig}-\pageref{h-tortile1/fig},
  {\sl p}/\pageref{ribbon-tot/fig},
  {\sl s}/\pageref{pr-antipode/fig}].}
\centerline{\fig{}{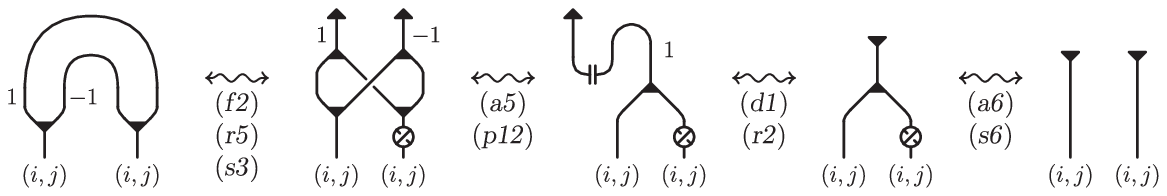}}\vskip-4pt
\end{Figure}

\begin{Figure}[htb]{threeman2/fig}{}
 {Proof of \(P) and \(P') in $\partial\H^r_n$ ($i < j$)
 [{\sl d}/\pageref{defn-boundH/fig},
  {\sl r}/\pageref{ribbon1/fig},
  {\sl s}/\pageref{antipode/fig}].}
\vskip3pt\centerline{\fig{}{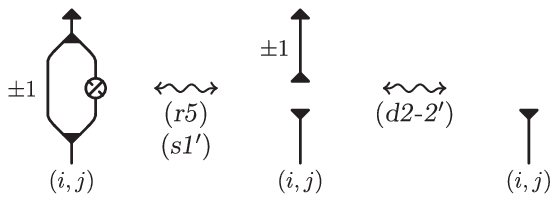}}\vskip-4pt
\end{Figure}

\section{Appendix: proof of Proposition \ref{surfcat/theo}%
\label{appendix1/sec}}

We start with Proposition \ref{1isot/theo} which expresses 1-isotopy of ribbon surfaces
through moves of planar diagrams. Namely, up to isotopy of planar diagrams, we have
three types of moves: those which change the ribbon surface by 3-dimensional diagram
isotopy while preserving the ribbon intersections and the core graph, described in
Figure \ref{ribbsurf15/fig} (p. \pageref{ribbsurf15/fig}); those which change the core
graph, depicted in Figure \ref{ribbsurf16/fig} (p. \pageref{ribbsurf16/fig}); those
which interpret the 1-isotopy moves of Figure \ref{ribbsurf5/fig} in terms of planar
diagrams, as shown Figure \ref{ribbsurf17/fig} (p. \pageref{ribbsurf17/fig}). Of
course, in this context we disregard the moves of Figure \ref{ribbsurf14/fig} (p.
\pageref{ribbsurf14/fig}), since isotopy of planar diagrams is not required to preserve
$y$-coordinate. These last moves will become relevant when we will switch to the
categorical point of view of Proposition \ref{surfcat/theo}.

\begin{proposition}\label{1isot/theo}
 The ribbon surfaces represented by two planar diagrams are 1-isotopic if and only if
their diagrams are related by a finite sequence of planar diagram isotopies (induced 
by ambient isotopies of the projection plane) and moves as in Figures
\ref{ribbsurf15/fig}, \ref{ribbsurf16/fig} and \ref{ribbsurf17/fig} 
(all considered up to planar diagram isotopy).
\end{proposition}

\begin{proof}
 The ``if'' part is trivial, since all the moves in Figures \ref{ribbsurf15/fig},
\ref{ribbsurf16/fig} represent special 3-dimensional diagram isotopies, while
the moves in Figure \ref{ribbsurf17/fig} are equivalent to the 1-isotopy moves in
Figure \ref{ribbsurf5/fig} up to 3-dimensional diagram isotopy.

In order to prove the ``only if'' part, we consider two planar diagrams representing 
ribbon surfaces $F_0$ and $F_1$ as regular neighborhoods of their core graphs $G_0$ and
$G_1$, such that there is a 3-dimensional diagram isotopy $H :(F,G) \times [0,1] \to
R^3$ taking $(F_0,G_0)$ to $(F_1,G_1)$. Contrary to our general convention, here we
think of $F$ as a 3-dimensional diagram, i.e a singular surface with ribbon
self-intersections in $R^3$. Notice that the intermediate pairs $(F_t,G_t) =
H((F,G),t)$ with $0 < t < 1$ do not necessarily project suitably into $R^2$ to give
planar diagrams.

Of course, we can assume that $H$ is smooth, as a map defined on a pair of smooth
stratified spaces, and that the graph $G_t$ regularly projects to a diagram in $R^2$
for every $t \in [0,1]$, except a finite number of $t$'s corresponding to extended
Reidemeister moves for graphs. For such exceptional $t$'s, the lines tangent to $G_t$
at its vertices are assumed not to be vertical.

We define $\Gamma \subset G \times [0,1]$ as the subspace of pairs $(x,t)$ for which
the plane $T_{x_t}F_t$ tangent to $F_t$ at $x_t = H(x,t)$ is vertical (if $x \in G$ is
a singular vertex, there are two such tangent planes and we require that one of them is
vertical).

By a standard transversality argument, we can perturb $H$ in such a way that:
\begin{itemize}\itemsep0pt
\item[{\sl a}\/)]\vskip-\lastskip 
$\Gamma$ is a graph imbedded in $G \times [0,1]$ as a smooth stratified subspace
of constant codimension 1 and the restriction $\eta: \Gamma \to [0,1]$ of the height
function $(x,t) \mapsto t$ is a Morse function on each edge of $\Gamma$;
\item[{\sl b}\/)] the edges of $\Gamma$ locally separate regions consisting of points
$(x,t)$ for which the projection of $F_t$ into $R^2$ has opposite local orientations
at $x_t$;
\item[{\sl c}\/)] the two planes tangent to any $F_t$ at a singular vertex of $G_t$ are
not both vertical, and if one of them is vertical then it does not contain both the
lines tangent to $G_t$ at that vertex.
\end{itemize}\vskip-\lastskip

As a consequence of {\sl b}\/), for each flat vertex $x \in G$ of valency one (resp.
three) there are finitely many points $(x,t) \in \Gamma$, all of which have the same 
valency one (resp. three) as vertices of $\Gamma$. Similarly, as a consequence of
{\sl c}\/), for each singular vertex $x \in G$ there are finitely many points $(x,t)
\in \Gamma$, all of which have valency one or two as vertices of $\Gamma$. Moreover,
the above mentioned vertices of $\Gamma$ of valency one or three are the only 
vertices of $\Gamma$ of valency $\neq 2$.

Let $0 < t_1 < t_2 < \dots < t_k < 1$ be the critical levels $t_i$ at which one of
the following facts happens:
\begin{itemize}\itemsep0pt
\item[{\sl 1}\/)]\vskip-\lastskip 
$G_{t_i}$ does not project regularly in $R^2$, since
there is one point $x_i$ along an edge of $G$ such that the line tangent to $G_{t_i}$
at $H(x_i,t_i)$ is vertical;
\item[{\sl 2}\/)] $G_{t_i}$ projects regularly in $R^2$, but the projection of
$G_{t_i}$ is not a graph diagram, due to a multiple tangency or crossing;
\item[{\sl 3}\/)] there is one point $(x_i,t_i) \in \Gamma$ with $x_i$ a uni-valent or
a singular vertex of $G$;
\item[{\sl 4}\/)] there is one critical point $(x_i,t_i)$ for the function $\eta$ along
an edge of $\Gamma$.
\end{itemize}\vskip-\lastskip

Without loss of generality, we assume that only one of the four cases above occurs
for each critical level $t_i$. Notice that the points $(x,t)$ of $\Gamma$ such that
$x \in G$ is a flat tri-valent vertex represent a subcase of {\sl 2}\/) and for this
reason they are not included in {\sl 3}\/).

\begin{Figure}[b]{ribbsurf6/fig}{}
 {Reversing a wrong projection by an auxiliary half-twist.}
\centerline{\fig{}{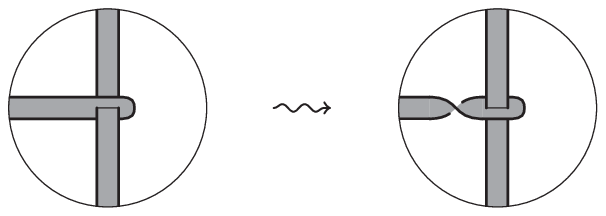}}
\end{Figure}

For $t \in [0,1] - \{t_1, t_2, \dots, t_k\}$, there exists a sufficiently small
regular neighbor\-hood $N_t$ of $G_t$ in $F_t$, such that the pair $(N_t,G_t)$ projects
to a planar diagram, except for the possible presence of some ribbon intersections
projecting in the wrong way, as in the left side of Figure \ref{ribbsurf6/fig} (notice
the difference with the ribbon intersection in Figure \ref{ribbsurf2/fig}).
We fix this problem by inserting an auxiliary positive half-twists along the tounges
containing those ribbon intersections, as shown in the right side of Figure
\ref{ribbsurf6/fig}. The resulting ribbon surfaces, still denoted by $N_t$, projects
to a true planar diagram.

Actually, we modify the $N_t$'s all together to get a new isotopy where no wrong
projection of ribbon intersection occurs, so that $N_t$ projects to a true planar
diagram for each $t \in [0,1] - \{t_1, t_2, \dots, t_k\}$. Namely, at each critical
level when a wrong projection of a ribbon intersection is going to appear in the
original isotopy, we insert an auxiliary half-twis, to prevent the projection from
becoming wrong. Such half-twist remains close to the ribbon intersection until the
first critical level when the projection becomes good again in the original isotopy
(remember that 3-dimensional diagram isotopy preserves ribbon intersections). At that
critical level we remove the auxiliary half-twist. We remark that the second part of
condition {\sl c}\/) is violated when inserting/removing an auxiliary half-twist
at critical points of type {\sl 2}\/), as it can be seen by looking at Figure
\ref{ribbsurf6a/fig}.

\begin{Figure}[htb]{ribbsurf6a/fig}{}{}
\vskip6pt\centerline{\fig{}{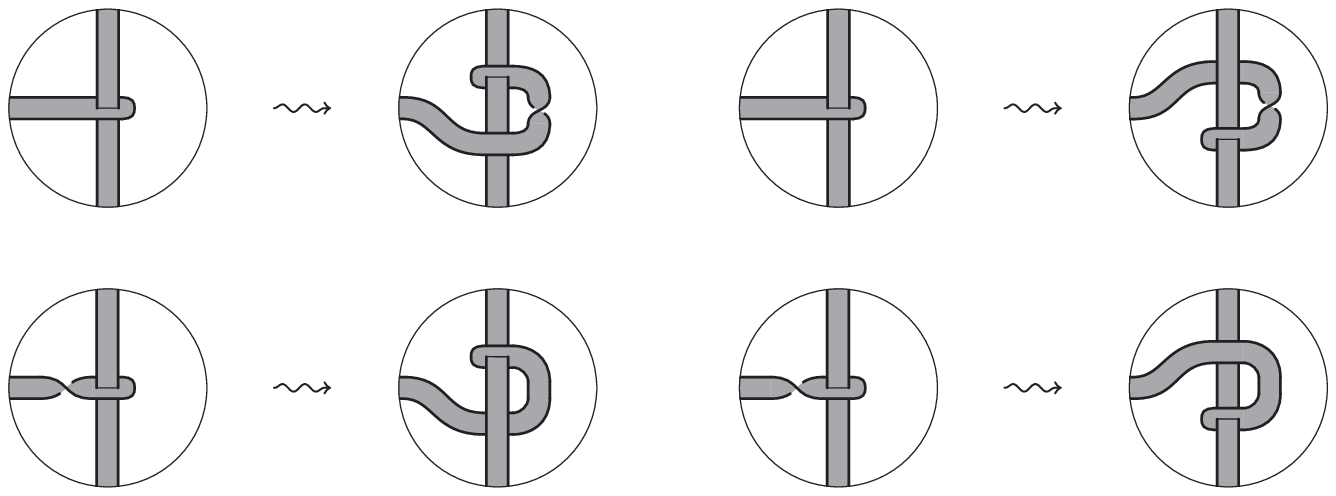}}
\end{Figure}

We observe that the planar diagram of $N_t$ is uniquely determined up to diagram
isotopy by that of its core $G_t$ and by the tangent planes of $F_t$ at $G_t$. 
In fact, the half-twists of $N_t$ along the egdes of $G_t$ correspond to the
transversal intersections of $\Gamma$ with $G \times \{t\}$ and their signs, depend
only on the local behaviour of the tangent planes of $F_t$. In particular, the planar
diagrams of $(N_0,G_0)$ and $(N_1,G_1)$ coincide, up to diagram isotopy, with the
original ones of $(F_0,G_0)$ and $(F_1,G_1)$.

If the interval $[t',t'']$ does not contain any critical level $t_i$, then each single
half-twist persists between the levels $t'$ and $t''$, and hence the planar isotopy
relating the diagrams of $G_{t'}$ and $G_{t''}$ also relate the diagrams of $N_{t'}$
and $N_{t''}$, except for possible slidings of half-twists along ribbons over/under
crossings. Therefore the planar diagrams of $(N_{t'}, G_{t'})$ and $(N_{t''}, G_{t''})$
are equivalent up to diagram isotopy and moves \(I8-9) in Figure \ref{ribbsurf15/fig}.

On the other hand, if the interval $[t',t'']$ is a sufficiently small neighborhood of
a  critical level $t_i$, then the planar diagrams of $N_{t'}$ and $N_{t''}$ are related
by the moves in Figure \ref{ribbsurf15/fig}, depending on the type of $t_i$ as follows.

If $t_i$ is of type {\sl 1}\/), then a positive/negative kink is appearing (resp.
disappearing) along an edge of the core graph. When the kink is positive and $(x_i,
t_i)$ is a local maximum (resp. minimum) point for $\eta$, i.e. two positive
half-twists along the ribbon corresponding to the edge are being converted into a kink
(resp. viceversa), the diagrams of $N_{t'}$ and $N_{t''}$ are directly related by move
\(I11). The cases when $(x_i,t_i)$ is not a local maximum (resp. minimum) point for
$\eta$, that is one or two negative half-twists appear (disappear) together with the
kink, can be reduced to the previous case by means of move \(I12).
On the other hand, by using the regular isotopy moves \(I7-9) in order to create or
delete in the usual way a pair of canceling kinks (without introducing any half-twist)
along the ribbon, we can reduce the case of a negative kink to that of a positive one.

If $t_i$ is of type {\sl 2}\/), then either a regular isotopy move is occurring between
$G_{t'}$ and $G_{t''}$ or two tangent lines at a tri-valent vertex $x_i$ of the graph
project to the same line in the plane. In the first case, the regular isotopy move
occurring between $G_{t'}$ and $G_{t''}$, trivially extends to one of the moves
\(I7-9). In the second case, $x_i$ may be either a flat or a singular vertex. If $x_i$
is a flat vertex, then the tangent plane to $F_t$ at $H(x_i,t)$ is vertical for $t =
t_i$ and its projection reverses the orientation when $t$ passes from $t'$ to $t''$.
Move \(I15) (modulo moves \(I7) and \(I12)) describes the effect on the diagram of
such a reversion of the tangent plane. If $x_i$ is a singular vertex, then $N_{t'}$
changes into $N_{t''}$ in one of the four ways shown in Figure \ref{ribbsurf6/fig}. In
the top (resp. bottom) line auxiliary half-twists are inserted (resp. removed)
according to what we have said above, while left and right columns differ for the
edges which present the common tangency at the critical level. All these modifications
can be reduced to moves \(I14-14') by using moves \(I7) and \(I12).

If $t_i$ is of type {\sl 3}\/), then either a half-twist is appearing/disappearing
at the tip of the tongue of surface corresponding to a uni-valent vertex or one of the
two bands at the ribbon intersection corresponding to a singular vertex is being
reversed in the plane projection. The first case corresponds to move \(I13) (here we
have a positive half-twist, for dealing with a negative one we combine this move with
\(I12)). The second case may happen in two different ways, depending on which band is
being reversed. If such band is the one passing through the other in the ribbon
intersection, then, we can transform $N_{t'}$ into $N_{t''}$ by applying moves \(I10)
and \(I12). Otherwise, the projection of the ribbon intersection is changing from good
to wrong one or viceversa, and the appearing/disappearing half-twist is compensated by
the auxiliary one up to move \(I12).

Finally, if $t_i$ is of type {\sl 4}\/), a pair of canceling half-twists is appearing
or disappearing along an egde of the graph. This is just move \(I12).

\smallskip

At this point, in order to conclude that the moves in the statement of the theorem
suffice to realize 3-dimensional diagram isotopy between any two planar diagrams of a
given ribbon surface $F$, it is left to prove that, given two different core graphs
$G',G''$ of $F$ as above, the  planar diagrams $F'$ and $F''$ determined respectively
by $G'$ and $G''$, are related by those moves. This is quite straightforward. In fact,
by cutting $F$ along the ribbon intersection arcs, we get a new surface $\hat F$ with
some marked arcs. This operation also makes the graphs $G'$ and $G''$ into two simple
spines $T'$ and $T''$ of $\hat F$ relative to those marked arcs (Figure
\ref{ribbsurf11/fig} shows the effect of the cut at the ribbon intersections in Figure
\ref{ribbsurf4/fig}).

\begin{Figure}[htb]{ribbsurf11/fig}{}{}
\vskip6pt\centerline{\fig{}{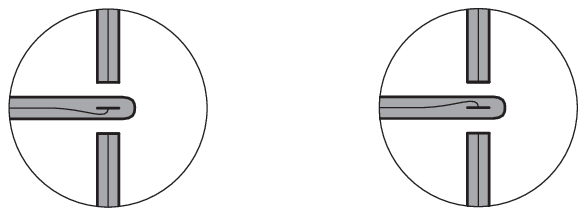}}
\end{Figure}

From intrinsic point of view, that is considering $\hat F$ as an abstract surface and
forgetting its inclusion in $R^3$, the theory of simple spines tells us that the moves
in Figures \ref{ribbsurf4/fig} and \ref{ribbsurf16/fig} suffice to transform $T'$ into
$T''$. In particular, the first and last moves in Figure \ref{ribbsurf16/fig}
correspond to the well known moves for simple spines of surfaces, while the one in the
midle relates the different positions of the spine with respect to the marked arcs in
the interior of $\hat F$. It remains only to observe that, up to a 3-dimensional
diagram isotopy preserving the core graph, hence up to the moves in Figure
\ref{ribbsurf15/fig}, the portion of the surface involved in each single spine
modification can be isolated in the planar diagram as in Figure \ref{ribbsurf16/fig}.

An analogous observation holds for the moves in Figure \ref{ribbsurf17/fig}, which
realize in terms of planar diagrams the 1-isotopy moves of Figure \ref{ribbsurf5/fig},
up to 3-dimensional diagram isotopy. 
\end{proof}

Now we need to restrict ourselves to planar diagrams which are oriented in a special
way with respect to the $y$-axis. This requires an explicit description of the ambient
isotopy in the projection plane in terms of moves.

\smallskip

A planar diagram of an \rs-tangle is said to be in {\sl normal position} with respect
to the $y$-axis if its core graph satisfies the following properties: 
\begin{itemize}\itemsep0pt
\item[{\sl a}\/)]\vskip-\lastskip
 each edge projects to a regular smooth arc immersed in $R \times [0,1]$, such that
the $y$-coordinate restricts to a Morse function on it;
\item[{\sl b}\/)] vertices, half-twists, crossings and local minimum/maximum points for
the $y$-coor\-dinate along the edges have all different $y$-coordinate (in particular,
there are no horizontal tangencies at vertices, half-twists and crossings).
\end{itemize}\vskip-\lastskip

\smallskip

Figure \ref{ribbsurf10/fig} shows the different ways, up to plane isotopy
preserving $y$-coordinate, to put the spots of Figures \ref{ribbsurf2/fig} (p.
\pageref{ribbsurf2/fig}) and \ref{ribbsurf3/fig} (p. \pageref{ribbsurf3/fig}) in
normal position with respect to the $y$-axis, by planar diagram isotopies wich do not
introduce any local minimum/maximum for the $y$-coordinate along the edges of the core
graph.

\begin{Figure}[htb]{ribbsurf10/fig}{}{}
\vskip6pt\centerline{\fig{}{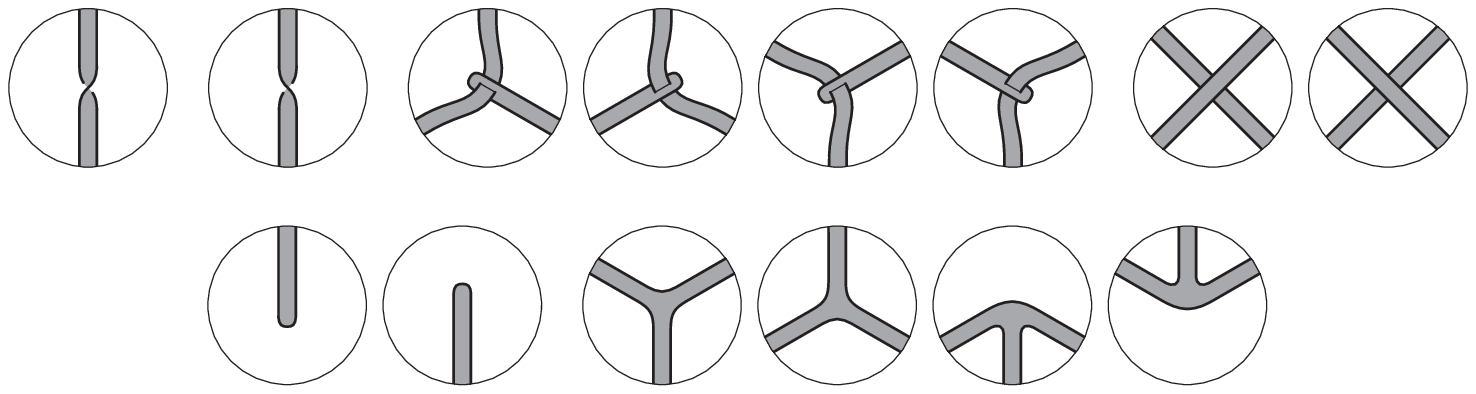}}
\end{Figure}

We notice that all such local configurations appear among the elementary \rs-tangle
diagrams in Figure \ref{ribbsurf13/fig} (p. \pageref{ribbsurf13/fig}), except for some
of those at tri-valent vertices of the core graph. Namely, only the first two of those
at a singular vertex (that means at a ribbon intersection) and the first of those at a
flat tri-valent vertex are considered as elementary \rs-tangle diagrams. The others
can be expressed in terms of them as in Figure \ref{ribbsurf10a/fig}.

\begin{Figure}[htb]{ribbsurf10a/fig}{}{}
\vskip4pt\centerline{\fig{}{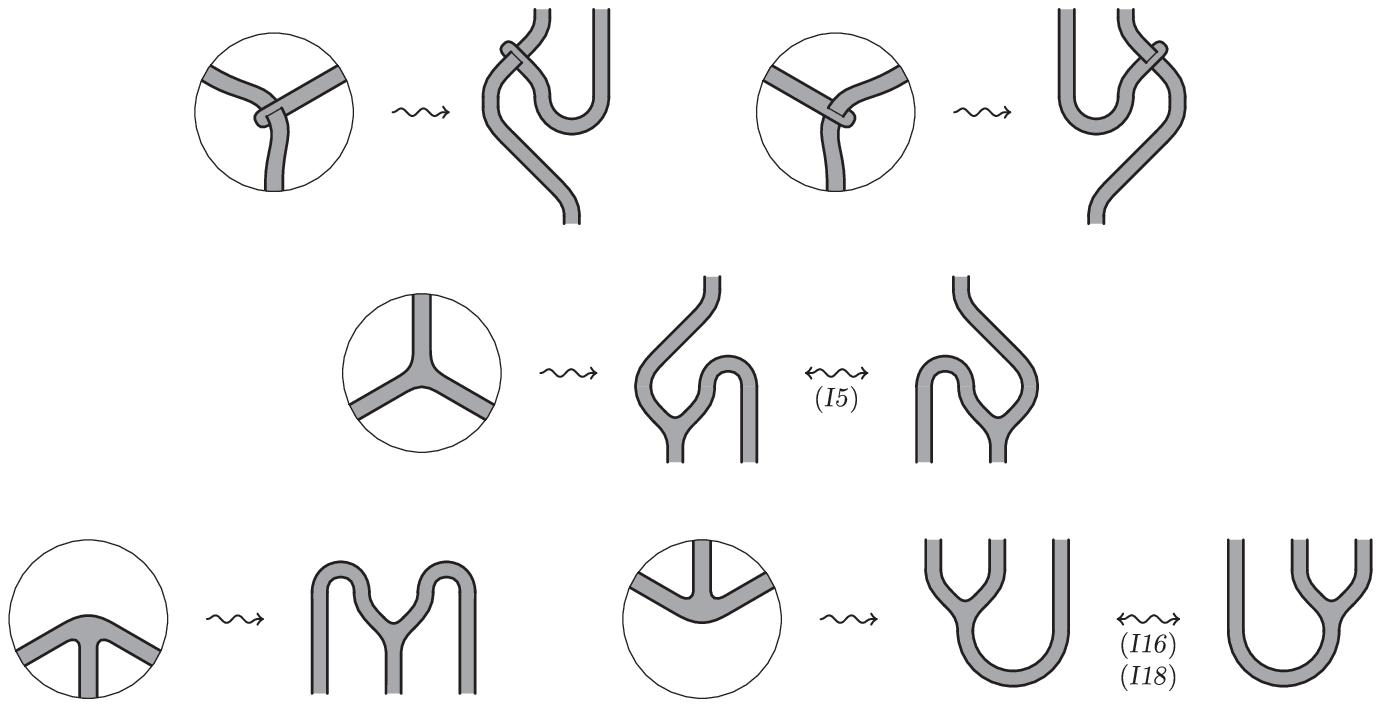}}
\end{Figure}

On the other hand, all the elementary \rs-tangle diagrams in Figure
\ref{ribbsurf13/fig} are in normal position with respect to the $y$-axis. 
Hence, this is also true for any iterated product/composition of them, up to vertical
perturbations to get property {\sl b}\/).

\begin{lemma}\label{vertstat/theo} 
 Through planar diagram isotopy, any planar diagram of an \rs-tangle can be 
 presented as an iterated product/composition of the elementary ones of Figure
\ref{ribbsurf13/fig} in normal position with respect to the $y$-axis. Moreover, any
two such presentations of the same \rs-tangle are related by a finite sequence of
plane isotopies preserving the $y$-coordinate, moves as in Figures
\ref{ribbsurf14/fig} and \ref{ribbsurf10b/fig} (all considered up to plane isotopy
preserving the $y$-coordinate).
\end{lemma}

\begin{Figure}[htb]{ribbsurf10b/fig}{}{}
\centerline{\fig{}{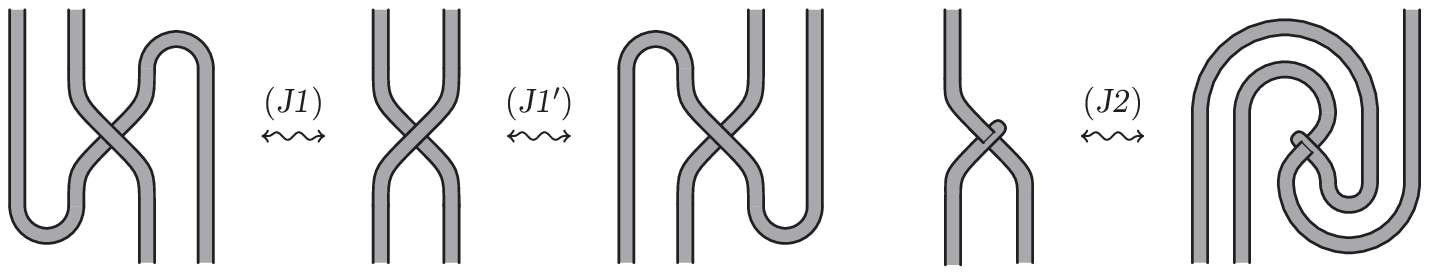}}
\end{Figure}

\begin{proof}
 Here we think of a planar diagram of the core graph of an \rs-tangle as a planar
graph in itself, whose vertices, other than the uni- and tri-valent flat ones and the
singular tri-valent ones, also include bi- and four-valent vertices respectively at the
half-twists and at the crossings.

In the light of the discussion above, the first part of the statement is essentially
trivial. In fact, any planar diagram can be perturbed to get normal position with
respect to the $y$-axis and then made into an iterated product/composition of the
elementary \rs-tangle diagrams by local isotopies as in Figure \ref{ribbsurf10a/fig}.

The proof of the second part is analogous to the proof of the previous Proposition. 
We start with an arbitrary smooth planar diagram isotopy relating any two given
presentions as in the stament. Then we use transversality to perturb this isotopy in
such a way the diagram is in normal position with respect to the $y$-axis at all the
levels except a finite number of critical ones. At these critical levels the planar
diagram of the core graph presents exactly one of the following properties:
\begin{itemize}\itemsep0pt
\item[{\sl 1}\/)]\vskip-\lastskip 
the $y$-coordinate on one edge is not a Morse function;
\item[{\sl 2}\/)] one of the vertices (including half-twists and crossings) has a
horizontal tangent;
\item[{\sl 3}\/)]
two points among the extremal ones along edges and the vertices (including half-twists
and crossings) have the same $y$-coordinate.
\end{itemize}\vskip-\lastskip

Away from critical levels the isotopy can be assumed to preserve the $y$-coordinate,
while, as it is argued below, the moves of Figures \ref{ribbsurf14/fig} and
\ref{ribbsurf10b/fig} allow us to realize all the changes occurring in the planar
diagram when passing through one critical level.

In fact, the cases of critical levels of types {\sl 1}\/) and {\sl 3}\/) are
respectively covered by moves \(I2-2') and \(I1). At critical levels of type {\sl 2}\/)
the vertex with horizontal tangency is swiching from one to another of its normal
positions depicted in Figure \ref{ribbsurf10/fig}. At the same time one extremal point
(resp. one pair of canceling extremal points) for the $y$-coordinate is
appearing/disappearing along the edge (resp. the opposite edges) presenting the
horizontal tangency. The cases when the vertex we are considering is a uni-valent flat
vertex, a half-twists or a crossing correpond respectively to moves \(I3), \(I4) and
\(J1-1'), modulo moves \(I1) and \(I2-2'). On the other hand, in order to deal with
singular and tri-valent flat vertices, we need to replace the normal positions which
are missing in the elementary diagrams as indicated in Figure \ref{ribbsurf10a/fig}.
After that, modulo moves \(I1) and \(I2-2'), all the cases reduce to moves \(I6) and
\(J2) for singular vertices, and to move \(I5) for tri-valent vertices.
\end{proof}

{\sl Proof of Proposition \ref{surfcat/theo}.}\kern1ex
First of all we observe that in Proposition \ref{1isot/theo} the moves of Figures
\ref{ribbsurf15/fig}--\ref{ribbsurf17/fig} are considered up to plane isotopy,
disregarding their normal position with respect to the $y$-axis. On the contrary, here
those moves need to be interpreted in a more restrictive way, assuming that they are
in the preferred normal position given in the figures, up to plane isotopy preserving
the $y$-coordinate. However, as a consequence of Lemma \ref{vertstat/theo}, the two
points of view coincide in the presence of the moves of Figures \ref{ribbsurf14/fig}
and \ref{ribbsurf10b/fig}.

\begin{Figure}[b]{ribbsurf10c/fig}{}{}
\vskip-3pt\centerline{\fig{}{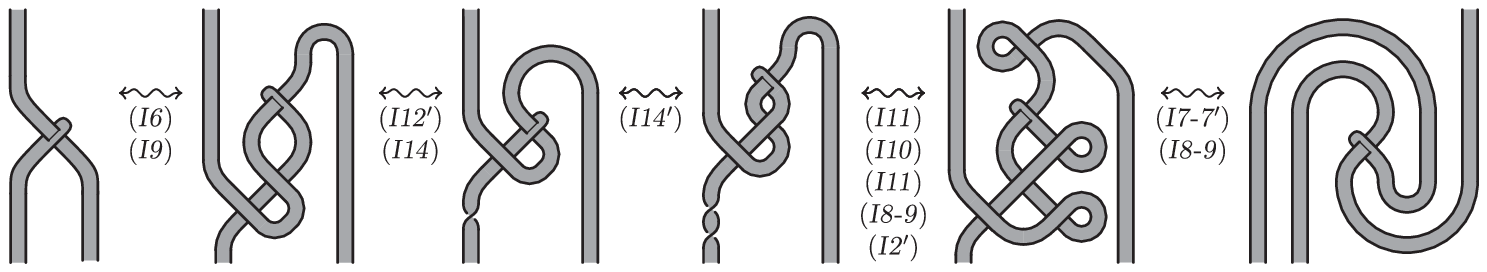}}
\end{Figure}

Then Proposition \ref{surfcat/theo} immediately follows from Proposition
\ref{1isot/theo} and Lemma \ref{vertstat/theo}, once the auxiliary moves in Figure
\ref{ribbsurf10b/fig} are shown to be consequences of those in Figures 
\ref{ribbsurf14/fig}--\ref{ribbsurf17/fig}.
Actually, move \(J1-1') can be easily derived from the moves \(I2-2'), \(I7-7) and
\(I8-9), while Figure \ref{ribbsurf10c/fig} describes how to get \(J2).
$\;\square$

\section{Appendix: proof of some relations in $\H^r(\G)$%
\label{appendix/sec}}

Here we prove the relations in Figures \ref{ribbon-isot/fig} and \ref{ribbon-tot/fig}
and show that \(r8') is a consequence of the rest of the axioms for $\H^r(\G)$.

We start with the ``easy'' isotopy moves in Figure \ref{ribbon-isot/fig}, postponing
the proof of \(p3) to the end.  \(p4) is a direct consequence of the definition \(f2)
in Figure \ref{form/fig} and of the axioms \(r4) and \(r5-5') in Figure
\ref{ribbon1/fig}. Concerning \(p5) and \(p5'), they can be seen to be equivalent by
using \(f6-6') in Figure \ref{h-tortile1/fig} and \(p1) in Figure \ref{coform-s/fig}.
On the other hand, they are trivial when $i \neq j$, due to the definition \(r7) and
the duality moves in Figure \ref{form-uni/fig}, while the proof of \(r5') for $i = j$
is given in Figure \ref{proof-p5/fig}.

\begin{Figure}[htb]{proof-p5/fig}{}
 {Proof of \(p5') for $i = j$
 [{\sl f}/\pageref{h-tortile1/fig},
  {\sl r}/\pageref{ribbon2/fig}, 
  {\sl s}/\pageref{antipode/fig}-\pageref{pr-antipode/fig}].}
\centerline{\fig{}{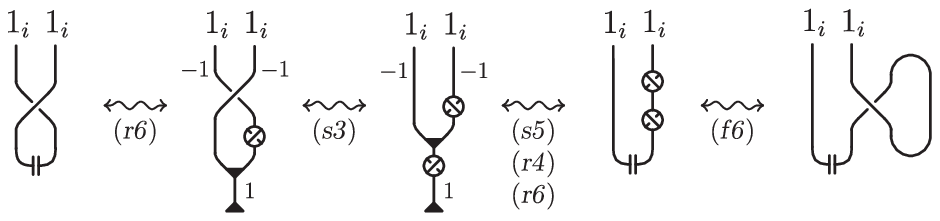}}
\end{Figure}

As already mentioned, \(p5) and \(p5') are particular cases of a more general set moves
which allow to interchange the positions of any coform and copairing appearing on the
same string connecting two polarized vertices. We want to show that \(p5) and \(p5')
actually imply all such moves. Since the braid axioms in Figure \ref{isot/fig} allow
coforms and copairings to pass crissings when we slide them along strings, we can
assume that the coform and the copairing which we want to exchange are contiguous as 
in Figure \ref{proof-p5gen/fig}. Now, the first move in this figure is a trivial
consequence of move \(f3-3') in Figure \ref{form/fig}, while the other two are
equivalent to \(p5) and \(p5') up to \(f3-3') and the braid axioms.

\begin{Figure}[htb]{proof-p5gen/fig}{}{}
\centerline{\fig{}{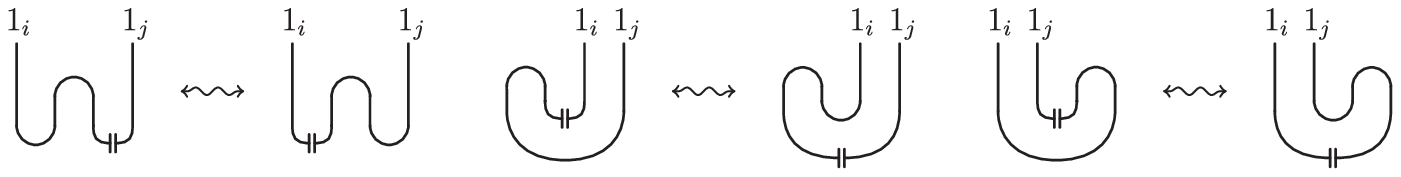}}
\end{Figure}

\begin{Figure}[b]{proof-p6-7/fig}{}
 {Proof of \(p6) and \(p7)
 [{\sl p}/\pageref{ribbon-tot/fig},
  {\sl r}/\pageref{ribbon1/fig}-\pageref{ribbon2/fig},
  {\sl s}/\pageref{antipode/fig}].}
\vskip-8pt\centerline{\fig{}{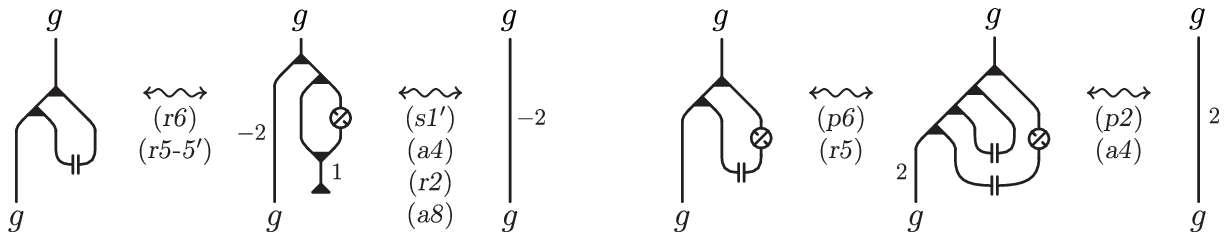}}\vskip-4pt
\end{Figure}

\begin{Figure}[htb]{proof-p8/fig}{}
 {Proof of \(p8)
 [{\sl a}/\pageref{algebra/fig},
  {\sl p}/\pageref{coform-s/fig},
  {\sl r}/\pageref{ribbon1/fig}-\pageref{ribbon4/fig},
  {\sl s}/\pageref{antipode/fig}-\pageref{pr-antipode/fig}].}
\centerline{\fig{}{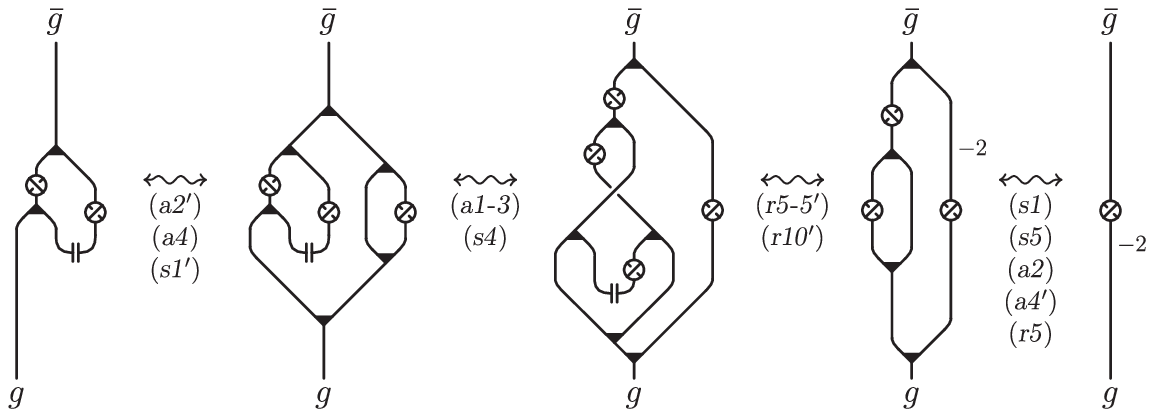}}\vskip-4pt
\end{Figure}

\begin{Figure}[htb]{proof-p10/fig}{}
 {Proof of \(p10) 
 [{\sl a}/\pageref{algebra/fig},
  {\sl r}/\pageref{ribbon1/fig}-\pageref{ribbon2/fig},
  {\sl s}/\pageref{antipode/fig}].}
\centerline{\fig{}{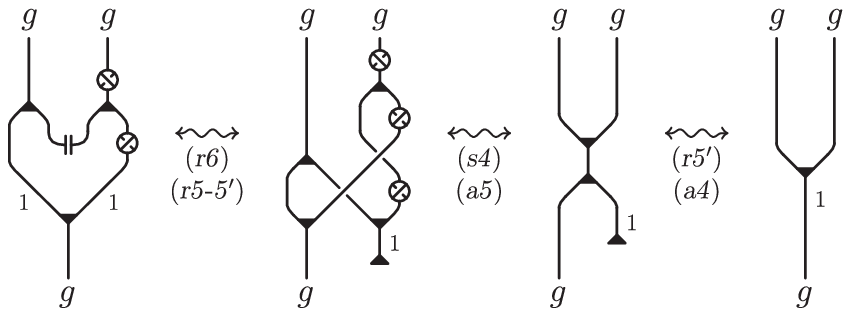}}\vskip-4pt
\end{Figure}

\begin{Figure}[htb]{proof-p11/fig}{}
 {Proof of \(p11) 
 [{\sl a}/\pageref{algebra/fig},
  {\sl f}/\pageref{form/fig},
  {\sl i}/\pageref{unimodular/fig},
  {\sl s}/\pageref{antipode/fig}-\pageref{pr-antipode/fig}, 
  {\sl p}/\pageref{coform-s/fig}-\pageref{ribbon-tot/fig}].}
\centerline{\fig{}{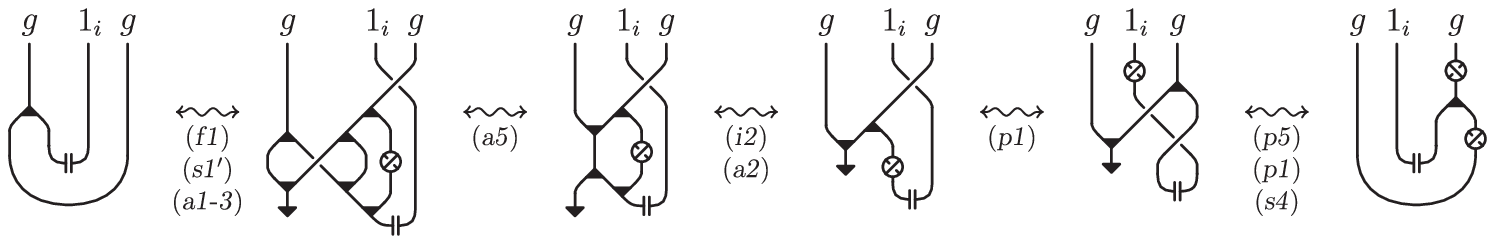}}\vskip-4pt
\end{Figure}

\begin{Figure}[htb]{proof-p12/fig}{}
 {Proof of \(p12)
 [{\sl a}/\pageref{algebra/fig},
  {\sl f}/\pageref{h-tortile1/fig},
  {\sl i}/\pageref{unimodular/fig},
  {\sl p}/\pageref{coform-s/fig},
  {\sl r}/\pageref{ribbon1/fig}-\pageref{ribbon2/fig},
  {\sl s}/\pageref{pr-antipode/fig}].}
\centerline{\fig{}{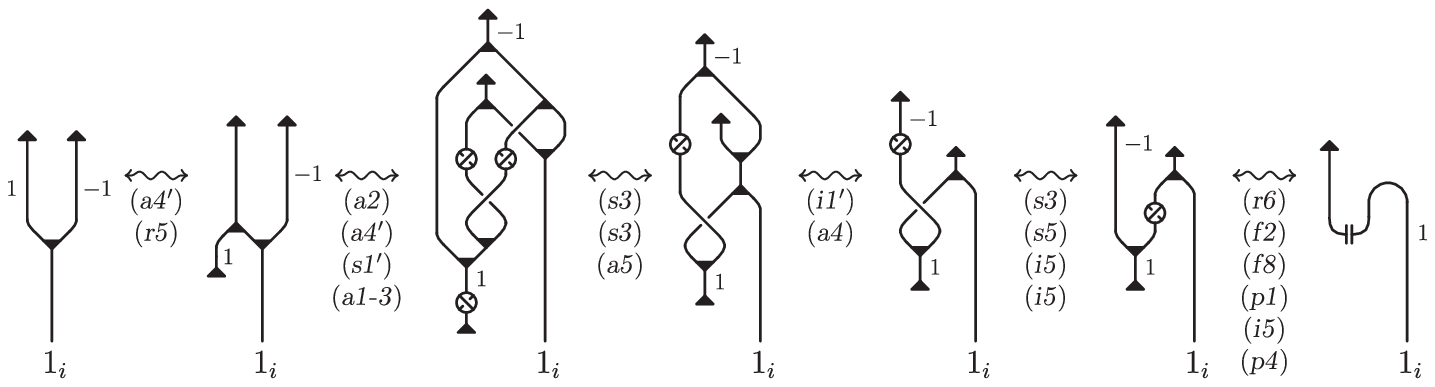}}\vskip-4pt
\end{Figure}

\begin{Figure}[htb]{proof-p3/fig}{}
 {Proof of \(p3) 
 [{\sl f}/\pageref{form/fig},
  {\sl p}/\pageref{ribbon-tot/fig},
  {\sl r}/\pageref{ribbon1/fig}].}
\centerline{\fig{}{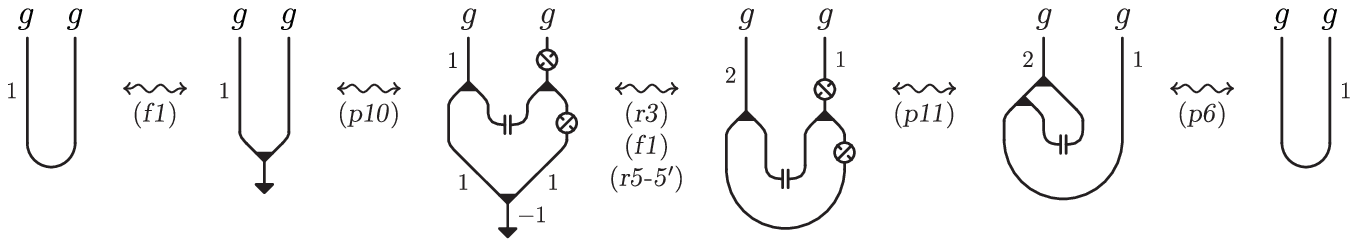}}\vskip-4pt
\end{Figure}

Observe that the proof of \(p5') and therefore of the moves interchanging coforms
copairings appearing on the same string, uses only two of the ribbon axioms: the
definition of copairing (\(r6) in Figure \ref{ribbon2/fig}) and the commutativity
of the ribbon morphism with the antipode (\(r4) in Figure \ref{ribbon1/fig}).
It is easy to see that \(r8') can be obtained from \(r8) through such interchangings 
and the braid axioms. Therefore \(r8') is a consequence of the rest of the axioms of
$\H^r(\G)$. 

The proof of the moves in Figure \ref{ribbon-tot/fig} is given in Figures
\ref{proof-p6-7/fig}--\ref{proof-p12/fig}. We limits ourselves to consider the moves
\(p6)--\(p10), being the proof of the corresponding moves \(p6')--\(p10') analogous.
Observe that \(p9) follows from \(p8) in the same way as \(p7) follows from
\(p6) (cf. Figure \ref{proof-p6-7/fig}) and we leave the proof to the reader.

Finally, we prove \(p3) in Figure \ref{proof-p3/fig}.


\thebibliography{[00]}

\bibitem{BM03} I. Bobtcheva, M.G. Messia, {\sl HKR-type invariants of 4-thickenings
of 2-dimen\-sional CW-complexes}, Algebraic and Geometric Topology {\bf 3} (2003),
33--87.\label{references/sec}

\bibitem{BP} I. Bobtcheva, R. Piergallini, {\sl Covering Moves and Kirby Calculus},
preprint ArXiv:math.GT/0407032.

\bibitem{BQ} I. Bobtcheva, F. Quinn, {\sl The reduction of quantum invariants of
4-thickenings}, Fundamenta Mathematicae {\bf 188} (2005), 21--43.

\bibitem{BKLT97} Y. Bespalov, T. Kerler, V. Lyubashenko, V. Turaev, {\sl Integrals for
braided Hopf algebras}, preprint ArXiv:q-alg/9709020.

\bibitem{FR79} R. Fenn, C. Rourke, {\sl On Kirby's calculus of links}, Topology
{\bf 18} (1979), 1--15.

\bibitem{GS99} R.E. Gompf, A.I. Stipsicz, {\sl 4-manifolds and Kirby calculus},
Grad. Studies in Math. {\bf 20}, Amer. Math. Soc. 1999. 

\bibitem{IP02} M. Iori, R. Piergallini, {\sl 4-manifolds as covers of $S^4$ branched
over non-singular surfaces}, Geometry \& Topology {\bf 6} (2002), 393--401.

\bibitem{H00} K. Habiro, {\sl Claspers and finite type invariants of links},
Geometry \& Topology, {\bf 4} (2000), 1--83.

\bibitem{He96} M. Hennings, {\sl Invariants from links and 3-manifolds obtained from
Hopf algebras}, J. London Math. Soc. (2) {\bf 54} (1996), 594--624.

\bibitem{Ke97} T. Kerler, {\sl Genealogy of nonpertrubative quantum invariants of
3-manifolds -- The surgical family}, in ``Geometry and Physics'', Lecture Notes in
Pure and Applied Physics {\bf 184}, Marchel Dekker 1997, 503--547.

\bibitem{Ke99} T. Kerler, {\sl Bidged links and tangle presentations of cobordism
categories}, Adv. Math {\bf 141} (1999), 207--281.

\bibitem{KL01} T. Kerler, V.V. Lyubashenko, {\sl Non-semisimple topological quantum
field theories for 3-manifolds with corners}, Lecture Notes in Mathematics {\bf 1765}, 
Springer Verlag 2001.

\bibitem{Ke02} T. Kerler, {\sl Towards an algebraic characterization of 3-dimensional
cobordisms}, Contemporary Mathematics {\bf 318} (2003), 141--173.

\bibitem{Ki89} R. Kirby, {\sl The topology of 4-manifolds}, Lecture Notes in
Mathematics {\bf 1374}, Springer-Verlag 1989.

\bibitem{Ku94} G. Kuperberg, {\sl Non-involutory Hopf algebras and 3-manifold
invariants}, Duke Math. J. {\bf 84} (1996), 83--129.

\bibitem{L93} G. Lusztig, {\sl Introduction to quantum groups}, Progress in
Mathematics {\bf 110}, Birkh\"auser 1993.

\bibitem{McL} S. MacLane, {\sl Natural associativities and commutativities}, Rice
Univ. Studies {\bf 49} (1963), 28--46.

\bibitem{MaP92} S. Matveev, M. Polyak, {\sl A geometrical presentation of the surface
mapping class group and surgery}, Comm. Math. Phys. {\bf 160} (1994), 537--550.

\bibitem{M78} J.M. Montesinos, {\sl 4-manifolds, 3-fold covering spaces and ribbons},
Trans. Amer. Math. Soc. {\bf 245} (1978), 453--467.

\bibitem{MP98} M. Mulazzani, R. Piergallini, {\sl Lifting braids}, Rend. Ist. Mat.
Univ. Trieste {\bf XXXII} (2001), Suppl. 1, 193--219.

\bibitem{O02} T. Ohtsuki, {\sl Problems on invariants of knots and 3-manifolds},
Geom. Topol. Monogr. {\bf 4}, in ``Invariants of knots and 3-manifolds (Kyoto, 2001)'',
Geom. Topol. Publ. 2002, 377--572.

\bibitem{P91} R. Piergallini, {\sl Covering Moves}, Trans Amer. Math. Soc. {\bf 325}
(1991), 903--920.

\bibitem{P95} R. Piergallini, {\sl Four-manifolds as $4$-fold branched covers of
$S^4$}, Topology {\bf 34} (1995), 497--508.

\bibitem{PZ03} R. Piergallini, D. Zuddas, {\sl A universal ribbon surface in $B^4$},
Proc. London Math. Soc. {\bf 90} (2005), 763--782.

\bibitem{RT91} N.Yu. Reshetikhin, V.G. Turaev, {\sl Invariants of 3-manifold via
link polynomials and quantum groups}, Invent. Math. {\bf 103} (1991), 547--597.

\bibitem{Ru85} L. Rudolph, {\sl Special position for surfaces bounded by closed
braids}, Rev. Mat. Ibero-Americana {\bf 1} (1985), 93--133; revised version: preprint
2000.

\bibitem{Sh94} M.C. Shum, {\sl Tortile tensor categories}, Journal of Pure and
Applied Algebra {\bf 93} (1994), 57--110.

\bibitem{S96} T. Standford, {\sl Finite type invariants of knots, links and graphs},
Topology {\bf 35} (1996), 1027--1050.

\bibitem{V01} A. Virelizier, {\sl Alg\'ebres de Hopf gradu\'ees et fibr\'es plats
sur les 3-vari\'et\'es}, PhD thesis, Institut de recherche math\'ematique
avanc\'ee, Universit\'e Louis Pasteur et CNRS 2001.
\end{document}